\documentclass[leqno]{book}
\usepackage[english]{babel}
\usepackage[cp1251]{inputenc}
\usepackage{verbatim}
\usepackage{amsmath}
\usepackage{amssymb}
\usepackage{amsfonts}
\usepackage{amsthm}
\usepackage[mathscr]{eucal}
\usepackage{graphicx}
\usepackage{url}
\usepackage{mathtools}
\usepackage{manfnt}
\usepackage[usenames,dvipsnames,svgnames,table]{xcolor}
\usepackage[dvipdfmx]{hyperref}
\usepackage{makeidx}

\def\ig#1#2#3#4{\begin{figure}[!ht]\begin{center}%
\includegraphics[height=#2\textheight]{#1.eps}\caption{#4}\label{#3}%
\end{center}\end{figure}}

\def\thtext#1{
  \catcode`@=11
  \gdef\@thmcountersep{. #1}
  \catcode`@=12
}

\def\threst{
  \catcode`@=11
  \gdef\@thmcountersep{.}
  \catcode`@=12
}

\theoremstyle{plain}
\newtheorem{thm}{Theorem}[chapter]
\newtheorem{prop}[thm]{Proposition}
\newtheorem{cor}[thm]{Corollary}

\newtheorem{lem}[thm]{Lemma}

\theoremstyle{definition}
\newtheorem{agree}[thm]{Agreement}

\newtheorem{constr}{\color{red}Construction}[chapter]

\newtheorem{dfn}[thm]{Definition}

\newtheorem{examp}[thm]{Example}
\newtheorem{exe}{Exercise}[chapter]
\newtheorem*{notation}{Notation}
\newtheorem*{hint}{Hint}
\newtheorem{prb}{\color{blue}Problem}[chapter]
\newtheorem{rk}[thm]{Remark}

\newenvironment{plan}{\scriptsize\trivlist\item[\hskip \labelsep{\bf Schedule.}]}{\endtrivlist}
\newenvironment{sol}{\trivlist\item[\hskip \labelsep{\it Solution.}]}{\endtrivlist}

 \catcode`@=11
 \def\.{.\spacefactor\@m}
 \catcode`@=12

\def\overln#1{
\setbox0=\hbox{$#1$}
\overline{\hbox to 0.8\wd0{%
\vphantom{\hbox{$#1$}}}}
\hskip-0.85\wd0\hbox{$#1$}
}

\def\underln#1{
\setbox0=\hbox{$#1$}
\underline{\hbox to 0.8\wd0{%
\vphantom{\hbox{$#1$}}}}
\hskip-0.85\wd0\hbox{$#1$}
}

\def\dotfill#1{\cleaders\hbox to #1{.}\hfill}
\newcommand\dotline[2][.5em]{\leavevmode\hbox to #2{\dotfill{#1}\hfil}}

\def\N{{\mathbb N}}
\def\Q{{\mathbb Q}}
\def\R{\mathbb R}

\def\Z{{\mathbb Z}}

\def\a{\alpha}

\def\b{\beta}
\def\e{\varepsilon}
\def\dl{\delta}
\def\D{\Delta}
\def\g{\gamma}
\def\G{\Gamma}
\def\l{\lambda}

\def\om{\omega}
\def\Om{\Omega}

\def\r{\rho}
\def\s{\sigma}

\def\t{\tau}

\def\v{\varphi}

\def\0{\emptyset}
\def\:{\colon}
\def\<{\langle}
\def\>{\rangle}

\def\c{\circ}
\def\d{\partial}

\def\fsim{\!\!\sim\,}
\def\hsb#1{\overbracket[0.5pt][0.4ex]{\,#1\,}}
\def\imply{\Rightarrow}

\def\lra{\Leftrightarrow}

\def\mat(#1,#2,#3,#4){\left(\begin{array}{cc}#1&#2\\#3&#4\end{array}\right)}

\def\rom#1{\emph{#1}}

\def\({\rom(}
\def\){\rom)}

\def\simind#1{\mathop{\sim}\nolimits_#1}
\def\sm{\setminus}
\def\ss{\subset}
\def\sp{\supset}
\def\toGH{\xrightarrow{\GH}}
\def\toH{\xrightarrow{d_H}}

\def\Vec(#1,#2){\left(\begin{array}{c}#1\\ #2\end{array}\right)}
\def\x{\times}

\def\bA{{\bar A}}

\def\bB{{\bar B}}

\def\bd{\overline{d}}

\def\bga{{\bar\gamma}}
\def\bG{{\bar G}}

\def\bR{{\bar R}}

\def\bsg{{\bar\s}}

\def\bx{{\bar x}}
\def\by{{\bar y}}
\def\bX{\overline{\strut X}}
\def\bY{\overline{\strut Y}}

\def\ud{\underline{d}}

\def\CL{\operatorname{CL}}
\def\const{\operatorname{const}}

\def\cov{\operatorname{cov}}

\def\diam{\operatorname{diam}}
\def\dil{\operatorname{dil}}
\def\dis{\operatorname{dis}}

\def\Fin{\operatorname{Fin}}
\def\GH{\operatorname{GH}}
\def\GHlim{\operatorname{GH-lim}}
\def\GL{\operatorname{GL}}

\def\id{\operatorname{id}}
\def\im{\operatorname{im}}

\def\Iso{\operatorname{Iso}}

\def\Int{\operatorname{Int}}
\def\mf{\operatorname{mf}}
\def\MF{\operatorname{MF}}

\def\mst{\operatorname{mst}}
\def\MST{\operatorname{MST}}

\def\O{\operatorname{O}}
\def\opt{{\operatorname{opt}}}
\def\pack{\operatorname{pack}}

\def\SL{\operatorname{SL}}

\def\smt{\operatorname{smt}}
\def\SMT{\operatorname{SMT}}

\def\SO{\operatorname{SO}}

\def\U{\operatorname{U}}

\def\cA{\mathcal{A}}

\def\cB{\mathcal{B}}
\def\cC{\mathcal{C}}
\def\cCL{\mathcal{C\!L}}
\def\cD{\mathcal{D}}

\def\cF{\mathcal{F}}

\def\cH{\mathcal{H}}

\def\cK{\mathcal{K}}
\def\cL{\mathcal{L}}
\def\cM{\mathcal{M}}

\def\cP{\mathcal{P}}
\def\cR{\mathcal{R}}

\def\cT{\mathcal{T}}
\def\cU{\mathcal{U}}
\def\cV{\mathcal{V}}

\def\cY{\mathcal{Y}}

\begin{document}
 \title{Lectures on Hausdorff and Gromov--Hausdorff Distance Geometry\\
 The course was given at Peking University, Fall 2019}
 \author{Alexey A. Tuzhilin}
 \date{}

 \maketitle

 \tableofcontents

 \chapter{Elements of general topology.}\label{ch:GenTop}
 \markboth{\chaptername~\thechapter.~Elements of general topology.}%
          {\chaptername~\thechapter.~Elements of general topology.}

\begin{plan}
Definition of topology and topological space, induced topology, subspace of topological space, discrete and anti-discrete topologies, metric spaces and metric topology, standard topology on Euclidean space, base of topology, cover of set and subset, Zariski topology, Sorgenfrey topology, subbase of topology, disjoint union of topological spaces, Cartesian product of topological spaces, Tychonoff or product topology, quotient topology and quotient space, Vietoris topology, continuous mapping, homeomorphism, embedding, convergence of sequences, closure, interior, boundary, dense subsets, separability, separated or Hausdorff topological space,  connected and disconnected topological spaces, connected components, path-connected topological spaces, open cover, subcover, compact and sequentially compact topological spaces, bounded metric spaces, hyperspaces.
\end{plan}

In this chapter we present an introduction to general topology.

\section{Basic facts of general topology}
\markright{\thesection.~Basic facts of general topology}

For a set $X$, let $2^X$ denote the collection of all subsets of $X$. If $\cA\ss2^X$ is a family of subsets of $X$, then $\cup\cA$ and $\cap\cA$ denote the union and the intersection of the elements from $\cA$, respectively. If $\cA$ is an indexed family, i.e., $\cA=\{A_i\}_{i\in I}$, then we use $\cup_{i\in I}A_i$ and $\cap_{i\in I}A_i$ for the union and the intersection. If different elements of $\cA$ do not intersect each other (such family $\cA$ is called \emph{disjoint\/}), then we write $\sqcup\cA$ instead of $\cup\cA$, to emphasize that $\cA$ is disjoint; similarly, we write $\sqcup_{i\in I}A_i$ instead of $\cup_{i\in I}A_i$ for indexed families. We can define $\sqcup_{i\in I}A_i$ also in the case when some differen $A_i$ intersect each other, in particular, when they coincide. In this situation we simply consider $A_i$ for different $i$ as nonintersecting sets. This can be done in a formal way if we change the elements $a_i\in A_i$ to $(a_i,i)$ and identify $A_i$ with the set $\bigl\{(a_i,i)\bigr\}_{a_i\in A_i}$.

\begin{dfn}
A set $\t=\{U_\a\}_{\a\in A}\ss2^X$ is called \emph{a topology on $X$} if $\t$ satisfies the following properties:
\begin{enumerate}
\item $\0,\,X\in\t$;
\item for any $\cA\ss\t$ we have $\cup\cA\in\t$ (the union of arbitrary collection of elements from $\t$ belongs to $\t$);
\item for any \textbf{finite\/} $\cA\ss\t$ it holds $\cap\cA\in\t$ (the intersection of arbitrary finite collection of elements from $\t$ belongs to $\t$).
\end{enumerate}
\end{dfn}

\begin{dfn}
A set $X$ with a given topology $\t\ss2^X$ is called \emph{a topological space}. It is also convenient to denote the topological space $X$ as the pair $(X,\t)$. Also, speaking about the topological space $X$, we will often denote the topology defined on it by $\t_X$, without specifically mentioning it every time.

The elements of $X$ are usually called \emph{points}, and the elements of $\t$ are called \emph{open sets}. A set $F\ss X$ is called \emph{closed\/} if its complement is open.
\end{dfn}

\begin{prb}
Show that the family of all closed subsets of a topological space $X$ contains $\0$ and $X$, and that the intersection of any collection of closed subsets, as well as the union of any finite collection of closed subsets are some closed sets.
\end{prb}

Take an arbitrary $Y\ss X$ and consider the family $\t_Y:=\{U\cap Y: U\in\t_X\}$.

\begin{prb}
Prove that $\t_Y$ is a topology on $Y$.
\end{prb}

\begin{dfn}
The $\t_Y$ is called \emph{the topology on $Y$ induced from $X$}. The $Y$ with the topology $\t_Y$ is usually called \emph{a subspace\/} of the topological space $X$.
\end{dfn}

Generally speaking, there are many different topologies on each set $X$, and the inclusion relation generates a partial order on the set of all such topologies. The smallest topology in this order consists of two elements: $\t_a=\{\0,X\}$. It is called \emph{anti-discrete}. The largest topology consists of all subsets: $\t_d=2^X$. It is called \emph{discrete}. If $\cT$ is a collection of topologies defined on the same set $X$, then $\cap\cT$ is a topology as well; it is smaller than each topology $\t\in\cT$; for any topology $\t'$ on $X$ that is smaller than all topologies from $\cT$ it holds $\t'\ss\cap\cT$.

The most important  for us example of topology  will be generated by metrics. Namely, suppose that for a set $X$ a function $\r\:X\x X\to\R$ is given, which has the following properties:
\begin{enumerate}
\item for any $x,y\in X$ it holds $\r(x,y)\ge0$, and $\r(x,y)=0$ iff $x=y$ (positive definiteness);
\item for any $x,y\in X$ we have $\r(x,y)=\r(y,x)$ (symmetry);
\item for any $x,y,z\in X$ it holds $\r(x,y)+\r(y,z)\ge\r(x,z)$ (triangle inequality).
\end{enumerate}
Then $\r$ is called \emph{a metric}, and the set $X$ with the metric $\r$ is called \emph{a metric space}. It is also convenient to denote the metric space $X$ as the pair $(X,\r)$. Each $Y\ss X$ endowed with the restriction of $\r$ to $Y\x Y$ is called \emph{a subspace of $X$}.

\begin{examp}\label{examp:EuclideanMetric}
In calculus, the standard example of a metric space is the arithmetic space $\R^n$ with the Euclidean metric defined on it: for $x=(x_1,\ldots,x_n)$ and $y=(y_1,\ldots,y_n)$ it holds $\r(x,y)=\sqrt{\sum_{i=1}^n(x_i-y_i)^2}$. This metric is called \emph{Euclidean}. We will also call \emph{Euclidean\/} the arithmetic space itself, endowed with the Euclidean metric.
\end{examp}

Let $X$ be a metric space with a metric $\r$. For every $x\in X$ and $r>0$ we put
$$
U_r(x):=\bigl\{y\in X:\r(x,y)<r\bigr\}
$$
and call \emph{an open ball of radius $r>0$ and center $x$}. Using the metric $\r$, we construct the natural topology $\t_\r$, called \emph{the metric topology\/}: we assign the subset $U\ss X$ to open sets of the metric topology $\t_\r$ if and only if $U$ is either empty, or $U$ is a union of open balls. Equivalent definition: $U\in\t_\r$ if and only if for any point $x\in U$ there exists $r>0$ such that $U_r(x)\ss U$ (together with each point the set $U$ contains some open ball with the center at this point).

\begin{prb}
Prove that the family $\t_\r$ is a topology.
\end{prb}

\begin{rk}\label{rk:StandTop}
Unless otherwise stated, on the real line $\R$ and, more generally, on the arithmetic space $\R^n$, we consider the topologies generated by the Euclidean metric (see\. Example~\ref{examp:EuclideanMetric}). This topology is called \emph{standard}.
\end{rk}

\section{Base and subbase}
\markright{\thesection.~Base and subbase}
The construction of metric topology described above leads to the following important notion. Similar to linear algebra, where to describe a linear space it is enough to choose a family of vectors that, using linear combinations, generates the whole space, to define a topology, one can also select a subfamily of open sets and generate the topology by means of set-theoretic operations.

\begin{dfn}
A family $\b\ss\t$ is called \emph{a base of the topology $\t$} if every nonempty open set $U\in\t$ is representable as a union of some elements from $\b$.
\end{dfn}

Thus, by the definition of metric topology, its base is the family of all open balls.

\begin{rk}
Note that a given topology can have many different bases. For example, not all balls can be selected as the base of the metric space, but only, say, of radii not exceeding $1$, or of only rational radii, or of only radii of the form $1/n$, etc. On the Euclidean line, for example, only rational numbers can be selected as centers.
\end{rk}

We note two important properties of the topology base on the set $X$:
\begin{enumerate}
\item each point $x\in X$ is contained in some element from the base (otherwise the set $X$ cannot be obtained as the union of some elements from the base);
\item a nonempty intersection of any two elements of the base is representable as the union of some elements from the base (otherwise this intersection will not belong to the topology).
\end{enumerate}

It turns out that these two properties completely characterize the families that are the bases of some topologies. Before formulating the corresponding criterion, we introduce a definition of cover, which will be useful to us both here and hereinafter.

\begin{dfn}
A family $\cA\ss2^X$ is called \emph{a cover of the set $X$} if $X=\cup\cA$. A family $\cA\ss2^X$ is called \emph{a cover of $Y\ss X$} if $Y\ss\cup\cA$.
\end{dfn}

It is clear that each base of a topological space $X$ is a cover of $X$.

\begin{prb}\label{prb:criterion-base}
Prove that a family $\b\ss2^X$ is a base of some topology $\t$ on $X$ if and only if $\b$ is a cover of $X$, and for any intersecting $B_1,B_2\in\b$ their intersection $B_1\cap B_2$ is the union of some elements from $\b$. Moreover, each family satisfying these properties, generates a unique topology.
\end{prb}

Notice that a collection $\b$ of open sets in a topological space $X$ which satisfies the conditions of Problem~\ref{prb:criterion-base} may generate a topology $\t$ different from $\t_X$. What do we need to add for to be sure that $\t=\t_X$? The answer can be obtained from the following more general result that is often used in proving the coincidence of topologies.

\begin{prb}\label{prb:topologies-equality}
Let some topologies $\t_1$ and $\t_2$ with bases $\b_1$ and $\b_2$ be given on a set $X$. Then $\t_1=\t_2$ if and only if for any $x\in X$ the following condition is fulfilled: for any $B_2\in\b_2$, $x\in B_2$, there is $B_1\in\b_1$ for which $x\in B_1\ss B_2$, and vice versa, for any $B_1\in\b_1$, $x\in B_1$, there exists $B_2\in\b_2$ for which $x\in B_2\ss B_1$. In particular, for a topological space $X$, a collection $\b$ of open sets satisfying the conditions of Problem~\ref{prb:criterion-base} is a base of the topology $\t_X$ if and only if for each open set $U\in\t_X$ and any point $x\in U$ there exists some $B\in\b$ such that $x\in B\ss U$.
\end{prb}

Recall that two sets are called \emph{equivalent\/} if there exists a bijection between them. The equivalence classes of sets are called \emph{cardinalities\/} or \emph{cardinal numbers}. The cardinality of a set $X$ will be denoted by $\#X$.

\begin{examp}\label{examp:Zaris}
Let $X$ be an infinite set of cardinality $n$, and $m$ be an infinite cardinal number, with $m\le n$. Consider a family $\cF$ of all $F\ss X$ such that $\#F<m$, and let $\b_m=\{X\sm F:F\in\cF\}$. Then $\b_m$ is a base of some topology $\t$, which we call \emph{the Zariski topology of the weight $m$}.
\end{examp}

\begin{prb}
Prove that the family $\b_m$ from Example~\ref{examp:Zaris} is a base of some topology.
\end{prb}

\begin{examp}\label{examp:Sorgenfrey-line}
As we noted above, on the Euclidean line we can take the family of all intervals as the base of the standard topology. Another interesting example of topology is obtained if, instead of intervals, we take all possible half-intervals of the form $[a,b)$. The corresponding topology is called \emph{the arrow topology\/} or \emph{the Sorgenfrey topology}.
\end{examp}

\begin{rk}
The Sorgenfrey topology contains the standard topology of the line, since each interval $(a,b)$ can be represented as a union of half-intervals $[a+1/n,b)$, $n\in\N$.
\end{rk}

\begin{prb}
Show that the collection of all possible half-intervals of the form $[a,b)\ss\R$ form a base of some topology that contains the standard topology.
\end{prb}

If, to generate a topology, we allowed also to use finite intersections, then the generating family can, generally speaking, be reduced.

\begin{dfn}
A family $\s\ss\t$ is called \emph{a subbase of the topology $\t$} if the set of all finite intersections of elements from $\s$ forms a base of the topology $\t$.
\end{dfn}

It is clear that, like the base, each subbase of a topological space $X$ is a cover of $X$.

\begin{prb}
Prove that a family $\s\ss2^X$ is a subbase of some topology on $X$ if and only if $\s$ is a cover of $X$. Moreover, each cover of $X$ generates a unique topology.
\end{prb}

\begin{examp}
The family of all subsets of the real line $\R$, each of which is an open ray, forms a subbase of the standard topology and is not a base of this topology.
\end{examp}

\section{Standard constructions of topologies}
\markright{\thesection.~Standard constructions of topologies}

This section provides examples of standard constructions that allow to build new examples of topological spaces from existing ones.

\begin{constr}
Let $\s$ be an arbitrary family of subsets of a set $X$, and $\cT_\s$ be the family of all topologies on $X$ containing $\s$. Then $\t:=\cap\cT_\s$ is the smallest topology containing $\s$. If $\s$ is a cover of $X$, then $\s$ is a subbase of $\t$.
\end{constr}

\begin{constr}[Disjoint union]
Let $\cA=\{X_i\}_{i\in I}$ be a family of topological spaces. We define a topology on $\sqcup\cA=\sqcup_{i\in I}X_i$, setting its base to be equal to $\sqcup_{i\in I}\t_{X_i}$. The set $\sqcup_{i\in I}X_i$ with the corresponding topology is called \emph{the disjoint union of the topological spaces $X_i$}.
\end{constr}

\begin{constr}[Cartesian product]\label{constr:Tychonoff}
Let $\cA=\{X_i\}_{i\in I}$ be an arbitrary family of topological spaces. The set of all mappings $w\:I\to\sqcup_{i\in I}X_i$ such that $w(i)\in X_i$ for every $i\in I$  is called \emph{the Cartesian product\/} of the spaces $X_i$ and is denoted by $W:=\prod\cA=\prod_{i\in I}X_i$. In particular, if all $X_i$ are equal to the same space $X$, then $W=X^I$,  where the latter, recall, denotes the set of all mappings from $I$ to $X$. For convenience, we will often write $w_i$ instead of $w(i)$, and we will call this value \emph{the $i$-th coordinate of the point $w\in\prod_{i\in I}X_i$}.

We define a topology on $W$, choosing as its subbase the family of all products $\prod_{i\in I}U_i$, $U_i\in\t_{X_i}$, in which only one $U_i$ can differ from $X_i$. The corresponding base consists of $\prod_{i\in I}U_i$, $U_i\in\t_{X_i}$, in which only a finite number of $U_i$ can differ from $X_i$. This topology is called \emph{the product topology\/} or \emph{the Tychonoff topology}.

In the case when $I$ is a finite set, say, $I=\{1,\ldots,n\}$, then the Cartesian product of the spaces $X_i$ is denoted by $X_1\x\cdots\x X_n$. In particular, in this way one can define the standard topology on the $n$-dimensional arithmetic space $\R^n$.
\end{constr}

\begin{prb}
Show that the standard topology of the Euclidean space $\R^n$ coincides with the topology of the Cartesian product $\R\x\cdots\x\R$ of real lines endowed with the standard topology.
\end{prb}

\begin{examp}\label{examp:Sorgenfrey-space}
\emph{Sorgenfrey space\/} is the Cartesian product of the Sorgenfrey lines from Example~\ref{examp:Sorgenfrey-line} and is used in general topology to illustrate numerous exotic possibilities.
\end{examp}

\begin{constr}[Quotient topology]\label{constr:factor_top}
Let $X$ be a topological space, and $\nu$ be some equivalence relation on $X$. Denote by $X/\nu$ the set of classes of this equivalence. For each $x\in X$ denote by $[x]\in X/\nu$ the $\nu$-equivalence class containing $x$, and let $\pi\:X\to X/\nu$, $\pi\:x\mapsto[x]$, be the canonical projection. Then the family of all $U\ss X/\nu$ such that $\pi^{-1}(U)\in\t_X$ forms a topology called \emph{the quotient topology}. The set $X/\nu$ endowed with the quotient topology is called \emph{the quotient space}.
\end{constr}

\begin{constr}[Vietoris topology]\label{constr:Vietoris_top}
Let $X$ be an arbitrary topological space. For each finite collection of open sets $U_1,\ldots,U_n\in\t_X$ we put
$$
\<U_1,\ldots,U_n\>=\{Y\ss X:Y\ss\cup_{i=1}^nU_i\ \text{and $Y\cap U_i\ne\0$ for all $i=1,\ldots,n$}\}.
$$
Note that if at least one of $U_i$ is empty, then $\<U_1,\ldots,U_n\>=\0$.

\begin{prb}\label{prb:Vietoris}
Show that the families
$$
\s=\bigl\{\<U\>:U\in\t_X\bigr\}\cup\bigl\{\<X,U\>:U\in\t_X\bigr\}\ \ \text{and}\ \
\b=\bigl\{\<U_1,\ldots,U_n\>:U_1,\ldots,U_n\in\t_X\bigr\}
$$
form respectively a subbase and the corresponding base of some topology on $2^X$.
\end{prb}

The topology on $2^X$ defined in Problem~\ref{prb:Vietoris} is called \emph{the Vietoris topology}.
\end{constr}

\begin{rk}
Usually, Vietoris topology is defined on the family of all nonempty closed subsets of a topological space.
\end{rk}

\begin{prb}
Prove that each construction of these section provides a topology.
\end{prb}

\section{Continuous mappings}
\markright{\thesection.~Continuous mappings}

As a rule, all considered mappings between topological spaces are consistent with topologies. Such mappings are called continuous. We give three equivalent definitions of continuity.

\begin{dfn}
\emph{A neighborhood of a point $x\in X$} of a topological space $X$ is an arbitrary open set $U\ss X$ containing $x$. \emph{A neighborhood of a subset $Z$} of a topological space $X$ is an arbitrary open set $U\ss X$ containing $Z$.
\end{dfn}

\begin{rk}
For convenience, we denote an arbitrary neighborhood of a point $x\in X$ by $U^x$.
\end{rk}

\begin{dfn}\label{dfn:cont1}
A mapping $f\:X\to Y$ between topological spaces \emph{is continuous at $x\in X$} if for any neighborhood $U^{f(x)}$ there exists a neighborhood $U^x $ such that $f(U^x)\ss U^{f(x)}$. The mapping $f$, continuous at all points $x\in X$, is called \emph{continuous}.
\end{dfn}

\begin{dfn}\label{dfn:cont2}
A mapping $f\:X\to Y$ between topological spaces \emph{is continuous} if for any open set $U\ss Y$ its preimage $f^{- 1}(U)\ss X$ is open (the preimage of each open set is open).
\end{dfn}

\begin{dfn}\label{dfn:cont3}
A mapping $f\:X\to Y$ between topological spaces \emph{is continuous} if for any closed set $F\ss Y$ its preimage $f^{- 1}(F)\ss X$ is closed (the preimage of each closed set is closed).
\end{dfn}

\begin{prb}
Prove that the definitions~\ref{dfn:cont1}, \ref{dfn:cont2}, and \ref{dfn:cont3} are equivalent.
\end{prb}

\begin{prb}
Let $f\:X\to Y$ be a mapping of topological spaces and $\s$ a subbase of the topology on the space $Y$. Prove that $f$ is continuous if and only if $f$-preimage of each element from the subbase $\s$ is open in $X$.
\end{prb}

\begin{rk}
When we speak of the continuity of a function $f\:X\to\R$ or, more generally, of a vector-valued mapping $f\:X\to\R^n$ from a topological space $X$, then, unless otherwise stated, on $\R$ and $\R^n$ we consider the standard topologies (see\. Remark~\ref{rk:StandTop}).
\end{rk}

\begin{prb}
Show that the identity mapping and the composition of continuous mappings are continuous.
\end{prb}

\begin{prb}
Let $X$ be a topological space, and $Z\ss X$ be its subspace. Show that the inclusion mapping $i\:Z\to X$, $i(z)=z$ for each point $z\in Z$, is continuous.
\end{prb}

\begin{prb}
Let $X$, $Y$ be topological spaces, $W\ss Y$ be a subspace of $Y$, and $f\:X\to W$ be a continuous mapping. Let $g\:X\to Y$ be a mapping coinciding with $f$: for each $x\in X$ it holds $f(x)=g(x)$. Prove that the mapping $g$ is continuous.
\end{prb}

Let $f\:X\to Y$ be an arbitrary mapping of sets. Choose arbitrary subsets $Z\ss X$ and $W\ss Y$ such that $f(Z)\ss W$. Then \emph{the restriction $f|_{Z,W}$ of the mapping $f$ to $Z$ and $W$} is the mapping $g\:Z\to W$ that coincides on the domain with the mapping $f$, i.e., for any $x\in Z$ it holds $f(x)=g(x)$.

\begin{prb}\label{prb:continuous-restriction}
Let $f\:X\to Y$ be a continuous mapping of topological spaces, $Z\ss X$, $W\ss Y$, $f(Z)\ss W$. Then the restriction $f|_{Z,W}\:Z\to W$ is also continuous as the mapping of the topological spaces $Z$ and $W$ with the topologies induced on them from $X$ and $Y$, respectively.
\end{prb}

\begin{prb}
Let $\{X_i\}_{i\in I}$ be a cover of a topological space $X$ by open subsets $X_i$, and $f\:X\to Y$ a mapping to a topological space $Y$. Show that $f$ is continuous if and only if all the restrictions $f|_{X_i}$ are continuous. In particular, this holds when $X=\sqcup_{i\in I}X_i$ is the disjoint union of some topological spaces. Will this result remain true if we replace $\{X_i\}$ with a cover of $X$ by arbitrary sets?
\end{prb}

\begin{prb}
Let $X=\sqcup_{i\in I}X_i$ be the disjoint union of some topological spaces and $f\:X\to Y$ be a map into a topological space $Y$. Show that $f$ is continuous if and only if all its restrictions $f|_{X_i}$ are continuous.
\end{prb}

\begin{prb}
Let $\{X_i\}_{i\in I}$ be a family of topological spaces and $X=\prod_{i\in I}X_i$. We define the canonical projection $\pi_i\:X\to X_i$ by setting $\pi_i(x)=x_i$. Prove that the product topology on $X$ is the smallest of those topologies in which all the projections $\pi_i$ are continuous.
\end{prb}

\begin{prb}
Let $\{Y_i\}_{i\in I}$ be a family of topological spaces, $Y=\prod_{i\in I}Y_i$, and $f_i\:X\to Y_i$ be mappings from some topological space $X$. We construct the mapping $F:=\prod_{i\in I}f_i\:X\to Y$  by associating with each point $x\in X$ the element $y\in Y$ defined as follows: $y_i=f_i(x)$. Prove that the mapping $F$ is continuous if and only if all $f_i$ are continuous.
\end{prb}

\begin{prb}
Let $A\ss\R^n$ be an arbitrary subset, $(x^1,\ldots,x^n)$ the Cartesian coordinates on $\R^n$, $f\:A\to\R^m$ a continuous mapping, $(y^1,\ldots,y^m)$ the Cartesian coordinates on $\R^m$, and $y^i=y^i(x^1,\ldots,x^n)$ the coordinate functions of the mapping $f$. Prove that the mapping $f$ is continuous if and only if all the coordinate functions $y^i=y^i(x^1,\ldots,x^n)$ are continuous.
\end{prb}

\begin{prb}
Describe all continuous functions on a topological space with Zariski topology.
\end{prb}

\section{Homeomorphisms, embeddings}
\markright{\thesection.~Homeomorphisms, embeddings}

An important particular case of continuous mapping is a homeomorphism.

\begin{dfn}
A mapping $f\:X\to Y$ of topological spaces is called \emph{a homeomorphism\/} if it is bijective, and both the maps $f$ and $f^{-1}$ are continuous. Topological spaces between which there is a homeomorphism are called \emph{homeomorphic}.
\end{dfn}

\begin{rk}
A homeomorphism, being a bijection, identifies not only points of spaces, but also identifies topologies, establishing a one-to-one correspondence between them. For clarity, we can imagine that the homeomorphism $f:X\to Y$ is a replacement for the names of points in the space $X$: the point $x\in X$ is ``renamed'' to $f(x)$, without changing the topology. From these considerations it follows that all topological properties of homeomorphic spaces are the same.
\end{rk}

The next exercise follows directly from Problem~\ref{prb:continuous-restriction}.

\begin{prb}
Let $f\:X\to Y$ be a homeomorphism, and $Z\ss X$, $W=f(Z)$. Prove that the restriction $f|_{Z,W}\:Z\to W$ is also a homeomorphism.
\end{prb}

An injective mapping $f\:X\to Y$ of topological spaces is called \emph{an embedding of $X$ into $Y$} if the restriction $f|_{X,f(X)}$ is a homeomorphism.

\begin{prb}
Show that every embedding is continuous. Give an example of a continuous injective mapping of topological spaces that is not an embedding.
\end{prb}

\section{Convergence of sequences}
\markright{\thesection.~Convergence of sequences}

\emph{A sequence\/} in a set $X$ is an arbitrary mapping $x\:\N\to X$ from the set of natural numbers $\N=\{1,2,\ldots\}$. For convenience, the points $x(n)$ are usually denoted by $x_n$. Also, for brevity, it is customary to say that \emph{a sequence of points $x_n$} is given.

\begin{dfn}
A sequence of points $x_n$ in a topological space $X$ is called \emph{convergent\/} if, for some $x\in X$, called \emph{a limit of this sequence}, the following holds: for any neighborhood $U^x$ there exists $N\in\N$ such that for all $n\ge N$ we have $x_n\in U^x$. If the sequence is not convergent, then it is called \emph{divergent}.
\end{dfn}

\begin{prb}
Let $\om$ be a character not contained in $\N$. We define a topology on the set $\bar{N}=\N\cup\{\om\}$, taking as a base all points from $\N$, as well as all sets of the form $\{n\ge N\}\cup\{\om\}$, $N\in\N$. Show that a sequence $x\:\N\to X$ converges if and only if the mapping $x$ can be extended to a continuous mapping on $\bar{\N}$.
\end{prb}

\begin{prb}
Show that a continuous mapping $f\:X\to Y$ of topological spaces takes convergent sequences to convergent ones. Show that if $X$ is  a metric space, then every mapping $g\:X\to Y$ that takes convergent sequences into convergent ones is continuous. Give an example of a topological space $X$ and a mapping $h\:X\to Y$ into a topological space $Y$, which takes convergent sequences into convergent ones, but is not continuous nonetheless.
\end{prb}

\begin{prb}\label{prb:seqInMetricWIthoutConverge}
Let $x_1,x_2,\ldots$ be a sequence of points in a metric space $X$. Suppose that for some point $x\in X$ each neighborhood of $x$ intersects the set $\{x_i\}_{i=1}^\infty\sm\{x\}$. Prove that the sequence $x_1,x_2,\ldots$ contains a convergent subsequence. Prove that if a sequence of points in a metric space does not contain any convergent subsequence, then for each $x\in X$ there exists $r>0$ such that the open ball $U_r(x)$ does not contain points of this sequence other than $x$.
\end{prb}

\section{Closure, interior, boundary, dense subsets, separability}
\markright{\thesection.~Closure, interior, boundary, everywhere dense subsets, separability}

Let $Y$ be a subset of a topological space $X$. A point $x\in X$ is called \emph{an adherent point}, or \emph{a closure point}, or \emph{a contact point for $Y$} if every neighborhood of $x$ intersects $Y$. The set of all adherent points of the set $Y$ is called its \emph{closure\/} and is denoted by $\bY$.

\begin{prb}
Prove that the closure of a set $Y\ss X$ is the smallest closed subset of $X$ containing $Y$, i.e., $\bY$ is the intersection of all closed sets containing $Y$.
\end{prb}

A subset of $Y$ of a topological space $X$ is called \emph{everywhere dense in $X$} if $\bY=X$.

\begin{examp}
The set of all rational numbers, like the set of all irrational numbers, are everywhere dense in the real line.
\end{examp}

\begin{prb}
Let the topology of Zariski of weight $m$ be given on an infinite set $X$. Then a subset $Y\ss X$ is everywhere dense in $X$ if and only if $\#Y\ge m$.
\end{prb}

A topological space is called \emph{separable\/} if it contains an everywhere dense sequence.

\begin{examp}
Each finite space is separable. The Euclidean space $\R^n$ is also separable: as an everywhere dense sequence we can take arbitrary numbered set of all points with rational coordinates. Each space with a countable Zariski topology is separable. Sorgenfrey space (Example~\ref{examp:Sorgenfrey-space}) is separable.
\end{examp}

\begin{prb}\label{prb:SeparabilityAndCountableBase}
Show that in metric space, separability is equivalent to having a countable base. Extract from this that every subset of a separable metric space is separable. Show that an open subset of an arbitrary separable topological space is separable. Give an example of a separable topological space containing an non-separable subset (use the Sorgenfrey plane).
\end{prb}

A point $x$ from a subset $Y$ of a topological space $X$ is called \emph{interior for $Y$} if some neighborhood of $x$ is contained in $Y$. The family of all interior points of the set $Y$ is called its \emph{interior\/} and is denoted by $\Int Y$.

\begin{prb}
Show that the interior $\Int Y$ is the largest open subset of $X$ contained in $Y$.
\end{prb}

\begin{prb}
Prove that a subset $Y$ of the topological space $X$ is closed if and only if $Y=\bY$, and is open if and only if $Y=\Int Y$.
\end{prb}

A point $x\in X$ of a topological space $X$ is called \emph{a boundary point for a subset $Y\ss X$} if each neighborhood of $x$ intersects both $Y$ and its complement $X\sm Y$. The set of all boundary points of the set $Y$ is called its \emph{boundary\/} and is denoted by $\d Y$.

\begin{prb}
Prove that the boundary $\d Y$ is a closed subset of $X$, and
$$
\d Y=\bY\sm\Int Y=\bY\cap\overline{\strut X\sm Y}.
$$
\end{prb}

\section{Separated spaces}
\markright{\thesection.~Separated spaces}

There are a number of separation axioms that generate various classes of topological spaces. We will not dwell on this in detail here, but formulate only one axiom of separation, which will be useful to us in the future.

A topology on a set $X$, as well as the topological space $X$ itself, is called \emph{Hausdorff\/} or \emph{separated\/} if any two points of $X$ have disjoint neighborhoods.

\begin{examp}
Each discrete topology is Hausdorff. Each metric space is Hausdorff. If the set $X$ consists of more than one point, then the anti-discrete topology is not Hausdorff. Also, the Zariski topology is not a Hausdorff topology (see Example~\ref{examp:Zaris}).
\end{examp}

\begin{prb}
Show that in a Hausdorff topological space every point is closed. Give an example of a non-Hausdorff topological space in which all points are closed.
\end{prb}

\begin{prb}
Show that the disjoint union and the Cartesian product of Hausdorff topological spaces are also Hausdorff.
\end{prb}

\begin{prb}
Show that in a Hausdorff topological space the limit of a convergent sequence is uniquely defined. Give an example of a topological space in which each sequence converges to each point.
\end{prb}

\begin{prb}
Describe what sequences in a space with Zariski topology are convergent, and what limits each convergent sequence has.
\end{prb}

\section{Connected spaces}
\markright{\thesection.~Connected spaces}

We say that a set $X$ \emph{is partitioned into subsets $\{X_i\}_{i\in I}$} if $X=\sqcup_{i\in I}X_i$.

A topological space (its topology) is called \emph{disconnected} if it can be partitioned into two nonempty open (equivalently, closed) sets. If such a partition does not exist, then the topological space is called \emph{connected}. A subset of a topological space is \emph{connected\/} (\emph{disconnected\/}), if such is the topology induced on it. In other words, a subset $Y$ of a topological space $X$ is \emph{disconnected\/} if there exist $U,V\in\t_X$ such that $Y\ss U\cup V$, and both intersections $Y\cap U$ and $Y\cap V$ are nonempty and do not intersect each other.

\begin{prb}
Prove that each segment $[a,b]\ss\R$ is connected.
\end{prb}

\begin{prb}
Prove that the closure of a connected subset of a topological space is connected.
\end{prb}

\begin{prb}
Let $\{A_i\}_{i\in I}$ be a family of connected pairwise intersecting subsets of a topological space $X$, then the set $\cup_{i\in I}A_i$ is connected.
\end{prb}

\begin{prb}
Show that the image of a connected topological space under a continuous mapping is also connected.
\end{prb}

\begin{prb}
Prove that every continuous function on a connected topological space takes all intermediate values.
\end{prb}

The maximum (by inclusion) connected  subset of a topological space is called a \emph{connected component\/} of this space.

\begin{prb}
Show that each connected component is closed, and that each topological space is uniquely partitioned into its connected components. If such a partition is finite, then connected components are also open. Give an example of a topological space in which some connected components are not open.
\end{prb}

\section{Path-connected spaces}
\markright{\thesection.~Path-connected spaces}

\emph{A curve\/} in a topological space $X$ is an arbitrary continuous map $\g\:[a,b]\to X$. It is said that the curve $\g$ \emph{joins the points $\g(a)$ and $\g(b)$}.

A topological space $X$ is called \emph{path-connected\/} if any two of its points can be connected by a curve.

\begin{prb}
Prove that a path-connected topological space is connected. Give an example of a connected space that is not path-connected.
\end{prb}

\section{Compact and sequentially compact spaces}
\markright{\thesection.~Compact and sequentially compact spaces}

\emph{A subcover of a cover\/} is a subfamily of a cover, which itself is a cover. A cover of a topological space composed of open sets is called \emph{open}.

\begin{dfn}\label{dfn:compact}
A topological space $X$ is called \emph{compact\/} if a finite subcover can be found in any of its open covers.
\end{dfn}

\begin{rk}
To define a cover for a subset $Y$ of a topological space $X$, it is more convenient to modify Definition~\ref{dfn:compact}, rather than reduce it to the corresponding concept for the induced topology. Namely, an open cover of $Y$ is a family $\{U_a\}_{a\in A}$ of open subsets of $X$ such that $Y\ss\cup_{a\in A}U_a$. Other definitions do not change.
\end{rk}

\begin{prb}
Show that a finite union of compact subsets of a topological space is compact.
\end{prb}

\begin{prb}
Prove the following statements:
\begin{enumerate}
\item the image under a continuous mapping from a compact topological space is compact;
\item a closed subset of a compact topological space is compact;
\item a compact subset of a Hausdorff topological space is closed;
\item a continuous bijective mapping from a compact topological space to a Hausdorff space is a homeomorphism;
\item give an example of an infinite topological space in which all subsets are compact. Note that in such a space there are compact subsets that are not closed;
\item give an example of a continuous bijective mapping of topological spaces that is not a homeomorphism.
\end{enumerate}
\end{prb}

\begin{prb}[Alexander subbase theorem]
Let $X$ be a topological space and $\s$ its subbase. Prove that $X$ is compact if and only if each cover of $X$ by elements of the subbase $\s$ has a finite subcover.
\end{prb}

\begin{prb}[Tychonoff's theorem]
Prove that the Cartesian product $\prod_{i\in I}X_i$ of topological spaces $X_i$, endowed with Tychonoff topology, is compact if and only if all $X_i$ are compact.
\end{prb}

\begin{prb}
Prove that each segment $[a,b]\ss\R$ is compact.
\end{prb}

A subset of a metric space is called \emph{bounded\/} if it is contained in some ball.

\begin{prb}
Prove that a subset of a Euclidean space is compact if and only if it is closed and bounded.
\end{prb}

\begin{prb}
Prove that every compact metric space is bounded. Prove that a continuous function on a compact topological space is bounded and takes its largest and smallest values.
\end{prb}

\begin{dfn}
A topological space is called \emph{sequentially compact\/} if every sequence of its points has a convergent subsequence.
\end{dfn}

\begin{prb}\label{prb:MinMaxSequentialCompact}
Prove that every sequentially compact metric space is bounded. Prove that a continuous function on a sequentially compact topological space is bounded and takes its largest and smallest values.
\end{prb}

\begin{rk}
Note that compactness and sequential compactness in the case of general topological spaces are not related to each other: there are spaces that have one of these properties and do not have the other. Since the examples are quite complicated, we will omit them here. However, everything is much simpler for metric spaces, see Theorem~\ref{thm:ComcactSequentialCompact} in Chapter~\ref{ch:TopMetr}.
\end{rk}

\section{Hyperspaces}\label{sec:hyperspaces}
\markright{\thesection.~Hyperspaces}
A family of various subsets of a topological or metric space endowed with a certain topology or metric is called \emph{a hyperspace}. In Construction~\ref{constr:Vietoris_top} above, we defined the Vietoris topology on the set $2^X$ of all subsets of a topological space $X$. Thus, $2^X$ is a special case of hyperspace.

Here are a few more examples of hyperspaces (in all these spaces, the topology is induced from $2^X$ in the standard way):
\begin{itemize}
\item by $\cP_0(X)\ss2^X$ we denote the set of all nonempty subsets of $X$;
\item by $\cCL(X)\ss\cP_0(X)$ we denote the set of all nonempty closed subsets of $X$;
\item by $\cC(X)\ss\cCL(X)$ the set of all nonempty closed connected subsets of $X$;
\item by $\cC_n(X)\ss\cCL(X)$ the set of all nonempty closed subsets of $X$ having at most $n$ components;
\item by $\cC_\infty(X)\ss\cCL(X)$ the set of all nonempty closed subsets of $X$, each of which has finitely many components;
\item by $\cK(X)\ss\cP_0(X)$ the set of all nonempty compact subsets of $X$;
\item by $\cF_n(X)\ss\cK(X)$ the set of all nonempty at most $n$-point subsets of $X$;
\item by $\cF_\infty(X)\ss\cK(X)$ the set of all nonempty finite subsets of $X$.
\end{itemize}

There are numerous, usually obvious, connections between these spaces. For example, if the space $X$ is Hausdorff, then $\cK(X)\ss\cCL(X)$; if $X$ is Hausdorff and compact, then $\cK(X)=\cCL(X)$.

A connected nonempty compact Hausdorff topological space is called \emph{a continuum}. If the space $X$ is Hausdorff, then $\cK(X)\cap\cC(X)$ is the set of all continua. In some literature, for example, in~\cite{Nadler1}, this space is denoted by $\cC(X)$.

\begin{prb}
Let $X=\{a,b\}$. We define the following topology on $X$ as follows: $\t=\bigl\{\0,X,\{a\}\bigr\}$. Find out what the space $\CL(X)$ is.
\end{prb}

\begin{dfn}
A topological space is called \emph{a space of class $T_0$} if, for any two different points of this space, at least one of them has a neighborhood that does not contain the second point.
\end{dfn}

\begin{prb}
Prove that the space $\CL(X)$ is always a space of class $T_0$.
\end{prb}

\begin{prb}
A topological space is called \emph{a space of class $T_1$} if, for any two different points of this space, each of them has a neighborhood that does not contain the remaining point.
\end{prb}

\begin{prb}
Prove that if $X$ is a space of class $T_1$, then $\CL(X)$ is also a space of class $T_1$. Give an example that demonstrates that the converse is not true.
\end{prb}

\begin{prb}
Prove that the space $\cP_0(X)$ belongs to the class $T_1$ if and only if the space $X$ is discrete.
\end{prb}


\phantomsection
\renewcommand\bibname{References to Chapter~\thechapter}
\addcontentsline{toc}{section}{\bibname}
\renewcommand{\refname}{\bibname}

\vfill\eject

\section*{\Huge Exercises to Chapter~\thechapter}
\markright{Exercises to Chapter~\thechapter.}

\begin{exe}
Show that the family of all closed subsets of a topological space $X$ contains $\0$ and $X$, and that the intersection of any collection of closed subsets, as well as the union of any finite collection of closed subsets are some closed sets.
\end{exe}

\begin{exe}
Let $X$ be a topological space and $Y\ss X$. Consider the family $\t_Y:=\{U\cap Y: U\in\t_X\}$. Prove that $\t_Y$ is a topology on $Y$.
\end{exe}

\begin{exe}
For a metric space $(X,\r)$ define $\t_\r\ss2^X$ as the collection consisting of the empty set and all possible unions of open balls. Prove that the family $\t_\r$ is a topology.
\end{exe}

\begin{exe}\label{exe:criterion-base}
Prove that a family $\b\ss2^X$ is a base of some topology $\t$ on $X$ if and only if $\b$ is a cover of $X$, and for any intersecting $B_1,B_2\in\b$ their intersection $B_1\cap B_2$ is the union of some elements from $\b$. Moreover, each family satisfying these properties, generates a unique topology.
\end{exe}

\begin{exe}
Let some topologies $\t_1$ and $\t_2$ with bases $\b_1$ and $\b_2$ be given on a set $X$. Then $\t_1=\t_2$ if and only iff for any $x\in X$ the following condition is fulfilled: for any $B_2\in\b_2$, $x\in B_2$, there is $B_1\in\b_1$ for which $x\in B_1\ss B_2$, and vice versa, for any $B_1\in\b_1$, $x\in B_1$, there exists $B_2\in\b_2$ for which $x\in B_2\ss B_1$. In particular, for a topological space $X$, a collection $\b$ of open sets satisfying the condition of Exercise~\ref{exe:criterion-base} is a base of the topology $\t_X$ if and only if for each open set $U\in\t_X$ and any point $x\in U$ there exists some $B\in\b$ such that $x\in B\ss U$.
\end{exe}

\begin{exe}
Let $X$ be an infinite set of cardinality $n$, and $m$ be an infinite cardinal number, with $m\le n$. Consider a family $\cF$ of all $F\ss X$ such that $\#F<m$, and let $\b_m=\{X\sm F:F\in\cF\}$. Prove that the family $\b_m$ is a base of some topology.
\end{exe}

\begin{exe}
Show that the collection of all possible half-intervals of the form $[a,b)\ss\R$ form a base of some topology that contains the standard topology.
\end{exe}

\begin{exe}
Prove that a family $\s\ss2^X$ is a subbase of some topology on $X$ if and only if $\s$ is a cover of $X$. Moreover, each cover of $X$ generates a unique topology.
\end{exe}

\begin{exe}
Show that the standard topology of the Euclidean space $\R^n$ coincides with the topology of the Cartesian product $\R\x\cdots\x\R$ of real lines endowed with the standard topology.
\end{exe}

\begin{exe}
Let $X$ be an arbitrary topological space. For each finite collection $U_1,\ldots,U_n\in\t_X$ we put
$$
\<U_1,\ldots,U_n\>=\{Y\ss X:Y\ss\cup_{i=1}^nU_i,\ \text{and $Y\cap U_i\ne\0$ for all $i=1,\ldots,n$}\}.
$$
Show that the families
$$
\s=\bigl\{\<U\>:U\in\t_X\bigr\}\cup\bigl\{\<X,U\>:U\in\t_X\bigr\}\ \ \text{and}\ \
\b=\bigl\{\<U_1,\ldots,U_n\>:U_1,\ldots,U_n\in\t_X\bigr\}
$$
form respectively a subbase and the corresponding base of some topology on $2^X$.
\end{exe}

\begin{exe}
Prove that each construction from the section ``Standard constructions of topologies'' provides a topology.
\end{exe}

\begin{exe}
Prove that the definitions~\ref{dfn:cont1}, \ref{dfn:cont2}, and \ref{dfn:cont3} are equivalent.
\end{exe}

\begin{exe}
Let $f\:X\to Y$ be a mapping of topological spaces and $\s$ be a subbase of the topology on the space $Y$. Prove that $f$ is continuous if and only if $f$-preimage of each element from the subbase $\s$ is open in $X$.
\end{exe}

\begin{exe}
Show that the identity mapping and the composition of continuous mappings are continuous.
\end{exe}

\begin{exe}
Let $X$ be a topological space, and $Z\ss X$ be its subspace. Show that the inclusion mapping $i\:Z\to X$, $i(z)=z$ for each point $z\in Z$, is continuous.
\end{exe}

\begin{exe}
Let $X$, $Y$ be topological spaces, $W\ss Y$ be a subspace of $Y$, and $f\:X\to W$ be a continuous mapping. Let $g\:X\to Y$ be a mapping coinciding with $f$: for each $x\in X$ it holds $f(x)=g(x)$. Prove that the mapping $g$ is continuous.
\end{exe}

\begin{exe}
Let $f\:X\to Y$ be a continuous mapping of topological spaces, $Z\ss X$, $W\ss Y$, $f(Z)\ss W$. Then the restriction $f|_{Z,W}\:Z\to W$ is also continuous as the mapping of the topological spaces $Z$ and $W$ with the topologies induced on them from $X$ and $Y$, respectively.
\end{exe}

\begin{exe}
Let $\{X_i\}_{i\in I}$ be a cover of a topological space $X$ by open subsets $X_i$, and $f\:X\to Y$ a mapping to a topological space $Y$. Show that $f$ is continuous if and only if all the restrictions $f|_{X_i}$ are continuous. In particular, this holds when $X=\sqcup_{i\in I}X_i$ is the disjoint union of some topological spaces. Will this result remain true if we replace $\{X_i\}$ with a cover of $X$ by arbitrary sets?
\end{exe}

\begin{exe}
Let $X=\sqcup_{i\in I}X_i$ be the disjoint union of some topological spaces and $f\:X\to Y$ be a map into the topological space $Y$. Show that $f$ is continuous if and only if all its restrictions $f|_{X_i}$ are continuous.
\end{exe}

\begin{exe}
Let $\{X_i\}_{i\in I}$ be a family of topological spaces and $X=\prod_{i\in I}X_i$. We define the canonical projection $\pi_i\:X\to X_i$ by setting $\pi_i(x)=x_i$. Prove that the product topology on $X$ is the smallest of those topologies in which all the projections $\pi_i$ are continuous.
\end{exe}

\begin{exe}
Let $\{Y_i\}_{i\in I}$ be a family of topological spaces, $Y=\prod_{i\in I}Y_i$, and $f_i\:X\to Y_i$ are mappings from some topological space $X$. We construct the mapping $F:=\prod_{i\in I}f_i\:X\to Y$  by associating with each point $x\in X$ the element $y\in Y$ defined as follows: $y_i=f_i(x)$. Prove that the mapping $F$ is continuous if and only if all $f_i$ are continuous.
\end{exe}

\begin{exe}
Let $A\ss\R^n$ be an arbitrary subset, $(x^1,\ldots,x^n)$ the Cartesian coordinates on $\R^n$, $f\:A\to\R^m$ a continuous mapping, $(y^1,\ldots,y^m)$ the Cartesian coordinates on $\R^m$, and $y^i=y^i(x^1,\ldots,x^n)$ the coordinate functions of the mapping $f$. Prove that the mapping $f$ is continuous if and only if all the coordinate functions $y^i=y^i(x^1,\ldots,x^n)$ are continuous.
\end{exe}

\begin{exe}
Describe all continuous functions on a topological space with Zariski topology.
\end{exe}

\begin{exe}
Let $f\:X\to Y$ be a homeomorphism, and $Z\ss X$, $W=f(Z)$. Prove that the restriction $f|_{Z,W}\:Z\to W$ is also a homeomorphism. Show that the letters \textsf{b, c, f, g, i, h, o} are pairwise non-homeomorphic.
\end{exe}

\begin{exe}
Show that every embedding is continuous. Give an example of a continuous injective mapping of topological spaces that is not an embedding.
\end{exe}

\begin{exe}
Let $\om$ be a character not contained in $\N$. We define a topology on the set $\bar{N}=\N\cup\{\om\}$, taking as a base all points from $\N$, as well as all sets of the form $\{n\ge N\}\cup\{\om\}$, $N\in\N$. Show that a sequence $x\:\N\to X$ converges if and only if the mapping $x$ can be extended to a continuous mapping on $\bar{\N}$.
\end{exe}

\begin{exe}
Show that a continuous mapping $f\:X\to Y$ of topological spaces takes convergent sequences to convergent ones. Show that if $X$ is  a metric space, then every mapping $g\:X\to Y$ that takes convergent sequences into convergent ones is continuous. Give an example of a topological space $X$ and a mapping $h\:X\to Y$ into a topological space $Y$, which takes convergent sequences into convergent ones, but is not continuous nonetheless.
\end{exe}

\begin{exe}
Let $x_1,x_2,\ldots$ be a sequence of points in a metric space $X$. Suppose that for some point $x\in X$ each neighborhood of $x$ intersects the set $\{x_i\}_{i=1}^\infty\sm\{x\}$. Prove that the sequence $x_1,x_2,\ldots$ contains a convergent subsequence. Extract from this that if a sequence of points in a metric space does not contain any convergent subsequence, then for each $x\in X$ there exists $r>0$ such that the open ball $U_r(x)$ does not contain points of this sequence other than $x$.
\end{exe}

\begin{exe}
Prove that the closure of a set $Y\ss X$ is the smallest closed subset of $X$ containing $Y$, i.e., $\bY$ is the intersection of all closed sets containing $Y$.
\end{exe}

\begin{exe}
Let the topology of Zariski of weight $m$ be given on an infinite set $X$. Then a subset $Y\ss X$ is everywhere dense in $X$ if and only if $\#Y\ge m$.
\end{exe}

\begin{exe}
Show that in metric space, separability is equivalent to having a countable base. Extract from this that every subset of a separable metric space is separable. Show that an open subset of an arbitrary separable topological space is separable. Give an example of a separable topological space containing an non-separable subset (use the Sorgenfrey plane).
\end{exe}

\begin{exe}
Show that the interior $\Int Y$ is the largest open subset of $X$ contained in $Y$.
\end{exe}

\begin{exe}
Prove that a subset $Y$ of the topological space $X$ is closed if and only if $Y=\bY$, and is open if and only if $Y=\Int Y$.
\end{exe}

\begin{exe}
Prove that the boundary $\d Y$ is a closed subset of $X$, and
$$
\d Y=\bY\sm\Int Y=\bY\cap\overline{\strut X\sm Y}.
$$
\end{exe}

\begin{exe}
Show that in a Hausdorff topological space every point is closed. Give an example of a non-Hausdorff topological space in which all points are closed.
\end{exe}

\begin{exe}
Show that the disjoint union and the Cartesian product of Hausdorff topological spaces are also Hausdorff.
\end{exe}

\begin{exe}
Show that in a Hausdorff topological space the limit of a convergent sequence is uniquely defined. Give an example of topological space in which each sequence converges to each point.
\end{exe}

\begin{exe}
Describe what sequences in a space with Zariski topology are convergent, and what limits each convergent sequence has.
\end{exe}

\begin{exe}
Prove that each segment $[a,b]\ss\R$ is connected.
\end{exe}

\begin{exe}
Prove that the closure of a connected subset of a topological space is connected.
\end{exe}

\begin{exe}
Let $\{A_i\}_{i\in I}$ be a family of connected pairwise intersecting subsets of a topological space $X$, then the set $\cup_{i\in I}A_i$ is connected.
\end{exe}

\begin{exe}
Show that the image of a connected topological space under a continuous mapping is also connected.
\end{exe}

\begin{exe}
Prove that every continuous function on a connected topological space takes all intermediate values.
\end{exe}

\begin{exe}
Show that each connected component is closed, and that each topological space is uniquely partitioned into its connected components. If such a partition is finite, then connected components are also open. Give an example of a topological space in which some connected components are not open.
\end{exe}

\begin{exe}
Prove that a path-connected topological space is connected. Give an example of a connected space that is not path-connected.
\end{exe}

\begin{exe}
Show that a finite union of compact subsets of a topological space is compact.
\end{exe}

\begin{exe}
Prove the following statements:
\begin{enumerate}
\item the image under a continuous mapping from a compact topological space is compact;
\item a closed subset of a compact topological space is compact;
\item a compact subset of a Hausdorff topological space is closed;
\item a continuous bijective mapping from a compact topological space to a Hausdorff space is a homeomorphism;
\item give an example of an infinite topological space in which all subsets are compact. Note that in such a space there are compact subsets that are not closed;
\item give an example of a continuous bijective mapping of topological spaces that is not a homeomorphism.
\end{enumerate}
\end{exe}

\begin{exe}[Alexander subbase theorem]
Let $X$ be a topological space and $\s$ its subbase. Prove that $X$ is compact if and only if each cover of $X$ by elements of the subbase $\s$ has a finite subcover.
\end{exe}

\begin{exe}[Tychonoff's theorem]
Prove that the Cartesian product $\prod_{i\in I}X_i$ of topological spaces $X_i$, endowed with Tychonoff topology, is compact if and only if all $X_i$ are compact.
\end{exe}

\begin{exe}
Prove that each segment $[a,b]\ss\R$ is compact.
\end{exe}

\begin{exe}
Prove that a subset of a Euclidean space is compact if and only if it is closed and bounded.
\end{exe}

\begin{exe}
Prove that every compact metric space is bounded. Prove that a continuous function on a compact topological space is bounded and takes its largest and smallest values.
\end{exe}

\begin{exe}
Prove that every sequentially compact metric space is bounded. Prove that a continuous function on a sequentially compact topological space is bounded and takes its largest and smallest values.
\end{exe}

\begin{notation}
The following matrix groups consist of real matrices of size $n\x n$ and are considered as subsets of $\R^{n^2}$ with the induced topology (their rows or columns are written out one after another and form vectors): $\O(n)$ consists of all orthogonal matrices (\emph{orthogonal group\/}); $\SO(n)$ consists of all orthogonal matrices with determinant $1$ (\emph{special orthogonal group\/}); $\GL(n)$ consists of all nondegenerate matrices (\emph{general linear group\/}); $\SL(n)$ consists of all matrices with determinant $1$ (\emph{special linear group\/}).
\end{notation}

\begin{exe}
Find out which of the following matrix groups are connected, which are compact:
$$
\O (n),\ \SO (n),\ \GL (n),\ \SL (n).
$$
\end{exe}

\begin{exe}
Let $X=\{a,b\}$. We define the following topology on $X$: $\t=\bigl\{\0,X,\{a\}\bigr\}$. Find out what the space $\CL(X)$ is.
\end{exe}

\begin{dfn}
A topological space is called \emph{a space of class $T_0$} if, for any two different points of this space, at least one of them has a neighborhood that does not contain the second point.
\end{dfn}

\begin{exe}
Prove that the space $\CL(X)$ is always a space of class $T_0$.
\end{exe}

\begin{dfn}
A topological space is called \emph{a space of class $T_1$} if, for any two different points of this space, each of them has a neighborhood that does not contain the remaining point.
\end{dfn}

\begin{exe}
Prove that if $X$ is a space of class $T_1$, then $\CL(X)$ is also a space of class $T_1$. Give an example that demonstrates that the converse is not true.
\end{exe}

\begin{exe}
Prove that the space $\cP_0(X)$ belongs to the class $T_1$ if and only if the space $X$ is discrete.
\end{exe}

\vfill\eject
 \chapter{Introduction to Metric Spaces.}\label{ch:TopMetr}
 \markboth{\chaptername~\thechapter.~Introduction to Metric Spaces.}%
          {\chaptername~\thechapter.~Introduction to Metric Spaces.}

\begin{plan}
Distance, pseudometric or semimetric, metric, open ball, closed ball, sphere, distance to a nonempty subset, open and closed $r$-neighborhoods of a nonempty subset, diameter of a subset, natural quotient of pseudometric space, Lipschitz mappings, Lipschitz constant, dilatation of a Lipschitz mapping, uniform continuity, bi-Lipschitz mappings, isometric mapping, isometry, isometry group, standard constructions of metrics, multiplying a metric by a number, adding to metric a constant, induced distance, semidirect product of metric spaces, examples, Levenshtein distance, elements of graph theory, metric construction for graphs, distance on a connected graph, distance on a connected weighted graph, Cayley graph of a group, quotient pseudometric and metric, generalized distance, disjoint union of generalized spaces, maximal pseudometric, the relation between quotient and maximal pseudometrics and metrics, isometries groups actions and quotient spaces, metrized graphs, polyhedron spaces, convergence of sequences and completeness, fundamental sequences, completion of a metric space, equivalence of compactness and sequential compactness for metric spaces, completeness and total boundedness equivalent to compactness for metric spaces, canonical isometric embeddings of metric spaces, Frechet-Kuratowski embedding to the space of bounded continuous functions, Frechet embedding of a separable metric space to the space of bounded sequences.
\end{plan}

In this section, we discuss some basic facts from metric spaces theory.

\section{Distance}
\markright{\thesection.~Distance}

Let $X$ be an arbitrary set. A function $\r\:X\x X\to\R$ is called \emph{a distance}, if it is non-negative, it is equal to zero on all pairs of the form $(x,x)$, and it is \emph{symmetric\/}: $\r(x,y)=\r(y,x)$ for any $x,y\in X$.

If the distance function $\r$ satisfies \emph{the triangle inequality}, namely, $\r(x, y)+\r(y,z)\ge\r(x,z)$, then such $\r$ is called \emph{a pseudometric\/} or \emph{a semimetric}, and the set $X$ with the introduced pseudometric on it is \emph{a pseudometric space\/} or \emph{a semimetric space}.

Thus, the metric defined above is a pseudometric not equal to zero on the pairs $(x,y)$ with $x\ne y$.

\begin{rk}
It will be convenient for us to introduce a universal notation for the distance function defined on an arbitrary set, namely, the distance between points $x$ and $y$ from this set will be denoted by $|xy|$. Even if several sets are considered at the same time with the distances given on them, we will denote these distances in the same way, while understanding what distance is used according to which set the corresponding points belong to.
\end{rk}

For each distance given on the set $X$, a number of subsets arise that play an important role in the study of geometry:
\begin{itemize}
\item\emph{an open ball of radius $r>0$ and center $x\in X$}: $U_r(x):=\bigl\{y\in X:|xy|<r\bigr\}$ (above we defined an open ball for a metric);
\item\emph{a closed ball of radius $r\ge0$ and center $x\in X$}: $B_r(x):=\bigl\{y\in X:|xy|\le r\bigr\}$;
\item\emph{a sphere of radius $r\ge0$ and center $x\in X$}: $S_r(x):=\bigl\{y\in X:|xy|=r\bigr\}$;
\end{itemize}
Note that in metric space, for $r=0$, both a closed ball and a sphere of radius $r$ degenerate to a point and are called \emph{degenerate}.

For any point $x\in X$ and a nonempty set $A$, define \emph{the distance from $x$ to $A$} by setting
$$
|xA|=|Ax|=\inf\bigl\{|xa|:a\in A\bigr\}.
$$
This concept gives rise to a number of objects:
\begin{itemize}
\item for $r>0$ \emph{an open $r$-neighborhood of the set $A$} is $U_r(A)=\bigl\{y\in X:|Ay|<r\bigr\}$;
\item for $r\ge0$ \emph{a closed $r$-neighborhood of the set $A$} is $B_r(A)=\bigl\{y\in X:|Ay|\le r\bigr\}$.
\end{itemize}

In addition, for a nonempty $A\ss X$, a numerical characteristic of the set $A$ is defined --- its \emph{diameter}
$$
\diam A=\sup\bigl\{|aa'|:a,a'\in A\bigr\}.
$$
Moreover, the diameter can also be naturally determined for the empty set by setting $\diam\0=0$.

\begin{rk}
Sometimes an object $A$ under consideration can simultaneously belong to different metric spaces, for example, in one space it is a subset, and in another it is a point. Then, in the notation introduced above, $U_r(A)$, $B_r(A)$, etc., as an upper index, we will add either the name of the space to which $A$ belongs or the metric of this space. A typical case: $A$ can be considered both as a nonempty subset of the space $X$ with the metric $d$, and as a point in the space $\cH(X)$ of all nonempty closed bounded subsets of $X$ with some metric $d_H$. In this situation, instead of $U_r(A)$ we will write $U_r^X (A)$ or $U_r^d (A)$ in the first case, and $U_r^{\cH (X)}(A)$ or $U_r^{d_H}(A)$ in the second one.
\end{rk}

Each pseudometric on the set $X$ defines a natural equivalence relation: $x\sim y$ if and only if $|xy|=0$. Let $X/\!\!\sim$ be the set of classes of this equivalence, and for each $x\in X$ denote by $[x]$ the equivalence class containing $x$.

\begin{prb}\label{prb:pseudoTOmetric}
Prove that for any $x,y\in X$ and $x'\in[x]$ and $y'\in[y]$ it is true that $|x'y'|=|xy|$. Thus, on the set $X/\!\!\sim$ the corresponding distance function is correctly defined: $\bigl|[x][y]\bigr|=|xy|$. Show that this distance function is a metric.
\end{prb}

\begin{rk}
If $\r$ is the pseudometric on $X$, then the quotient space $X/\!\!\sim$ from Problem~\ref{prb:pseudoTOmetric} is sometimes denoted by $X/\r$.
\end{rk}

\begin{prb}\label{prb:metric-simple-props}
Let $X$ be an arbitrary metric space, $x,y\in X$, $r\ge0$, $s,t>0$, and $A\ss X$ be nonempty. Verify that
\begin{enumerate}
\item\label{prb:metric-simple-props:01} $U_s\bigl(\{x\}\bigr)=U_s\bigl(x)$ and $B_r\bigl(\{x\}\bigr)=B_r\bigl (x)$;
\item\label{prb:metric-simple-props:02} the functions $y\mapsto|xy|$, $y\mapsto|yA|$ are continuous;
\item\label{prb:metric-simple-props:03} an open neighborhood $U_s(A)$ is an open subset of $X$, and a closed neighborhood $B_r(A)$ is a closed subset of $X$;
\item\label{prb:metric-simple-props:1} $U_t\bigl(U_s(A)\bigr)\ss U_{s+t}(A)$ and construct an example demonstrating that the left-hand side can be different from the right-hand side;
\item\label{prb:metric-simple-props:2} $B_t\bigl(B_s(A)\bigr)\ss B_{s+t}(A)$ and construct an example demonstrating that the left-hand side can be different from the right-hand side;
\item\label{prb:metric-simple-props:3} $\d U_s(x)$, $\d B_s(x)$ are not related by any inclusion; $\d U_s(x)\ss S_s(x)$ and $\d B_r(x)\ss S_r(x)$; the both previous inclusions can be strict;
\item\label{prb:metric-simple-props:4} $\diam U_s(x)\le\diam B_s(x)\le2s$;
\item\label{prb:metric-simple-props:5} $\diam U_s(A)\le\diam B_s(A)\le\diam A+2s$.
\end{enumerate}
\end{prb}

\section{Lipschitz mappings and isometries}\label{sec:LipschitzMappins}
\markright{\thesection.~Lipschitz mappings and isometries}

Mappings of metric spaces that distort distances no more than a certain finite number of times are called Lipschitz. More formally, a mapping $f\:X\to Y$ of metric spaces is called \emph{Lipschitz}, if there exists $C\ge0$ such that for any $x,x'\in X$ the inequality $\bigl|f(x)f(x')\bigr|\le C|xx'|$ holds. Each such $C$ is called \emph{a Lipschitz constant}. Sometimes, for brevity, a Lipschitz mapping with Lipschitz constant $C$ is called \emph{$C$-Lipschitz}. For $1$-Lipschitz mapping is reserved the term \emph{nonexpanding}.

\begin{prb}
Let $\cL(f)\ss\R$ be the set of all Lipschitz constants for a Lipschitz mapping $f$. Prove that $\inf\cL(f)$ is also a Lipschitz constant.
\end{prb}

The smallest Lipschitz constant for a Lipschitz mapping $f$ is called \emph{the dilatation of the mapping $f$} and is denoted by $\dil f$.

A mapping $f\:X\to Y$ of metric spaces is called \emph{uniformly continuous\/} if for any $\e>0$ there exists $\dl>0$ such that for any points $x,x'\in X$, $|xx'|<\dl$, we have $|f(x)f(x')|<\e$.

\begin{prb}
Show that each Lipschitz mapping is uniformly continuous, and each uniformly continuous mapping is continuous.
\end{prb}

A bijective mapping $f$ between metric spaces such that $f$ and $f^{- 1}$ are Lipschitz is called \emph{bi-Lipschitz}. It is clear that each bi-Lipschitz mapping is a homeomorphism.

A mapping of metric spaces $f\:X\to Y$ is called \emph{isometric\/} if it preserves distances: $\bigl|f(x)f(x')\bigr|=|xx'|$ for any $x,x'\in X$. A bijective isometric mapping is called \emph{an isometry}.

\begin{prb}\label{prb:IsomHomeo}
Show that each isometry is a homeomorphism, in particular, each isometric mapping from one metric space into another one is an embedding, i.e., we recall, it is a homeomorphism with its image.
\end{prb}

\begin{prb}
Verify that the identity map, the composition of isometries, and the inverse mapping to an isometry are also isometries, i.e., the set of all isometries of an arbitrary metric space forms a group.
\end{prb}

\begin{dfn}
The set of all isometries $f\:X\to X$ of the metric space $X$ endowed with the composition operation is called \emph{the isometry group of the space $X$} and is denoted by $\Iso(X)$.
\end{dfn}


\begin{prb}
Let $X$ be an arbitrary metric space, $x,y\in X$ and $A\ss X$ be nonempty. Prove that $|Ax|+|xy|\ge|Ay|$, so that the function $\r_A(x)=|Ax|$ is $1$-Lipschitz and, therefore, uniformly continuous.
\end{prb}

\section{Standard constructions of metrics}\label{sec:StandCOnstrMetric}
\markright{\thesection.~Standard constructions of metrics}

This section provides examples of standard constructions that build various metric spaces.

\begin{constr}[Multiplying a metric by a number]
If $d$ is a distance function on the set $X$, then for each real $\l>0$ the function $\l\,d$ is also a distance on $X$; moreover, if $d$ was a metric (pseudometric), then $\l\,d$ is the distance of the same type. The corresponding space will be denoted by $\l\,X$.
\end{constr}

\begin{constr}[Adding a constant]
Let $d$ be the distance function on $X$. We define an analogue of the Kronecker symbols for $x,y\in X$ by setting $\dl_{xy}=0$ for any $x\ne y$, and $\dl_{xx}=1$ for any $x$. Then for every real $c\ge0$ the function $(d+c)(x,y)=d(x,y)+c\,(1-\dl_{xy})$ is a distance. Moreover, if $d$ is a pseudometric (metric) and $c>0$, then $d+c$ is a metric. Note that this construction can also be extended to some negative numbers $c$.
\end{constr}

\begin{prb}\label{prb:metric-add-const}
Let $d$ be a metric. Find the least possible $c$ such that $d+c$ is a pseudometric. Verify that for such $c$ and any $c'>c$ the function $d+c'$ is a metric.
\end{prb}

\begin{constr}[Induced distance] In fact, we have already implicitly used this obvious construction. Let $d$ be a distance function on the set $X$, and $Y$ be a nonempty subset of $X$. We define a distance function on $Y$ by setting $|yy'|=d(y,y')$ for any $y,y'\in Y$, and call it \emph{induced from $d$} or \emph{the restriction of $d$ to $Y$}. Note that the restriction of a metric (pseudometric) is always a metric (pseudometric). To emphasize that $Y$ is equipped  with the induce metric, we call such $Y$ \emph{a subspace of $X$}. However, in what follows, each subset of a metric space will be considered as a subspace, unless otherwise stated.
\end{constr}

\begin{constr}[Product]\label{constr:metricProd} Let $\{X_i\}_{i\in I}$ be some family of sets, and $d_i$ be a distance function on $X_i$. We call each $X_i$ equipped with $d_i$ \emph{a space}. We set $X=\prod_{i\in I}X_i$, and for any $x,x'\in X$ we consider the function $d_I(x, x')\:I\to\R$, $d_I(x,x')(i)=d_i\bigl(x_i,x'_i\bigr)$. Put
$$
D_I=\cup_{x,x'\in X}\{d_I(x,x')\},
$$
and let $\cV_I\sp D_I$ be some linear subspace of the linear space $\cF_I$ consisting of all real-valued functions defined on the set $I$, and $\cV_I^+\ss\cV_I$ consists of all functions with non-negative values. We call the subset $\cV_+\ss\cV$ by \emph{the first orthant of $\cV$}. Denote by $0\in\cF_I$ the zero function: $0(i)=0$ for any $i\in I$. Let $\r\:\cV_I^+\to[0,\infty)$ be an arbitrary function such that $\r(0)=0$. This function $\r$ generates a distance function $d_\r$ on $X$: $d_\r(x,x')=\r\bigl(d_I(x,x')\bigr)$. The function $\r$ is sometimes called \emph{a binder}, and the function $d_\r$ is called \emph{the distance generated by $\r$}; the space $X$ with the distance function $d_\r$ is called \emph{the semidirect product of the spaces $X_i$ w.r.t\. $\r$}. Below we will show a few examples in which all $d_i$ are metrics, and the resulting $d_\r$ is a metric as well.
\end{constr}

\begin{examp}
Consider the standard partial ordering on the first orthant $\cV_+$ of $\cV$: $v\le w$ if and only if $v(i)\le w(i)$ for all $i\in I$.

We call the function $\r$ from Construction~\ref{constr:metricProd}
\begin{itemize}
\item \emph{positively definite\/} if it vanishes only at the origin;
\item \emph{subadditive\/} if for any $v,w\in\cV_+$ it holds $\r(v+w)\le\r(v)+\r(w)$;
\item \emph{monotone\/} if it is monotone w.r.t\. the standard partial order described above.
\end{itemize}

\begin{prop}\label{prop:subadditiveMonotone}
Under notations of Construction~$\ref{constr:metricProd}$, let all $X_i$ be metric spaces, and $\r$ be positively definite, subadditive and monotone. Then the distance function $d_\r$ on $\prod_{i\in I}X_i$ is a metric.
\end{prop}

\begin{proof}
Positive definiteness of $d_\r$ follows from the one of metrics and of the function $\r$. It remains to verify the triangle inequality. Take arbitrary $x,x',x''\in X$, then
\begin{multline*}
d_\r(x,x'')=\r\bigl(d_I(x,x'')\bigr)\le\r\bigl(d_I(x,x')+d_I(x',x'')\bigr)\le\r\bigl(d_I(x,x')\bigr)+\r\bigl(d_I(x',x'')\bigr)=
d_\r(x,x')+d_\r(x',x'').
\end{multline*}
Here the first inequality follows from monotonicity of $\r$ and triangle inequalities for metrics, and the second one from subadditivity of $\r$.
\end{proof}
\end{examp}

\begin{examp}
Let $I$ be a finite set, say, $I=\{1,\ldots, n\}$, and $\cV_I=\cF_I=\R^n$. For $\r$ we choose the corresponding restriction of one of the following standard norms on $\cV_I$:
$$
\bigl\|(v_1,\ldots,v_n)\bigr\|_p=\sqrt[p]{\sum_{i=1}^n|v_i|^p}\ \ \text{for $1\le p<\infty$\ \ and}\ \
\bigl\|(v_1,\ldots,v_n)\bigr\|_\infty=\max\bigl\{|v_i|:i=1,\ldots,n\bigr\}.
$$
Then the distance function $d_\r$ is a metric because $\r$ is positively definite, subadditive and monotone (verify), thus we can apply Proposition~\ref{prop:subadditiveMonotone}. For $p=2$ (when the norm is Euclidean) we get \emph{the direct product\/} of the spaces $X_i$.
\end{examp}

\begin{examp}
On the first orthant $\R^2_+$ of the plane $\R^2$ with coordinates $(x,y)$ consider the following function $\r(x,y)=x+\sqrt y$. This function is obviously positively definite and monotonic. To prove it is subadditive, make the following calculation:
$$
\r(x+x',y+y')=x+x'+\sqrt{y+y'}\le x+\sqrt y+x'+\sqrt{y'}=\r(x,y)+\r(x',y').
$$
Let $X_1=X_2=\R$, then, by Proposition~\ref{prop:subadditiveMonotone}, the distance function $d_\r$ is a metric. In explicit form
$$
d_\r\bigl((x,y),(x',y')\bigr)=|x-x'|+\sqrt{|y-y'|}.
$$
\end{examp}

\begin{examp}
Let $I$ be an arbitrary set, and the family $\{X_i\}_{i\in I}$ consists of spaces with diameters bounded by the same constant, i.e., the set of numbers $\{\diam X_i\}_{i\in I}$ is bounded. Then, for $\cV$, we can choose the linear space of all bounded functions $v\:I\to\R$, and for $\r$ we can choose the norm $\|v\|_\infty=\sup_{i\in I}\Bigl\{\bigl|v(i)\bigr|\Bigr\}$. If all $X_i$ are metric spaces, then $d_\r$ is a metric.
\end{examp}

The previous construction can be generalized if we consider not all the Cartesian product $X=\prod_{i\in I} X_i$, but only a part of it extracted by some condition.

\begin{examp}\label{examp:ellpANDinfty}
We choose the set of natural numbers $\N$ as $I$, and the real line $\R$ as $X_i$, then each element of $\prod_{i\in\N}X_i$ can be represented as a sequence $(x_1,x_2,\ldots)$, where $x_i\in X_i=\R$.

Consider a subset $\ell_p\ss\prod_{i\in\N}X_i$ consisting of all sequences $x=(x_1,x_2,\ldots)$ for which
$$
\|x\|_p:=
\begin{cases}
\sqrt[p]{\strut\sum_{i=1}^\infty|x_i|^p}&\text{for $1\le p<\infty$, and}\\
\sup\bigl\{|x_i|:i\in\N\bigr\}&\text{for $p=\infty$}
\end{cases}
$$
is finite. Then $\|\cdot\|_p$ is a norm on $\ell_p$. The resulting metric space is denoted by $\ell_p$. In particular, the space $\ell_\infty$ consists of all bounded sequences.
\end{examp}

\begin{examp}
We choose an arbitrary set of indices $I$, and consider $I$ as topological space with discrete topology. We put all $X_i$ equal to $\R$, and as a subset of $\prod_{i\in I}X_i$ we take the family $\Fin_I$ of all maps $x\:I\to\R$ with compact support, i.e., such that $x(i)\ne0$ only for a finite number of indices $i$. It is clear that $\Fin_I$ forms a linear space, and each function
$$
\|x\|_p=\sqrt[p]{\sum_{i\in I}|x_i|^p}\ \ \text{for $1\le p<\infty$~~~and}\ \
\|x\|_\infty=\max_{i\in I}\bigl\{|x_i|\bigr\}
$$
is a norm on the space $\Fin_I$. If all $X_i$ are metric space, and $\r$ is a norm described above, then $d_\r$ is a metric.
\end{examp}

\begin{constr}[Levenshtein distance] Let $A$ be some set, and $A^*$ be the family of all finite sequences of elements from $A$, as well as an empty sequence $\l$. We will interpret $A$ as \emph{an alphabet\/} of some language. Then elements from $A$ are naturally called \emph{letters}, and the elements from $A^*$ are \emph{words}. \emph{An editorial operation on $A^*$} is a word transformation consisting either in the exclusion of one of the letters from the word (\emph{deletion\/}), or in the insertion of a letter into the word (\emph{insertion\/}), or in replacing one letter with another one (\emph{substitution\/}). The smallest number of editorial operations needed to move from one word to another is called \emph{the Levenshtein distance\/} and generates a metric on $A^*$. This distance plays an important role in linguistics and bioinformatics.
\end{constr}

Now we need to recall the notion of graph. We consider more general graphs permitting infinite numbers of vertices and edges.

\subsection{Elements of graph theory, and metric constructions for graphs}
Consider a triple $G=(V,E,i)$, where $V$ and $E$ are arbitrary sets, and $i\:E\to\cF_2(V)$ a mapping to the set of all at most $2$-point nonempty subsets of $V$, see Section~\ref{sec:hyperspaces}. The set $V$ is also denoted by $V(G)$ and called \emph{the set of vertices\/}; the set $E$ is also denoted by $E(G)$ and called \emph{the set of edges\/}; the mapping $i\:E\to\cF_2(V)$ is also denoted by $i_G$ and called \emph{the incidence mapping}, and it defines how edges are glued to vertices; the triple $G$ is called \emph{a non-oriented graph}. An edge $e$ is called \emph{a loop}, if $\#i(e)=1$. The number $\#i^{-1}\bigl(i(e)\bigr)$ is called \emph{multiplicity of the edge $e$}. An edge $e$ is called \emph{multiple}, if its multiplicity is more that $1$. Graph without loops and multiple edges is called \emph{simple}. In a simple graph each edge $e$ can be considered as the pair $\{v,w\}=i(e)$ of different vertices $v$ and $w$, that is why we can reduce the definition by setting $G=(V,E)$, where $E$ is a subset of $V^{(2)}$, where $V^{(2)}$ is the set of all 2-point subsets of $V$. \emph{A subgraph $H$ of the graph $G$} is a triple $(V',E',i')$ such that $V'\ss V$, $E'\ss E$, and $i'=i|_{E'}$.

If instead $\cF_2(V)$ we consider $V\x V$, then we obtain definition of \emph{an oriented graph}.

Let $G=(V,E,i)$ be an arbitrary non-oriented graph. For each edge $e\in E$ and each vertex $v\in i(e)$ we say that $e$ and $v$ are \emph{incident\/}; if $i(e)=\{v,w\}$, the vertices $v$ and $w$ are called \emph{adjacent}, and we say that $e$ \emph{joins $v$ and $w$}; if $i(e)=\{v\}$, we say that $v$ is \emph{adjacent to itself}, that $e$ \emph{joins $v$ with $v$}. The number of non-loop edges plus the double number of loops, incident to a vertex $v\in V$, is called \emph{the degree of $v$} and is denoted by $\deg v$ or $\deg_G v$.

\emph{A walk of length $n$ joining some vertices $v$ and $w$} is a sequence $v=v_0,\,e_1,\,v_1,\,e_2,\ldots,e_n,\,v_n=w$ of alternating vertices and edges such that each edge $e_k$ joins $v_{k-1}$ and $v_k$. If $G$ is a simple graph, we do not need to indicate the edges since each pair of adjacent vertices defines uniquely the edge incident to these vertices. The walt is called \emph{closed\/} if $v_0=v_n$, and it is called \emph{open\/} otherwise. \emph{A trail\/} is a walk with no repeated edges. \emph{A path\/} is an open trail with no repeated vertices. \emph{A circuit\/} is a closed trail. \emph{A cycle\/} is a circuit with no repeated vertices.

In the case of an oriented graph $G$, if $i(e)=(v,w)$, then we say that \emph{$e$ starts at $v$ and ends at $w$}. We redefine the walk just demanding that each $e_k$ starts at $v_{k-1}$ and ends at $v_k$. All other definitions remain the same.

A graph $G$ is called \emph{connected}, if each pair of its vertices are joined by a walk. A graph without cycles is called \emph{a forest}, and a connected forest is called \emph{a tree}. Clearly that each forest is a simple graph.

\emph{A weighted graph\/} is a graph $G=(V,E,i)$ equipped with \emph{a weight function $\om\:E\to[0,\infty)$} (sometimes it is useful to consider more general weight functions, for instance, with possibility of negative values). Sometimes we denote the weighted graph as $(V,E,i,\om)$ or $(G,\om)$, and, in the case of simple graph $G$, by $(V,E,\om)$. \emph{The weight of a subgraph\/} of a weighted graph is the sum of the weights of edges from this subgraph. We can extend this definition to paths and cycles considering them as subgraphs of $G$. In the case of the walk, its weight is defined as the sum of weights of its consecutive edges. For graphs without weight functions these notions are defined as well by assigning the weight $1$ to each edge by default.

\begin{constr}[Distance on a connected graph]\label{constr:graphMetric}
Let $\G=(V,E,i)$ be a connected graph. We define a distance function on $V$, setting it equal to the infimum of length of walks joining a given pair of vertices (we can always change the walks to paths with the same result). This function is a metric on $V$ (verify).
\end{constr}

\begin{constr}[Distance on a connected weighted graph]\label{constr:DistOnWeightedGraph}
If $\G=(V,E,i,\om)$ is a weighted connected graph with the weight function $\om\:E\to[0,\infty)$, then we define the distance function by setting it equal to the infimum of the weights of all walks (or paths) joining a given pair of vertices. The resulting function is a pseudometric on $V$. Note that even if $\om$ is everywhere positive, the constructed distance function does not have to be a metric (consider a graph with a positive weight function in which some pair of vertices are connected by an infinite number of paths with weights tending to zero).
\end{constr}

\begin{constr}[Cayley graph of a group]
Let $G$ be an arbitrary group with a set $S$ of generators. \emph{The Cayley graph of the pair $(G,S)$} is the directed graph $\G(G,S)=(G,E,i)$, in which $(g,h)=i(e)$ for some $e\in E$ if and only if $h=gs$ for some $s\in S$.

We note that in geometric group theory the set $S$ usually satisfies the following properties: the neutral element is not contained in $S$, and $S=S^{-1}$, i.e., if $s\in S$, then $s^{-1}\in S$. In this case, the graph $\G(G,S)$ has no loops, and if $(g,h)=i(e)$ for some $e\in E$, then $(h,g)=i(e')$ for some $e'\in E$, since $h=gs$ implies $g=hs^{-1}$. Under these assumptions, by Cayley graph we mean a simple non-oriented graph in which each pair of mutually opposite oriented edges were replaced by one (non-oriented) edge. For such Cayley graph, Construction~\ref{constr:graphMetric} defines a metric on $G$, which is used to determine the growth rate of the group $G$.
\end{constr}

\begin{rk}
In what follows, we will always consider only those $S$ that satisfies the both conditions: $S=S^{-1}$ and $S$ does not contain the neutral element of $G$. However, we will not list in $S$ all inverses of its elements supposing that this holds by default.
\end{rk}

\begin{prb}
Describe the Cayley graphs for the following groups $G$ and generating sets $S$:
\begin{enumerate}
\item $G=\Z$ and $S=\{1\}$;
\item $G=\Z_m$ and $S=\{1\}$;
\item $G=\Z^2$ and $S=\bigl\{(1,0)\,(0,1)\bigr\}$;
\item $G=\Z^2$ and $S=\bigl\{(1,0)\,(0,1),\,(1,1)\bigr\}$;
\item $G$ is a free group with generators $a$ and $b$.
\end{enumerate}
\end{prb}

\begin{constr}[Quotient spaces]
Let $(X,\r)$ be a metric (pseudometric) space, and $\sim$ an equivalence relation on $X$. Define the following \emph{quotient distance function\/} on $X$:
$$
\r_\sim(x,y)=\inf\Bigl\{\sum_{i=0}^n\r(p_i,q_i):p_0=x,\,q_n=y,\,n\in\N,\,q_i\sim p_{i+1}\ \text{for all $i$}\Bigr\}.
$$
The terms in the right-hand side can be naturally visualized as $p_0-q_0\sim p_1-q_1\sim\cdots\sim p_n-q_n$, where $p_i-q_i$ indicates that we calculate the distance between $p_i$ and $q_i$, and $q_i\sim p_{i+1}$ says that these points are taken from the same equivalency class, and thus the distance between them is zero. We call such sequences $p_0-q_0\sim p_1-q_1\sim\cdots\sim p_n-q_n$ an \emph{admissible\/} one joining the classes $[x]$ and $[y]$, and the value $\sum_{i=0}^n\r(p_i,q_i)$ \emph{the weight of this admissible sequence}.

\begin{prb}\label{prb:roSimIsPseudometric}
Prove that $\r_\sim$ is a pseudometric on $X$.
\end{prb}

Now, let $x\sim y$, then we can put $p_0=x=q_0$ and $p_1=q_1=y$, so $\r_\sim(x,y)\le\r(x,x)+\r(y,y)=0$, thus $\r_\sim(x,y)=0$. Therefore, each equivalence class of $\sim$ consists of points on $\r_\sim$-zero distance from each other, thus $\r_\sim$ generates correctly a pseudometric on $X/\!\!\sim$ which we also denote by $\r_\sim$ and call in \emph{the quotient pseudometric w.r.t\. $\sim$}; the space $X/\!\!\sim$ we will call the \emph{the quotient pseudometric space w.r.t\. $\sim$}. The next step is to consider the space $X/\r_\sim$ instead of the $X/\!\!\sim$: these two spaces are different when the distance between some distinct classes of $\sim$ vanishes. As above,  we denote by the same $\r_\sim$ the corresponding metric on $X/\r_\sim$, and the $\r_\sim$ and $X/\r_\sim$ we call \emph{the quotient metric space and the quotient metric space w.r.t\. $\sim$}, respectively.

To work effectively with the $\r_\sim$, let us define the following notions. The admissible sequence $p_0-q_0\sim p_1-q_1\sim\cdots\sim p_n-q_n$ from definition of $\r_\sim$ is called \emph{reducible\/} if it is possible to delete a part of it in such a way that its ends remains to belong to the classes $[x],[y]\in X/\!\!\sim$, and its weight did not increase. Otherwise, its is called \emph{irreducible}. It is evident that to calculate the distance $\r_\sim$, it is sufficient to consider irreducible admissible sequences only. If $a$ and $b$ not equivalent, then we write it as $a\not\sim b$.

\begin{prb}\label{prb:admissibleSeq}
Let $\xi=(p_0-q_0\sim p_1-q_1\sim\cdots\sim p_n-q_n)$ be an irreducible admissible sequence. Prove that
\begin{enumerate}
\item for any $i<j$ we have $p_i\not\sim p_j$, $q_i\not\sim q_j$;
\item $q_i\sim p_j$ if and only if $j=i+1$;
\item for any $i$ we have $q_i\not= p_{i+1}$.
\end{enumerate}
\end{prb}

\begin{examp}
Let us show how Problem~\ref{prb:admissibleSeq} can be applied. Let $X=[0,a]\ss\R$ be a segment of real line. Identify its ends. This means, that $x,y\in X$ are equivalent if and only if either $x=y$, or $\{x,y\}=\{0,a\}$. Describe all irreducible admissible sequences $\xi$. The are two possibilities (verify):
\begin{enumerate}
\item $\xi$ consists of two different points, say $x$ and $y$, and its weight equals $|xy|$;
\item $\xi$ consists of $4$ points, its ends are distinct from $0$ and $a$, but the middle points are these $0$ and $a$; w.l.o.g., $\xi=(x-0\sim a-y)$, thus its weight equals $|0x|+|ya|$.
\end{enumerate}
Thus, on the topological circle $[0,a]/\sim$ we introduce the distance in the same way as we do for the standard circle when we choose the length of the shortest arc between the points.
\end{examp}
\end{constr}

\begin{prb}
Let $\sim$ be the trivial equivalence on a metric space $(X,\r)$, i.e., $x\sim y$ if and only if $x=y$. Prove that $\r_\sim=\r$.
\end{prb}

\begin{prb}\label{prb:bUpperBound}
Let $\sim$ be an equivalence on a pseudometric space $(X,\r)$. Prove that for any $x,y\in X$ it holds $\r_\sim(x,y)\le\r(x,y)$. Thus, if we define the function $b\:X\to X\to\R$ such that $b(x,y)=0$ for $x\sim y$, and $b(x,y)=\r(x,y)$ otherwise, then $\r_\sim\le b$.
\end{prb}

\begin{constr}[Generalization of distance]
It is useful to allow infinite distances. A function $\r\:X\x X\to[0,\infty]$ that satisfies the axioms of distance (pseudometric, metric) will be called \emph{generalized\/} one. The corresponding spaces $X$ with such distances we call \emph{generalized\/} as well.
\end{constr}

\begin{rk}
In some monographs the authors work with generalized distances from the very beginning and because of that they call such distances without the word ``generalized''. However, the distances with values in $[0,\infty)$ they call \emph{finite\/} ones.
\end{rk}

The generalization of distance gives rise more metric constructions.

\begin{constr}[Disjoint union of spaces]
Let $\bigl\{(X_i,\r_i)\bigr\}_{i\in I}$ be an arbitrary family of generalized spaces. Consider on $\sqcup_{i\in I}X_i$ the distance function that is equal to $\r_i$ on $X_i$, and to $\infty$ for any pair of points from different $i$. The resulting space is called \emph{the disjoint union of spaces $X_i$}.

Evidently, if $X_i$ are generalized pseudometric (metric) spaces, then $\sqcup_{i\in I}X_i$ is the space of the same type.
\end{constr}

Let $X$ be a generalized pseudometric space. So, we can define two equivalence relations: $x\simind1y$ if and only if $|xy|<\infty$; and $x\simind2y$ if and only if $|xy|=0$.

\begin{prb}\label{prb:doubleQuotient}
Prove that each class of equivalence $\simind1$ is a pseudometric space (with finite distance), and that the distance between points from different classes equals $\infty$. Thus, if we denote by $X_i$ the classes of equivalence $\simind1$, then $X=\sqcup X_i$. Prove that the space $X/\!\!\simind2$ equals the disjoint union $\sqcup(X_i/\!\!\simind2)$ of metric spaces $X_i/\!\!\simind2$.
\end{prb}

Now we combine disjoint union and quotient operation.

\begin{examp}\label{examp:GluingSpace}
Let $X$ and $Y$ be metric spaces, $Z\ss X$, and $f\:Z\to Y$ a mapping. Consider on $X\sqcup Y$ with generalized metric $\r$ the equivalence relation $\sim$ which identifies each $a\in f(Z)$ with all $b\in f^{-1}(a)$. The quotient space $X\sqcup_fY:=(X\sqcup Y)/\r_\sim$ is called the result of \emph{gluing the spaces $X$ and $Y$ over the mapping $f$}.
\end{examp}

\begin{prb}\label{prb:GluingSpace}
Suppose that $f$ is isometric. Prove that the restrictions of $\r_\sim$ to $X$ and $Y$ coincides with the initial metrics of $X$ and $Y$, respectively.
\end{prb}

\begin{prb}
Let $y_0\in Y$ and $f(X)=y_0$. Prove that $X\sqcup_fY$ is isometric to $Y$.
\end{prb}

\begin{constr}[Maximal pseudometric]\label{constr:MaxPseudometric}
Consider an arbitrary function $b\:X\x X\to[0,\infty]$, and denote by $\cD_b$ the set of all generalized pseudometrics $d\:X\x X\to[0,\infty]$ such that $d(x,y)\le b(x,y)$ for any $x,y\in X$.

\begin{lem}\label{lem:MaxPseudometricExistsAndUnique}
There exists and unique $d_b\in\cD_b$ such that $d_b\ge d$ for all $d\in\cD_b$.
\end{lem}

\begin{proof}
We put $d_b(x,y)=\sup_{d\in\cD_b}d(x,y)$. It is evident that $d_b$ is nonnegative, symmetric, $d_b(x,x)=0$ for all $x\in X$, and $d_b\le b$. It remains to verify the triangle inequality:
$$
d_b(x,z)=\sup_{d\in\cD_b}d(x,z)\le\sup_{d\in\cD_b}\bigl(d(x,y)+d(y,z)\bigr)\le\sup_{d\in\cD_b}d(x,y)+\sup_{d\in\cD_b}d(y,z)= d_b(x,y)+d_b(y,z).
$$
The uniqueness of $d_b$ is evident.
\end{proof}

We denote $d_b$ by $\sup\cD_b$ and call it \emph{maximal pseudometric w.r.t\. the function $b$}.
\end{constr}

\begin{thm}\label{thm:QuotirntInTermsMaxFunct}
Let $(X,\r)$ be a pseudometric space, and $\sim$ an equivalence relation. We put
$$
b(x,y):=b_\sim(x,y)=
\begin{cases}
0,&\text{if $x\sim y$},\\
\r(x,y),&\text{otherwise}.
\end{cases}
$$
Then $d_b=\sup\cD_b=\r_\sim$, thus $X/\r_\sim=X/d_b$.
\end{thm}

\begin{proof}
By Problem~\ref{prb:bUpperBound}, we have $\r_\sim\le b$, thus $\r_\sim\in\cD_b$. It remains to prove that $\r_\sim\ge d$ for each $d\in\cD_b$. To do that, we write down
\begin{multline*}
\r_\sim(x,y)=\inf\Bigl\{\sum_{i=0}^n\r(p_i,q_i):p_0=x,\,q_n=y,\,n\in\N,\,q_i\sim p_{i+1}\ \text{for all $i$}\Bigr\}\ge\\ \ge
\inf\Bigl\{d(p_0,q_0)+d(q_0,p_1)+\cdots+d(q_{n-1},p_n)+d(p_n,q_n):n\in\N,\,q_i\sim p_{i+1}\ \text{for all $i$}\Bigr\}\ge d(x,y).
\end{multline*}
\end{proof}

Now we apply the quotation technique to obtain a few more important classes of metric spaces.

\begin{constr}[Groups actions]
Let $X$ be a metric space and $G\ss\Iso(X)$ a subgroup of its isometry group. Consider the action of $G$ on $X$, i.e., the mapping $\v\:G\x X\to X$ such that $\v(g,x)=:g(x)$ satisfies the following conditions:
\begin{enumerate}
\item $e(x)=x$ for the neutral element $e\in G$ and any $x\in X$;
\item $(hg)(x)=h\bigl(g(x)\bigr)$ for any $g,h\in G$ and $x\in X$.
\end{enumerate}

Consider the following equivalence relation on $X$: $x\sim y$ if and only if $g(x)=y$ for some $g\in G$. We say that the equivalence $\sim$ \emph{is generated by the action of $G$ on $X$}. The corresponding quotient space $X/\!\!\sim$ is usually denoted by $X/G$. The sets $G(x)$ are called \emph{orbits}, they are elements of $X/G$; the set $X/G$ is called \emph{the orbit-space}.
\end{constr}

\begin{prb}
Let $X$ be a metric space and $G\ss\Iso(X)$ a subgroup of its isometry group. For each two elements $G(x),G(y)\in X/G$ we set $d\bigl(G(x),G(y)\bigr)=\inf\bigl\{|x'y'|:x'\in G(x),\, y'\in G(y)\bigr\}$. Prove that $d=\r_\sim$, where the equivalence $\sim$ is generated by the action of $G$ on $X$.
\end{prb}

\begin{prb}\label{prb:StandardTorus}
Let $S^1$ be the standard unit circle in the Euclidean plane. As a distance between $x,y\in S^1$ we take the length of the shortest arc of $S_1$ between $x$ and $y$. By \emph{the standard torus\/} we mean the direct product $T^2=S^1\x S^1$ (with the Euclidean binder). We describe the points on the both $S^1$ by their polar angles $\v_1$ and $\v_2$, defined up to $2\pi$. So, \emph{the shifts\/} $s_{a,b}\:(\v_1,\v_2)\mapsto(\v_1+a,\v_2+b)$ are isometries of $T^2$. Fix some $(a,b)\in\R^2$ and consider a subgroup $G_{a,b}\ss\Iso(T^2)$ consisting of all shifts $s_{ta,tb}$, $t\in\R$. For different $a$ and $b$, find the corresponding pseudometric and metric quotient spaces.
\end{prb}

\begin{constr}[Metrized graphs]
Take an arbitrary collection of segments $\bigl\{[a_k,b_k]\bigr\}_{k\in K}$, and on the set $\sqcup\{a_k,b_k\}$ consider an equivalence $\sim$. The generalized pseudometric metric space $\bigl(\sqcup[a_k,b_k]\bigr)/\sim$ is called a \emph{metrized graph}. To obtained its combinatorial structure, we represent it as we did before, namely, as a triple $G=(V,E,i)$. To do that, we put $V:=\bigl(\sqcup\{a_k,b_k\}\bigr)/\sim$, $E:=\bigl\{[a_k,b_k]\bigr\}_{k\in K}$, and $\pi\:\sqcup\{a_k,b_k\}\to\bigl(\sqcup\{a_k,b_k\}\bigr)/\sim$ be the canonical projection, then $i\bigl([a_i,b_i]\bigr):=\pi\bigl(\{a_i,b_i\}\bigr)$.
\end{constr}

\begin{constr}[Polyhedral spaces]
By \emph{a polyhedron of dimension $n$} we mean an intersection of a finite number of half-spaces in $\R^n$ that has nonempty interior. For each point $x$ of a polyhedron $W$ we define its \emph{dimension $\dim_Wx$} as follows: consider all affine subspaces $L$ containing $x$ such that $x$ is an interior point of $W\cap L$, and put $\dim_Wx$ to be the maximal dimension of such $L$. Denote by $W^k$ the subset in $W$ consisting of all points of dimension $k$. One can show that $W^k$ consists of connected component each of which belongs of an affine subspace $L$ of dimension $k$. The closures of these connected components are called \emph{the faces of $W$ of dimension $k$}.

Now, consider an arbitrary collection $\{W_k\}_{k\in K}$ of polyhedra (probably, of different dimensions), and for some pairs $(W_k,W_l)$, $k\ne l$, choose some faces $F_k^p\ss W_k$ and $F_l^q\ss W_l$ such that there exists an isometry $\v\:F_k^p\to F_l^q$. Consider the generalized metric space $\sqcup_{k\in K}W_k$ and the equivalence relation generated by the isometries $\v$: we put in one class each point $x\in F_k^p$ and $\v(x)\in F_l^q$, for all isometries $\v$ (we glue the faces by the isometry $\v$). The resulting quotient generalized space is called \emph{a polyhedron space}.
\end{constr}

\begin{prb}
Represent the standard torus from Problem~\ref{prb:StandardTorus} as a polyhedron space.
\end{prb}

\begin{prb}
Verify that the above constructions actually, as declared, define (pseudo-)metrics.
\end{prb}

A huge collection of metric spaces can be found in~\cite{Deza2}.

\section{Convergence and completeness}
\markright{\thesection.~Convergence and completeness}

Since each metric space is Hausdorff, the following result holds.

\begin{prop}
If a sequence converges in a metric space, then its limit is uniquely determined.
\end{prop}

A sequence $x_1,\,x_2,\ldots$ of points of a metric space $X$ is called \emph{fundamental\/} if for any $\e>0$ there exists $N$ such that for all $m,\ n\ge N$ the inequality $|x_mx_n|<\e$ holds. If every fundamental sequence in a metric space $X$ is convergent, then such $X$ is called \emph{complete}.

\begin{prb}
Show that a subspace of a complete metric space is complete if and only if it is closed.
\end{prb}

\begin{thm}\label{thm:completion}
Each metric space $X$ is an everywhere dense subspace of some complete space called a \textbf{completion of $X$}. The completion of the space $X$ is uniquely defined: for any two completions $X'\sp X$ and $X''\sp X$ there is an isometry $f\:X'\to X''$ that is identical on $X$.
\end{thm}

\begin{constr}[Completion]
The standard construction of completion of a metric space $X$ is as follows:
\begin{itemize}
\item the set $W$ of all fundamental sequences in the space $X$ is considered;
\item on $W$, a pseudometric is set equal to the limit of the distances between the points of two fundamental sequences (the existence of the limit follows from the fundamentality);
\item $X$ can be realized as a subset of $W$ by associating each point $x\in X$ with the constant sequence $x,\,x,\ldots$.
\item the desired completion is obtained by factorization as in Problem~\ref{prb:pseudoTOmetric}.
\end{itemize}
\end{constr}

The uniqueness of completion implies the following result.

\begin{prb}
Let $X$ be an arbitrary subspace of a complete metric space. Then the closure $\bX$ of the set $X$ is a completion of the space $X$.
\end{prb}

\begin{prb}
Let $f\:X\to Y$ be a bi-Lipschitz mapping of metric spaces. Prove that $X$ is complete if and only of $Y$ is complete. Construct a homeomorphism of metric spaces that does not preserves completeness.
\end{prb}

\begin{prb}
Show that a metric space is complete if and only if the following condition holds: for any sequence of closed subsets $X_1\sp X_2\sp X_3\sp\cdots$ such that $\diam X_n\to0$ as $n\to0$, the intersection $\cap_{i=1}^\infty X_i$ is not empty (in fact, it consists of unique element). Show that the condition $\diam X_n\to0$ is essential.
\end{prb}

\begin{prb}[Fixed-point theorem]
Let $f\:X\to X$ be a $C$-Lipschitz mapping of a complete metric spaces $X$. Prove that for $C<1$ there exists and unique a point $x_0$ such that $f(x_0)=x_0$ (it is called \emph{the fixed point of the mapping $f$}).
\end{prb}

\section{Compactness and sequential compactness}
\markright{\thesection.~Compactness and sequential compactness}

Recall that a topological space is called \emph{compact} if a finite subcover can be extracted from any of its open cover, and \emph{sequentially compact} if any sequence of its points has a convergent subsequence. As noted in Section~\ref{ch:GenTop}, in the case of general topological spaces, compactness and sequential compactness are different concepts. However, for metric spaces this is not so.

\begin{thm}\label{thm:ComcactSequentialCompact}
A metric space is compact if and only if it is sequentially compact.
\end{thm}

\begin{proof}
Let $X$ be a compact metric space and $x_1, x_2,\ldots$ be an arbitrary sequence of points from $X$. We must show that there is a convergent subsequence in this sequence.

Suppose this is not so, then
\begin{itemize}
\item the sequence $x_i$ contains an infinite number of different points;
\item for each point $x\in U:=X\sm\{x_1,\, x_2,\ldots\}$ there is an (open) neighborhood $U^x$ that does not contain points of the sequence $x_i$, therefore $U=\cup_{x\in U}U^x$ and, thus, $U$ is an open set;
\item for each $x_i$ there is a neighborhood $U^{x_i}$ for which $U^{x_i}\cap\{x_1,\, x_2,\ldots\}=\{x_i\}$.
\end{itemize}
But then the family $\{U,\,U^{x_1},\,U^{x_2},\ldots\}$ is an open cover of $X$ from which it is impossible to choose a finite subcover (each finite subcover contains only finitely many points of the sequence $x_i$). The obtained contradiction shows that $X$ is sequentially compact.

Now let the space $X$ be sequentially compact. Suppose that there exists an open covering $\cU=\{U_\a\}$ of the space $X$ that does not contain finite subcovers. On $X$ we define a function $\r\:X\to\R$ as follows:
$$
\r(x)=\sup\{r\in\R\mid\exists U_\a\in\cU:\U_r(x)\ss U_a\}.
$$
This function is everywhere finite, since the sequentially compact space is bounded, and everywhere positive by virtue of the definition of metric topology.

\begin{lem}\label{lem:CoveringsInscribe1Lipschitz}
The function $\r$ defined above is $1$-Lipschitz and, therefore, continuous.
\end{lem}

\begin{proof}
Assume the contrary, i.e., that for some $x,\,y\in X$ the inequality $|\r (y)-\r (x)|>|xy|$ holds. Without loss of generality, we assume that $\r(y)>\r(x)$, then $\r(y)>\r(x)+|xy|$. Increase the number $\r(x)$ a little to $\r'$ and slightly decrease the number $\r(y)$ to $\r''$ so that $\r''>\r'+|xy|$ is still true, then
\begin{itemize}
\item $U _{\r'}(x)\ss U _{\r''}(y)$, because for an arbitrary point $z\in U _{\r'}(x)$ we have
$$
|zy|\le|zx|+|xy|<\r'+|xy|<\r'';
$$
\item there is $U_\a\in\cU$ for which $U_{\r''}(y)\ss U_\a$.
\end{itemize}
But then $U _{\r'}(x)\ss U_\a$, therefore $\r(x)\ge\r'$ by the definition of the function $\r$, a contradiction.
\end{proof}

By virtue of Problem~\ref{prb:MinMaxSequentialCompact}, the function $\r$ achieves the smallest value $\r_0$, which is, therefore, strictly positive. Put $r=\r_0/2$. Then for each point $x\in X$ there exists $U_\a\in\cU$ such that $U_r(x)\ss U_\a$.

Choose an arbitrary point $x_1\in X$, and let $U_1\in\cU$ be such that $U_r(x_1)\ss U_1$. There is $x_2\in X\sm U_1$. Choose $U_2\in\cU$ such that $U_r(x_2)\ss U_2$. In general, if $x_1,\ldots,x_n$ and $U_1,\ldots,U_n$ are selected, then there is $x_{n+1}\in X\sm\cup_{i=1}^nU_i$ and $U_{n+1}\in\cU$ such that $U_r(x_{n+1})\ss U_{n+1}$. Since $\cU$ does not have a finite subcover, we construct an infinite sequence $x_1,\,x_2,\ldots$, and it is clear that every point $x_{n+1}$ lies outside $\cup_{i=1}^nU_r(x_i)$, so for any $x_i$ and $x_j$, $i\ne j$, we have $|x_ix_j|\ge r$. But such a sequence does not contain a convergent subsequence, which contradicts the sequential compactness of $X$.
\end{proof}

\begin{prb}
Give an example of a topological space that is
\begin{enumerate}
\item compact, but not sequentially compact;
\item sequentially compact, but not compact.
\end{enumerate}
\end{prb}

\begin{prb}[Lebesgue's lemma]
Let $X$ be a compact metric space. Prove the following statement: for any open cover $\{U_i\}_{i\in I}$ of $X$ there exists $\r>0$ such that for any $x\in X$ one can find $U_i$ with $B_\r(x)\ss U_i$.
\end{prb}

\begin{prb}
Show that each continuous mapping $f\:X\to Y$ from a compact metric space to an arbitrary metric space is uniformly continuous.
\end{prb}

\begin{prb}
Prove that the diameter $\diam X$ of a compact metric space $X$ if finite, and that there exist $x,y\in X$ such that $\diam X=|xy|$.
\end{prb}

\subsection{Completeness and compactness}

Theorem~\ref{thm:ComcactSequentialCompact} easily yields the following result.

\begin{cor}\label{cor:CompactComplete}
A compact metric space is complete.
\end{cor}

\begin{proof}
Let $X$ be a compact metric space. In the space $X$ we choose an arbitrary fundamental sequence $x_1,\, x_2,\ldots$. By Theorem~\ref{thm:ComcactSequentialCompact}, there exists a subsequence $x_{i_1},\,x_{i_2},\ldots$, converging to some point $x\in X$. Due to fundamentality of $x_1,\,x_2,\ldots$, the whole this sequence converges to $x$, therefore the space $X$ is complete.
\end{proof}

\begin{dfn}
For $\e>0$, a subset $S$ of a metric space $X$ is called \emph{an $\e$-net} if for any point $x\in X$ there exists $s\in S$ such that $|xs|<\e$. A metric space is called \emph{totally bounded} if for any $\e>0$ there exists a finite $\e$-net in it.
\end{dfn}

\begin{thm}\label{thm:CompleteCompact}
A metric space is compact if and only if it is complete and totally bounded.
\end{thm}

\begin{proof}
Let $X$ be a compact metric space. Then, by virtue of Corollary~\ref{cor:CompactComplete}, it is complete. We choose an arbitrary $\e>0$, then $\{U_\e(x)\}_{x\in X}$ is an open cover of $X$. Since $X$ is compact, there exists a finite subcover $\{U_\e(x_i)\}_{i=1}^n$. But then $\{x_1,\ldots,x_n\}$ is a finite $\e$-net. Thus, the space $X$ is totally bounded.

Now let $X$ be a complete and totally bounded metric space. We prove that $X$ is sequentially compact and apply Theorem~\ref{thm:ComcactSequentialCompact}. Consider an arbitrary sequence $x_i$ in $X$. For each $n\in\N$, consider a finite $1/n$-net $S_n$. The balls $\{U_1(s)\}_{s\in S_1}$ cover $X$, so there is a ball $U_1$ among them that contains infinitely many elements of the sequence $x_i$. The balls $\{U_{1/2}(s)\}_{s\in S_2}$ cover $U_1$, so there is a ball $U_2$ among them for which $C_2=U_1\cap U_2$ contains infinitely many $x_i$. If $U_1,\,U_2,\ldots,U_k$ are already selected so that $C_k=\cap_{j=1}^kU_j$ contains infinitely many $x_i$, then the family $\{U_{1/(k+1)}(s)\}_{s\in S_{k+1}}$, since it covers $C_k$, contains a ball $U_{k+1}$ for which there are infinitely many $x_i$ in $C_{k+1}=\cap_{j=1}^{k+1}U_j$.

Now, we choose an arbitrary $x_{i_1}\in C_1$. Since there are infinitely many points of our sequence in $C_2$, there exists $i_2>i_1$ such that $x_{i_2}\in C_2$. Continuing this process, we construct a subsequence $x_{i_1},x_{i_2},\ldots$ such that $x_{i_k},x_{i_{k+1}},\ldots\in C_k$ for each $k$, and since $\diam C_k\to\0$ as $k\to\infty$, this subsequence is fundamental and, therefore, converges because $X$ is complete. Thus, $X$ is sequentially compact.
\end{proof}

\begin{prb}
Prove that every compact metric space is separable.
\end{prb}

A subset of a topological space is called \emph{nowhere dense\/} if its closure has empty interior.

\begin{prb}[Baire's theorem]
Prove that a complete metric space cannot be covered by at most countably many nowhere dense subsets. Moreover, the complement of the union of at most countably many nowhere dense subsets is everywhere dense.
\end{prb}

\begin{prb}
Prove that a compact metric space $X$ cannot be isometrically mapped to a subspace $Y\ss X$ such that $Y\ne X$. In other words, each isometric mapping $f\:X\to X$ for a compact metric space $X$ is surjective.
\end{prb}

\begin{prb}
Let $X$ be a compact metric space and $f\:X\to X$ be a mapping. Prove that
\begin{enumerate}
\item if $f$ is surjective and nonexpanding, then $f$ is an isometry\/\rom;
\item if $\bigl|f(x)f(y)\bigr|\ge|xy|$ for all $x,y\in X$, then $f$ is an isometry.
\end{enumerate}
\end{prb}

\section{Canonical isometric embeddings of metric spaces}
\markright{\thesection.~Canonical isometric embeddings of metric spaces}

Let $X$ be an arbitrary metric space. We denote by $C(X)$ the vector space of functions continuous on $X$, and by $C_b(X)$ the subspace of $C (X)$ composed of all bounded functions, and consider on $C_b(X)$ the norm $\|f\|_\infty=\sup_{x\in X}\bigl|f(x)\bigr|$ and the corresponding metric $|fg|_\infty=\|f-g\|_\infty$.

We define a mapping $\nu\:X\to C_b (X)$ as follows. For each point $x\in X$, by $d_x$ we denote the function $d_x\:X\to[0,\infty)$ defined by the rule $d_x(y)=|xy|$. This function is $1$-Lipschitz and, therefore, $d_x\in C(X)$. Now we fix some point $p\in X$ and consider the function $d_x-d_p$. Then for each $y\in X$ we have $\bigl|d_x(y)-d_p (y)\bigr|\le|xp|$, so the function $d_x-d_p$ is bounded and, thus, belongs to $C_b(X)$.

\begin{thm} [Frechet, Kuratowski]\label{thm:Frechet-Kuratowski}
The mapping $\nu\:X\to C_b (X)$ defined by the formula $\nu\:x\mapsto d_x-d_p$ is an isometric embedding.
\end{thm}

\begin{proof}
By virtue of Problem~\ref{prb:IsomHomeo}, it suffices to verify that the mapping $\nu$ is isometric.

Choose arbitrary $x,y\in X$, then
$$
\bigl|(d_x-d_p)(d_y-d_p)\bigr|_\infty=\sup_{z\in X}\bigl|d_x (z)-d_y (z)\bigr|\le|xy|
$$
by triangle inequality. On the other hand, if $z=y$, then $\bigl|d_x(z)-d_y(z)\bigr|=|xy|$, which implies the required result.
\end{proof}

Let us prove a finer result. Recall that by $\ell_\infty$ we denote the space of all bounded sequences with the metric defined by the norm $\|\cdot\|_\infty$.

\begin{thm} [Frechet]\label{thm:Frechet-separable}
Let $X$ be a separable metric space, then $X$ can be isometrically embedded into $\ell_\infty$.
\end{thm}

\begin{proof}
Consider an everywhere dense sequence $x_1,x_2,\ldots$ in $X$, which exists due to separability. Choose an arbitrary point $p\in X$ and associate with each point $x\in X$ a sequence $\nu(x)$ obtained by restricting the function $d_x-d_p$ to the sequence $x_1,x_2,\ldots$, namely, $\nu(x)(i)=d_x(x_i)-d_p(x_i)$. Since the functions $d_x-d_p$ are bounded, then $\nu(x)\in\ell_\infty$. The inequality $\bigl|\nu(x)\nu(y)\bigr|_\infty\le|xy|$ is verified in the same way as in the proof of Theorem~\ref{thm:Frechet-Kuratowski}. The converse inequality follows from the fact that the subset $\{x_i\}$ is everywhere dense in $X$: consider $x_{i_k}\to x$, then
$$
\bigl|\nu(x)\nu(y)\bigr|_\infty=\sup_{z\in X}\bigl|d_x(z)-d_y(z)\bigr|\ge\bigl|d_x( x_{i_k})-d_y(x_{i_k})\bigr|\to d_y(x)=|xy|\ \text{as $k\to\infty$}.
$$
The theorem is proved.
\end{proof}


\phantomsection
\renewcommand\bibname{References to Chapter~\thechapter}
\addcontentsline{toc}{section}{\bibname}
\renewcommand{\refname}{\bibname}

\vfill\eject

\section*{\Huge Exercises to Chapter~\thechapter}
\markright{Exercises to Chapter~\thechapter.}

\begin{exe}
Let $X$ be a pseudometric space and $\sim$ is the natural equivalence relation: $x\sim y$ if and only if $|xy|=0$. For each $x\in X$ denote by $[x]$ the equivalence class containing $x$. Prove that for any $x,y\in X$, $x'\in[x]$, and $y'\in[y]$ it is true that $|x'y'|=|xy|$. Thus, on the set $X/\!\!\sim$ the corresponding distance function is correctly defined: $\bigl|[x][y]\bigr|=|xy|$. Show that this distance function is a metric.
\end{exe}

\begin{exe}
Let $X$ be an arbitrary metric space, $x,y\in X$, $r\ge0$, $s,t>0$, and $A\ss X$ be nonempty. Verify that
\begin{enumerate}
\item $U_s\bigl(\{x\}\bigr)=U_s\bigl(x)$ and $B_r\bigl(\{x\}\bigr)=B_r\bigl (x)$;
\item the functions $y\mapsto|xy|$, $y\mapsto|yA|$ are continuous;
\item an open neighborhood $U_s(A)$ is an open subset of $X$, and a closed neighborhood $B_r(A)$ is a closed subset of $X$;
\item $U_t\bigl(U_s(A)\bigr)\ss U_{s+t}(A)$ and construct an example demonstrating that the left-hand side can be different from the right-hand side;
\item $B_t\bigl(B_s(A)\bigr)\ss B_{s+t}(A)$ and construct an example demonstrating that the left-hand side can be different from the right-hand side;
\item $\d U_s(x)$, $\d B_s(x)$ are not related by any inclusion; $\d U_s(x)\ss S_s(x)$ and $\d B_r(x)\ss S_r(x)$; the both previous inclusions can be strict;
\item $\diam U_s(x)\le\diam B_s(x)\le2s$;
\item $\diam U_s(A)\le\diam B_s(A)\le\diam A+2s$.
\end{enumerate}
\end{exe}

\begin{exe}
Let $\cL(f)\ss\R$ be the set of all Lipschitz constants for a Lipschitz mapping $f$. Prove that $\inf\cL(f)$ is also a Lipschitz constant.
\end{exe}

\begin{exe}
Show that each Lipschitz map is uniformly continuous, and each uniformly continuous map is continuous.
\end{exe}

\begin{exe}
Show that each isometry is a homeomorphism, in particular, each isometric mapping of one metric space into another one is an embedding, i.e., we recall, it is a homeomorphism with an image.
\end{exe}

\begin{exe}
Verify that the identity map, the composition of isometries, and the inverse mapping to an isometry are also isometries, i.e., the set of all isometries of an arbitrary metric space forms a group.
\end{exe}


\begin{exe}
Let $X$ be an arbitrary metric space, $x,y\in X$ and $A\ss X$ be nonempty. Prove that $|Ax|+|xy|\ge|Ay|$, so that the function $\r_A(x)=|Ax|$ is $1$-Lipschitz and, therefore, uniformly continuous.
\end{exe}

\begin{exe}
Describe the Cayley graphs for the following groups $G$ and generating sets $S$:
\begin{enumerate}
\item $G=\Z$ and $S=\{1\}$;
\item $G=\Z_m$ and $S=\{1\}$;
\item $G=\Z^2$ and $S=\bigl\{(1,0)\,(0,1)\bigr\}$;
\item $G=\Z^2$ and $S=\bigl\{(1,0)\,(0,1),\,(1,1)\bigr\}$;
\item $G$ is a free group with generators $a$ and $b$.
\end{enumerate}
\end{exe}

\begin{exe}
Let $(X,\r)$ be a metric (pseudometric) space, and $\sim$ an equivalence relation on $X$. Define the following \emph{quotient distance function\/} on $X$:
$$
\r_\sim(x,y)=\inf\Bigl\{\sum_{i=0}^n\r(p_i,q_i):p_0=x,\,q_n=y,\,n\in\N,\,q_i\sim p_{i+1}\ \text{for all $i$}\Bigr\}.
$$
Prove that $\r_\sim$ is a pseudometric on $X$.
\end{exe}

\begin{exe}
Let $\xi=(p_0-q_0\sim p_1-q_1\sim\cdots\sim p_n-q_n)$ be an irreducible admissible sequence. Prove that
\begin{enumerate}
\item for any $i<j$ we have $p_i\not\sim p_j$, $q_i\not\sim q_j$;
\item $q_i\sim p_j$ if and only if $j=i+1$;
\item for any $i$ we have $q_i\not= p_{i+1}$.
\end{enumerate}
\end{exe}

\begin{exe}
Let $\sim$ be the trivial equivalence on a metric space $(X,\r)$, i.e., $x\sim y$ if and only if $x=y$. Prove that $\r_\sim=\r$.
\end{exe}

\begin{exe}
Let $\sim$ be an equivalence on a pseudometric space $(X,\r)$. Prove that for any $x,y\in X$ it holds $\r_\sim(x,y)\le\r(x,y)$. Thus, if we define the function $b\:X\x X\to\R$ such that $b(x,y)=0$ for $x\sim y$, and $b(x,y)=\r(x,y)$ otherwise, then $\r_\sim\le b$.
\end{exe}

\begin{exe}
Let $X$ be a generalized pseudometric space. We can define two equivalence relations: $x\simind1y$ if and only if $|xy|=\infty$, and $x\simind2y$ if and only if $|xy|=0$. Prove that each class of equivalence $\simind1$ is a pseudometric space (with finite distance), and that the distance between points from different classes equals $\infty$. Thus, if we denote by $X_i$ the classes of equivalence $\simind1$, then $X=\sqcup X_i$. Prove that the space $X/\!\!\simind2$ equals the disjoint union $\sqcup(X_i/\!\!\simind2)$ of metric spaces $X_i/\!\!\simind2$.
\end{exe}

\begin{exe}
Let $X$ and $Y$ be metric spaces, $Z\ss X$, and $f\:Z\to Y$ is an isometric mapping. Let $\r$ be the metric on the $X\sqcup_fY$. Prove that the restrictions of $\r$ to $X$ and $Y$ coincides with the initial metrics of $X$ and $Y$, respectively.
\end{exe}

\begin{exe}
Let $y_0\in Y$ and $f(X)=y_0$. Prove that $X\sqcup_fY$ is isometric to $Y$.
\end{exe}

\begin{exe}
Let $X$ be a metric space and $G\ss\Iso(X)$ a subgroup of its isometry group. For each two elements $G(x),G(y)\in X/G$ we set $d\bigl(G(x),G(y)\bigr)=\inf\bigl\{|x'y'|:x'\in G(x),\, y'\in G(y)\bigr\}$. Prove that $d=\r_\sim$, where the equivalence $\sim$ is generated by the action of $G$ on $X$.
\end{exe}

\begin{exe}\label{exe:StandardTorus}
Let $S^1$ be the standard unit circle in the Euclidean plane. As a distance between $x,y\in S^1$ we take the length of the shortest arc of $S_1$ between $x$ and $y$. By \emph{the standard torus\/} we mean the direct product $T^2=S^1\x S^1$ (with the Euclidean binder). We describe the points on the both $S^1$ by their polar angles $\v_1$ and $\v_2$, defined up to $2\pi$. So, \emph{the shifts\/} $s_{a,b}\:(\v_1,\v_2)\mapsto(\v_1+a,\v_2+b)$ are isometries of $T^2$. Fix some $(a,b)\in\R^2$ and consider a subgroup $G_{a,b}\ss\Iso(T^2)$ consisting of all shifts $s_{ta,tb}$, $t\in\R$. For different $a$ and $b$, find the corresponding pseudometric and metric quotient spaces.
\end{exe}

\begin{exe}
Represent the standard torus from Exercise~\ref{exe:StandardTorus} as a polyhedron space.
\end{exe}

\begin{exe}
Verify that the constructions given in Section~\ref{sec:StandCOnstrMetric} do define (pseudo-)metrics, as declared.
\end{exe}

\begin{exe}
Let $d$ be a metric. Find the least possible $c$ such that $d+c$ is a pseudometric. Verify that for such $c$ and any $c'>c$ the function $d+c'$ is a metric.
\end{exe}

\begin{exe}
Show that a subspace of a complete metric space is complete if and only if it is closed.
\end{exe}

\begin{exe}
Let $X$ be an arbitrary subspace of a complete metric space. Then the closure $\bX$ of the set $X$ is a completion of the space $X$.
\end{exe}

\begin{exe}
Let $f\:X\to Y$ be a bi-Lipschitz mapping of metric spaces. Prove that $X$ is complete if and only of $Y$ is complete. Construct a homeomorphism of metric spaces that does not preserves completeness.
\end{exe}

\begin{exe}
Show that a metric space is complete if and only if the following condition holds: for any sequence of closed subsets $X_1\sp X_2\sp X_3\sp\cdots$ such that $\diam X_n\to0$ as $n\to0$, the intersection $\cap_{i=1}^\infty X_i$ is not empty (in fact, it consists of unique element). Show that the condition $\diam X_n\to0$ is essential.
\end{exe}

\begin{exe}[Fixed-point theorem]
Let $f\:X\to X$ be a $C$-Lipschitz mapping of a complete metric spaces $X$. Prove that for $C<1$ there exists and unique a point $x_0$ such that $f(x_0)=x_0$ (it is called \emph{the fixed point of the mapping $f$}).
\end{exe}

\begin{exe}
Give an example of a topological space that is
\begin{enumerate}
\item compact, but not sequentially compact;
\item sequentially compact, but not compact.
\end{enumerate}
\end{exe}

\begin{exe}[Lebesgue's lemma]
Let $X$ be a compact metric space. Prove the following statement: for any open cover $\{U_i\}_{i\in I}$ of $X$ there exists $\r>0$ such that for any $x\in X$ one can find $U_i$ with $B_\r(x)\ss U_i$.
\end{exe}

\begin{exe}
Show that each continuous mapping $f\:X\to Y$ from a compact metric space to an arbitrary metric space is uniformly continuous.
\end{exe}

\begin{exe}
Prove that the diameter $\diam X$ of a compact metric space $X$ if finite, and that there exist $x,y\in X$ such that $\diam X=|xy|$.
\end{exe}

\begin{exe}
Prove that every compact metric space is separable.
\end{exe}

\begin{exe}[Baire's theorem]
A subset of a topological space is called \emph{nowhere dense\/} if its closure has empty interior. Prove that a complete metric space cannot be covered by at most countably many nowhere dense subsets. Moreover, the complement of the union of at most countably many nowhere dense subsets is everywhere dense.
\end{exe}

\begin{exe}
Prove that a compact metric space $X$ cannot be isometrically mapped to a subspace $Y\ss X$ such that $Y\ne X$. In other words, each isometric mapping $f\:X\to X$ for a compact metric space $X$ is surjective.
\end{exe}

\begin{exe}
Let $X$ be a compact metric space and $f\:X\to X$ be a mapping. Prove that
\begin{enumerate}
\item if $f$ is surjective and nonexpanding, then $f$ is an isometry\/\rom;
\item if $\bigl|f(x)f(y)\bigr|\ge|xy|$ for all $x,y\in X$, then $f$ is an isometry.
\end{enumerate}
\end{exe}

\vfill\eject
 \chapter{Curves in Metric Spaces.}
 \markboth{\chaptername~\thechapter.~Curves in Metric Spaces.}%
          {\chaptername~\thechapter.~Curves in Metric Spaces.}

\begin{plan}
Curves in a topological space, parameter of a curve, reparametrization, polygonal line in a metric space, its edges, the length of the edge, the length of the polygonal line, the length of a curve in a metric space, rectifiable curves, properties of the length functional,  intrinsic metric, generalized intrinsic pseudometric, maximal and minimal generalized pseudometrics, minimum of generalized intrinsic pseudometrics, quotients of generalized intrinsic pseudometric spaces, Hopf-Rinow condition, Hopf-Rinow theorem Part 1, convergence and uniform convergence in terms of the corresponding product spaces, limits of sequences of Lipschitz mappings, arc-length and uniform curves, reparametrizations, uniform reparametrizations, Arzela-Ascoli theorem, shortest curves and geodesics, existences theory for shortest curves, geodesic metric space, midpoints and $\e$-midpoints, existence of shortest curves in term of midpoints, intrinsic metrics and $\e$-midpoints.
\end{plan}

In this section we discuss some results related to the geometry of curves in metric spaces.

Recall that \emph{a curve\/} in a topological space $X$ is any continuous mapping $\g\:[a,b]\to X$ from a segment $[a,b]\ss\R$ with the standard topology; the variable $t\in[a,b]$ is called \emph{the parameter of the curve $\g$}, and the curve $\g$ is sometimes written in the form $\g(t)$.

Each homeomorphism $\v\:[c,d]\to[a,b]$ generates a new curve $\g\c\v\:[c,d]\to X$, about which we say that it is obtained from $\g$ \emph{by replacement $\v$ of the parameter $t$ with the parameter $s\in[c,d]$}. Moreover, if there is no misunderstanding, instead of the curve $(\g\c\v)(s)$ we simply write $\g(s)$. Such replacement $\v$ is also called \emph{a reparametrization}.

Note that each reparametrization is a strictly monotonic continuous function. If the function $\v$ grows, then we say that $\v$ \emph{reverses the direction}, otherwise that it \emph{changes the direction}.

\section{Rectifiable curves}\label{sec:reacifible_curves}
\markright{\thesection.~Rectifiable curves}
Let $X$ be a metric space. A finite sequence $L=(A_0,\ldots,A_n)$ of points in the space $X$ we called \emph{a polygonal line in $X$}; moreover, the pairs $(A_{i-1},A_i)$ will be called \emph{edges of the polygonal line $L$}, and the numbers $|A_{i-1} A_i|$ \emph{the lengths\/} of these edges. The sum of the lengths of all these edges we call \emph{the length of the polygonal line $L$} and denote by $|L|$.

Let $\g\:[a,b]\to X$ be an arbitrary curve. For each partition $\xi=(a=t_0<t_1<\cdots<t_m=b)$, consider the corresponding polygonal line $L_\g(\xi)=\bigl(\g(t_0),\ldots,\g(t_m)\bigr)$ (such polygonal lines will be called \emph{inscribed in the curve $\g$\/}), then the value
$$
|\g|=\sup\Bigl\{\bigl|L_\g(\xi)\bigr|:\text{$\xi$ is a partition of the segment $[a,b]$}\Bigr\}
$$
is called \emph{the length of the curve $\g$}. A curve $\g$ is called \emph{rectifiable\/} if $|\g|<\infty$.

Let us give some examples of rectifiable curves.

\begin{examp}\label{examp:LipschitzCurveIsRectifiable}
Each $C$-Lipschitz curve $\g\:[a,b]\to X$ is rectifiable, because for any partition $\xi$ of the segment $[a,b]$ we have $|L_\g(\xi)|\le C(b-a)$ and, therefore, $|\g|\le C(b-a)<\infty$.
\end{examp}

Denote by $\Om(X)$ the family of all curves in a metric space $X$, and by $\Om_0(X)\ss\Om(X)$ the subfamily of all rectifiable curves. Note that on $\Om(X)$ and on $\Om_0(X)$ there are defined
\begin{enumerate}
\item \emph{the restriction\/} of each curve $\g\:[a,b]\to X$ to each subsegment $[c,d]\ss[a,b]$;
\item \emph{the gluing $\g_1\cdot\g_2$} of those pairs of curves $\g_1\:[a,b]\to X$, $\g_2\:[b,c]\to X$ for which $\g_1(b)=\g_2(b)$, namely, $(\g_1\cdot\g_2)\:[a,c]\to X$ is the curve whose restrictions to $[a,b]$ and $[b,c]$ coincide with $\g_1$ and $\g_2$, respectively;
\item \emph{the reparametrization\/} and equivalence identifying curves that differ by parameterization.
\end{enumerate}

The following proposition describes some properties of the length of a curve.

\begin{prop}\label{prop:length-functional}
Let $X$ be an arbitrary metric space. Then
\begin{enumerate}
\item{\bf generalized triangle inequality:}\label{prop:length-functional:0} if $\g\in\Om(X)$ joins the points $x,y\in X$, then $|\g|\ge|xy|$\rom;
\item{\bf additivity:}\label{prop:length-functional:1} if $\g=\g_1\cdot\g_2$ is the gluing of curves $\g_1,\g_2\in\Om(X)$ then $|\g|=|\g_1|+|\g_2|$\rom;
\item{\bf continuity:}\label{prop:length-functional:2} for any $\g\in\Om_0(X)$, $\g\:[a,b]\to X$, the function $f(t)=\bigl|\g|_{[a, t]}\bigr|$ is continuous\/\rom;
\item{\bf independence from parameter:}\label{prop:length-functional:3} for each $\g\in\Om(X)$, $\g\:[a,b]\to X$, and reparametrization $\v\:[c,d]\to[a,b]$, it holds $|\g|=|\g\c\v|$\rom;
\item{\bf consistency with topology:}\label{prop:length-functional:4} for each $x\in X$, $\e> 0$, $y\in X\sm U_\e (x)$, and a curve $\g\in\Om(X)$ joining $x$ and $y$, it holds $|\g|\ge\e$\rom;
\item{\bf lower semicontinuity:}\label{prop:length-functional:5} for any sequence $\g_n\in\Om_0(X)$ that converges pointwise to some $\g\in\Om_0(X)$, we have
$$
|\g|\le\liminf_{n\to\infty}|\g_n|.
$$
\end{enumerate}
\end{prop}

\begin{proof}
Only the items~(\ref{prop:length-functional:2}) and~(\ref{prop:length-functional:5}) are nontrivial, we prove them.

(\ref{prop:length-functional:2}) Choose an arbitrary $t\in[a,b]$ and show that for any $\e>0$ there is $\dl>0$ such that for all $s\in[a,b]\cap(t-\dl,t+\dl)$ the inequality $\bigl|f(t)-f(s)\bigr|<\e$ holds. Put $\ell=|\g|$. By definition, there exists a partition $\xi$ of the segment $[a,b]$ such that $\ell-\e/2<\bigl|L_\g(\xi)\bigr|\le\ell$. If $t\not\in\xi$, add it to $\xi$ (we denote the resulting partition by the same letter). It is clear that for the resulting partition, $\ell-\e/2 <\bigl|L_\g(\xi)\bigr|\le\ell$ is still satisfied.

For $\dl_1$ we take the distance from $t$ to the nearest element of the partition $\xi$, other than $t$. Since subdivisions of the partition $\xi$ can change the length of the polygonal line $L_\g(\xi)$ only within $(\ell-\e/2,\ell]$, then for each $s\in[a,b]\cap(t-\dl_1,t+\dl_1)$ the length $\ell_{ts}=\bigl|f(t)-f(s)\bigr|$ of the fragment of the curve $\g$ between points $\g(t)$ and $\g(s)$ differs from $\bigl|\g(t)\g(s)\bigr|$ by less than $\e/2$. On the other hand, since the map $\g$ is continuous, there exists $\dl_2>0$ such that for all $s\in[a,b]\cap(t-\dl_2,t+\dl_2)$ we have $\bigl|\g(t)\g(s)\bigr|<\e/2$. It remains to put $\dl=\min\{\dl_1,\dl_2\}$.

(\ref{prop:length-functional:5}) Choose an arbitrary $\e>0$ and show that for sufficiently large $n$ the inequality $|\g|\le|\g_n|+\e$ holds, thus $|\g|\le\liminf_{n\to\infty}|\g_n|+\e$ and, due to the arbitrariness of $\e$, we get what is required.

So, let $\e>0$ be fixed. Choose a partition $\xi=(a=t_0<t_1<\cdots<t_m=b)$ of the segment $[a,b]$ such that $|\g|-\bigl|L_\g(\xi)\bigr|<\e/2$. There is $N$ such that for any $n>N$ and all $i$ the inequality $\bigl|\g (t_i)\g_n(t_i)\bigr|<\frac{\e}{4m}$ holds. This immediately implies that
$$
\bigl|\g(t_{i-1})\g(t_i)\bigr|<\bigl|\g_n(t_{i-1})\g_n(t_i)\bigr|+\frac{\e}{2m},
$$
therefore $\bigl|L_\g(\xi)\bigr|<\bigl|L_{\g_n}(\xi)\bigr|+\e/2$. Thus,
$$
|\g|<\bigl|L_\g(\xi)\bigr|+\e/2 <|L_{\g_n}(\xi)|+\e/2+\e/2\le|\g_n|+\e,
$$
as required.
\end{proof}

\begin{prb}
Prove the remaining items of Proposition~\ref{prop:length-functional}.
\end{prb}

\begin{prb}
Will the items (\ref{prop:length-functional:2}) and (\ref{prop:length-functional:5}) of Proposition~\ref{prop:length-functional} remain true if we change $\Om_0(X)$ to $\Om(X)$?
\end{prb}

\begin{prb}
Show that a piecewise smooth curve in $\R^n$ is Lipschitzian with a Lipschitz constant equal to the maximum modulus of the velocity vector of the curve, therefore every such curve is rectifiable.
\end{prb}

Let $X$ be a metric space in which any two points are connected by a rectifiable curve. Then for any $x,y\in X$ the quantity
$$
d_{in}(x,y)=\inf\bigl\{|\g|:\text{$\g$ is a curve joining $x$ and $y$}\bigr\}
$$
is finite.

\begin{prb}
Let $X$ be a metric space in which any two points are connected by a rectifiable curve.
\begin{enumerate}
\item Prove that $d_{in}$ is a metric.
\item Denote by $\t$ the metric topology of $X$ w.r.t\. the initial metric on $X$, by $\t_{in}$ the metric topology w.r.t\. $d_{in}$, by $X_{in}$ the set $X$ with metric $d_{in}$ and topology $\t_{in}$. Show that $\t\ss\t_{in}$. In particular, if a mapping $\g\:[a,b]\to X_{in}$ is continuous, then the mapping $\g\:[a,b]\to X$ is continuous as well.
\item Construct an example when $\t\ne\t_{in}$.
\item Prove that for each rectifiable curve $\g\:[a,b]\to X$ the mapping $\g\:[a,b]\to X_{in}$ is continuous.
\item Denote by $|\g|_{in}$ the length of a curve $\g\:[a,b]\to X_{in}$. Show that for each curve $\g\:[a,b]\to X$ which is also a curve in $X_{in}$, it holds $|\g|=|\g|_{in}$. Thus, the sets of rectifiable curves for $X$ and $X_{in}$ coincide, and each non-rectifiable curve in $X$ is either a non-rectifiable one in $X_{in}$, or the mapping $\g\:[a,b]\to X_{in}$ is discontinuous.
\item Construct an example of continuous mapping $\g\:[a,b]\to X$ such that the mapping $\g\:[a,b]\to X_{in}$ is not continuous. Notice that the curve $\g\:[a,b]\to X$ can not be rectifiable.
\end{enumerate}
\end{prb}

\begin{dfn}
If $d_{in}$ coincides with the original metric, then the original metric is called \emph{intrinsic}. A metric space with an intrinsic metric is also called \emph{intrinsic}.
\end{dfn}

\begin{prb}
Let $X$ be a metric space in which any two points are connected by a rectifiable curve. Prove that the metric $d_{in}$ is intrinsic.
\end{prb}

\begin{examp}\label{examp:TwoMetricsOnCircle}
Let $S^1$ be the standard circle on the Euclidean plane $\R^2$.
\begin{enumerate}
\item If for $x,y\in S^1$ we put $d(x,y)$ equal to the distance in $\R^2$ between these points, then the metric $d$ on $S^1$ is not intrinsic.
\item If, for $d(x,y)$, we choose the length of the smaller of the two arcs of the circle $S^1$ into which $x$ and $y$ divide it, then the resulting metric will be intrinsic.
\end{enumerate}
\end{examp}

\begin{rk}
If we allow generalized pseudometrics, then we can define $d_{in}$ and intrinsic metric not only for those metric spaces where each pair of points is connected by a rectifiable curve, but for generalized pseudometric spaces too. To do that, we need only put $\inf\0=\infty$, where $\inf$ is applied to subsets of $[0,\infty]$. Thus, if $x,y\in X$ cannot be joined by a curve, then we get $d_{in}(x,y)=\infty$. The same holds for $x$ and $y$ joined by non-rectifiable curves only.

Recall that in the previous chapter we introduced the equivalence relation $\sim_1$ such that $x\sim_1 y$ if and only if the generalized distance between these points equals $\infty$. Applying this equivalence to the generalized $d_{in}$, we get a partition of the space $X$ into metrics subspaces with finite $d_{in}$.
\end{rk}

\section{Maximal and minimal pseudometrics, quotients}
\markright{\thesection.~Maximal and minimal pseudometrics, quotients}

In Chapter~\ref{ch:TopMetr} we defined maximal pseudometric, see Construction~\ref{constr:MaxPseudometric}, and we demonstrated its relation with quotient distance. Now, let us define maximal and minimal pseudometrics for arbitrary families of generalized pseudometrics. Also, we apply these generalizations to investigation of intrinsic distances.

\begin{constr}[Maximal and minimal pseudometric for arbitrary family of pseudometrics]\label{constr:MaxMinPseudometricArbitrary}
Let $X$ be a set, and $\cD$ an arbitrary nonempty family of generalized pseudometrics on $X$. Then we consider the standard partial order on the set of all generalized pseudometrics on $X$, namely, $d_1\le d_2$ if $d_1(x,y)\le d_2(x,y)$ for all $x,y\in X$. By means of this partial order, we define $\inf\cD$ and $\sup\cD$ in the standard way.

Put $\bd_\cD(x,y)=\sup_{d\in\cD}d(x,y)$, then the same proof as for Lemma~\ref{lem:MaxPseudometricExistsAndUnique} can be carried out to obtain that $\bd_\cD$ is a generalized pseudometric and $\bd_\cD=\sup\cD$, thus, $\sup\cD$ exists for any nonempty $\cD$. We call $\bd_\cD$ \emph{the maximal generalized pseudometric for the family $\cD$}. Notice that Construction~\ref{constr:MaxPseudometric} is a particular case of the present one: $d_b=\sup\cD_b$.

Now, let us put $b(x,y)=\inf_{d\in\cD}d(x,y)$. This $b$ is not necessarily a generalized pseudometric (construct an example), however, the $d_b=\sup\cD_b$ is, and we denote it by $\ud_\cD$.

\begin{prop}
Under above notations, $\ud_\cD=\inf\cD$, and, thus, $\inf\cD$ exists for any nonempty $\cD$.
\end{prop}

\begin{proof}
Since for any $d'\in\cD_b$ and $d\in\cD$ we have $d'\le b\le d$, then $\ud_\cD\le\cD$, i.e., $\ud_\cD$ is a lower bound for $\cD$.

Now, let $d'$ be an arbitrary lower bound for $\cD$, then $d'\le b$ and, therefore, $d'\in\cD_b$. Thus $d'\le\ud_\cD$.
\end{proof}

We call $\bd_\cD$ \emph{the minimal generalized pseudometric for the family $\cD$}.

Also, we can define $\sup\cD$ and $\inf\cD$ for empty family $\cD$: the first is the zero pseudometric, and the last is the generalized pseudometric equal $\infty$ for any pair of distinct points.
\end{constr}

Below we will use the following result.

\begin{prb}\label{prb:CompareMetricsAndCurves}
Let $\r_1\le\r_2$ be generalized pseudometrics on a set $X$, and $Y$ be a topological space. Prove that each mapping $f\:Y\to X$, continuous w.r.t\. $\r_2$, is also continuous w.r.t\. $\r_1$, in particular, if $\g$ is a curve in $(X,\r_2)$, then $\g$ is also a curve in $(X,\r_1)$; moreover, if $\r'_1$ and $\r'_2$ denote the corresponding generalized intrinsic pseudometrics, then $\r'_1\le\r'_2$.
\end{prb}

\begin{prop}\label{prop:InfIntrinsicAsWell}
Let $X$ be a set, and $\cD$ an arbitrary family of intrinsic generalized pseudometrics on $X$. Then $\ud_\cD=\inf\cD$ is intrinsic.
\end{prop}

\begin{proof}
We put $b(x,y)=\inf_{d\in\cD}d(x,y)$, $\r=\ud_\cD=d_b$, and denote by $\r'\ge\r$ the generalized intrinsic pseudometric corresponding to $\r'$. Since $\r\le d$ for all $d\in\cD$, and all $d$ are intrinsic, then $\r'\le d$ for all $d\in\cD$ due to Problem~\ref{prb:CompareMetricsAndCurves}. Thus, $\r'\le b$ and, therefore, $\r'\in\cD_b$. However, $\r$ is maximal for the class $\cD_b$, so $\r\ge\r'$ and we get $\r=\r'$.
\end{proof}

In the previous chapter we considered the construction of quotient spaces. What can we say about the quotient space if the distance function of the initial one is intrinsic?

\begin{prop}
Let $X$ be a generalized pseudometric space whose distance function is intrinsic, and $\sim$ an arbitrary equivalence on $X$. Then the generalized pseudometric of $X/\!\!\sim$ is intrinsic as well.
\end{prop}

\begin{proof}
Define $d\:X\x X\to[0,\infty]$ by setting $d(x,y)=0$ if $x\sim y$, and $d(x,y)=\infty$ otherwise. It is easy to see that $d$ is a generalized pseudometric. Thus, the space $(X,d)$ is partitioned into subspaces $X_i$ such that $|X_iX_j|=\infty$ for $i\ne j$, and the distance between points of each $X_i$ vanishes. Thus, the topology on each $X_i$ is anti-discrete, so each mapping $[a,b]\to X_i$ is continuous, i.e., it is a curve, and such curve has zero length. This implies that the distance function of $X$ is intrinsic.

Let $\r$ be the original distance function of $X$. Let us put $\cD=\{\r,d\}$, then $\r_\sim=\inf\cD$, because $\min\{\r,d\}$ equals the function $b_\sim$ from the definition of the quotient space. It remains to apply Proposition~\ref{prop:InfIntrinsicAsWell}.
\end{proof}

\begin{prb}
Let $X$ be an arbitrary set covered by a family $\{X_i\}_{i\in I}$ of generalized pseudometric spaces. Denote the distance function on $X_i$ by $\r_i$, and consider the set $\cD$ of all generalized pseudometrics $d$ on $X$ such that for any $i$ and $x,y\in X_i$ it holds $d(x,y)\le\r_i(x,y)$. Extend each $\r_i$ to the whole $X$ by setting $\r'_i(x,y)=\infty$ if at least one of $x$, $y$ does not belong to $X_i$, and $\r'_i(x,y)=\r(x,y)$ otherwise (it is easy to see that each $\r'_i$ is a generalized pseudometric). Denote by $\cD'$ the set of all such $\r'_i$. Prove that $\sup\cD=\inf\cD'$, and if all $\r_i$ are intrinsic, then $\sup\cD$ is intrinsic as well.
\end{prb}

\begin{prb}
Let $\cD$ be a collection of generalized pseudometrics defined on the same set $X$, and $X_d$ for $d\in\cD$ denote the generalized pseudometric space $(X,d)$. Put $W=\sqcup_{d\in\cD}X_d$ and denote by $\r$ the generalized pseudometric of $W$. Define on $W$ an equivalence relation $\sim$ by identifying those points $x_d\in X_d$ and $x_{d'}\in X_{d'}$ which correspond to the same point $x$ of the set $X$. The equivalence class of these points $x_d$ and $x_{d'}$ we denote by $[x]$. Denote by $\r_\sim$ the quotient generalized pseudometric on $W/\!\!\sim$.  Define the mapping $\v\:W/\!\!\sim\to X$ as $\v\:[x]\to x$, then $\v$ is bijective, and $\r_\sim$ can be considered as a generalized pseudometric on $X$. Prove that $\r_\sim=\inf\cD$.
\end{prb}

\begin{prb}
Let $\r_1$ and $\r_2$ be intrinsic metrics on a set $X$. Suppose that these metrics generate the same topology, and that each $x\in X$ has a neighborhood $U^x$ such that the restrictions of $\r_1$ and $\r_2$ to $U^x$ coincide. Prove that $\r_1=\r_2$. Show that the condition ``$\r_1$ and $\r_2$ are intrinsic'' is essential.
\end{prb}

\section{Hopf--Rinow condition}
\markright{\thesection.~Hopf--Rinow condition}

General metric spaces can be geometrically very different from $\R^n$. For example, in discrete spaces, balls of nonzero radius can coincide with their centers. In particular, the distance from an arbitrary point to such a ball will be equal to the distance from this point to the center. In spaces with an internal metric, this does not occur.

\begin{thm}[Hopf--Rinow condition]\label{thm:DistToDistInnerMetric}
Let $X$ be a space with an intrinsic metric, $x,\,y\in X$, $x\ne y$, and $0<r\le|xy|$. Then
$$
\bigl|yU_r(x)\bigr|=|xy|-r.
$$
\end{thm}

\begin{rk}
For general metric spaces $X$, Theorem~\ref{thm:DistToDistInnerMetric} does not hold. For example, if $X=\{x,y\}$, $|xy|=1$, and $r=0.5$, then $U_r(x)=\{x\}$, $\bigl|yU_r(x)\bigr|=1\ne|xy|-r=0.5$.
\end{rk}

\begin{proof} [Proof of Theorem~$\ref{thm:DistToDistInnerMetric}$]
For any point $z\in U_r(x)$ we have $|yz|\ge|yx|-|zx|>|xy|-r$, therefore $\bigl|yU_r (x)\bigr|\ge|xy|-r$. Let us prove that the converse inequality also holds.

For each $0<\e<r$ we consider a rectifiable curve $\g\:[0,1]\to X$, $x=\g(0)$ and $y=\g(1)$, for which $|\g|\le|xy|+\e$. We define a continuous function $f(t)=\bigl|x\g(t)\bigr|$, $f(0)=0$, $f(1)=|xy|$, and choose an arbitrary $t_0$ such that $f(t_0)=r-\e$. We denote by $\g_1$ the part of the curve $\g$ between $0$ and $t_0$, and by $\g_2$ the remaining part of the curve $\g$. Then $|\g_1|\ge r-\e$ by Item~(\ref{prop:length-functional:0}) of Proposition~\ref{prop:length-functional}, so $|\g_2|\le|xy|-r+2\e$ and, by the same proposition, $\bigl|\g(t_0)y\bigr|\le|\g_2|\le|xy|-r+2\e$. However, $\g(t_0)\in U_r(x)$, therefore $\bigl|yU_r(x)\bigr|\le|xy|-r+2\e$. Since $\e$ is arbitrary, we obtain what is required.
\end{proof}

\begin{rk}
The Hopf--Rinow condition can also be satisfied in spaces whose metric is not intrinsic, for example, in the metric space $\Q$ of all rational numbers (with the standard distance function).
\end{rk}

We give some corollaries from Theorem~\ref{thm:DistToDistInnerMetric}. First we give a necessary definition.

Let $X$ be a metric space, $x\in X$, $r\ge0$. Note that a closed ball $B_r(x)$ is a closed set, but, generally speaking, different from the closure of the open ball $U_r(x)$: if, as in the above example, $X$ consists of two points $x$ and $y$ at the distance $1$, then $U_1(x)=\{x\}$, $B_1(x)=\{x,y\}$, $\overline{U_1(x)}=\{x\}\ne B_1(x)$. However, if the metric of the space $X$ is intrinsic, then Theorem~\ref{thm:DistToDistInnerMetric} immediately implies the following result.

\begin{cor}\label{cor:InnerMetricCLosureOfOpenBall}
Let $X$ be a space with an intrinsic metric. Then $B_r(x)=\overline{U_r(x)}$.
\end{cor}

\begin{proof}
A point $y$ is adherent for a ball $U_r(x)$ if and only if $|yU_r(x)|=0$, thus $|xy|\le r$, i.e., $y\in B_x(r)$ and, therefore, $\overline{U_r (x)}\ss B_r(x)$. Let us prove the reverse inclusion.

Let $y\in B_r(x)$. If $|xy|<r$, then $y\in U_r(x)\ss\overline{U_r(x)}$. If $|xy|=r$, then, by Theorem~\ref{thm:DistToDistInnerMetric}, we have $|yU_r(x)|=|xy|-r=0$, thus $y\in\overline{U_r(x)}$.
\end{proof}

The following result will be used in the proof of the first part of Hopf--Rinow Theorem.

\begin{cor}\label{cor:DiskExtantionAtNet}
Let $X$ be a space with intrinsic metric and $\e>0$. Then for each $\e$-net $S$ in the ball $B_r(x)\ss X$ and any $\dl'>\dl>0$ we have $B_{r+\dl}(x)\ss\cup_{s\in S}U_{\e+\dl'}(s)$, i.e., $S$ is $(\e+\dl')$-net for $B_{r+\dl}(x)$.
\end{cor}

\begin{proof}
For any point $y\in B_{r+\dl}(x)$ we have $|xy|\le r+\dl$, therefore either $y\in U_r(x)$ and, thus, $|yU_r(x)|=0$, or, by Theorem~\ref{thm:DistToDistInnerMetric}, $\bigl|yU_r(x)\bigr|\le\dl$ holds. Thus, for any $\dl'>\dl$ there exists $z\in U_r(x)\ss B_r(x)$ such that $|yz|<\dl'$. On the other hand, there exists $s\in S$ for which $U_\e(s)\ni z$, whence $|sy|\le|sz|+|zy|<\e+\dl'$, therefore $y\in U_{\e+\dl'}(s)$, as required.
\end{proof}

\section{Local compactness}
\markright{\thesection.~Local compactness}

\begin{dfn}
A metric space $X$ is called \emph{locally compact} if for every point $x\in X$ there exists $\e>0$ such that the closed ball $B_\e(x)$ is compact.
\end{dfn}

\begin{prb}
Prove that a metric space $X$ is locally compact if and only if for each point $x\in X$ there exists a neighborhood with compact closure.
\end{prb}

\begin{rk}
Unlike compactness, local compactness, even in combination with the intrinsic metric, does not guarantee the completeness of the metric space. An obvious example is an open ball in Euclidean space. Another example is the Euclidean space with a point removed.
\end{rk}

\begin{thm}[Hopf-Rinow, Part 1]\label{thm:Hopf-Rinow:1}
Let $X$ be a locally compact space with intrinsic metric. Then the space $X$ is complete if and only if every closed ball in $X$ is compact.
\end{thm}

\begin{proof}
Suppose first that each closed ball is compact. We prove the completeness. Consider an arbitrary fundamental sequence $x_1,\,x_2,\ldots$. Then there exists $r$ such that all $x_n$ are contained in $B_r(x_1)$. By Theorem~\ref{thm:CompleteCompact}, the ball $B_r(x_1)$ is complete, therefore the sequence $x_1,\,x_2,\ldots$ converges to some point $x\in B_r(x)\ss X$, as required.

Now let the space $X$ be complete. On $X$ we define a function $\r\:X\to[0,\infty]$ as follows:
$$
\r(x)=\sup\{r>0:\text{the ball $B_r(x)$ is compact}\}.
$$

\begin{lem}\label{lem:RhoEqualINfty}
Suppose that there exists a point $x_0\in X$ such that $\r(x_0)=\infty$. Then each ball $B_r(x)$ is compact and, therefore, $\r$ is identically equal to $\infty$.
\end{lem}

\begin{proof}
For every $x$ and $r>0$, the ball $B_r(x)$ is contained in some compact ball $B_{r'}(x_0)$, therefore, since the set $B_r(x)$ is closed, the ball $B_r(x)$ is also compact.
\end{proof}

Thus, it suffice to prove that there exists a point $x_0\in X$ such that $\r(x_0)=\infty$. Assume the contrary, i.e., that the function $\r$ is everywhere finite.

\begin{lem}
The function $\r$ is $1$-Lipschitz and, therefore, continuous.
\end{lem}

\begin{proof}
Otherwise, there exists $x,y\in X$ such that $\bigl|\r(x)-\r(y)\bigr|>|xy|$. To be definite, assume that $\r(x)\ge\r(y)$, thus $\r(x)>\r(y)+|xy|$, and if $\e>0$ is chosen in such a way that $\r(x)>\r(y)+2\e+|xy|$, then $B_{\r(y)+\e}(y)\ss B_{\r(x)-\e}(x)$, however, $B_{\r(x)-\e}(x)$ is compact, thus $B_{\r(y)+\e}(y)$ is compact as well, a contradiction with definition of $\r(y)$.
\end{proof}

\begin{lem}\label{lem:bBRxIsComplact}
Under the assumptions made, the ball $B_{\r(x)}(x)$ is compact for every $x$.
\end{lem}

\begin{proof}
Since the ball $B_{\r(x)}(x)$ is a closed subset of the complete space $X$, this ball is also complete. Therefore, by Theorem~\ref{thm:CompleteCompact}, it suffices to prove that for each $\e>0$ this ball contains a finite $\e$-net.

To do that, we choose $0<r<\r(x)$ such that $\dl:=\r(x)-r<\e/2$, then the ball $B_r(x)$ is compact, thus it contains a finite $(\e/2)$-net $S$. By Corollary~\ref{cor:DiskExtantionAtNet}, the set $S$ is $(\e/2+\e/2)$-net for the ball $B_{r+\dl}(x)=B_{\r(x)}(x)$.
\end{proof}

Since $\r$ is a continuous function, its restriction to the compact set $B_{\r(x)}(x)$ attains its minimum and, thus, this minimum is positive. We denote this minimum by $\e$, then, by Lemma~\ref{lem:bBRxIsComplact}, all the balls $B_\e(y)$, $y\in B_{\r(x)}(x)$ are compact. Let $S$ be a finite $(\e/2)$-net in $B_{\r(x)}(x)$, and $0<\dl<\e/2$, then, by Corollary~\ref{cor:DiskExtantionAtNet}, the set $S$ is $(\e/2+\e/2)$-net for the ball $B_{\r(x)+\dl}(x)$. In particular, $B_{\r(x)+\dl}(x)$ is contained in the set $\cup_{s\in S}B_\e(s)$ which is compact as a finite union of compact sets. Therefore, $B_{\r(x)+\dl}(x)$ is compact that contradicts to definition of the function $\r$.
\end{proof}

\begin{dfn}
A metric space in which every closed ball is compact is called \emph{proper} or \emph{boundedly compact}.
\end{dfn}

\begin{cor}
A metric space with an intrinsic metric is boundedly compact if and only if it is locally compact and complete.
\end{cor}

\begin{prb}
Show that a metric spaces is boundedly compact if and only if its compact subsets are exactly those subsets that are closed and bounded.
\end{prb}

\section{Lipschitz, convergence and uniform convergence}\label{sec:LipConvEquiConv}
\markright{\thesection.~Lipschitz, convergence and uniform convergence}

In this section, we state and prove a few useful technical results regarding the convergence of Lipschitz mappings.

Let $f_n\:X\to Y$ be a family of arbitrary mappings from a set $X$ to a metric space $Y$. We say that the sequence $f_n$ \emph{converges pointwise\/} to a mapping $f\:X\to Y$ if for each $x\in X$ the sequence $f_n(x)$ converges to $f(x)$. The sequence $f_n$ \emph{converges uniformly\/} to a mapping $f\:X\to Y$ if for any $\e>0$ there exists $N$ such that for every $n\ge N$ the inequality $\bigl|f(x)f_n(x)\bigr|<\e$ holds for all $x\in X$.

\begin{rk}
Let us introduce the convergences described above by representing the mappings $f_n$ as points of the set $\prod_{x\in X}Y$, which, recall, we defined as the family $Y^X$ of all mappings from $X$ to $Y$. Now we model the pointwise and uniform convergences by means of some topologies.

We start with the case of pointwise convergence. We define on $Y^X=\prod_{x\in X}Y$ the Tychonoff's topology, see Construction~\ref{constr:Tychonoff}.

\begin{prb}
Show that convergence in the Tychonoff topology of points $f_n$ to a point $f$ is equivalent to pointwise convergence of the mappings $f_n$ to the mapping $f$.
\end{prb}

To model uniform convergence, we first give some definitions. A mapping $f\:X\to Y$ is called \emph{bounded\/} if its image $f (X)$ is a bounded subset of $Y$. The family of all bounded mappings from $X$ to $Y$ we denote by $\cB(X,Y)$. We define the following distance function on $\cB(X,Y)$: $|fg|=\sup_{x\in X}\bigl|f(x)g(x)\bigr|$.

\begin{prb}
Prove that the distance function defined above on $\cB(X,Y)$ is a metric, and that the convergence in this metric of a sequence $f_n\in\cB(X,Y)$ to a point $f\in\cB(X,Y)$ is equivalent to the uniform convergence of the mappings $f_n$ to the mapping $f$.
\end{prb}

For arbitrary $f_n$ and $f$ we can also define such convergence by considering $|fg|=\sup_{x\in X}\bigl|f(x)g(x)\bigr|$ as generalized metric.

\begin{prb}
Prove that the generalized distance function $|fg|=\sup_{x\in X}\bigl|f(x)g(x)\bigr|$ defined on $Y^x$ is a generalized metric, and that the convergence in this generalized metric of a sequence $f_n\in Y^X$ to a point $f\in Y^X$ is equivalent to the uniform convergence of the mappings $f_n$ to the mapping $f$.
\end{prb}
\end{rk}

\begin{prop}\label{prop:ConvergenceLipschizMaps}
Let $X$ be compact, and $Y$ be arbitrary metric spaces, and let $f_n\:X\to Y$ be a sequence of $C$-Lipschitz mappings converging pointwise to some mapping $f\: X\to Y$. Then $f$ is a $C$-Lipschitz mapping, and the sequence $f_n$ converges to $f$ uniformly.
\end{prop}

\begin{proof}
To verify that the mapping $f$ is $C$-Lipschitz, it is sufficient to pass to the limit in the inequality $\bigl|f_n(x)f_n(x')\bigr|\le C\cdot|xx'|$ for arbitrary fixed $x,\,x'\in X$.

We now prove uniform convergence. Choose an arbitrary $\e>0$ and show that there exists $N$ such that for all $n> N$ and all $x\in X$ we have $\bigl|f(x)f_n(x)\bigr|<\e$.

Put $\dl=\e/(3C)$, and let $\{x_i\}\ss X$ be a finite $\dl$-net. We choose $N$ such that for all $n>N$ and all $i$ the inequality $\bigl|f (x_i) f_n (x_i)\bigr|<\e/3$ holds.

Fix an arbitrary $x\in X$. There is $i$ such that $|xx_i|<\dl$. Since $f_n$ and $f$ are $C$-Lipschitz, we conclude that $\bigl|f_n(x)f_n(x_i)\bigr|\le C\cdot|xx_i|<\e/3$ and, similarly, $\bigl|f(x)f(x_i)\bigr|<\e/3$, therefore
$$
\bigl|f(x)f_n(x)\bigr|\le
\bigl|f(x)f(x_i)\bigr|+\bigl|f(x_i)f_n(x_i)\bigr|+\bigl|f_n(x_i)f_n(x)\bigr|<\e/3+\e/3+\e/3=\e,
$$
as required.
\end{proof}

The following version of the previous statement is useful in studying curves.

\begin{cor}\label{cor:ConvergenceLipschizCurves}
Let $X$ be a metric space, and $\g_n\:[a, b]\to X$ be a sequence of $C$-Lipschitz curves converging pointwise to a mapping $\g\:[a,b]\to X$. Then $\g$ is a $C$-Lipschitz curve, and the sequence $\g_n$ converges to $\g$ uniformly.
\end{cor}

The proposition below can be proved similarly to Proposition~\ref{prop:ConvergenceLipschizMaps}.

\begin{prop}\label{prop:PartialConvergenceLipschizMaps}
Let $X$ be compact, $Y$ be an arbitrary metric spaces, and $f_n\:X\to Y$ be a sequence of $C$-Lipschitz mappings. Suppose that for some everywhere dense subset $Z\ss X$ the sequence $f_n|_Z$ converges pointwise. Then the sequence $f_n$ converges pointwise to some mapping $f\: X\to Y$ and, therefore, by Proposition~$\ref{prop:ConvergenceLipschizMaps}$, this convergence is uniform, and the mapping $f$ is $C$-Lipschitz.
\end{prop}

\begin{cor}\label{cor:PartialConvergenceLipschiCurves}
Let $\g_n\:[a, b]\to X$ be a sequence of $C$-Lipschitz curves in a metric space $X$. Suppose that for some everywhere dense subset $Z\ss[a,b]$ the sequence of mappings $\g_n|_Z$ converges pointwise. Then the sequence of curves $\g_n$ converges pointwise to some curve $\g\:[a,b]\to X$ and, therefore, by virtue of Corollary~$\ref{cor:ConvergenceLipschizCurves}$, this convergence is uniform, and the curve $\g$ is $C$-Lipschitz.
\end{cor}

\section{Arc-length and uniform curves}
\markright{\thesection.~Arc-length and uniform curves}

\begin{dfn}\label{dfn:ArcLengthAndUninform}
A curve $\g(s)$ and its parameter $s\in[a,b]$ are called \emph{natural\/} or \emph{arc-length}, if for any $a\le s_1\le s_2\le b$ it holds $\bigl|\g|_{[s_1,s_2]}\bigr|=s_2-s_1$. A curve $\g(t)$ and its parameter $t\in[a,b]$ are called \emph{uniform}, if there exists $\l\ge0$ such that for any $a\le t_1\le t_2\le b$ it holds $\bigl|\g|_{[t_1,t_2]}\bigr|=\l(t_2-t_1)$; the value $\l$ is called \emph{the velocity\/} or \emph{the speed of uniform $\g$}.
\end{dfn}

\begin{rk}
Let $\g(t)$, $t\in[a,b]$, be uniform curve.
\begin{enumerate}
\item If $a\ne b$, then its velocity $\l$ is uniquely determined by the equation $|\g|=\l(b-a)$. In particular, $|\g|=0$, i.e\., $\g$ is a constant mapping if and only if $\l=0$.
\item If $a=b$, then $\l$ can be arbitrary. However, it is natural to consider this case as the limiting one for constant mappings $\g$, where $\l=0$. So, to be definite, we make the following agreement: \textbf{if $a=b$ then $\l=0$}.
\end{enumerate}
\end{rk}

Thus, under the above agreement, the velocity of a uniform curve $\g$ vanishes if and only if $\g$ is a constant mapping. We call such curves \emph{degenerate}, and all the remaining curves \emph{nondegenerate}. So, each degenerate curve is uniform.

\begin{rk}
The following simple observations concern the relations between arc-length and uniform curves.
\begin{enumerate}
\item Each arc-length curve is uniform.
\item A degenerate curve $\g(s)$, $s\in[a,b]$, is arc-length if and only if $a=b$; in this case the arc-length curve has zero velocity.
\item A degenerate curve $\g(t)$, $t\in[a,b]$, with $a<b$ is uniform but not arc-length.
\item A nondegenerate uniform curve is arc-length if and only if its velocity equals $1$.
\end{enumerate}
\end{rk}

\begin{prop}\label{prop:EquiParLip}
If a curve $\g(t)$, $t\in[a,b]$, is uniform with the velocity $\l$, then the mapping $\g$ is $\l$-Lipschitz.
\end{prop}

\begin{proof}
For any $a\le t_1\le t_2\le b$ we have
$$
\bigl|\g(t_1)\g(t_2)\bigr|\le\bigl|\g|_{[t_1,t_2]}\bigr|=\l(t_2-t_1).
$$
\end{proof}

\begin{rk}
The following types of curves $\g\:[a,b]\to X$ cannot be reparameterized to arc-length curves:
\begin{enumerate}
\item not rectifiable $\g$ (otherwise, $b=\infty$);
\item $\g=\const$ when $a\ne b$;
\item more general, $\g$ containing \emph{stops}, i.e., when there exists $[\a,\b]\ss[a,b]$, $\a\ne\b$, such that $\g|_{[\a,\b]}=\const$.
\end{enumerate}
\end{rk}

Curves that do not contain stops are called \emph{non-stop\/} ones.

\begin{prb}
Let $\g$ be a curve in a metric space. Prove that
\begin{enumerate}
\item nondegenerate $\g$ can be reparameterized to an arc-length or, more generally, to a uniform one if and only if $\g$ is rectifiable and non-stop;
\item degenerate $\g$ can be reparameterized to an arc-length one if and only if its domain is singleton (indeed, such $\g$ is arc-length itself and, thus, it need not a reparametrization);
\item degenerate $\g$ is always uniform.
\end{enumerate}
\end{prb}

In~\cite{BurBurIvaEng3} it is proposed to extend the class of reparametrizations, namely, to consider monotone (not necessarily strictly monotone) surjective mappings between the domains of the curves. It turns out that with this definition of reparametrization, it is possible to introduce an arc-length parameter on any rectifiable curve.

\begin{dfn}
We say that curves $\g\:[a,b]\to X$ and $\bga\:[c,d]\to X$ are obtained from each other by \emph{a monotone reparametrization} if either there exists a monotone surjective mapping $\v\:[c,d]\to[a,b]$ such that $\bga=\g\c\v$, or there exists a monotone surjective mapping $\psi\:[a,b]\to[c,d]$ such that $\g=\bga\c\psi$.
\end{dfn}

\begin{rk}\label{rk:MonotoneRepLength}
It is easy to see that monotone reparametrization does not change the length.
\end{rk}

\begin{prb}\label{prb:ArcLengthParameter}
Prove that a curve $\g$ in a metric space can be monotonically reparameterized to an arc-length or, more generally, a uniform one if and only if $\g$ is rectifiable.

The reparameterized curve is unique upto the choice of its domain and direction. In arc-length case one can choose any segment of the length $|\g|$. In the uniform case the domain can be arbitrary nondegenerate segment for nondegenerate $\g$, and arbitrary segment for degenerate $\g$.
\end{prb}

\begin{rk}
An instructive example is a curve that is a parametrization of the segment $[0,1]\ss\R$ by the Cantor staircase. The Cantor staircase is a graph of a function $f\:[0,1]\to[0,1]$, the construction of which we will now describe.

At the points $0$ and $1$, we set the value of the function $f$ equal to $0$ and $1$, respectively. Next, we divide the segment $[0,1]$ into three equal parts and on the middle interval we set $f$ equal to $1/2$. The remaining two segments are again divided into three equal parts each, and on the middle intervals we assume that the function $f$ is equal to the arithmetic mean of its values at the nearest intervals where it is defined. Thus, on the left-most interval, the function $f$ is equal to $1/4$, and on the right-most interval it is equal to $3/4$. Continuing this process to infinity, we define a function $f$ on an everywhere dense subset of the segment $[0,1]$, which is the complement to the Cantor set. Extend $f$ to the remaining points of the segment $[0,1]$ by continuity (make sure that this can be done).

We now consider the Cantor staircase as a curve $f\:[0,1]\to\R$ on the Euclidean line. Note that the subset of the segment $[0,1]$, on which the point of this curve changes its position, is the Cantor set that has measure zero. So, this curve stops almost everywhere, however, its length equals $1$ and it can be reparameterized to an arc-length curve.
\end{rk}

\section{Arzela-Ascoli Theorem}
\markright{\thesection.~Arzela-Ascoli Theorem}

Developing the ideas from Section~\ref{sec:LipConvEquiConv}, we formulate and prove a variant of the famous Arzela--Ascoli theorem. First we give necessary definitions.

\begin{dfn}
Let $\g_n\:[a_n,b_n]\to X$ be a sequence of curves in a metric space $X$. We say that this sequence \emph{converges\/} (\emph{uniformly converges\/}) to a curve $\g\:[a,b]\to X$ if there exist curves $\bga_n\:[c,d]\to X$ and $\bga\:[c,d]\to X$ obtained from $\g_n$ and $\g$, respectively, by monotone reparametrization, such that the mappings $\bga_n$ converge (uniformly converge) to the mapping $\bga$.
\end{dfn}

\begin{thm}[Arzela--Ascoli]\label{thm:ArzelaAscoli}
Let $X$ be a compact metric space, and $\g_n$ be a sequence of curves in $X$. Suppose that the lengths of the curves $\g_n$ are uniformly bounded, i.e., there exists a real number $C$ such that $|\g_n|\le C$ for all $n$. Then in this sequence there is a subsequence that converges uniformly to a curve whose length is at most $C$.
\end{thm}

\begin{proof}
In virtue of Problem~\ref{prb:ArcLengthParameter}, the curves $\g_n$ can be monotonically reparameterized to uniform curves $\bga_n\:[0,1]\to X$ with speeds at most $C$. It follows from Remark~\ref{rk:MonotoneRepLength} and Proposition~\ref{prop:EquiParLip} that all the curves $\bga_n$ are $C$-Lipschitz.

Choose a countable everywhere dense subset $Z\ss[0,1]$, $Z=\{z_i\}_{i=1}^\infty$. The sequence $\bigl(\bga_n(z_1)\bigr)_{n=1}^\infty$ has a convergent subsequence $\bigl(\g^1_n(z_1)\bigr)_{n=1}^\infty$; the sequence $\bigl(\g^1_n(z_2)\bigr)_{n=1}^\infty$ has a convergent subsequence $\bigl(\g^2_n(z_2)\bigr)_{n=1}^\infty$, etc. Then the sequence $\bigl(\g^n_n(z_k)\bigr)_{n=1}^\infty=\bigl(\bga_{n_i}(z_k)\bigr)_{i=1}^\infty$ convergence for any $k$ (Cantor diagonal process). Let us put $f(z_k)=\lim_{i\to\infty}\bga_{n_i}(z_k)$, then $\bga_{n_i}|_Z\to f$.

By Corollary~\ref{cor:PartialConvergenceLipschiCurves}, the mappings $\bga_{n_i}\:[0,1]\to X$ converge uniformly to some $C$-Lipschitz curve $\bga\:[0,1]\to X$. By Item~(\ref{prop:length-functional:5}) of Proposition~\ref{prop:length-functional}, we have $|\bga|\le\liminf_{n_i\to\infty}|\bga_{n_i}|\le C$, as required.
\end{proof}

\section{Existence of shortest curves}
\markright{\thesection.~Existence of shortest curves}

We apply the previous results to investigation of curves of smallest length.

\begin{dfn}
A rectifiable curve in a metric space is called \emph{shortest\/} if its length is equal to the infimum of the lengths of all the curves joining its ends.
\end{dfn}

\begin{rk}
If $X$ is a space with an intrinsic metric, then a curve $\g$ in $X$ joining $x$ and $y$ is \emph{shortest} if and only if $|xy|=|\g|$.
\end{rk}

The following proposition is obvious.

\begin{prop}
A curve in a metric space is shortest if and only if each of its parts is a shortest curve.
\end{prop}

\begin{prb}
Prove that an arc-length curve $\g\:[a,b]\to X$ in a space $X$ with an intrinsic metric is shortest if and only if $\g$ is an isometric embedding.
\end{prb}

\begin{dfn}
A curve $\g\:[a,b]\to X$ in a metric space $X$ is called \emph{locally shortest} if for each $t\in[a,b]$ there exists an interval $(\a,\b)\ss\R$ containing $t$ such that $\g|_{[a,b]\cap [\a,\b]}$ is a shortest curve.
\end{dfn}

\begin{dfn}
A uniform locally shortest curve is called \emph{a geodesic}.
\end{dfn}

Arzela-Ascoli theorem, together with a few other previous propositions, implies the following result.

\begin{cor}\label{cor:CompactSpaceAndShortestCurves}
Any two points $x$ and $y$ of a compact metric space $X$ that are joined by a rectifiable curve are also joined by a shortest curve.
\end{cor}

\begin{proof}
Let $\ell$ be the infimum of the lengths of the curves joining $x$ and $y$. There is a sequence $\g_n$ for which $|\g_n|\to\ell$ and, thus, the lengths of $\g_n$ are uniformly bounded. Theorem~\ref{thm:ArzelaAscoli} implies that the sequence $\g_n$ contains a subsequence $\g_{n_i}$ which uniformly converging to some curve $\g$. By Item~(\ref{prop:length-functional:5}) of Proposition~\ref{prop:length-functional}, we have $|\g|\le\liminf_{i\to\infty}|\g_{n_i}|=\ell$, however, by the minimality of $\ell$, it holds $|\g|\ge\ell$, therefore $|\g|=\ell$ and, thus, $\g$ is a shortest curve.
\end{proof}

\begin{rk}
Corollary~\ref{cor:CompactSpaceAndShortestCurves} remains true if we change compact $X$ to a boundedly compact one (verify this).
\end{rk}

\begin{dfn}
A metric on $X$ is called \emph{strictly intrinsic\/} if any two points in $X$ are joined by a curve whose length is equal to the distance between these points. A metric space with strictly intrinsic metrics is called \emph{strictly intrinsic\/} or \emph{geodesic}.
\end{dfn}

Taking into account the Hopf--Rinow theorem, we obtain the following

\begin{cor}\label{cor:LocalCompactInnerCompleteIsStrictlyInner}
Each complete locally compact space with an intrinsic metric is a geodesic space.
\end{cor}

\section{Shortest curves and midpoints}
\markright{\thesection.~Shortest curves and midpoints}

\begin{dfn}
A point $z$ of a metric space is called \emph{a midpoint between or for points $x$ and $y$} of this space if $|xz|=|yz|=\frac12|xy|$.
\end{dfn}

\begin{thm}\label{thm:MiddlePoints}
Let $X$ be a complete metric space. Suppose that for each pair of points $x,\,y\in X$ there is a midpoint. Then $X$ is a geodesic space.
\end{thm}

\begin{proof}
Choose two arbitrary points $x$ and $y$ from $X$. We show that these points can be joined by a curve $\g\:[0,1]\to X$, for which $|\g|=|xy|$.

We will sequentially determine the map $\g$ for various points of the segment $[0,1]$. Put $\g(0)=x$ and $\g(1)=y$. Next, let $\g(1/2)$ be a midpoint between $x$ and $y$; $\g(1/4)$ be a midpoint between $\g(0)$ and $\g(1/2)$, and $\g (3/4)$ be a midpoint between $\g(1/2)$ and $\g(1)$. Continuing this process, we define $\g$ at all binary rational points of the segment $[0,1]$, i.e., at all points of the form $m/2^n$, where $0\le m\le2^n$ is an integer, and $n=0,1,\ldots$. Note that the set of all binary rational points of the segment $[0,1]$ is everywhere dense in $[0,1]$. In addition, it is easy to show that the constructed mapping $\g$ is $|xy|$-Lipschitz. The proof of the following technical lemma is left as an exercise.

\begin{lem}\label{lem:ContinuousExtensionOfLipschitz}
Let $Z$ be an everywhere dense subset of a metric space $X$, and $f\:Z\to Y$ be some $C$-Lipschitz mapping into a complete metric space $Y$. Then there exists a unique continuous mapping $F\:X\to Y$ extending $f$. Moreover, the mapping $F$ is also $C$-Lipschitz.
\end{lem}

\begin{prb}
Prove Lemma~\ref{lem:ContinuousExtensionOfLipschitz}.
\end{prb}

So, using Lemma~\ref{lem:ContinuousExtensionOfLipschitz}, we extend by continuity the mapping $\g$ onto the entire segment $[0,1]$, and we again denote the resulting $|xy|$-Lipschitz curve by $\g$. As noted in Example~\ref{examp:LipschitzCurveIsRectifiable}, it holds $|\g|\le|xy|(1-0)=|xy|$, from where, by virtue of Item~(\ref{prop:length-functional:0}) of Proposition~\ref{prop:length-functional}, we have $|\g|=|xy|$ and, therefore, $\g$ is a shortest curve.
\end{proof}

\begin{rk}
In a complete metric space, the property of a metric to be intrinsic is not sufficient for midpoints and shortest curves between any points to exist. Consider a countable family of segments $[0,1+1/n]$, $n\in\N$, each with the standard metric, and glue all their zeros at one point $A$, and at another point $B$ we glue all the other ends $1+1/n$. If $x$ and $y$ belong to different segments, say to $[0,1+1/n]$ and $[0,1+1/m]$, then we set the distance between $x$ and $y$ equal to $\min(x+y,1-x+1/n+1-y+1/n)$ (i.e., the intrinsic circle metric is considered on each pair of glued segments). Then the distance between $A$ and $B$ is $1$ and is not reached on any curve. In addition, there is no midpoints between $A$ and $B$.
\end{rk}

\begin{dfn}
A point $z$ of a metric space $X$ is called \emph{an $\e$-midpoint between or for points $x$ and $y$} of this space if $\bigl||xz|-\frac12|xy|\bigr|\le\e$ and $\bigl||yz|-\frac12|xy|\bigr|\le\e$.
\end{dfn}

\begin{thm}\label{thm:E-MiddlePoints}
Let $X$ be a complete metric space. Suppose that for each pair of points $x,\,y\in X$ and each $\e>0$, there is an $\e$-midpoint. Then the metric of $X$ is intrinsic.
\end{thm}

\begin{proof}
The proof is similar with the one of Theorem~\ref{thm:MiddlePoints}, however, now we find not strict midpoints, but approximate ones, making sure that the total ``spread'' is not large (we use the fact that $\sum_{i=1}^\infty\e/2^i=\e$).
\end{proof}

There are also converse obvious statements, even without assuming the completeness of the ambient space.

\begin{prop}\label{prop:e-middle-inner-metric}
In a space with an intrinsic\/ \(strictly intrinsic\/\) metric, for any two points and any $\e>0$ there is an $\e$-midpoint\/  \(a midpoint\/\), respectively.
\end{prop}

\begin{prb}
Prove Proposition~\ref{prop:e-middle-inner-metric}.
\end{prb}


\phantomsection
\renewcommand\bibname{References to Chapter~\thechapter}
\addcontentsline{toc}{section}{\bibname}
\renewcommand{\refname}{\bibname}

\vfill\eject

\section*{\Huge Exercises to Chapter~\thechapter}
\markright{Exercises to Chapter~\thechapter.}

\begin{exe}\label{exe:length-functional}
Let $X$ be an arbitrary metric space and $\Om(X)$ the family of all curves in $X$. Verify that
\begin{enumerate}
\item\label{exe:length-functional:1} if $\g\in\Om(X)$ joins the points $x,y\in X$, then $|\g|\ge|xy|$;
\item\label{exe:length-functional:2} if $\g=\g_1\cdot\g_2$ is the gluing of curves $\g_1,\g_2\in\Om(X)$ then $|\g|=|\g_1|+|\g_2|$;
\item\label{exe:length-functional:3} for each $\g\in\Om(X)$, $\g\:[a,b]\to X$, and reparametrization $\v\:[c,d]\to[a,b]$, it holds $|\g|=|\g\c\psi|$;
\item\label{exe:length-functional:4} for each $x\in X$, $\e>0$, $y\in X\sm U_\e (x)$ and the curve $\g\in\Om (X)$ joining $x$ and $y$, $|\g|\ge\e$ holds;
\item\label{exe:length-functional:5} is it true that for any $\g\in\Om(X)$, $\g\:[a,b]\to X$, the function $f(t)=\bigl|\g|_{[a, t]}\bigr|$ is continuous?
\item\label{exe:length-functional:6} is it true that for any sequence $\g_n\in\Om(X)$ converging pointwise to some $\g\in\Om(X)$, we have
$$
|\g|\le\liminf_{n\to\infty}|\g_n|?
$$
\end{enumerate}
\end{exe}

\begin{exe}
Show that the piecewise smooth curve in $\R^n$ is Lipschitzian with a Lipschitz constant equal to the maximum modulus of the velocity vector of the curve, therefore each such curve is rectifiable.
\end{exe}

\begin{exe}
Let $X$ be a metric space in which any two points are connected by a rectifiable curve.
\begin{enumerate}
\item Prove that $d_{in}$ is a metric.
\item Denote by $\t$ the metric topology of $X$ w.r.t\. the initial metric on $X$, by $\t_{in}$ the metric topology w.r.t\. $d_{in}$, by $X_{\in}$ the set $X$ with metric $d_{in}$ and topology $\t_{in}$. Show that $\t\ss\t_{in}$. In particular, if a mapping $\g\:[a,b]\to X_{in}$ is continuous, then the mapping $\g\:[a,b]\to X$ is continuous as well.
\item Construct an example when $\t\ne\t_{in}$.
\item Prove that for each rectifiable curve $\g\:[a,b]\to X$ the mapping $\g\:[a,b]\to X_{in}$ is continuous.
\item Denote by $|\g|_{in}$ the length of a curve $\g\:[a,b]\to X_{in}$. Show that for each curve $\g\:[a,b]\to X$ which is also a curve in $X_{in}$, it holds $|\g|=|\g|_{in}$. Thus, the sets of rectifiable curves for $X$ and $X_{in}$ coincide, and each non-rectifiable curve in $X$ is either a non-rectifiable one in $X_{in}$, or the mapping $\g\:[a,b]\to X_{in}$ is discontinuous.
\item Construct an example of continuous mapping $\g\:[a,b]\to X$ such that the mapping $\g\:[a,b]\to X_{in}$ is not continuous. Notice that the curve $\g\:[a,b]\to X$ can not be rectifiable.
\end{enumerate}
\end{exe}

\begin{exe}
Let $X$ be a metric space in which any two points are connected by a rectifiable curve. Prove that the metric $d_{in}$ is intrinsic.
\end{exe}

\begin{exe}
Let $\r_1\le\r_2$ be generalized pseudometrics on a set $X$, and $Y$ be a topological space. Prove that each mapping $f\:Y\to X$, continuous w.r.t\. $\r_2$, is also continuous w.r.t\. $\r_1$, in particular, if $\g$ is a curve in $(X,\r_2)$, then $\g$ is also a curve in $(X,\r_1)$; moreover, if $\r'_1$ and $\r'_2$ denote the corresponding generalized intrinsic pseudometrics, then $\r'_1\le\r'_2$.
\end{exe}

\begin{exe}
Let $X$ be an arbitrary set covered by a family $\{X_i\}_{i\in I}$ of generalized pseudometric spaces. Denote the distance function on $X_i$ by $\r_i$, and consider the set $\cD$ of all generalized pseudometrics $d$ on $X$ such that for any $i$ and $x,y\in X_i$ it holds $d(x,y)\le\r_i(x,y)$. Extend each $\r_i$ to the whole $X$ by setting $\r'_i(x,y)=\infty$ if at least one of $x$, $y$ does not belong to $X_i$, and $\r'_i(x,y)=\r(x,y)$ otherwise (it is easy to see that each $\r'_i$ is a generalized pseudometric). Denote by $\cD'$ the set of all such $\r'_i$. Prove that $\sup\cD=\inf\cD'$, and if all $\r_i$ are intrinsic, then $\sup\cD$ is intrinsic as well.
\end{exe}

\begin{exe}
Let $\cD$ be a collection of generalized pseudometrics defined on the same set $X$, and $X_d$ for $d\in\cD$ denote the generalized pseudometric space $(X,d)$. Put $W=\sqcup_{d\in\cD}X_d$ and denote by $\r$ the generalized pseudometric of $W$. Define on $W$ an equivalence relation $\sim$ by identifying those points $x_d\in X_d$ and $x_{d'}\in X_{d'}$ which correspond to the same point $x$ of the set $X$. The equivalence class of these points $x_d$ and $x_{d'}$ we denote by $[x]$. Denote by $\r_\sim$ the quotient generalized pseudometric on $W/\!\!\sim$.  Define the mapping $\v\:W/\!\!\sim\to X$ as $\v\:[x]\to x$, then $\v$ is bijective, and $\r_\sim$ can be considered as a generalized pseudometric on $X$. Prove that $\r_\sim=\inf\cD$.
\end{exe}

\begin{exe}
Let $\r_1$ and $\r_2$ be intrinsic metrics on a set $X$. Suppose that these metrics generate the same topology, and that each $x\in X$ has a neighborhood $U^x$ such that the restrictions of $\r_1$ and $\r_2$ to $U^x$ coincide. Prove that $\r_1=\r_2$. Show that the condition ``$\r_1$ and $\r_2$ are intrinsic'' is essential.
\end{exe}

\begin{exe}
Prove that a metric space $X$ is locally compact if and only if for each point $x\in X$ there exists a neighborhood with compact closure.
\end{exe}

\begin{exe}
Show that a metric spaces is boundedly compact if and only if its compact subsets are exactly those subsets that are closed and bounded.
\end{exe}

\begin{exe}
Let $X$ be an arbitrary set and $Y$ an arbitrary metric space. Consider the collection of sets of the form $\prod_{x\in X}V(x)\ss\prod_{x\in X} Y$, where $\bigl\{V(x)\bigr\}_{x\in X }$ is the family of nonempty open subsets of $Y$ such that for all $x\in X$, except for their finite number, $V(x)=Y$. Show that the family defined in this way forms a basis of a topology, and the convergence in this topology of points $f_n$ to a point $f$ is equivalent to pointwise convergence of the mappings $f_n$ to the mapping $f$.
\end{exe}

\begin{exe}
Let $X$ be an arbitrary set and $Y$ an arbitrary metric space. A mapping $f\: X\to Y$ is called \emph{bounded\/} if its image $f(X)$ is a bounded subset of $Y$. The family of all bounded mappings from $X$ to $Y$ we denote by $\cB(X,Y)$. We define the following distance function on $\cB(X,Y)$: $|fg|=\sup_{x\in X}\bigl|f(x)g(x)\bigr|$. Prove that the distance function defined above is a metric, and that the convergence in this metric of a sequence $f_n\in\cB(X,Y)$ to some $f\in\cB(X,Y)$ is equivalent to uniform convergence of the mappings $f_n$ to the mapping $f$.
\end{exe}

\begin{exe}
Let $X$ be an arbitrary set and $Y$ an arbitrary metric space. Define the following generalized distance function on $Y^X$: $|fg|=\sup_{x\in X}\bigl|f(x)g(x)\bigr|$. Prove that the generalized distance function defined above is a generalized metric, and that the convergence in this generalized metric of a sequence $f_n\in Y^X$ to some $f\in Y^X$ is equivalent to uniform convergence of the mappings $f_n$ to the mapping $f$.
\end{exe}

\begin{exe}
Let $\g$ be a curve in a metric space. Prove that
\begin{enumerate}
\item nondegenerate $\g$ can be reparameterized to an arc-length or, more generally, to a uniform one if and only if $\g$ is rectifiable and non-stop;
\item degenerate $\g$ can be reparameterized to an arc-length one if and only if its domain is singleton (indeed, such $\g$ is arc-length itself and, thus, it need not a reparametrization);
\item degenerate $\g$ is always uniform.
\end{enumerate}
\end{exe}

\begin{exe}
Prove that a curve $\g$ in a metric space can be monotonically reparameterized to an arc-length or, more generally, a uniform one if and only if $\g$ is rectifiable.

The reparameterized curve is unique upto the choice of its domain and direction. In arc-length case one can choose any segment of the length $|\g|$. In the uniform case the domain can be arbitrary nondegenerate segment for nondegenerate $\g$, and arbitrary segment for degenerate $\g$.
\end{exe}

\begin{exe}
Prove that an arc-length curve $\g\:[a,b]\to X$ in a space $X$ with an intrinsic metric is shortest if and only if $\g$ is an isometric embedding.
\end{exe}

\begin{exe}
Let $Z$ be an everywhere dense subset of a metric space $X$, and $f\:Z\to Y$ be some $C$-Lipschitz map into a complete metric space $Y$. Then there exists a unique continuous mapping $F\:X\to Y$ extending $f$. Moreover, the mapping $F$ is also $C$-Lipschitz.
\end{exe}

\begin{exe}
Show that in a space with an intrinsic (strictly intrinsic) metric, for any two points and any $\e>0$ there is an $\e$-midpoint (a midpoint), respectively.
\end{exe}

\vfill\eject
 \chapter{Extreme graphs and networks.}
 \markboth{\chaptername~\thechapter.~Extreme graphs and networks.}%
          {\chaptername~\thechapter.~Extreme graphs and networks.}

\begin{plan}
Simple graphs, finite graphs, vertices, edges, isomorphism of graphs, adjacency, incidence, neighborhood of a vertex, subgraph, spanning subgraph, complete graph, subgraph generated by vertices, subgraph generated by edges, walk, degenerate and non-degenerated walks, open and closed walks, trail, path, circuit, cycle, connected graph, components of a graph, forest, tree, weighted graph, the weight of a subgraph, the weight of trail, the weight of walk, operations on graphs, union, disjoint union, intersection, difference, deleting edges, deleting vertices, quotient graphs, quotient by an edge, splitting a vertex, splitting off a vertex, graphs in metric spaces, the length of an edge, the length of a graph, the length of minimum spanning tree, minimum spanning tree, the length of Steiner minimal tree, Steiner minimal tree, the length of minimal filling, minimal filling, $\mst$-spectrum of finite metric space, calculation of $\mst$-spectrum in terms of partitions, graphs with boundaries, boundary (fixed) vertices, interior (movable) vertices, networks, parameterizing graphs of networks, boundary of a network, the length of a network, splitting and splitting off for networks, full Steiner tree, Steiner minimal trees existence in boundedly compact metric spaces.
\end{plan}

In this section, we collect information about various kinds of extreme graphs and networks. We will consider two types of such graphs: minimal spanning trees and shortest trees, also called Steiner minimal trees.

\section{Necessary information from graph theory}\label{sec:graphPrelim}
\markright{\thesection.~Necessary information from graph theory}
We will consider only simple graphs, so in what follows by a \emph{graph\/} we mean a pair $G=(V,E)$ consisting of two sets $V$ and $E$, respectively called \emph{the set of vertices\/} and \emph{the set of edges\/} of the graph $G$; the elements from $V$ are called \emph{vertices}, and from $E$ are called \emph{edges\/} of the graph $G$. The set $E$ is a subset of the family of two-element subsets of $V$. If $V$ and $E$ are finite sets then the graph $G$ is called \emph{finite}.

It is convenient to use the following notation:
\begin{itemize}
\item if $\{v,w\}\in E$ is an edge of the graph $G$, then we will write it in the form $vw$ or $wv$; we will also say that the edge $vw$ \emph{joins the vertices $v$ and $w$}, and that $v$ and $w$ are \emph{the vertices of the edge $vw$};
\item if the sets $V$ and $E$ are not explicitly indicated, and only the notation for the graph $G$ is introduced, then the set of vertices of this graph is usually written as $V(G)$, and the set of edges is denoted by $E(G)$.
\end{itemize}

Recall some concepts from the graph theory. Graphs $G=(V,E)$ and $H=(W,F)$ are called \emph{isomorphic\/} if there exists a bijective map $f\:V\to W$ such that $uv\in E$ if and only if $f(u)f(v)\in F$. Such a mapping $f$ is called an \emph{isomorphism of the graphs $G$ and $H$}. Isomorphic graphs are often identified and, therefore, are not distinguished.

Two vertices $v,w\in V(G)$ are called \emph{adjacent\/} if $vw\in E(G)$. Two different edges $e_1,e_2\in E(G)$ are called \emph{adjacent\/} if they have a common vertex, i.e., if $e_1\cap e_2\ne\0$. Each edge $vw\in E(V)$ and its vertex, i.e., $v$ or $w$, \emph{are incident\/} to each other. The set of vertices of a graph $G$ adjacent to a vertex $v\in V$ is called the \emph{neighborhood of the vertex $v$} and denoted by $N_v$. The cardinal number of edges incident to a vertex $v$ is called the \emph{degree of the vertex $v$} and is denoted by $\deg v$, so $\deg v=\#N_v$.

A \emph{subgraph\/} of a graph $G=(V,E)$ is each graph $H=(W,F)$ provided that $W\ss V$ and $F\ss E$. The fact that a graph $H$ is a subgraph of a graph $G$ will be written as $H\ss G$. If $W=V$ then the subgraph $H\ss G$ is called \emph{spanning}.

On the set of all graphs whose vertex sets lie in a given set $V$, the inclusion relation $\ss$ defines a partial order. The smallest element in this order is the empty graph $(\0,\0)$; the greatest one is called \emph{the complete graph on $V$}, which we denote by $K(V)$: in $K(V)$ each pair of vertices are joined by an edge. This partial order induces the one on the set of all subgraphs of a graph $G=(V,E)$: now the smallest element is again the empty graph $(\0,\0)$, but the greatest one is $G$.

For each $W\ss V$ we define \emph{the subgraph $G(W)$ of the graph $G$ generated by $W$}: its set of vertices coincides with $W$, and its set of edges consists of all $e\in E$ that connect the vertices from $W$. In other words, $G(W)$ is maximal among subgraphs of $G$ whose vertex sets coincides with $W$.

Now, we define a similar construction by interchanging vertices and edges. Namely, for $F\ss E$ we define \emph{the subgraph $G(F)$ of the graph $G$ generated by $F$}: its set of edges coincides with $F$, and its set of vertices is the collection of all vertices of $G$ incident to edges from $F$. In what follows we also apply this construction.

A finite sequence $\g=(v_0=v,v_1,\ldots,v_k=w)$ of vertices of a graph $G$ is called a \emph{walk of length $k$ joining $v$ and $w$} if for every $i=1,\ldots,k$ the vertices $v_{i-1}$ and $v_i$ are adjacent, and the edges $e_i=v_{i-1}v_i$ are called the \emph{edges of the walk $\g$}. A walk containing at least one edge is called \emph{nondegenerate}, and not containing is called \emph{degenerate}. The walt is called \emph{closed\/} if $v_0=v_n$, and it is called \emph{open\/} otherwise. \emph{A trail\/} is a walk with no repeated edges. \emph{A path\/} is an open trail with no repeated vertices. \emph{A circuit\/} is a closed trail. \emph{A cycle\/} is a circuit with no repeated vertices.

A graph $G$ is called \emph{connected}, if each pair of its vertices are joined by a walk. Maximal (by inclusion) connected subgraphs of a graph $G$ are called \emph{components\/} of $G$. A graph without cycles is called \emph{a forest}, and a connected forest is called \emph{a tree}.

\emph{A weighted graph\/} is a graph $G=(V,E)$ equipped with \emph{a weight function $\om\:E\to[0,\infty)$} (sometimes it is useful to consider more general weight functions, for instance, with possibility of negative values or $\infty$). Sometimes we denote such weighted graph as $(V,E,\om)$ or $(G,\om)$. \emph{The weight $\om(H)$ of a subgraph $H\ss G$} is the sum of the weights of edges from this subgraph: $\om(H)=\sum_{e\in E(H)}\om(e)$. We can extend this definition to trails, in particular, to  paths, circuits and cycles, considering them as subgraphs of $G$. In the case of the walk $\g=(v_0=v,v_1,\ldots,v_k=w)$, its weight is defined as the sum of weights of its consecutive edges: $\om(\g)=\sum_{i=1}^n\in(v_{i-1}v_i)$. For graphs without weight functions these notions are defined as well by assigning the weight $1$ to each edge by default.

\begin{rk}
As in the case of metric spaces, we sometimes won't explicitly denote the weight function. Instead of that, when we speak about weighted graph $G$, the weights of all the objects $x$ related to such $G$ we denote by $|x|$, for example, for $e\in E$ by $|e|$ we mean the weight of this edge, and $H\ss G$ by $|H|$ we mean the weight of $H$, etc.
\end{rk}

\subsection{Some operations on graphs}
Let $H_1=(W_1,F_1)$ and $H_2=(W_2,F_2)$ be subgraphs of a graph $G$. Then the following subgraphs are defined:
\begin{itemize}
\item \emph{the union\/} $H_1\cup H_2=(W_1\cup W_2,F_1\cup F_2)$;
\item \emph{the disjoint union\/}: if $W_1\cap W_2=\0$, then to emphasize this fact, instead of $H_1\cup H_2$ we write $H_1\sqcup H_2$;
\item \emph{the intersection\/} $H_1\cap H_2=(W_1\cap W_2,F_1\cap F_2)$;
\item \emph{the difference\/} $H_1\sm H_2=H_1(W_1\sm W_2)$.
\end{itemize}

\begin{rk}
We can define the operations described above on any graphs $H_i$, not only on subgraphs of a graph. To reduce these definitions to the previous ones, we consider $H_i$ as subgraphs of the graph $K(W_1\cup W_2)$.
\end{rk}

\begin{prb}
Show that each forest is the disjoint union of trees that are components of this forest.
\end{prb}

If $G=(V,E)$ is a graph, and $F$ is a set possibly intersecting $E$, then the operation of \emph{deleting the set of edges $F$ from the graph $G$} produces the graph $G\sm^e F:=(V,E\sm F)$. If $F=\{e\}$ then instead of $G\sm^e\{e\}$ we will write $G\sm^e e$. The operation of \emph{deleting a set of vertices $W$ from the graph $G$} produces the graph $G\sm^v W:=G(V\sm W)$. If $W=\{w\}$ then instead of $G\sm^v\{w\}$ we will write $G\sm^v w$. If it is clear that $F$ refers to edges, or $W$ refers to vertices, we write simplified $G\sm F$ or $G\sm W$, respectively.

Using the operation of deleting edges, we define \emph{the complement of a graph $G=(V,E)$} or, in other words, \emph{the graph dual to $G$} to be the graph $\bG=K(V)\sm E$. Thus, the dual graph $\bG$ has the same set of vertices $V$, and its edges are exactly those edges of the complete graph on $V$ that were absent in the original graph $G$.

Another useful operation for us produces \emph{a quotient graph $G=(V,E)$}: let $\sim$ be an equivalence relation on $V$, and $V=\sqcup_{i\in I}V_i$ the partition into classes of this equivalence. We put $V/\fsim=\{V_i\}$, and as $E/\fsim$ we take the set of pairs $V_iV_j$, $V_i\ne V_j$ for which there exist $v_i\in V_i$, $v_j\in V_j$, such that $v_iv_j\in E$. By \emph{the quotient graph $G/\fsim$} we call the graph $(V/\fsim,\,E/\fsim)$. An important particular case of this operation creates \emph{the quotient of $G$ by an edge $e=vw\in E$}: the result is the graph $G/\fsim$ for the equivalence relation identifying the vertices $v$ and $w$. We denote this quotient graph by $G/e$.

The following notation and concepts are also useful: the equivalence class containing a given vertex $v$ will be denoted by $[v]$; the mapping $\pi\:V\to V/\fsim$, $\pi\:v\mapsto[v]$, is called \emph{the canonical projection}.

\begin{prb}
Let $G=(V,E)$ be a connected graph, and $\sim$ an arbitrary equivalence relation on the set $V$. Show that the graph $G/\fsim$ is connected.
\end{prb}

\begin{prb}\label{prb:factorTreeIsTree}
Let $G=(V,E)$ be an arbitrary tree, and $\sim$ be an equivalence relation on the set $V$ such that for each class $V_i$ of this equivalence the subgraph $G(V_i)$ is a tree. Show that then $G/\fsim$ is a tree.
\end{prb}

In some cases, the following operations are inverse to the quotient by an edge. We define two such operations: splitting a vertex of degree greater than or equal to $4$, and splitting off some vertex of degree $1$ from a vertex of degree greater than or equal to $2$.

So, let $G=(V,E)$ be a graph, $v\in V$, $\deg v\ge4$. We partition the neighborhood $N_v$ of the vertex $v$ into two sets $V_1$ and $V_2$, each of which contains at least two vertices. Consider the graph $G\sm v$, add to its vertex set $V\sm\{v\}$ two elements $w_1$ and $w_2$ not contained in $V\sm\{v\}$, and, to the set of edges, all pairs of the form $w_1v_1$, $v_1\in V_1$, $w_2v_2$, $v_2\in V_2$, as well as the pair $w_1w_2$. We call the obtained graph \emph{the result of splitting the vertex $v$}, and the edge $w_1w_2$ \emph{the splitting edge}. It is clear that the graph obtained from $G$ by the composition of splitting a vertex and the quotient by the corresponding splitting edge, is isomorphic to $G$ (just in this sense, the splitting is inverse to the quotient operation).

To determine the splitting off a vertex of degree $1$ from $v$, add to $V$ an element $w$ not contained in $V$, and add the edge $vw$ to $E$. The obtained graph is called \emph{the result of splitting off the vertex $w$ from the vertex $v$}, and the edge $wv$ is called \emph{the splitting edge}. It is clear that the quotient by the splitting edge is isomorphic to the original graph (in this sense, splitting off is also inverse to the quotient operation).

\section{Graphs and optimization problems}
\markright{\thesection.~Graphs and optimization problems}

Let $G=(V,E)$ be an arbitrary graph. We say that the graph $G$ is defined \emph{in a metric space $X$} if $V\ss X$. For every such graph, \emph{the length $|e|$ of its edge $e=vw$} is defined as the distance $|vw|$ between the ending vertices $v$ and $w$ of these edge, as well as \emph{the length $|G|$ of the graph $G$} itself as the sum of the lengths of all its edges. More generally, one can replace the metric space $X$ with a weighted complete graph $\bigl(K(X),\om\bigr)$; another possibility --- to consider some weighted graph $(H,\om)$ with $V(H)=X$ (not necessarily the complete one), such that $G$ is a subgraph of $H$. Let us note that each metric space $(X,\r)$ can be considered as a weighted complete graph, namely, as $\bigl(K(X),\r\bigr)$.

\subsection{Minimum spanning tree problem}
Let $M$ be a metric space. We consider $M$ as a weighted complete graph $K(M)$, and denote by $\cT(M)$ the set of all spanning trees in $K(M)$. Then we put
$$
\mst(M)=\inf_{T\in\cT(M)}|T|
$$
and call it \emph{the length of minimum spanning tree on $M$}. Each $T\in\cT(M)$ with $|T|=\mst(M)$ is call \emph{a minimum spanning tree on $M$}. The set of all minimum spanning trees on $M$ is denoted by $\MST(M)$.

\begin{rk}
If $M$ is finite, them $\MST(M)\ne\0$. For infinite $M$ the situation is rather more difficult, see~\cite{IvaNikTuz} and~\cite{IvaTuzInfTrees}.
\end{rk}

\begin{examp}
If all nonzero distances in $M$ are the same, then every spanning tree in $K(M)$ is minimal, so $\MST(M)=\cT(M)$.

If $\#M=3$, then each minimum spanning tree is obtained from the complete graph $K(M)$ by deleting the longest edge (if there are several, then any of them).
\end{examp}

\begin{prb}
Let $M$ be a finite metric space. Partition $M$ into nonempty subsets $M_1$ and $M_2$, and let $v_i\in M_i$ were chosen in such a way that $|v_1v_2|=|M_1M_2|$. Prove that there exists a minimum spanning tree $T\in\MST(M)$ such that $v_1v_2\in E(T)$.
\end{prb}

\begin{rk}
The problem of finding a minimum spanning tree can be naturally generalized. Let $M$ be a set. Consider a connected weighted graph $H$ with $V(H)=M$, and denote by $\cT(H)$ the set of all spanning subtrees of $H$. Then we put
$$
\mst(H)=\inf_{T\in\cT(H)}|T|
$$
and call it \emph{the weight of minimum spanning tree in $H$}. If there exists $T\in\cT(H)$ such that $|T|=\mst(H)$ then we call such $T$ \emph{a minimum spanning tree in $H$}. The set of all minimum spanning trees in $H$ is denoted by $\MST(H)$. If $M$ is a metric space, and $H=K(M)$ the corresponding weighted complete graph, then $\mst(H)=\mst(M)$ and $\MST(H)=\MST(M)$.

Note that there are a number of fast algorithms that solve the problem of finding a minimum spanning tree in a finite weighted connected graph. The most popular of them are Kruskal~\cite{Kruscal4} and Prim~\cite{Prim4} algorithms.
\end{rk}

\subsection{Steiner minimal tree problem}
Now we generalize the notion of minimum spanning tree. To do that, we consider $M$ as a subset of another metric spaces $X$, then we will minimize $\mst(V)$ over all $M\ss V\ss X$. Namely, we put
$$
\smt_X(M)=\inf\bigl\{\mst(V):M\ss V\ss X\bigr\}
$$
and call it \emph{the length of Steiner minimal tree on $M$}. Each $T\in\cT(V)$ for $M\ss V\ss X$ is called \emph{a shortest tree on $M$\/} or \emph{a Steiner minimal tree on $M$} if $|T|=\smt_X(M)$. The set of all Steiner minimal trees on $M$ is denoted by $\SMT_X(M)$. If it is clear or not important what $X$ the set $M$ belongs to, we simply write $\smt(M)$ and $\SMT(M)$ omitting $X$.

\begin{rk}
The following terminology is convenient when we study Steiner minimal trees or minimum spanning trees: if $G=(V,E)$ is a graph such that $V\ss X$, then we say that $G$ is a graph \emph{in the space $X$}. If $M\ss V$ then we say that $G$ \emph{joins $M$}; if $M=V$ then we say that $G$ \emph{spans $M$}. Thus, looking for minimum spanning trees we minimize the length of the trees spanning $M$, and for Steiner minimal trees we deal with the trees in the space $X$ joining $M$.
\end{rk}

\begin{rk}
The classical problem of finding a shortest tree is formulated for the Euclidean plane $X=\R^2$. The case $\#M=3$ arose as early as 1643 in works of Fermat~\cite{Fermat4}. For an arbitrary finite number of points on the Euclidean plane, the problem was posed by Jarn\'{\i}k and K\"ossler in 1934~\cite{JK4}. Courant and Robbins~\cite{CourRobb4} mistakenly called the problem of finding a shortest tree on the Euclidean plane the Steiner problem. Due to popularity of the book~\cite{CourRobb4}, this title has been fixed. The Steiner problem can be solved by Melzak's algorithm~\cite{Melzak4} and its many improvements, see for example~\cite{Hwang4} and~\cite{WWZ4}. As shown in~\cite{GGJ4}, the Steiner problem is algorithmically complex ($NP$-complete).
\end{rk}

\begin{rk}\hfill
\begin{enumerate}
\item Generally speaking, the set $\SMT(M)$ can be empty, also for finite $M$, however, the value $\smt(M)$ is always defined.
\item The set $\SMT_X(M)$ and the value $\smt_X(M)$ depend not only on the distances between points from $M$, but also on the geometry of the ambient space $X$: isometric $M$ lying in different metric spaces $X$ can be joined by Steiner minimal trees of different lengths. Some details on the theory of Steiner minimal trees can be found, for example, in~\cite{ITCRC4} or~\cite{HwRiW4}.
\end{enumerate}
\end{rk}

\begin{prb}
Find all Steiner minimal trees for $3$-point boundaries in the Euclidean plane. How many such trees exist for different boundaries?
\end{prb}

\begin{prb}
Find all Steiner minimal trees for the vertices of a square in the Euclidean plane. How many such trees exist?
\end{prb}

\begin{prb}
Find all Steiner minimal trees for $3$-point boundaries in the plane with $\ell_1$-metric defined by the norm $\bigl\|(x,y)\bigr\|=|x|+|y|$. How many such trees exist for different boundaries?
\end{prb}

\begin{prb}\label{prb:CompleteBut SMTnot}
Construct an example of a complete metric space and of some its finite subset $M$, such that there is no a Steiner minimal tree joining $M$.
\end{prb}

\subsection{One-dimensional minimal filling problem}
We now fix a finite metric space $M$. We will embed it isometrically into various metric spaces $X$, and minimize $\smt_X(M)$ over all such embeddings. To overcome the Cantor paradox, we put
$$
\mf(M)=\inf\Bigl\{r\mid\text{there exists an isometric embedding $\nu\:M\to X$ with $\smt_X\bigl(\nu(M)\bigr)\le r$}\Bigr\}
$$
and call it \emph{the length of minimal filling of $M$}. Each tree $G\in\SMT_X\bigl(\nu(M)\bigr)$ such that $|G|=\mf(M)$ is called \emph{a minimal filling of $M$}. The set of all minimal fillings of $M$ is denoted by $\MF(M)$.

\begin{rk}
Each graph $G$ in a metric space $X$ can be naturally considered as a weighted graph with the weight function assigning to the edges their lengths. Thus, each minimal filling is a weighted graph. The triangle inequality in $X$ leads to the fact that the distances between points in $M$ are majorized by the length of pathes in $G$ connecting these points.
\end{rk}

All this motivates an alternative equivalent definition of minimal fillings. Namely, let $G=(V,E,\om)$ be a weighted connected graph. Recall that in Construction~\ref{constr:DistOnWeightedGraph} we introduced the corresponding pseudometric $d_\om$ on $V$ as follows: for arbitrary $v,w\in W$ we put
$$
d_\om(v,w)=\inf\bigl\{\om(\g):\text{$\g$ is a walk joining $v$ and $w$}\bigr\}.
$$
Let $M$ be a metric space. A connected weighted graph $G=(V,E,\om)$ joining $M$ is called \emph{a filling of $M$} if for any $v,w\in M$ we have $|vw|\le d_\om(v,w)$.

\begin{prb}\label{prb:MinFilAlt}
Prove that for any metric space $M$ it holds
$$
\mf(M)=\inf\bigl\{\om(G):\text{$G$ is a filling of $M$}\bigr\}.
$$
\end{prb}

The following results were obtained in~\cite{IvaTuzMinFil4}.

\begin{prb}
Prove that for any finite metric space $M$ there exists a minimal filling.
\end{prb}

\begin{prb}
Let $M$ be a finite metric space with equal non-zero distances. Describe all minimal filling of $M$.
\end{prb}

\begin{rk}
The multidimensional problem on minimal fillings was formulated by M.Gromov~\cite{GromovMinFil4}. One-dimensional minimal filling as a stratified version of the Gromov's problem was studied by Ivanov and Tuzhilin~\cite{IvaTuzMinFil4}.
\end{rk}

\section{$\mst$-spectrum of a finite metric space}\label{sec:mst-spectrum}
\markright{\thesection.~$\mst$-spectrum of a finite metric space}
In this section we consider only finite metric spaces $M$, i.e., $\#M<\infty$.

To start with, we note that the minimum spanning tree, generally speaking, is not uniquely defined. For $G\in\MST(M)$, by $\s(G)$ we denote the vector whose elements are the lengths of the edges of the tree $G$ sorted in descending order. The following result is well known, however, we present its proof for completeness.

\begin{prop}\label{prop:mst-spect}
For any $G_1,G_2\in\MST(M)$ it holds $\s(G_1)=\s(G_2)$.
\end{prop}

\begin{proof}
Recall the standard algorithm for converting one minimum spanning tree to another~\cite{Kruscal4}.

Let $G_1\ne G_2$, $G_i=(M,E_i)$, then $E_1\ne E_2$ and $\#E_1=\#E_2$, therefore, there exists $e\in E_2\sm E_1$. The graph $G_1\cup e$ has a cycle $C$ containing the edge $e$. There is no longer edge in the $C$ cycle than $e$, because otherwise $G_1\not\in\MST(M)$. The forest $G_2\sm e$ consists of two trees whose vertex sets we denote by $V'$ and $V''$. Clearly, $M=V'\sqcup V''$. The cycle $C$ contains an edge $e'\ne e$ joining a vertex from $ V'$ with a vertex from $ V''$. This edge does not lie in $E_2$, otherwise $G_2$ would contain a loop. Therefore, $e'\in E_1\sm E_2$.

The graph $G_2\cup e'$ also contains some cycle $C'$. By the choice of $e'$, the cycle $C'$ also has the edge $e$. Similarly to the above, the length of the edge $e$ is less than or equal to the length of the edge $e'$, otherwise $G_2\not\in\MST(M) $. Therefore, $|e|=|e'|$.

Replacing the edge $e'$ in $G_1$ with $e$, we get a tree $G_1'$ of the same length, i.e., it is a minimum spanning tree as well, and $G_1'$ and $G_2$ have one common edge more than the trees $G_1$ and $G_2$. Thus, in a finite number of steps, we rebuild the tree $G_1$ into the tree $G_2$, passing through minimum spanning trees. It remains to notice that $\s(G_1')=\s(G_1)$, therefore, $\s(G_1)=\s(G_2)$.
\end{proof}

Proposition~\ref{prop:mst-spect} motivates the following definition.

\begin{dfn}
For any finite metric space $M$, by $\s(M)$ we denote $\s(G)$ for an arbitrary $G\in\MST(M)$ and call it \emph{the $\mst$-spectrum of the space $M$}.
\end{dfn}

\begin{constr}
For a set $M$ and a cardinal number $k\le\#M$, by $\cD_k(M)$ we denote the family of all possible partitions of the set $M$ into $k$ of its nonempty subsets. Now let $M$ be a metric space and $D=\{M_i\}_{i\in I}\in\cD_k(M)$. Put $\a(D)=\inf\bigl\{|M_iM_j|:i\ne j\bigr\}$.
\end{constr}

\begin{thm}\label{thm:spect-calc}
Let $M$ be a finite metric space and $\s(M)=(\s_1,\ldots,\s_{n-1})$. Then
$$
\s_k=\max\bigl\{\a(D):D\in\cD_{k+1}(M)\bigr\}.
$$
\end{thm}

\begin{proof}
Let $G=(M,E)\in\MST(M)$ and the set $E$ be ordered so that $|e_i|=\s_i$. Denote by $D=\{M_1,\ldots,M_{k+1}\}$ the partition of the set $M$ into the sets of vertices of the trees $G\sm\{e_i\}_{i=1}^k$.

\begin{lem}\label{lem:canon-part}
We have $\a(D)=|e_k|$.
\end{lem}

\begin{proof}
Indeed, we choose arbitrary $M_i$ and $M_j$, $i\ne j$, in them we take points $P_i$ and $P_j$, respectively, and let $\g$ be the unique path in $G$, joining $P_i$ and $P_j$. Then $\g$ contains some edge $e_p$, $1\le p\le k$. However, due to the minimality of the tree $G$, we have $|P_iP_j|\ge|e_p|\ge\min_i|e_i|=|e_k|$, thus $|M_iM_j|\ge|e_k|$, so $\a(D)\ge|e_k|$. On the other hand, if $i$ and $j$ are chosen so that $e_k$ joins $M_i$ and $M_j$, then we get $\a(D)\le|M_iM_j|=|e_k|$.
\end{proof}

Now consider an arbitrary partition $D'=\{M_1',\ldots,M_{k+1}'\}$.

\begin{lem}\label{lem:noncanon-part-less}
We have $\a(D')\le\a (D)$.
\end{lem}

\begin{proof}
By virtue of Lemma~\ref{lem:canon-part}, it suffices to show that $\a(D')\le|e_k|$. Denote by $E'$ the set consisting of all edges $e_p\in E$, for each of which there are $M_i'$ and $M_j'$, $i\ne j$, such that $e_p$ joins $M_i'$ and $M_j'$. Since $G$ is connected, the set $E'$ consists of at least $k$ edges; otherwise, the set of indices $\{1,\ldots,k+1\}$ is split into two nonempty subsets $I$ and $J$ such that the sets $\cup_{i\in I}M_i'$ and $\cup_{j\in J}M_j'$ that generate the partition $M$ are not joined by any edge from $E$. On the other hand, if some $M_i'$ and $M_j'$, $i\ne j$, are joined by an edge $e'\in E'$, then $|M_i'M_j'|\le|e'|$, hence $\a(D')=\min|M_i'M_j'|\le\min_{e'\in E'}|e'|\le|e_k|$.
\end{proof}

Lemma~\ref{lem:noncanon-part-less} completes the proof of the theorem.
\end{proof}

\section{Networks}
\markright{\thesection.~Networks}
To study Steiner minimal trees and minimal fillings, it is sometimes more convenient to work with so-called networks instead of the graphs in metric spaces. For example, if we investigate deformations of such graphs perturbing the positions of some their vertices, it may happen that after such a perturbation some vertices coincide, however we would like to preserve the structure of the graph by considering the coinciding vertices as different ones. To achieve this, we suppose that the graph is given apart of the metric space, and ``the positions of its vertices in the space'' are provided by a mapping from the vertex set of the graph to this space. Such mappings are called networks.

\begin{rk}\label{rk:ThreeObservations}
Let us discuss three additional observations.
\begin{enumerate}
\item In the optimization problems we usually deal with connected graphs, thus the domain of each network will be the vertex set of a connected graph, more often, of a tree.
\item We usually investigate boundary-value problems, that is why we need to partition the vertices into boundary and all remaining (nonboundary) ones.
\item We usually minimize the length of a graph. If a graph $G$ in a metric space contains a nonboundary vertex of degree $1$ or $2$, we can simplify $G$ preserving its connectedness, boundary, and not increasing its length. In the case of degree $1$ nonboundary vertex, we can simply remove the edge incident to this vertex. In the case of degree $2$ nonboundary vertex, we can change the both edges incident to this vertex by the unique edge joining the remaining vertices of these edges. That is why we usually assume that each boundary contains all the vertices of degree $1$ and $2$.
\end{enumerate}
\end{rk}

\begin{rk}
Indeed, we could also destroy possible cycles in the graph we optimize, and thus we might restrict ourselves with trees. However, in what follows we will see that such restriction leads us to some inconvenience, that is why we do not limit ourselves with trees only, but develop the corresponding theory for general connected graphs.
\end{rk}

Now we are ready to give formal definitions.

We will assume that in each graph $G$ there is a certain set of vertices $\d G\ss V(G)$ containing all vertices of degree $1$ and $2$ (the set $\d G$ can be empty), which is called \emph{the boundary of the graph $G$}, and the vertices from $\d G$ are called \emph{boundary\/} ones. The remaining vertices of the graph $G$ are called \emph{interior\/} ones. Sometimes the boundary vertices are also called \emph{fixed}, while the interior vertices are called \emph{movable}.

Let $G=(V,E)$ be a connected graph with some boundary $\d G$. \emph{A network of the type $G$ in a metric space $X$} is an arbitrary mapping $\G\:V\to X$. The $G$ is also called \emph{the parameterizing graph\/} of $\G$. In what follows, we will transfer to the networks all the terminology from the graph theory related to their parameterizing graphs.

Let $\G\:V\to X$ be a network parameterized by a connected graph $G=(V,E)$ with a boundary $\d G\ss V$. Then
\begin{enumerate}
\item the restrictions of $\G$ to the vertices and the edges of $G$ are called \emph{the vertices\/} and \emph{the edges\/} of $\G$, respectively;
\item the restriction of $\G$ to $\d G$ is called \emph{the boundary of $\G$} and is denoted by $\d\G$;
\item for each $vw\in E$ the value $\bigl|\G(v)\G(w)\bigr|$ is called \emph{the length of the edge $\G\:\{v,w\}\to X$};
\item the sum of the lengths of all edges of $\G$ is called \emph{the length of $\G$} and is denoted by $|\G|$.
\end{enumerate}

\begin{examp}\label{examp:factorForGraphInX}
Let $G=(V,E)$ be a connected graph in a metric space $X$ joining $M\ss X$, i.e., $M\ss V\ss X$, such that $M$ contains all the vertices of $G$ of degree $1$ and $2$. Put $\d G=M$, and define a network $\G\:V\to X$ as the inclusion mapping $\G\:v\mapsto v$. Thus, $\G$ is a network in $X$ of the type $G$, and $|\G|=|G|$.
\end{examp}

\begin{examp}\label{examp:InjectiveNet}
Let $X$ be a metric space, $G=(V,E)$ a connected graph with $\d G=M\ss X$. Let $\G\:V\to X$ be a network whose restriction to $M$ is the inclusion: $\G(v)=v$ for all $v\in M$. Suppose that the mapping $\G$ is injective. By identifying each vertex $v\in V$ with its image $\G(v)\in X$, we can consider $G$ as a graph in $X$, so the length $|G|$ of $G$ is defined. Then, with this identification, we have $|\G|=|G|$.
\end{examp}

\begin{dfn}\label{dfn:MainTypeNetwors}
As in Example~\ref{examp:InjectiveNet}, given a metric space $X$, let $G=(V,E)$ be a connected graph with $\d G=M\ss X$, and $\G\:V\to X$ be a network whose boundary $\d\G\:M\to X$ is the inclusion: $\d\G(v)=v$ for all $v\in M$. For such $\G$ we say that \emph{$\G$ joins the subset $M$ of the space $X$} (we do not assume here that $\G$ is injective outside $M$).
\end{dfn}

\subsection{Networks and quotients}

Now let $V'$ be a set, $\sim$ an equivalence relation on $V'$, and suppose that $V=V'/\fsim$. Denote by $\pi\:V'\to V$ the canonical projection. Let $G'=(V',E')$ be a connected graph such that $G=G'/\fsim$ and $\d G=\pi(\d G')$. Let $\G\:V\to X$ be a network of the type $G$ with the boundary $\d\G\:\d G\to X$, then the composition $\G'=\G\c\pi\:V'\to X$ is correctly defined and is a network in $X$ of the type $G'$ with the boundary $\d\G'=\d\G\c\pi|_{\d G'}$.

\begin{prb}\label{prb:twoNetsFactor}
Show that
\begin{enumerate}
\item\label{prb:twoNetsFactor:1} $|\G'|\ge|\G|$;
\item\label{prb:twoNetsFactor:2} if $G'$ is a tree, $V=\{V'_i\}$, and $G'(V'_i)\ss G'$ is a tree for each $i$, then $G$ is also a tree and $|\G'|=|\G|$;
\item\label{prb:twoNetsFactor:3} give an example in which $G'$ is a tree, $G$ is not a tree, and $|\G'|=|\G|$;
\item\label{prb:twoNetsFactor:4} give an example in which $G'$ and $G$ are trees, and $|\G'|>|\G|$.
\end{enumerate}
\end{prb}

\subsection{Splitting and splitting off for networks}
In this section we extend the operations of ``splitting a vertex'' and ``splitting off from a vertex'' which we defined above, to the case of graphs with boundaries and the corresponding networks. These operations enable us to simplify the structures of networks in consideration.

In Section~\ref{sec:graphPrelim} we defined splitting a vertex of degree greater than or equal to $4$, and splitting off a vertex of degree $1$ from a vertex of degree greater than or equal to $2$. For graphs with boundary, we refine these definitions.

We will split only interior vertices, while the resulting vertices will again be classified as interior; we will only split off from boundary vertices, and if we split off a vertex $w$ from the boundary vertex $v$, then we assign the vertex $v$ to interior one, and $w$ to boundary one.
The both these operations can be naturally defined for networks.

A graph $G$ with a boundary is called \emph{non-splittable\/} if no vertex can be split off from any boundary vertex, and no interior vertex can be split. For a finite graph $G$, we define \emph{the degree $p(G)$ of non-splitting}, setting it equal to the sum of $\deg v-3$ over all interior vertices $v$ of $G$. It is easy to see that the $G$, for which the degrees of all boundary vertices are $1$, is not splittable if and only if $p(G)=0$.

The next lemma will be the key point in the proof of Theorem~\ref{thm:smtInTermBinTrees}.

\begin{lem}\label{lem:SplittingOfAFiniteGraph}
For each finite graph $G=(V,E)$ with a boundary $\d G$ there exists a finite graph $G'=(V',E')$ with a boundary $\d G'$, and an equivalence relation $\sim$ on $V'$, such that $G=G'/\fsim$ and the following properties hold\/\rom:
\begin{enumerate}
\item all boundary vertices of $G'$ have degree $1$\rom;
\item all interior vertices of $G'$ have degree $3$\rom;
\item for connected $G$, the graph $G'$ is connected\/\rom;
\item for a tree $G$, the graph $G'$ is a tree\/\rom;
\item the $\sim$-class of each boundary vertex of $G'$ is a singleton\/\rom;
\item the canonical projection mapping $\pi\:V'\to V$ corresponding to $\sim$ is a bijection between $\d G'$ and $\d G$\rom;
\item if $G$ \(and $G'$\/\) is a tree, then for each $v\in V$, $W=\pi^{-1}(v)$, the graph $G'(W)\ss G'$ is a tree.
\end{enumerate}
\end{lem}

\begin{proof}
To start with, we split off a vertex of degree $1$ from each boundary vertex of $G$. As a result, we obtain a graph with all boundary vertices of degree $1$. In this graph, we will successively split all its interior vertices of degree greater than or equal to $4$. It is easy to see that the degree of non-splitting of this graph decreases by $1$ for each splitting, therefore, in a finite number of steps, we arrive to a non-splittable graph $G'=(V',E')$. We denote by $\sim$ the equivalence relation on $V'$, which is obtained from the trivial equivalence relation on $V$ (whose all equivalence classes are singletons) according to the following rule: for each splitting off and splitting, the resulting pair of vertices is equivalent to all those ones to which the original vertex was equivalent, and to each other. It is clear that $G'/\fsim=G$. It remains to notice that those splittings preserve the connectivity and do not create cycles.
\end{proof}

\section{Steiner minimal trees existence}\label{sec:STM}
\markright{\thesection.~Steiner minimal trees existence}
As we already mentioned in Problem~\ref{prb:CompleteBut SMTnot}, some metric spaces, also under assumption of completeness, may contain finite subsets which can not be joined by a Steiner minimal tree. In this section we prove that Steiner minimal trees always exist in boundedly compact metric spaces. To do that, we first need to reduce this problem to minimization of a finite number of continuous functions. To construct these functions, we show that it suffices to minimize the lengths of networks whose types can be chosen from some finite collection.

A tree $G$, and a network $\G$ of the type $G$ in a metric space $X$, both joining a set $M\ss X$, are called \emph{full Steiner trees\/} if all their boundary vertices have degree $1$ and all their interior vertices have degree $3$. Let us stress that the boundary $\d\G$ of this network is the inclusion $M\subset X$.

\begin{rk}
If the graph $G$ from Lemma~\ref{lem:SplittingOfAFiniteGraph} is a tree then the graph $G'$ from this lemma is a full Steiner tree.
\end{rk}

\begin{thm}\label{thm:smtInTermBinTrees}
Let $X$ be an arbitrary metric space and $M$ be a finite subset of $X$. Then $\smt(M)$ is equal to the infimum of the lengths of all full Steiner trees $\G$ joining $M$.
\end{thm}

\begin{proof}
Recall that
$$
\smt(M)=\inf\bigl\{|G|:\text{$G$ is a tree in $X$ with $\d G=M$}\bigr\}.
$$
As we mentioned in Remark~\ref{rk:ThreeObservations}, it suffices to consider only trees $G$ whose vertices of degree $1$ and $2$ belong to $M$. In what follows, we will minimize over such trees only.

\begin{prb}\label{prb:fullSteinerTreeInnerVertsNum}
Prove that any finite tree $G$ with $n\ge2$ boundary vertices contains at most $n-2$ interior vertices. The equality holds exactly in the case when $G$ is a full Steiner tree.
\end{prb}

Problem~\ref{prb:fullSteinerTreeInnerVertsNum} states that in calculation of $\smt(M)$ we can consider only finite trees $G$ with the boundary $M$. This enables us to use Lemma~\ref{lem:SplittingOfAFiniteGraph}, according to which for every $G$ there exists a full Steiner tree $G'$ and an equivalence relation $\sim$ on $V(G')$ such that $G=G'/\fsim$. Let $\pi\:V(G')\to V(G)$ be the corresponding canonical projection, then this $\pi$, being considered as a mapping from $V(G')$ to $X$, is a full Steiner tree joining $M$. Since $G$ and $G'$ are trees, then by Lemma~\ref{lem:SplittingOfAFiniteGraph} each equivalence class generates a subtree of $G'$, thus, by Problem~\ref{prb:twoNetsFactor}, we have $|\pi|=|G|$, that completes the proof.
\end{proof}

To list, up to a natural isomorphism, all the full Steiner trees from Theorem~\ref{thm:smtInTermBinTrees}, we construct a model set of such trees. By \emph{a model full Steiner tree\/} we mean a full Steiner tree $G=(V,E)$ with $V=\{1,2,\ldots,2n-2\}$ and $\d G=\{1,\ldots,n\}$. Two model full Steiner trees are called \emph{equivalent\/} if there is an isomorphism between them that is identical on the boundary $\{1,\ldots,n\}$. Thus, equivalent trees differ from each other by the numbering of their interior vertices. We denote by $\cB_n$ the set of all model full Steiner trees with $n$ boundary vertices considered up to the introduced equivalence. In other words, we construct $\cB_n$ by choosing in each equivalence class an arbitrary representative.

Now let $X$ be a metric space and $M$ be a finite subset of $X$ consisting of $n$ points. We enumerate the points from $M$ in an arbitrary way, i.e., consider some bijection $\v\:\{1,\ldots,n\}\to M$. Choose an arbitrary $G\in\cB_n$, and consider a network $\G$ of the type $G$ for which $\d\G=\v$. Then all such networks differ only in the ``positions'' of their interior vertices. The set of such networks is denoted by $[G,\v]$.

It is clear that all networks from Theorem~\ref{thm:smtInTermBinTrees} are obtained from the networks just described by identification, concordant with $\v$, of the set $\{1,\ldots,2n-2\}$ with the sets of vertices of the graphs parameterizing the former networks. Thus, we have proved the following result.

\begin{cor}\label{cor:smtInTermBinTrees}
Let $X$ be an arbitrary metric space and $M\ss X$ be a finite subset of $X$. Then
$$
\smt(M)=\inf\bigl\{|\G|:\G\in[G,\v],\,G\in\cB_n\bigr\}.
$$
\end{cor}

\begin{rk}
Note that the set $\cB_n$ by which we minimize in Corollary~\ref{cor:smtInTermBinTrees} is finite, and the set $[G,\v]$, by which minimization is also carried out, can be infinite. In some cases, it is easy to prove that the infimum in $[G,\v]$ is attained for every $G$, which immediately implies that $\SMT(M)\ne\0$ because $\cB_n$ is finite.
\end{rk}

Recall that a metric space is called boundedly compact if each of its closed balls is compact. Equivalent condition: a subset is compact if and only if it is closed and bounded.

We will present a technical result that is rather simple, but necessary in the future. Let $f\:X\to\R$ be some function defined on a metric space. We fix an arbitrary point $p\in X$ and for $r\ge0$ we set $F_p(r)=\inf_{x\in X\sm B_r(p)}f(x)$.

\begin{prb}\label{prb:functionGoesToInf}
Suppose that $F_p(r)\to\infty$ as $r\to\infty$. Prove that at any point $q\in X$ it holds $F_q(r)\to\infty$ as $r\to\infty$.
\end{prb}

From Problem~\ref{prb:functionGoesToInf}, the correctness of the following definition immediately follows. We say that a function $f\:X\to\R$, defined on a metric space, \emph{blows up at infinity\/} if the corresponding function $F_p(r)$ tends to infinity as $r\to\infty$ for some and, therefore, for every choice of the point $p\in X$.

\begin{prop}\label{prop:contFunInfAttain}
Let $f\:X\to\R$ be a continuous function defined on a boundedly compact metric space $X$. Suppose that $f$ blows up at infinity, then $f$ is bounded below and attains its infimum.
\end{prop}

\begin{proof}
Indeed, since the corresponding function $F_p(r)$ blows up at infinity, for some $r_0$ it holds $f(x)\ge0$ for all $x\in X\sm B_{r_0}(p)$. On the other hand, the function $f$ is bounded below on the closed ball $B_{r_0}(p)$ due to its compactness. Thus, the lower boundedness of the function $f$ is proved.

Let $f_0=\inf_{x\in X}f(x)$, and $r$ be such that for all $x\in X\sm B_r(p)$ we have $f(x)\ge f_0+1$. This means that the infimum of the function $f$ is attained on the ball $B_r(p)$, and, due to the compactness of this ball, there exists a point $x_0\in B_r(p)$ for which $f(x_0)=f_0$.
\end{proof}

\begin{thm}\label{thm:SMT-existence}
Let $X$ be a boundedly compact metric space. Then for every nonempty finite $M\ss X$ we have $\SMT(M)\ne\0$.
\end{thm}

\begin{proof}
We use Corollary~\ref{cor:smtInTermBinTrees}. If $\#M=n$, then we choose an arbitrary enumeration $\v\:\{1,\ldots,n\}\to M$, as well as an arbitrary model full Steiner tree $G\in\cB_n$. Then each network $\G\in[G,\v]$ is uniquely determined by the positions of its interior vertices, i.e., by the ``vector'' $z=\bigl(\G(n+1),\ldots,\G(2n-2)\bigr)\in X^{n-2}$. The function $\ell(z)=|\G|$ is continuous as the sum of continuous functions. In addition, this function blows up at infinity, therefore, by virtue of Proposition~\ref{prop:contFunInfAttain}, it attains its infimum. Also, there are a finite number of such functions in the formula from Corollary~\ref{cor:smtInTermBinTrees}, so the infimum from this formula is attained at a minimum point of one of these functions.
\end{proof}


\phantomsection
\renewcommand\bibname{References to Chapter~\thechapter}
\addcontentsline{toc}{section}{\bibname}
\renewcommand{\refname}{\bibname}

\vfill\eject

\section*{\Huge Exercises to Chapter~\thechapter}
\markright{Exercises to Chapter~\thechapter.}

\begin{exe}
Show that each forest is a disjoint union of trees that are components of this forest.
\end{exe}

\begin{exe}
Let $G=(V,E)$ be a connected graph, and $\sim$ an arbitrary equivalence relation on the set $V$. Show that the graph $G/\fsim$ is connected.
\end{exe}

\begin{exe}
Let $G=(V,E)$ be an arbitrary tree, and $\sim$ be an equivalence relation on the set $V$ such that for each class $V_i$ of this equivalence the subgraph $G(V_i)$ is a tree. Show that then $G/\fsim$ is a tree.
\end{exe}

\begin{exe}
Let $M$ be a finite metric space. Partition $M$ into nonempty subsets $M_1$ and $M_2$, and let $v_i\in M_i$ were chosen in such a way that $|v_1v_2|=|M_1M_2|$. Prove that there exists a minimum spanning tree $T\in\MST(M)$ such that $v_1v_2\in E(T)$.
\end{exe}

\begin{exe}
Find all Steiner minimal trees for $3$-point boundaries in the Euclidean plane. How many such trees exist for different boundaries?
\end{exe}

\begin{exe}
Find all Steiner minimal trees for the vertices of a square in the Euclidean plane. How many such trees exist?
\end{exe}

\begin{exe}
Find all Steiner minimal trees for $3$-point boundaries in the plane with $\ell_1$-metric defined by the norm $\bigl\|(x,y)\bigr\|=|x|+|y|$. How many such trees exist for different boundaries?
\end{exe}

\begin{exe}
Construct an example of a complete metric space and of some its finite subset $M$, such that there is no a Steiner minimal tree joining $M$.
\end{exe}

\begin{hint}
Consider on the set $X=\{0,1,2,\ldots\}$ the distance function $|mn|=1+\frac1{m+n}$, $m\ne n$. Prove that it is a complete metric. Consider the space $X^3$ with the complete metric generated by the norm $\|\cdot\|_\infty$. Put $M=\{(1,0,0),\,(0,1,0),\,(0,0,1)\}$. Prove that $\SMT(M)=\0$.
\end{hint}

\begin{exe}
Prove that for any metric space $M$ it holds
$$
\mf(M)=\inf\bigl\{\om(G):\text{$G$ is a filling of $M$}\bigr\}.
$$
\end{exe}

\begin{exe}
Prove that for any finite metric space $M$ there exists a minimal filling.
\end{exe}

\begin{exe}
Let $M$ be a finite metric space with equal non-zero distances. Describe all minimal fillings of $M$.
\end{exe}

\begin{exe}
Let $G=(V,E)$ be a connected graph in a metric space $X$ joining $M\ss X$, i.e., $M\ss V\ss X$. Put $\d G=M$, and define a network $\G\:V\to X$ as the embedding mapping $\G\:v\mapsto v$. Let $V'$ be a set, $\sim$ an equivalence relation on $V'$, and suppose that $V=V'\fsim$. Denote by $\pi\:V'\to V$ the canonical projection. Let $G'=(V',E')$ be a connected graph such that $G=G'/\fsim$. Then the composition $\G'=\G\c\pi\:V'\to X$ is a network in $X$ of the type $G'$. Show that
\begin{enumerate}
\item $|\G'|\ge|\G|$;
\item if $G'=(V',E')$ is a tree, and for each class $V'_i$ of the equivalence $\sim$ the subgraph $G'(V'_i)\ss G'$ is a tree, then $G$ is also a tree and $|\G'|=|\G|$;
\item give an example in which $G'$ is a tree, $G$ is not a tree, and $|\G'|=|\G|$;
\item give an example in which $G'$ and $G$ are trees, and $|\G'|>|\G|$.
\end{enumerate}
\end{exe}

\begin{exe}
Prove that any finite tree $G$ with $n\ge2$ boundary vertices contains at most $n-2$ interior vertices. The equality holds exactly in the case when $G$ is a full Steiner tree.
\end{exe}

\begin{exe}
Let $f\:X\to\R$ be some function defined on a metric space. We fix an arbitrary point $p\in X$ and for $r\ge0$ we set $F_p(r)=\inf_{x\in X\sm B_r (p)}f(x)$. Suppose that $F_p(r)\to\infty$ for $r\to\infty$. Prove that at any point $q\in X$ it holds $F_q(r)\to\infty$ as $r\to\infty$.
\end{exe}

\vfill\eject
 \chapter{Hausdorff distance.}
 \markboth{\chaptername~\thechapter.~Hausdorff distance.}%
          {\chaptername~\thechapter.~Hausdorff distance.}

\begin{plan}
Hausdorff distance, equivalence of three definitions, triangle inequality for Hausdorff distance, Hausdorff distance is a metric on the set of all closed bounded nonempty subsets, coincidence of Vietoris topology and metric topology generated by Hausdorff distance on the set of all compact subsets, limits theory for nonempty subsets, definition of $\limsup$ and some its properties (equivalent definitions), definition of $\liminf$ and some its properties (equivalent definitions), convergence w.r.t\. Hausdorff distance (Hausdorff convergence) and calculating $\liminf$, Hausdorff convergence of singletons, definition of $\lim$, Hausdorff convergence implies existence of $\lim$, the cases of decreasing and increasing sequences that are Hausdorff converging, equivalence of Hausdorff convergence and existence of $\lim$ in compact spaces, the cases of decreasing and increasing sequences in compact spaces, convergence in complete metric spaces, simultaneous completeness (total boundness, compactness) of the original space and the hyperspace of all closed bounded nonempty subsets, inheritance of the property to be geodesic for compact space and the hyperspace of all its closed nonempty subsets.
\end{plan}

Let $X$ be an arbitrary metric space. For arbitrary nonempty sets $A,B\ss X$ we put
\begin{align}
& d_H^1(A,B)=\max\Bigl(\sup\bigl\{|aB|:a\in A\bigr\},\,\sup\bigl\{|Ab|:b\in B\bigr\}\Bigr),\label{eq:Hminmax}\\
& d_H^2(A,B)=\inf\bigl\{r\in[0,\infty]:A\ss B_r(B)\ \&\ B_r(A)\sp B\bigr\},\label{eq:HB}\\
& d_H^3(A,B)=\inf\bigl\{r\in[0,\infty]:A\ss U_r(B)\ \&\ U_r(A)\sp B\bigr\}.\label{eq:HU}
\end{align}

\begin{prop}\label{prop:Hausd-Def-Euqiv}
For nonempty subsets $A,B\ss X$ of a metric space $X$ we have $d_H^1(A,B)=d_H^2(A B)=d_H^3(A,B)$.
\end{prop}

\begin{proof}
Put $r_i=d_H^i(A,B)$. First, let $r_1=\infty$. Without loss of generality, we assume that $\sup\bigl\{|aB|:a\in A\bigr\}=\infty$, but then neither $A\ss U_r(B)$, nor $A\ss B_r(B)$ holds for any finite $r>0$, therefore $r_2=r_3=\infty$. Now suppose that $r_1<\infty$, then for any $r>r_1$ all inclusions in the definitions of $r_2$ and $r_3$ take place, therefore $r_2$ and $r_3$ are also finite. So, we have shown that either all three $r_i$ are infinite at the same time, or all of them are finite.

It remains to analyze the case of finite $r_i$. First, we show that $r_1=r_2$ and $r_1=r_3$. Let us note that $r_2\le r_3$ because $A\ss U_r(B)$ implies $A\ss B_r(B)$ (and the same for $A$ and $B$ swaped).

By definition of $r_1$, we have $|aB|\le r_1$ for all $a\in A$, and $|Ab|\le r_1$ for all $b\in B$, therefore for all $r>r_1$ it holds $A\ss U_r(B)$ and $B\ss U_r(A)$, hence $r_3\le r$. Since $r>r_1$ is arbitrary, we have $r_2\le r_3\le r_1$. On the other hand, for any $r>r_2$ we have $A\ss B_r(B)$ and $B\ss B_r(A)$, so for any $a\in A$ and $b\in B$ it holds $|aB|\le r$ and $|Ab|\le r$, therefore, $r_1\le r$ and since $r>r_2$ is arbitrary, we get $r_1\le r_2\le r_3$.
\end{proof}

The value $d_H^i(A,B)$ from Proposition~\ref{prop:Hausd-Def-Euqiv} is denoted by $d_H(A,B)$. It is easy to see that $d_H$ is non-negative, symmetric, and $d_H(A,A)=0$ for any nonempty $A\ss X$, thus, $d_H$ is a generalized distance on the family $\cP_0(X)$ of all nonempty subsets of a metric space $X$. The function $d_H$ is called \emph{the Hausdorff distance}.

\begin{prop}\label{prop:HausdorffTriangle}
For an arbitrary metric space $X$, the function $d_H$ is a generalized pseudometric on $\cP_0(X)$.
\end{prop}

\begin{proof}
It remains to prove the triangle inequality. Choose arbitrary $A,B,C\in\cP_0(X)$ and set $c=d_H(A,B)$, $a=d_H(B,C)$, $b=d_H(A,C)$. We have to show that $b\le c+a$.

If either $c=\infty$, or $a=\infty$, then the equality holds. Suppose now that the both $c$ and $a$ are finite. Choose arbitrary finite $r>c$ and $s>a$, then $A\ss U_r(B)$ and $B\ss U_s(C)$ implies, by virtue of Item~(\ref{prb:metric-simple-props:1}) of Problem~\ref{prb:metric-simple-props}, that $A\ss U_r\bigl(U_s(C)\bigr)\ss U_{r+s}(C)$. Similarly, $U_{r+s}(A)\sp C$. Thus, $b\le r+s$. Since $r>c$ and $s>a$ are arbitrary, we obtain what is required.
\end{proof}

Denote by $\cH(X)\ss\cP_0(X)$ the set of all nonempty closed bounded subsets of a metric space $X$.

\begin{thm}
For an arbitrary metric space $X$, the generalized pseudometric $d_H$ is a metric on $\cH(X)$.
\end{thm}

\begin{proof}
Choose arbitrary $A,B\in\cH(A)$. Since they are bounded, for some $r>0$ we have $A\ss U_r(B)$ and $U_r(A)\sp B$, hence $d_H(A,B)<\infty$. Thus, $d_H$ is finite.

If $A\ne B$, then without loss of generality we can assume that there exists $a\in A\sm B$, but since the set $X\sm B$ is open, there exists $r>0$ such that $U_r(a)\cap B=\0$, in particular, $|aB|\ge r$ and, therefore, $d_H(A,B)\ge r$. Thus, $d_H$ is non-degenerate and, therefore, positively defined.
\end{proof}

Recall that by $\cK(X)$ we denoted the collection of all nonempty compact subsets of a topological space $X$. Since each compact subset of a metric space is closed and bounded, we have $\cK(X)\ss\cH(X)$ and, thus, we get

\begin{cor}
For an arbitrary metric space $X$, the function $d_H$ is a metric on $\cK(X)$.
\end{cor}

In what follows, when speaking about the distance in $\cH(X)$, we will always have in mind the Hausdorff metric, and for topology consider the corresponding metric one. Note that different authors use different notations for this hyperspace. We introduced the notation $\cH(X)$ by virtue of the fact that this is the largest natural set of subsets of a metric space on which the Hausdorff distance is defined.

We present a few Hausdorff distance properties in the next exercise.

\begin{prb}\label{prb:Hausdorff-evident-properties}
Prove the following statements for an arbitrary metric space $X$.
\begin{enumerate}
\item\label{prb:Hausdorff-evident-properties:1} Let $f\:X\to\cP_0(X)$ be given by the formula $f\:x\mapsto\{x\}$, then $f$ is an isometric embedding.
\item\label{prb:Hausdorff-evident-properties:3} For any $A,B\in\cP_0(X)$ we have $d_H(A, B)=d_H(A,\bB)=d_H(\bA,B)=d_H(\bA,\bB)$.
\item\label{prb:Hausdorff-evident-properties:4} For any $A,B\in\cP_0(X)$ we have $d_H(A,B)=0$ if and only if $\bA=\bB$.
\item\label{prb:Hausdorff-evident-properties:5} If $Y\ss X$ is an $\e$-net in $A\ss X$, then $d_H(A,Y)\le\e$.
\end{enumerate}
\end{prb}

\begin{prop}\label{prop:HausdAttainsForComp}
Let $X$ be an arbitrary metric space, and $A,B\in\cP_0(X)$, $r=d_H(A,B)$. Then $A\ss B_r(B)$, $B\ss B_r(A)$, and for less $r$ one of these inclusions fails. Thus, for $A,B\in\cP_0(X)$ we can change $\inf$ to $\max$ in equality~$(\ref{eq:HB})$.
\end{prop}

\begin{proof}
The fact that smaller $r$ do not fit follows directly from the definition of the Hausdorff distance. Let us prove the first part. Suppose the contrary, and let, say, $B\not\ss B_r(A)$. This means that there exists $b\in B$ for which $R:=|bA|>r$, thus for $r<s<R$ we have $B\not\ss B_s(A)$, which contradicts the definition of $d_H(A,B)$.
\end{proof}

\begin{prb}
Prove that for $A,B\in\cK(X)$ there exist $a\in A$ and $b\in B$ such that $d_H(A,B)=|ab|$. Is it possible to change $\cK(X)$ with $\cH(X)$?
\end{prb}

\begin{prb}
Let $X$ be an arbitrary metric space and $A,B,A',B'\in\cH(X)$ such that $A'\ss A$ and $B'\ss B$. Prove that $d_H(A\cup B',B\cup A')\le d_H(A,B)$.
\end{prb}

\begin{prb}
Let $X$ be an arbitrary metric space and $A,B,C\in\cH(X)$ such that $C\ss B$. Prove that $d_H(A,A\cup C)\le d_H(A,B)$.
\end{prb}

\section{Vietoris topology and Hausdorff metric}
\markright{\thesection.~Vietoris topology and Hausdorff metric}
In Construction~\ref{constr:Vietoris_top} we defined the Vietoris topology on the set of all nonempty subsets of a topological space $X$. Recall that a base of this topology is the family of sets
$$
\<U_1,\ldots,U_n\>=\{Y\ss X:Y\ss\cup_{i=1}^nU_i,\ \text{and $Y\cap U_i\ne\0$ for all $i=1,\ldots, n$}\}
$$
over all possible finite families $U_1,\ldots,U_n$ of open subsets of $X$.

\begin{thm}\label{thm:cHcompactHausdVietoris}
Let $(X,d)$ be an arbitrary metric space, then the metric topology on $\cK(X)$ defined by the Hausdorff metric $d_H$ coincides with the Vietoris topology. In particular, for a boundedly compact $X$ we have $\cH(X)=\cK(X)$, thus the above result holds if we change $\cK(X)$ with $\cH(X)$.
\end{thm}

\begin{proof}
We use Problem~\ref{prb:topologies-equality} to prove that each open set in metric topology generated by $d_H$ is also open in Vietoris topology. To do that, it suffices to take arbitrary $A\in\cK(X)$, $r>0$, and to construct an open neighborhood of $A$ in Vietoris topology that belongs to $U_r^{d_H}(A)$. Consider the family $\cC=\bigl\{U_{r/3}(a)\bigr\}_{a\in A}$, then $\cC$ is an open cover of the compact set $A$, thus we can extract from $\cC$ a finite subcover $\{U_1,\ldots,U_n\}$. Notice that $A\in\<U_1,\ldots,U_n\>$ because $A\ss\cup_{i=1}^nU_i$ and $A\cap U_i\ne\0$ for each $i=1,\ldots,n$. Further, we claim that $\<U_1,\ldots,U_n\>\ss U_r^{d_H}(A)$. Indeed, take an arbitrary $A'\in\<U_1,\ldots,U_n\>$, then $A'\ss U_{r/3}(A)$. Further, for every $a\in A$ there exist $U_i\ni a$ and $a'\in U_i\cap A'$, hence $|aa'|<2r/3$, so $A\ss U_{2r/3}(A')$, hence $d_H(A',A)\le2r/3<r$.

Now, consider an arbitrary $\<U_1,\ldots,U_n\>$ and any $A\in\<U_1,\ldots,U_n\>$. Let $C=X\sm\cup_{i=1}^nU_i$, then $C$ is a closed set that does not intersect $A$. By Item~(\ref{prb:metric-simple-props:02}) of Problem~\ref{prb:metric-simple-props}, the function $x\mapsto|xC|$ is continuous, thus its restriction onto the compact $A$ is bounded and attains its minimal value at some point $a\in A$. Since $C$ is closed and $a\not\in C$, we have $|aC|>0$ and, therefore, $r:=|AC|>0$.

Further, since $A\cap U_i\ne\0$ for each $i=1,\ldots,n$, we can choose some point $a_i$ in it. By the definition of the topology on $X$, for every $i$ there is $r_i$ such that $U_{r_i}(a_i)\ss U_i$. Put $\r=\min\{r,r_1,\ldots,r_n\}$. Then for each $B\in U_\r^{d_H}(A)$, that is, for each $B\in\cK(X)$ satisfying $d_H(A,B)<\r$, we have
\begin{enumerate}
\item $B\ss\cup_{a\in A} U_\r(a)\ss\cup_{i=1}^nU_i$, because $\r\le r$, and
\item for each $i=1,\ldots,n$ there is $b_i\in B$ such that $b_i\in U_\r(a_i)\ss U_i$ (because $\r\le r_i$), therefore $B$ intersects all $U_i$.
\end{enumerate}
It follows that $U_\r^{d_H}(A)\ss\<U_1,\ldots,U_n\>$, thus, $<U_1,\ldots,U_n\>$ is open in the metric topology generated by $d_H$. The proof is over.
\end{proof}

\begin{cor}
If two metrics on a set $X$ induce the same topology on $X$, then also the same metric topology is induced on the corresponding spaces $\cK(X)$. In other words, the metric topology of the space $\cK(X)$ does not depend on the specific form of the metric on $X$, but only on the topology defined by this metric on $X$. For boundedly compact $X$ the same is true for $\cH(X)$.
\end{cor}

The next example shows that the Vietoris and the metric topologies on $\cH(X)$ can be different also if $X$ is bounded and locally compact metric space.

\begin{examp}\label{examp:HasuNontVietoris}
Let $X=(-1,0)\cup(0,1)\ss\R$ endowed with the metric induced from the standard metric on $\R$. Then $A=(-1,0)$ is a closed bounded subset of $X$, thus $A\in\cH(X)$. Since $A$ is also an open set, then $\cU=\<A\>$ is a neighborhood of $A$ in Vietoris topology consisting of all closed nonempty subsets of $X$ which belong to $A$. However, for any $\e>0$ the ball $U^{d_H}_\e(A)\ss\cH(X)$ contains the closed subset $B=A\cup\{\e/2\}$ such that $B\not\in\cU$, therefore $U^{d_H}_\e(A)\not\ss\cU$ and, thus, $\cU$ is not open in the metric topology generated by the Hausdorff distance.
\end{examp}

\begin{examp}
In Example~\ref{examp:HasuNontVietoris} the space $X$ was not complete. Now we present an example of complete metric space where the Vietoris topology and the metric topology on $\cH(X)$ are different.

Let $X=\ell_2$ be the space of all sequences $(x_1,x_2,\ldots)$ of real numbers such that $\sum_{i=1}^\infty x^2<\infty$. For $\xi=(x_1,x_2,\ldots)\in\ell_2$ and $\xi'=(x'_1,x'_2,\ldots)\in\ell_2$ we put $d(\xi,\xi')^2={\sum_{i=1}^\infty(x_i-x'_i)^2}$.

\begin{prb}
Prove that $d$ is a metric on $\ell_2$, and that $(\ell_2,d)$ is a complete metric space.
\end{prb}

Denote by $e_i\in\ell_2$ the sequence $(x_1,x_2,\ldots)$ such that $x_i=1$ and $x_j=0$ for all $j\ne i$.

\begin{prb}
Prove that $A=\{e_i\}_{i=1}^\infty$ is a closed subset of $\ell_2$, thus $A\in\cH(\ell_2)$.
\end{prb}

Now, we put $U_i=U_{\frac1{2i}}(e_i)$ and let $U=\cup_{i=1}^\infty U_i$, then $U$ is an open neighborhood of $A$. Let $\cU=\<U\>\ss\cH(\ell_2)$ be the corresponding neighborhood of $A\in\cH(\ell_2)$ in the Vietoris topology. Choose an arbitrary $\e>0$. We show that $U_\e^{d_H}(A)\not\ss\cU$, thus, $\cU$ is not open in the metric topology generated by $d_H$. Consider an arbitrary $i$ such that $\frac1{2i}<\e$, and choose any $\frac1{2i}<\dl<\min\{1,\e\}$. Put $\xi=(1-\dl)e_i$ and $B=A\cup\{\xi\}\in\cH(\ell_2)$, then $d(\xi,e_i)=\dl<\e$, therefore, $d_H(A,B)<\e$. On the other hand, since $\frac1{2i}<\dl$, then $\xi\not\in U_i$; also, for any $j\ne i$ we have $d(\xi,e_j)=\sqrt{1+(1-\dl)^2}>1>\frac1{2j}$, thus $\xi\not\in U_j$, and, in account, $\xi\not\in\cU$.
\end{examp}

\begin{prb}
Let $X=\N$, and define two metrics on $X$: $d^1(x,y)=1$ for any $x\ne y$, and $d^2(x,y)=\bigl|\frac1x-\frac1y\bigr|$ for any $x,y$. Then the corresponding metric topologies are both discrete. Prove that the corresponding Hausdorff metric generates non-homeomorphic topologies on $\CL(X)$.
\end{prb}

\section{Limits Theory}
\markright{\thesection.~Limits Theory}
For the future, we need one technical result. Let $A_1,A_2,\ldots$ be a sequence of nonempty subsets of a metric space $X$.

\begin{dfn}\label{dfn:Flimsup}
Put
$$
\limsup A_k=\cap_{n=1}^\infty\overline{A_n\cup A_{n+1}\cup\cdots}
$$
and call it \emph{the upper limit\/} of the sequence $A_k$.
\end{dfn}

\begin{rk}\label{rk:FlimsupClosed}
Since $\limsup A_k$ is equal to the intersection of closed sets, it is always closed (possibly empty).
\end{rk}

\begin{prop}\label{prop:Flimsup}
We have
$$
\limsup A_k=\bigl\{x\in X:\forall\e>0\ \text{it holds $\#\{k:U_\e(x)\cap A_k\ne\0\}=\infty$}\bigr\}.
$$
\end{prop}

\begin{proof}
Put $B_n=\cup_{k=n}^\infty A_k$ and $A=\limsup A_k=\cap_{n=1}^\infty\bB_n$.

Choose an arbitrary $\e>0$. Let $x\in A$, then $x$ is adherent point of each set $B_n$, so $B_n\cap U_\e(x)\ne\0$. If the condition $A_k\cap U_\e(x)\ne\0$ were satisfied only for a finite number $k$, then for some $n$ they would have $B_n\cap U_\e(x)=\0$, so $x$ would not be the adherent point of the set $B_n$, a contradiction.

Conversely, suppose that for every $\e>0$ it holds $\#\{k:U_\e(x)\cap A_k\ne\0\}=\infty$. Then for every $n\in\N$ we have $U_\e(x)\cap B_n\ne\0$, so $x\in\bB_n$ for any $n$, therefore, $x\in A$.
\end{proof}

\begin{cor}\label{cor:Flimsup}
We have
$$
\limsup A_k=\{x\in X:\text{there exists a sequence $a_{i_k}\in A_{i_k}$ that converges to $x$}\}.
$$
\end{cor}

\begin{proof}
Let $x\in\limsup A_k$. By Proposition~\ref{prop:Flimsup}, for each $k\in\N$ one can find $a_{i_k}\in A_{i_k}$ such that $|a_{i_k}x|<1/k$, and there are infinitely many such $i_k$ (for fixed $k$). Thus, we can compose an increasing sequence $i_1<i_2 <\cdots$ for which $a_{i_k}\in A_{i_k}$, $|a_{i_k} x|<1/k$, and therefore, $a_{i_k}\to x$.

Conversely, let some sequence $a_{i_k}\in A_{i_k}$ converge to some $x\in X$. This means that for every $\e>0$ there is $n$ such that for any $k\ge n$ we have $U_\e (x)\cap A_{i_k}\ne\0$, therefore $x\in\limsup A_k$ by virtue of Proposition~\ref{prop:Flimsup}.
\end{proof}

\begin{rk}
The upper limit may be empty.
\end{rk}

\begin{examp}\label{examp:FlimsupEmpty}
Take as $X$ the interval $(0,1)$ with the standard metric, and put $A_k=\{1/k\}$, $k\in N$, then one cannot select a convergent subsequence $a_{i_k}$; therefore, $\limsup=\0$. A similar situation concerns the space $X=\R$ and the sequence $A_k=\{k\}$. Thus, the space $X$ can be totally bounded or complete, but it can contain sequences of subsets $A_k\ss X$ for which $\limsup A_k$ is empty.
\end{examp}

The construction of upper limit can be extended as follows.

\begin{dfn}\label{dfn:Fliminf}
Put
$$
\liminf A_k=\bigl\{x\in X:\forall\e>0\ \text{it holds $\#\{k:U_\e(x)\cap A_k=\0\}<\infty$}\bigr\}
$$
and call it \emph{the lower limit\/} of the sequence $A_k$.
\end{dfn}

\begin{rk}\label{rk:FliminfSSFlimsup}
Note that $\liminf A_k\ss\limsup A_k$.
\end{rk}

\begin{prop}\label{prop:Fliminf}
We have
$$
\liminf A_k=\{x\in X:\text{there exists a sequence $a_k\in A_k$ that converges to $x$}\}.
$$
\end{prop}

\begin{proof}
Let $x\in\liminf A_k$. By Definition~\ref{dfn:Fliminf}, for each $n\in\N$ one can find $k_n\in\N$ such that for all $k\ge k_n$ there exists $a_k^n\in A_k$ for which $|a_k^nx|<1/n$. Without loss of generality, we assume that the sequence $k_n$ is strictly monotone, and construct a sequence $a'_k\in A_k$, starting with arbitrary $a_1,\ldots,a_{k_1-1}$, and then adding $a_{k_1}^1,\ldots,a_{k_2-1}^1$,  $a_{k_2}^2,\ldots,a_{k_3-1}^2$, etc. It is clear that the sequence $a'_k\in A_k$ constructed in this way converges to $x$.

Conversely, let some sequence $a_k\in A_k$ converge to some $x\in X$. This means that for every $\e>0$ there is $n$ such that for any $k\ge n$ we have $U_\e(x)\cap A_k\ne\0$, therefore $x\in\liminf A_k$ by Definition~\ref{dfn:Fliminf}.
\end{proof}

\begin{prb}\label{prb:LiminfIsClosed}
Prove that $\liminf A_i$ is a closed subset of $X$.
\end{prb}

\begin{dfn}
If $X$ is a metric space, $A_k\in\cP_0(X)$, $k\in\N$, $A\in\cP_0(X)$, and $A_k$ converges to $A$ w.r.t\. the Hausdorff distance, then we write this as $A_k\toH A$.
\end{dfn}

\begin{prop}\label{prop:converge-eq-Fliminf}
If $A_k,A\in\cP_0(X)$ and $A_k\toH A$, then $\bA=\liminf A_k$, where $\bA$ is the closure of $A$.
\end{prop}

\begin{proof}
By Item~(\ref{prb:Hausdorff-evident-properties:3}) of Problem~\ref{prb:Hausdorff-evident-properties}, we have $A_k\to\bA$, thus, without loss of generality, we can assume that $A$ is a closed subset of $X$.

Now, we show that $A$ always contains the limit point of each convergent sequence $a_k\in A_k$, i.e., that $\liminf A_k\ss A$. Let $a_k\to a$. By definition, for any $\e>0$ there exists $n$ such that for all $k\ge n$ we have $d_H(A,A_k)<\e/2$ and $|aa_k|<\e/2$. Hence, for such $k$ we have $A_k\ss U_{\e/2}(A)$, in particular, $a_k\in U_{\e/2}(A)$. The latter means that for such $k$ there exist $a'_k\in A$ for which $|a'_ka_k|<\e/2$ and, therefore, $|a'_ka|<\e$. Thus, $a'_k\to a$ and therefore, since $A$ is closed, we have $a\in A$.

To complete the proof, it remains to verify that each point $a\in A$ is the limit point for some sequence $a_k\in A_k$, i.e., that $\liminf A_k\sp A$. For each $n$ there exists $k_n$ such that for all $k\ge k_n$ it holds $d_H(A_k,A)<1/n$. For such $k$ we have $A\ss U_{1/n}(A_k)$, in particular, $a\in U_{1/n}(A_k)$. It is clear that the sequence $k_n$ can be chosen strictly monotonic. For such $k_n$, we take $a_k$ as follows: for $1\le k\le k_1-1$ we choose $a_k\in A_k$ arbitrarily; for $k_1\le k<k_2$ we select $a_k\in A_k$ so that $|a_ka|<1$; for $k_2\le k <k_3$, select $a_k\in A_k$ so that $|a_ka|<1/2$; etc. Such sequence converges to $a$, as required.
\end{proof}

\begin{cor}\label{cor:convergeOnePointedSeq}
Let $A_k=\{a_k\}$ for all $k$. Then the sequence $A_k$ converges w.r.t\. $d_H$ if and only if the sequence $a_k$ converges. Moreover, if $a_k\to a$, then $A_k\toH\{a\}$.
\end{cor}

\begin{proof}
By virtue of Item~(\ref{prb:Hausdorff-evident-properties:1}) of Problem~\ref{prb:Hausdorff-evident-properties}, the mapping $f\:X\to\cH(X)$, $f\:x\mapsto\{x\}$, is isometric, therefore the convergence of $a_k\to a$ implies the convergence of $A_k\toH\{a\}$.

Conversely, if $A_k\toH A$, then $A\ne\0$ by definition. Also, $A$ cannot contain two different points. Indeed, if $a,a'\in A$, $a\ne a'$, then for $0<\e<|aa'|/2$ and sufficiently large $k$ we have $\{a,a'\}\ss A\ss U_\e(A_k)=U_\e(a_k)$, that is impossible. Thus, $A=\{a\}$ and for any $\e>0$ there exists $n$ such that for any $k\ge n$ we have $A_k\in U_\e(A)$, i.e., $a_k\in U_\e(a)$, thus $a_k\to a$.
\end{proof}

\begin{dfn}
If the upper and lower limits of a sequence $A_1,A_2,\ldots$ are nonempty and equal to each other, then the sequence $A_1,A_2,\ldots$ is said to have \emph{a limit}, which is denoted by $\lim A_k$. We write this as $A_k\to A$.
\end{dfn}

Let us discuss what the existence of the limit $\lim A_k$ means.

\begin{thm}\label{thm:convergeFlimExists}
If $A_k,A\in\cP_0(X)$ and $A_k\toH A$, then there exists $\lim A_k$ and $\bA=\lim A_k$.
\end{thm}

\begin{proof}
By Proposition~\ref{prop:converge-eq-Fliminf} and Remark~\ref{rk:FliminfSSFlimsup}, we have $\bA=\liminf A_k\ss\limsup A_k$. Therefore, it suffices to show that $\limsup A_k\ss\bA$.

Let $x\in\limsup A_k$. We show that $x$ is an adherent point for $A$ and, therefore, $x\in\bA$. To do this, in turn, it is enough to show that for an arbitrary $\e>0$ the ball $U_\e(x)$ intersects $A$.

Since $A_k\toH A$, there exists $n$ such that for every $k\ge n$ we have $A_k\ss U(A,\e/2)$. By Proposition~\ref{prop:Fliminf}, there exists $k\ge n$ for which $U_{\e/2}(x)\cap A_k\ne\0$, i.e., there exists $a_k\in A_k$ such that $|xa_k|<\e/2$. Since $A_k\ss U(A,\e/2)$, there exists $a\in A$ for which $|a_ka|<\e/2$, thus $|xa|<\e$ and, hence $U_\e(x)\cap A\ne\0$.
\end{proof}

Is the converse statement to that of Theorem~\ref{thm:convergeFlimExists} true as well? The following example demonstrates that this is not the case in general, also when $A_k,A\in\cH(X)$.

\begin{examp}\label{examp:FlimExistsNonConverge}
Let $X$ be the interval $(0,3)\ss\R$ with the standard distance function. Put $A_k=\{2,1/k\}$. Then $\limsup A_k=\liminf A_k=\{2\}$, so that $\lim A_k$ exists and is equal to $\{2\}$. If for some $A\in\cH(X)$ we have $A_k\toH A$, then, by virtue of Theorem~\ref{thm:convergeFlimExists}, we would have $A=\{2\}$. However, for $\e=1$ there is no $A_k$ that belongs to $U_\e(A)$, so the sequence $A_k$ diverges.

The case of $X=\R$ and $A_k=\{0,k\}$, for which $\liminf A_k=\limsup A_k=\{0\}$, but the sequence $A_k$ diverges, is similarly analyzed. Thus, the existence of a limit does not imply the convergence of the sequence $A_k\in\cH(X)$ for either totally bounded or complete $X$.
\end{examp}

We give some corollaries.

\begin{cor}\label{cor:convergeDecrease}
If $A_k,A\in\cP_0(X)$, $A_k\toH A$, and for all $k$ we have $A_k\sp A_{k+1}$, then $\bA=\cap_k\bA_k$.
\end{cor}

\begin{proof}
By Theorem~\ref{thm:convergeFlimExists}, we have
$$
\bA=\lim A_k=\limsup A_k=\cap_{k=1}^\infty\overline{A_k\cup A_{k+1}\cup\cdots}=\cap_{k=1}^\infty\bA_k.
$$
\end{proof}

\begin{cor}\label{cor:convergeIncrease}
If $A_k,A\in\cP_0(X)$, $A_k\toH A$, and for all $k$ we have $A_k\ss A_{k+1}$, then $\bA=\overline{\cup_kA_k}$.
\end{cor}

\begin{proof}
By Theorem~\ref{thm:convergeFlimExists}, we have
$$
\bA=\limsup A_k=\overline{A_1\cup A_2\cup\cdots},
$$
since all the sets $B_k=A_k\cup A_{k+1}\cup\cdots$ coincide.
\end{proof}

\subsection{Limits for compact $X$}

Example~\ref{examp:FlimsupEmpty} shows that a upper limit can be empty. However, for compact $X$ this is no longer the case.

\begin{prop}\label{prop:FlimsupCompNonEmpty}
For any compact metric space $X$ and any sequence of nonempty $A_k\ss X$ we have $\limsup A_k\ne\0$.
\end{prop}

\begin{proof}
Choose an arbitrary sequence $a_k\in A_k$, then it contains a subsequence that converges to some $x\in X$. By Corollary~\ref{cor:Flimsup}, we have $x\in\limsup A_k$, so the upper limit of this sequence is not empty.
\end{proof}

We will need the following technical result in the future.

\begin{prop}\label{prop:FlimsupCompHalfLim}
For any compact metric space $X$ and any sequence of nonempty $A_k\ss X$ the following statement holds\/\rom: for any $\e>0$ there exists $n$ such that for any $k\ge n$
$$
A_k\ss U_\e(\limsup A_k)\ \ \text{and}\ \ \liminf A_k\ss U_\e(A_k).
$$
\end{prop}

\begin{proof}
To start with, we prove the first inclusion. Put $A=\limsup A_k$, then, by Proposition~\ref{prop:FlimsupCompNonEmpty}, we have $A\ne\0$. Now, we suppose to the contrary that for some sequence $A_k$ there exists an $\e>0$ and a subsequence $A_{i_k}$ such that $A_{i_k}\not\ss U_\e(A)$. The latter means that in every $A_{i_k}$ there is a point $a_{i_k}$ for which $|a_{i_k}A|\ge\e$. Since $X$ is compact, there exists a subsequence in the sequence $a_{i_k}$ that converges to some $x\in X$. Then, by Corollary~\ref{cor:Flimsup}, we have $x\in A$, which implies that for some sufficiently large $k$ it holds $|a_{i_k}A|<\e$, a contradiction.

We now prove the second inclusion. Put $A=\liminf A_k$ and assume the contrary, i.e., that there exists $\e>0$ and a sequence $i_1<i_2<\cdots$ for which $A\not\ss U_\e(A_{i_k})$, in particular, $A\ne\0$. The latter is equivalent to the existence of $a'_{i_k}\in A$ such that $a'_{i_k}\not\in U_\e(A_{i_k})$. By Problem~\ref{prb:LiminfIsClosed}, the set $A$ is closed, thus $A$ is compact as a closed subset of a compact metric space, therefore, in the sequence $a'_{i_k}$ one can choose a subsequence converging to some $a\in A$. Without loss of generality, we assume that this subsequence coincides with the entire sequence $a'_{i_k}$.

Since $A=\liminf A_k$, by Proposition~\ref{prop:Fliminf} there exists a sequence $a_k\in A_k$ converging to $a$, thus $|a'_{i_k}a_{i_k}|\to0$, so for sufficiently large $k$ we have $a'_{i_k}\in U_\e(A_{i_k})$, a contradiction.
\end{proof}

\begin{cor}\label{cor:convergeCompandFlim}
For any compact metric space $X$ and any sequence of nonempty $A_k\ss X$, the $A_k$ converges w.r.t\. $d_H$ if and only if there exists $\lim A_k$. Moreover, if $A_k\toH A$, then $\bA=\lim A_k$.
\end{cor}

\begin{proof}
By Theorem~\ref{thm:convergeFlimExists}, the convergence of the sequence $A_k$ to some $A$ implies the existence of $\lim A_k$ and that $\bA=\lim A_k$. Thus, it remains to prove the converse.

So, let $\lim A_k=:A$ exist. We will show that $A_k\toH A$. By Proposition~\ref{prop:FlimsupCompHalfLim}, for any $\e>0$ there exists $n$ such that for any $k\ge n$ we have $A_k\ss U_\e(\limsup A_k)=U_\e(A)$ and $A=\liminf A_k\ss U_\e(A_k)$, hence $A_k\toH A$, as required.
\end{proof}

Let us demonstrate how the theory of limits works in the compact case.

\begin{cor}\label{cor:DecreasingStartetFromCompact}
For an arbitrary, not necessarily compact metric space $X$, and each decreasing sequence of nonempty sets $A_1\sp A_2\sp\cdots$ starting from compact $A_1$, we have $A:=\cap_k\bA_k\ne\0$ and $A_k\toH A$.
\end{cor}

\begin{proof}
Without loss of generality, we assume that $X=A_1$, then $X$ is a compact space. By Proposition~\ref{prop:FlimsupCompNonEmpty}, we have $\limsup A_k\ne\0$, however, by definition of $\limsup$, we have $\limsup A_k=A$, thus $A$ is not empty. On the other hand, each point $a\in A$ is the limit of the constant sequence $a_k=a\in\bA_k$, and for each $k$ there exists $a'_k\in A_k$ such that $|a_ka'_k|<1/k$, thus $a'_k\to a$ as well and, hence, $a\in\liminf A_k$. Since $a\in A$ is arbitrary then $\limsup A_k\ss\liminf A_k$. However, $\liminf A_k\ss\limsup A_k$, thus $\liminf A_k=A=\limsup A_k$, and it remains to use Corollary~\ref{cor:convergeCompandFlim}.
\end{proof}

\begin{rk}
In Corollary~\ref{cor:DecreasingStartetFromCompact} it may happen that $\cap_{k=1}^\infty A_k=\0$. Let $X=\R$, $A_1=[0,1]$, $A_k=(0,1/k)$ for $k>1$, then $\cap_{k=1}^\infty A_k=\0$, however,
$$
\limsup A_k=\cap_{k=1}^\infty\bA_k=\cap_{k=1}^\infty[0,1/k]=\{0\}\ne\0.
$$
\end{rk}

\begin{cor}
For an arbitrary, not necessarily compact, space $X$, if for an increasing sequence of nonempty sets $A_1\ss A_2\ss\cdots$ the set $A=\overline{\cup_{k=1}^\infty A_k}$ is compact, then $A_k\toH A$.
\end{cor}

\begin{proof}
Again, without loss of generality, we assume that $X=A$. By definition of $\limsup$, we have $\limsup A_k=A$, therefore, by Corollary~\ref{cor:Flimsup}, for every $a\in A$ there exists a sequence $a_{i_k}\in A_{i_k}$, $i_1<i_2<\cdots$, converging to $a$. Notice that each $a_{i_k}$ also belongs to all $A_p$ for $p>i_k$, therefore this sequence can be extended to a sequence $a_i\in A_i$ converging to $a$, so, by Proposition~\ref{prop:Fliminf}, we have $a\in\liminf A_k$, hence $\liminf A_k=\limsup A_k=A$, and it remains to use Corollary~\ref{cor:convergeCompandFlim}.
\end{proof}

\subsection{Limits for complete $X$}

We give a series of results for fundamental sequences of subsets of complete metric spaces. It will be convenient for us to use the following construction.

\begin{constr}
Let $Y$ be an arbitrary generalized pseudometric space and $\om=(y_1,y_2,\ldots)$ be some sequence in it. Put
$$
d_n(\om)=\sup_{p,q\ge n}|y_py_q|.
$$
Note that $d_n(\om)$ is a non-negative non-increasing sequence, therefore, it has a limit $d(\om)$.
\end{constr}

\begin{prop}\label{prop:fundANDdk}
The sequence $\om=(y_1,y_2,\ldots)$ is fundamental if and only if $d(\om)=0$.
\end{prop}

\begin{proof}
If $y$ is fundamental, then for every $\e>0$ there exists $n$ such that for any $p,q\ge n$ we have $|y_py_q|<\e$, so $d_n(\om)\le\e$ and, therefore, $d(\om)=0$. Conversely, if $d(\om)=0$, then for any $\e>0$ starting from some $n$ it holds $d_n(\om)<\e$, so for such $n$ and any $p,q\ge n$ we have $|y_py_q|<\e$, which means that the sequence $\om$ is fundamental.
\end{proof}

Now we present an important technical result on the fundamental sequences of an arbitrary metric space $X$.

\begin{prop}\label{prop:eachElA1inFund}
Let $\Om=(A_1,A_2,\ldots)$, $A_k\in\cP_0(X)$, be a fundamental sequence. Choose arbitrary $a_1\in A_1$ and $\e>0$. Then there exists a sequence $i_1=1<i_2<i_3 <\cdots$ of natural numbers and a fundamental sequence $\om=(a_{i_1},a_{i_2},a_{i_3},\ldots)$, $a_{i_k}\in A_{i_k}$, such that $d_1(\om)<d_1(\Om)+\e$.
\end{prop}

\begin{proof}
Proposition~\ref{prop:fundANDdk} implies $d_k(\Om)\to0$, so there exists $i_2\in\N$, $i_2>i_1$, such that $d_{i_2}(\Om)<\e/8$. Since $d_H(A_{i_1},A_{i_2})\le d_1(\Om)$, there exists $a_{i_2}\in A_{i_2}$ for which $|a_{i_1}a_{i_2}|<d_1(\Om)+\e/4$. Next, we find a natural $i_3>i_2$ such that $d_{i_3}(\Om)<\e/16$, as well as $a_{i_3}\in A_{i_3}$ for which $|a_{i_2}a_{i_3}|<\e/8$. Continuing this process, we construct a sequence $\om=(a_{i_1}=a_1,a_{i_2},a_{i_3},\ldots)$, $a_{i_k}\in A_{i_k}$, for which $|a_{i_{k-1}}a_{i_k}|<\e/2^k$ for $k\ge3$. By the triangle inequality, for any $p,q\in\N$, $2\le p\le q$, we have $|a_{i_p}a_{i_q}|<\e/2^p$, therefore, the sequence $\om$ is fundamental. On the other hand, for any $p,q\in\N$ we have $|a_{i_p}a_{i_q}|<d_1(\Om)+\e/2$, hence $d_1(\om)\le d_1(\Om)+\e/2<d_1(\Om)+\e$.
\end{proof}

Example~\ref{examp:FlimsupEmpty} demonstrates that for a complete $X$ a upper limit can be empty. However, for a fundamental sequence $A_k$ this is no longer the case.

\begin{prop}\label{prop:FlimsupCompleteNonEmpty}
Let $X$ be a complete space and $A_k\in\cP_0(X)$ be a fundamental sequence, then $\limsup A_k\ne\0$.
\end{prop}

\begin{proof}
We choose an arbitrary $a_1\in A_1$, then, by Proposition~\ref{prop:eachElA1inFund}, there exists a fundamental sequence $(a_{i_1}=a_1,a_{i_2},a_{i_3},\ldots)$ such that $a_{i_k}\in A_{i_k}$. Since $X$ is a complete space, this sequence converges to some $x\in X$. Proposition~\ref{prop:Flimsup} implies that $x\in\limsup A_k$, therefore $\limsup A_k\ne\0$.
\end{proof}

Now we will give an analogue of Proposition~\ref{prop:FlimsupCompHalfLim}.

\begin{thm}\label{thm:FlimsupCompleteHalfLim}
Let $X$ be an arbitrary metric space, $A_k\in\cP_0(X)$ a fundamental sequence, and $A=\limsup A_k$ \(possibly empty\/\). Then for any $\e>0$ there exists $n$ such that for any $k\ge n$
\begin{enumerate}
\item\label{thm:FlimsupCompleteHalfLim:1} $A\ss U_\e(A_k)$\rom;
\item\label{thm:FlimsupCompleteHalfLim:2} if the space $X$ is complete, then $A\ne\0$ and $A_k\ss U_\e(A)$.
\end{enumerate}
Therefore, in complete $X$ we have $A_k\toH A$ and, hence, there exists $\lim A_k$.
\end{thm}

\begin{proof}
Since the sequence $\Om'=(A_1, A_2,\ldots)$ is fundamental, Proposition~\ref{prop:fundANDdk} implies that $d_k(\Om')\to0$. Choose $n$ such that $d_n(\Om')<\e/2$. Consider an arbitrary $k\ge n$ and put $\Om=(A_k,A_{k+1},\ldots)$. Since the sequence $d_i(\Om')$ is monotonic, we have $d_1(\Om)<\e/2$.

(\ref{thm:FlimsupCompleteHalfLim:1}) If $A=\0$, then the inclusion is proved. Now let $A\ne\0$. We choose an arbitrary $a\in A$, then, by Corollary~\ref{cor:Flimsup}, there exists a sequence $\om=(a_{i_1},a_{i_2},\ldots)$ converging to $a$ so that for some $p\in\N$ we have $|a_{i_p}a|<\e/2$. On the other hand, since $d_1(\Om)<\e/2$, then $d_H(A_k,A_{i_p})<\e/2$, therefore there exists $a_k\in A_k$ for which $|a_ka_{i_p}|<\e/2$. By the triangle inequality, $|a_ka|<\e$, therefore, because $a$ is arbitrary, we have $A\ss U_\e(A_k)$.

(\ref{thm:FlimsupCompleteHalfLim:2}) Let us apply Proposition~\ref{prop:eachElA1inFund} to the sequence $\Om$. By this proposition, each $a_k\in A_k$ is included in some fundamental sequence $\om$ for which $d_1(\om)<d_1(\Om)+\e/2<\e$. Due to the completeness of the space $X$ and Proposition~\ref{prop:Flimsup}, the sequence $\om$ converges to some element $a\in A$, hence, taking into account the previous estimate on $d_1(\om)$, we get $|a_ka|<\e$, therefore $A_k\ss U_\e(A)$. The existence of the limit follows from Theorem~\ref{thm:convergeFlimExists}.
\end{proof}

\begin{cor}\label{cor:completenessInherit}
A metric space $X$ is complete if and only if $\cH(X)$ is complete.
\end{cor}

\begin{proof}
The completeness of $\cH(X)$ for complete $X$ is proved in Theorem~\ref{thm:FlimsupCompleteHalfLim}. Conversely, let $\cH(X)$ be a complete space. Consider an arbitrary fundamental sequence $a_k$ in $X$; then, due to the isometry of the mapping $x\mapsto\{x\}$, the sequence $A_k=\{a_k\}\in\cH(X)$ is also fundamental. Since $\cH(X)$ is complete, the sequence $A_k$ converges to some $A\in\cH(X)$. However, by Corollary~\ref{cor:convergeOnePointedSeq}, the sequence $a_k$ also converges.
\end{proof}

\subsection{Inheritance of total boundedness and compactness}

The purpose of this section is to prove the following theorem.

\begin{thm}\label{thm:Hausdorff-compactness}
Let $X$ be an arbitrary metric space. Then the following properties are simultaneously present or not in both $X$ and $\cH(X)$\/\rom:
\begin{enumerate}
\item\label{thm:Hausdorff-compactness:1} completeness\/ \(Khan\),
\item\label{thm:Hausdorff-compactness:2} total boundedness,
\item\label{thm:Hausdorff-compactness:3} compactness\/ \(Hausdorff, Blaschke\).
\end{enumerate}
\end{thm}

\begin{proof}
(\ref{thm:Hausdorff-compactness:1}) This is Corollary~\ref{cor:completenessInherit}.

(\ref{thm:Hausdorff-compactness:2}) Let $X$ be totally bounded. The total boundedness of $\cH(X)$ immediately follows from the following lemma and problem.

\begin{lem}
Let $W$ be an arbitrary metric space, $\e>0$, and $Y\ss W$ be some $\e$-net, then for any $\dl>\e$ the set $\cP_0(Y)$ is a $\dl$-net in $\cP_0(W)$.
\end{lem}

\begin{proof}
Indeed, we choose an arbitrary $F\in\cP_0(W)$ and put $M=\{y\in Y:|yF|<\e\}$, thus for every $y\in M$ there exists $w\in F$ such that $|yw|<\e$, therefore, $M\ss U_\e(F)$.

Further, since $Y$ is an $\e$-net for $W$, then for every $w\in F$ there exists $y\in Y$ for which $|yw|<\e$, hence $|yF|\le|yw|<\e$ and, therefore, $y\in M$. In particular, $M\ne\0$, i.e., $M\in\cP_0(W)$. Since $w$ is arbitrary, we get $F\ss U_\e(M)$, therefore, $d_H(F,M)\le\e<\dl$.
\end{proof}

\begin{prb}
If $W$ is a metric space, $Y\ss W$ is an $\e$-net, and $Z\ss W$ is not empty. Then $Z$ contains an $(2\e)$-net $S$ such that $\#S\le\# Y$.
\end{prb}

Conversely, let $\cH(X)$ be totally bounded. The total boundedness of $X$ follows from the following lemma.

\begin{lem}
Let $W$ be an arbitrary metric space, $\e>0$, and $\cY\ss\cH(W)$ be some $\e$-net. Denote by $M\ss W$ the set obtained by choosing at each element $Y\in\cY$ any one point $p(Y)$. Then $M$ is an $\e$-net in $W$ \(of the same cardinality as $\cY$\).
\end{lem}

\begin{proof}
We choose an arbitrary $w\in W$, then there exists $Y\in\cY$ such that $d_H\bigl(\{w\},Y\bigr)<\e$. Hence $Y\ss U_\e(w)$ and, therefore, $\bigl|p(Y)w\bigr|<\e$.
\end{proof}

(\ref{thm:Hausdorff-compactness:3}) This follows from the previous section and Theorem~\ref{thm:CompleteCompact}.
\end{proof}

\section{Inheritance of geodesic}
\markright{\thesection.~Inheritance of geodesic}
In this section we present the results from~\cite{Bryant-convexity5}.

\begin{thm}
Let $X$ be a compact metric space, then $\cH(X)$ is geodesic if and only if $X$ is geodesic.
\end{thm}

\begin{proof}
First, let $X$ be a geodesic space. To prove that $\cH(X)$ is geodesic, we use Theorems~\ref{thm:Hausdorff-compactness} and~\ref{thm:MiddlePoints}. The first of them claims that $\cH(X)$ is compact and, therefore, complete; the second is that $\cH(X)$ is geodesic if for any its elements there exists a midpoints. That is what we will prove.

Consider arbitrary $A,B\in\cH(X)$, $r:=d_H(A,B)$, put $C=B_{r/2}(A)\cap B_{r/2}(B)$, and show that $C$ is a midpoint between $A$ and $B$.

Indeed, $C$ is closed as the intersection of closed sets. We show that $C\ne\0$. Choose an arbitrary $a\in A$. By Proposition~\ref{prop:HausdAttainsForComp}, $A\ss B_r(B)$, therefore, since $B$ is compact, there exists $b\in B$ such that $|ab|\le r$. Since the space $X$ is geodesic, for the points $a$ and $b$ there is a midpoint $c$, then $c\in C$, which proves the nonemptiness of $C$ and, at the same time, that $A\ss B_{r/2}(C)$. In addition, $C\ss B_{r/2}(A)$ by construction, hence $d_H(A,C)\le r/2$. Similarly, $d_H(B,C)\le r/2$. Since $r=d_H(A,B)\le d_H(A,C)+d_H(C,B)$, then $d_H(A,C)=d_H(C,B)=r/2$, so $C$ is a midpoint between $A$ and $B$. We again use Theorem~\ref{thm:MiddlePoints} and conclude that the space $\cH(X)$ is geodesic.

Conversely, let the space $\cH(X)$ be geodesic. We choose arbitrary points $a,b\in X$, $r:=|ab|$, then, due to the geodesicity of $\cH (X)$, it contains $C$, which is a midpoint between $\{a\}$ and $\{b\}$. By Items~(\ref{prb:Hausdorff-evident-properties:1}) of Problem~\ref{prb:Hausdorff-evident-properties}, we have $d_H\bigl(\{a\},\{b\}\bigr)=r$, so $d_H\bigl(\{a\},C\bigr)=d_H\bigl(C,\{b\}\bigr)=r/2$, hence $C\ss B_{ r/2}(a)\cap B_{r/2}(b)$. We choose an arbitrary point $c\in C$, then $|ac|\le r/2$ and $|cb|\le r/2$; since $r=|ab|\le|ac|+|cb|$, then $|ac|=|cb|=r/2$ and, therefore, $c\in X$ is a midpoint between $a$ and $b$. It remains to use Theorem~\ref{thm:MiddlePoints}.
\end{proof}

We give a few examples to demonstrate that for the space $\cH(X)$ was geodesic, neither the total boundedness of $X$, nor its completeness is sufficient individually.

\begin{prb}
Denote by $X$ a subset of the space $\R^3$, which in the standard Cartesian coordinates has the form
$$
X=\{x^2+y^2<1,\,z\le1\}\cup\{x^2+y^2=1,\,z=-1,\,x\in\Q\}\cup\{x^2+y^2=1,\,z=1,\,x\in\R\sm\Q\}.
$$
Next, put $A=\{x^2+y^2=1,\,z=-1,\,x\in\Q\}$ and $B=\{x^2+y^2=1,\,z=1,\,x\in\R\sm\Q\}$. Show that $X$ is a geodesic space, but $A$ and $B$ are not connected by a shortest curve (there is no midpoint between them), so $\cH(X)$ is not geodesic.
\end{prb}

\begin{prb}
Consider the square $X=[-1,1]\x[-1,1]$ and introduce the following distance function on $X$:
$$
d\bigl((x,y),(x',y')\bigr)=|x-x'|+\min\{(1-y)+(1-y'),(y+1)+(y'+1)\}.
$$
Show that $(X,d)$ is a complete geodesic space, $A=[-1,1]\x{-1}$ and $B=[-1,1]\x\{1\}$ are compacts such that there is midpoint between them in $\cK(X)$, i.e., $\cK(X)$ is not geodesic.
\end{prb}

\begin{prb}
Consider two segments $I_0=[0,1]\x{0}$ and $I_1=[0,1]\x\{1\}$ in the Euclidean plane $\R^2$, and let $A\ss I_0$ be the subset of all points with rational abscissae, and $B\ss I_1$ the subset of all points with irrational abscissae. For each $a\in A$ and $b\in B$ we denote by $r(a,b)$ the distance between $a$ and $b$ in $\R^2$. Now, for each $a\in A$, $b\in B$ consider a circle $S(a,b)$ or the length $3$, and denote by the same $a$ and $b$ some points of $S(a,b)$ bounding an arc $\a(a.b)$ of the length $r(a,b)$. Let $C(a,b)\in S(a,b)$ be the middle point of the act complement to $\a(a,b)$. Denote by $X$ the quotient space obtained from all $S(a,b)$ by identifying their points $C(a,b)$. Also denote by $A$ and $B$ the subsets of $X$ consisting of all points $a$ and all points $b$, respectively. Prove that $X$ is complete and geodesic, $A,B\in\cH(X)$, however, there is no a midpoint between $A$ and $B$, thus, $\cH(X)$ is complete but not geodesic.
\end{prb}


\phantomsection
\renewcommand\bibname{References to Chapter~\thechapter}
\addcontentsline{toc}{section}{\bibname}
\renewcommand{\refname}{\bibname}

\vfill\eject

\section*{\Huge Exercises to Chapter~\thechapter}
\markright{Exercises to Chapter~\thechapter.}

\begin{exe}
Prove the following statements for an arbitrary metric space $X$.
\begin{enumerate}
\item Let $f\:X\to\cH(X)$ be given by the formula $f\:x\mapsto\{x\}$, then $f$ is an isometric embedding.
\item For any $A,B\in\cP_0(X)$ we have $d_H(A, B)=d_H(A,\bB)=d_H(\bA,B)=d_H(\bA,\bB)$.
\item For any $A,B\in\cP_0(X)$ we have $d_H(A,B)=0$ if and only if $\bA=\bB$.
\item If $Y\ss X$ is an $\e$-net in $A\ss X$, then $d_H(A,Y)\le\e$.
\end{enumerate}
\end{exe}

\begin{exe}
Prove that for $A,B\in\cK(X)$ there exist $a\in A$ and $b\in B$ such that $d_H(A,B)=|ab|$.  Is it possible to change $\cK(X)$ with $\cH(X)$?
\end{exe}

\begin{exe}
Let $X$ be an arbitrary metric space and $A,B,A',B'\in\cH(X)$ such that $A'\ss A$ and $B'\ss B$. Prove that $d_H(A\cup B',B\cup A')\le d_H(A,B)$.
\end{exe}

\begin{exe}
Let $X$ be an arbitrary metric space and $A,B,C\in\cH(X)$ such that $C\ss B$. Prove that $d_H(A,A\cup C)\le d_H(A,B)$.
\end{exe}

\begin{exe}
Prove that $\liminf A_i$ is a closed subset of $X$.
\end{exe}

\begin{exe}
Let $\ell_2$ denote the space of all sequences $(x_1,x_2,\ldots)$ of real numbers such that $\sum_{i=1}^\infty x^2<\infty$. For $\xi=(x_1,x_2,\ldots)\in\ell_2$ and $\xi'=(x'_1,x'_2,\ldots)\in\ell_2$ we put $d(\xi,\xi')^2={\sum_{i=1}^\infty(x_i-x'_i)^2}$.
Prove that $d$ is a metric on $\ell_2$, and that $(\ell_2,d)$ is a complete metric space.
\end{exe}

\begin{exe}
Denote by $e_i\in\ell_2$ the sequence $(x_1,x_2,\ldots)$ such that $x_i=1$ and $x_j=0$ for all $j\ne i$. Prove that $A=\{e_i\}_{i=1}^\infty$ is a closed subset of $\ell_2$, thus $A\in\cH(\ell_2)$.
\end{exe}

\begin{exe}
Let $X=\N$, and define two metrics on $X$: $d^1(x,y)=1$ for any $x\ne y$, and $d^2(x,y)=\bigl|\frac1x-\frac1y\bigr|$ for any $x,y$. Then the corresponding metric topologies are both discrete. Prove that the corresponding Hausdorff metric generates non-homeomorphic topologies on $\CL(X)$.
\end{exe}

\begin{exe}
Suppose that in a sequence $A_i\in\cH(X)$ all $A_i$ are connected, $A_i\to A\in\cH(X)$, and $\liminf A_i\ne\0$. Prove that $\limsup A_i$ is connected.
\end{exe}

\begin{exe}
Denote by $X$ a subset of the space $\R^3$, which in the standard Cartesian coordinates has the form
$$
X=\{x^2+y^2<1,\,z\le1\}\cup\{x^2+y^2=1,\,z=-1,\,x\in\Q\}\cup\{x^2+y^2=1,\,z=1,\,x\in\R\sm\Q\}.
$$
Next, put $A=\{x^2+y^2=1,\,z=-1,\,x\in\Q\}$ and $B=\{x^2+y^2=1,\,z=1,\,x\in\R\sm\Q\}$. Show that $X$ is a geodesic space, but $A$ and $B$ are not connected by a shortest curve (there is no midpoint between them), so $\cH(X)$ is not geodesic.
\end{exe}

\begin{exe}
Consider the square $X=[-1,1]\x[-1,1]$ and introduce the following distance function on $X$:
$$
d\bigl((x,y),(x',y')\bigr)=|x-x'|+\min\{(1-y)+(1-y'),(y+1)+(y'+1)\}.
$$
Show that $(X,d)$ is a complete geodesic space, $A=[-1,1]\x{-1}$ and $B=[-1,1]\x\{1\}$ are compacts such that there is no midpoint between them in $\cK(X)$, i.e., $\cK(X)$ is not geodesic.
\end{exe}

\begin{exe}
Consider two segments $I_0=[0,1]\x{0}$ and $I_1=[0,1]\x\{1\}$ in the Euclidean plane $\R^2$, and let $A\ss I_0$ be the subset of all points with rational abscissae, and $B\ss I_1$ the subset of all points with irrational abscissae. For each $a\in A$ and $b\in B$ we denote by $r(a,b)$ the distance between $a$ and $b$ in $\R^2$. Now, for each $a\in A$, $b\in B$ consider a circle $S(a,b)$ or the length $3$, and denote by the same $a$ and $b$ some points of $S(a,b)$ bounding an arc $\a(a.b)$ of the length $r(a,b)$. Let $C(a,b)\in S(a,b)$ be the middle point of the act complement to $\a(a,b)$. Denote by $X$ the quotient space obtained from all $S(a,b)$ by identifying their points $C(a,b)$. Also denote by $A$ and $B$ the subsets of $X$ consisting of all points $a$ and all points $b$, respectively. Prove that $X$ is complete and geodesic, $A,B\in\cH(X)$, however, there is no a midpoint between $A$ and $B$, thus, $\cH(X)$ is complete but not geodesic.
\end{exe}

\vfill\eject
 \chapter{Gromov--Hausdorff Distance.}\label{chap:Gromov_Hausdorff_metric}
 \markboth{\chaptername~\thechapter.~Gromov--Hausdorff Distance.}%
          {\chaptername~\thechapter.~Gromov--Hausdorff Distance.}

\begin{plan}
Realization of a pair of metric space, Gromov--Hausdorff distance, admissible metric on disjoin union of metric spaces, calculation of Gromov--Hausdorff distance in terms of admissible metrics, triangle inequality for Gromov--Hausdorff distance, positive definiteness of Gromov--Hausdorff distance for isometry classes of compact spaces, counterexample for boundedly compact spaces, Gromov--Hausdorff distance for separable spaces in terms of their isometric images in $\ell_\infty$, relations, distortion of a relation between metric spaces, correspondences, Gromov--Hausdorff distance in terms of correspondences, $\e$-isometries and Gromov--Hausdorff distance, irreducible correspondences between sets, existence of irreducible correspondences, irreducible correspondences as bijections of partitions of sets, Gromov--Hausdorff distance in terms of irreducible correspondences, examples: Gromov--Hausdorff distance between $2$- or $3$-point metric spaces, simple general properties of Gromov--Hausdorff distance, Gromov--Hausdorff convergence, inheritance of metric and topological properties while Gromov--Hausdorff convergence.
\end{plan}

In this section, we will study the Gromov--Hausdorff distance, see~\cite{BurBurIva6} and~\cite{Ghanaat6} for more details. All metric spaces are supposed to be nonempty.

\begin{constr}\label{constr:GHdistance}
Let $X$ and $Y$ be metric spaces. A triple $(X',Y',Z)$ consisting of a metric space $Z$ and its two subsets $X'$ and $Y'$, which are isometric respectively to $X$ and $Y$, will be called \emph{a realization of the pair $(X,Y)$}. We put
$$
d_{GH}(X,Y)=\inf\bigl\{r\in\R:\text{there exists a realization $(X',Y',Z)$ of $(X,Y)$ such that $d_H(X',Y')\le r$}\bigr\}.
$$
\end{constr}

\begin{rk}
The value $d_{GH}(X,Y)$ is evidently non-negative, symmetric, and $d_{GH}(X,X)=0$ for any metric space $X$. Thus, $d_{GH}$ is a distance function on each set of metric spaces.
\end{rk}

\begin{dfn}
The value $d_{GH}(X,Y)$ from Construction~\ref{constr:GHdistance} is called \emph{the Gromov--Hausdorff distance\/} between the metric spaces $X$ and $Y$.
\end{dfn}

\begin{prb}
Prove that for any metric spaces $X$ and $Y$ there exists a realization of $(X,Y)$.
\end{prb}

\begin{rk}
In some monographs the Gromov--Hausdorff distance $d_{GH}(X,Y)$ is defined as the infimum of the values $d_H(X',Y')$ over all realizations $(X',Y',Z)$ of the pair $(X,Y)$. However, we give such, at first glance, a more technically complicated definition of $d_{GH}$, to avoid the Cantor paradox, because the family of all realizations is no longer a set. Introducing $r$ and talking about \emph{the existence of a realization}, we thereby get rid of the need to consider all the realizations.
\end{rk}

\begin{notation}
In what follows, it will sometimes be necessary for us to explicitly indicate the space $X$ in which a particular metric is considered, as well as the corresponding Hausdorff metric. Thus, the distance between points $x,x'\in X$ will sometimes be denoted by $|xx'|_X$, and the corresponding Hausdorff distance between nonempty subsets $A$ and $B$ of the space $X$ by $d_H^X(A,B)$. In addition, if $\r$ is some metric on $X$, then the Hausdorff distance generated by this metric will sometimes be denoted by $\r_H$.
\end{notation}

It turns out that, to define the Gromov--Hausdorff distance, it suffices to consider only metric spaces of the form $(X\sqcup Y,\r)$, where $\r$ extends the original metrics of $X$ and $Y$. Such $\r$ will be called \emph{an admissible metric for $X$ and $Y$}, and the set of all admissible metrics for given $X$ and $Y$ will be denoted by $\cD(X,Y)$.

\begin{prb}
Prove that for any metric spaces $X$ and $Y$ there exists at least one admissible metric, i.e., the set $\cD(X,Y)$ is not empty.
\end{prb}

\begin{thm}\label{thm:GH_admiss_metrics}
For any metric spaces $X$ and $Y$, we have
\begin{equation}\label{eq:GH-alt-1}
d_{GH}(X,Y)=\inf\bigl\{\r_H(X,Y):\r\in\cD(X,Y)\bigr\}.
\end{equation}
\end{thm}

\begin{proof}
Denote by $d'_{GH}(X,Y)$ the right hand side of the equation~(\ref{eq:GH-alt-1}). Then $d_{GH}(X,Y)\le d'_{GH}(X,Y)$, because for every $\r\in\cD(X,Y)$ and $Z=(X\sqcup Y,\r)$, the triple $(X,Y,Z)$ is a realization of the pair $(X,Y)$ with the same Hausdorff distance $\r_H(X,Y)$. We now prove the inverse inequality.

By Construction~\ref{constr:GHdistance}, for any $\e>0$ there exists a realization $(X',Y',Z)$ of the pair $(X,Y)$ such that
$$
d_H(X',Y')\le d_{GH}(X,Y)+\e.
$$
If $X'$ and $Y'$ do not intersect each other, then we restrict the metric from $Z$ to $X'\sqcup Y'$, and, after identifying $X'$ and $Y'$ with $X$ and $Y$, respectively, we obtain an admissible metric $\r$ on $X\sqcup Y$ for which $\r_H(X,Y)\le d_{GH}(X,Y)+\e$. If $X'\cap Y'\ne\0$, we replace $Z$ by $Z\x\R$ with the metric $\bigl|(z,t)(z',s)\bigr|=|zz'|+|ts|$, also replace $X'$ and $Y'$ by the sets $X''=X'\x\{0\}$ and $Y''=Y'\x\{\e\}$, respectively, then $(X'',Y'',Z\x\R)$ is a realization of the $(X,Y)$ such that $X''\cap Y''=\0$, $d_H(X'',Y'')\le d_{GH}(X,Y)+2\e$, and aging we can restrict the metric of $Z\x\R$ to $X''\sqcup Y''$, identify $X$ with $X''$, $Y''$ with $Y$, and obtain an admissible metric $\r$ with $\r_H(X,Y)\le d_{GH}(X,Y)+2\e$. Since $\e$ is arbitrary, we get $d'_{GH}(X,Y)\le d_{GH}(X,Y)$.
\end{proof}

\begin{rk}
If $X$ and $Y$ are subsets of some metric space, then $d_{GH}(X,Y)\le d_H(X,Y)$. In particular, if $d_H(X,Y)=0$, then $d_{GH}(X,Y)=0$, so the Gromov--Hausdorff distance, like the Hausdorff distance, is not positively defined: for example, the Gromov--Hausdorff distance between the segment $[0,1]$ and the interval $(0,1)$ vanishes. However, if we restrict ourselves to compact metric spaces, then $d_{GH}(X,Y)=0$ if and only if $X$ and $Y$ are isometric (see below for the proof).
\end{rk}

\begin{prop}\label{prop:GH-triangle-inequality}
The function $d_{GH}$ satisfies the triangle inequality.
\end{prop}

\begin{proof}
Choose arbitrary metric spaces $X$, $Y$, and $Z$. We have to show that $d_{GH}(X,Z)\le d_{GH}(X,Y)+d_{GH}(Y,Z)$.

Choose any admissible metrics $\mu\in\cD(X,Y)$ and $\nu\in\cD(Y,Z)$. Recall that in Example~\ref{examp:GluingSpace} we have defined the gluing $U\sqcup_fV$ of metric spaces $U$ and $V$ over a mapping $f\:W\to V$, where $W\ss U$. Also, Problem~\ref{prb:GluingSpace} states that if $f$ is isometric, then the restrictions of the metric of $U\sqcup_fV$ to $U$ and $V$ coincides with the original ones. We apply this result to $U=X\sqcup Y$ and $V=Y\sqcup Z$ with the metrics $\mu$ and $\nu$, respectively, and to $W=Y$ and $f=\id\:Y\to Y\ss Y\sqcup Z$.

Consider the metric space $A=(X\sqcup Y)\sqcup_f(Y\sqcup Z)$ with the corresponding quotient metric $\r:=\r_\sim$. By Problem~\ref{prb:GluingSpace}, the restrictions of the metric $\r$ from $A$ to $X$ and $Z$ coincide with the original ones, thus $(X,Z,A)$ is a realization of $(X,Z)$. Also, by definition of $\r$, for any $x\in X$ and $z\in Z$ we have $\r(x,z)\le\mu(x,y)+\nu(y,z)$. Let us put $a=\mu_H(X,Y)$, $b=\nu_H(Y,Z)$, and choose some $\dl>0$. Then for any $x\in X$ there exists $y\in Y$ such that $\mu(x,y)<a+\dl/2$; also, for any $y\in Y$ there exists $z\in Z$ such that $\nu(y,z)<b+\dl/2$, thus for any $x\in X$ there exists $z\in Z$ and, similarly, for any $z\in Z$ there exists $x\in X$ such that $\r(x,z)\le\mu(x,y)+\nu(y,z)<a+b+\dl$, thus $\r_H(X,Z)\le a+b+\dl$. Since $\dl$ is arbitrary, we get $\r_H(X,Z)\le\mu_H(X,Y)+\nu_H(Y,Z)$. Thus,
\begin{multline*}
d_{GH}(X,Z)=\inf_{d\in\cD(X,Z)}d_H(X,Z)\le\inf_\r\r_H(X,Z)\le
\inf_{\mu\in\cD(X,Y),\ \nu\in\cD(Y,Z)}\bigl(\mu_H(X,Y)+\nu_H(Y,Z)\bigr)=\\ =\inf_{\mu\in\cD(X,Y)}\mu_H(X,Y)+\inf_{\nu\in\cD(Y,Z)}\nu_H(Y,Z)=d_{GH}(X,Y)+d_{GH}(Y,Z),
\end{multline*}
where the second equality holds because the values we minimize over $\mu$ and $\nu$ are independent.
\end{proof}

Thus, we have shown that on every set of metric spaces, the function $d_{GH}$ is a generalized pseudometric. If the diameters of all spaces in the family are bounded by the same number, then $d_{GH}$ is a (finite) pseudometric. As we already noted, $d_{GH}$ is not a metric in general. However, if we restrict ourselves to compact metric spaces considered upto isometry, then $d_{GH}$ will already be a metric.

\begin{prop}\label{prop:IsometricCompactsGHzero}
For compact metric spaces $X$ and $Y$, it holds $d_{GH}(X,Y)=0$ if and only if $X$ and $Y$ are isometric.
\end{prop}

\begin{proof}
If $X$ is isometric to $Y$ then we can take $Z=X'=Y'=X$, thus $(Z,X',Y')$ is a realization of $(X,Y)$, and $d_H(X',Y')=0$. Now we prove the converse statement.

Consider a sequence of admissible metrics $d^k\in\cD(X,Y)$ such that $d^k_H(X,Y)<1/k$, $k=1,\,2,\ldots$, then for each $x\in X\ss X\sqcup Y$ there exists $y\in Y\ss X\sqcup Y$ such that $d^k(x,y)<1/k$. We choose any such $y$ and put $I_k(x)=y$. Thus, we have constructed a mapping $I_k\:X\to Y$ (possibly discontinuous). Similarly, we define a mapping $J_k\:Y\to X$.

From the triangle inequality it follows that for any $x,x'\in X$ and $y,y'\in Y$ we have
\begin{equation}\label{eq:estimatesIJ}
\begin{gathered}
d^k\bigl(I_k(x),I_k(x')\bigr)<\frac2k+d^k(x,x'),\ \ d^k\bigl(J_k(y),J_k(y')\bigr)<\frac2k+d^k(y,y'),\\
d^k\bigl(x,J_k\c I_k(x)\bigr)<\frac2k,\ \ d^k\bigl(y,I_k\c J_k(y)\bigr)<\frac2k.
\end{gathered}
\end{equation}

Similarly to what was done in the proof of the Arzela--Ascoli theorem, we construct the ``limit'' mappings $I\:X\to Y$ and $J\:Y\to X$. Namely, we choose in $X$ a countable everywhere dense subset $S=\{x_1,x_2,\ldots\}$; using the Cantor diagonal process, we construct a subsequence $\{I_{k_1},I_{k_2},\ldots\}$ such that for every $i$ the sequence $I_{k_p}(x_i)$ converges to some $I(x_i)\in Y$ (here we use the compactness of $Y$). The mapping $I$ is $1$-Lipschitz because for every $x_i,x_j\in S$ we have
$$
\bigl|I(x_i)I(x_j)\bigr|=d^k\bigl(I(x_i),I(x_j)\bigr)\le\frac2k+d^k(x_i,x_j)=\frac2k+|x_ix_j|,
$$
and since $k$ is arbitrary, we have $\bigl|I(x_i)I(x_j)\bigr|\le|x_ix_j|$. Since $Y$ is complete, by Lemma~\ref{lem:ContinuousExtensionOfLipschitz} we can extend $I$ onto the whole $X$ by continuity to a $1$-Lipschitz mapping (we denote this mapping by the same symbol $I$). We proceed similarly with the sequence $J_k$. Passing to the limit in the last two inequalities~(\ref{eq:estimatesIJ}) as $k\to\infty$, we conclude that for any $x\in X$ and $y\in Y$ we have
\begin{equation}\label{eq:limitIJ_biject}
\bigl|x\,(J\c I)(x)\bigr|=0,\ \ \bigl|y\,(I\c J)(y)\bigr|=0.
\end{equation}
The relations~(\ref{eq:limitIJ_biject}) indicates that $I$ and $J$ are mutually inverse bijections, and the $1$-Lipschitz condition for the both $I$ and $J$ guarantees that them preserve the distance, i.e., $I$ and $J$ are isometric.
\end{proof}

The following exercise gives one possible generalization of Proposition~\ref{prop:IsometricCompactsGHzero}.

\begin{prb}\label{prb:GHcompactCompleteIsometric}
Prove that if a metric space $X$ is compact, a metric space $Y$ is complete, and $d_{GH}(X,Y)=0$, then $X$ is isometric to $Y$.
\end{prb}

\begin{sol}
By Theorem~\ref{thm:CompleteCompact}, it suffices to prove that $Y$ is a totally bounded space. Choose an arbitrary $\e>0$. Since $d_{GH}(X,Y)=0$ then there exists an admissible metric $d\in\cD(X,Y)$ such that $d_H(X,Y)<\e/3$. Since $X$ is compact, there exists a finite $(\e/3)$-net $S\ss X$. Since $d_H(X,Y)<\e/3$, for each $s\in S$ there exists $y_s\in Y$ for which $d(s,y_s)<\e/3$, and for each $y\in Y$ there exists $x\in X$ such that $|xy|<\e/3$. Since $S$ is an $(\e/3)$-net in $X$, there exists $s\in S$ for which $|xs|<\e/3$. Thus
$$
|yy_s|\le|yx|+|xs|+d(s,y_s)<\e,
$$
so the set $\{y_s\}_{s\in S}$ is a finite $\e$-net in $Y$.
\end{sol}

\begin{rk}\label{rk:properGHzeroNotIsometric}
Even if both metric spaces $X$ and $Y$ are boundedly compact, Proposition~\ref{prop:IsometricCompactsGHzero} may not hold. To describe the corresponding example, we denote by $X$ and $Y$ the subsets of the Euclidean plane $\R^2$ constructed as follows. Each of these spaces is obtained from the abscissa by adding vertical segments between $(m,0)$ and $(m,\ell_m)$, $m\in\Z$. In the case of $X$ we put $\ell_m=|\sin m|$, and in the case of $Y$ let $\ell_m=\bigl|\sin(m+1/2)\bigr|$. For the distances on $X$ and $Y$ we take the corresponding internal metrics. Notice that the spaces $X$ and $Y$ can be considered as graphs whose points $(m,0)$ correspond to the vertices of degree $3$, and the points $(m,\ell_m)$ to the vertices of degree $1$.

It is easy to see that for any $\e>0$ there exists $n\in\N$ such that for all $m\in\Z$ it holds
$$
\Bigl|\bigl|\sin(m+n)\bigr|-\bigl|\sin(m+1/2)\bigr|\Bigr|<\e.
$$
Therefore, the subset $X$ can be shifted along the abscissa by the vector $(n,0)$ in such a way that the resulting subset $X'$, which is isometric to $X$, satisfies $d_H(X',Y)<\e$. Thus, $d_{GH}(X,Y)=0$.

On the other hand, it is easy to see that the sets $\bigl\{|\sin m|\bigr\}_{m\in\Z}$ and $\Bigl\{\bigl|\sin(m+1/2)\bigr|\Bigr\}_{m\in\Z}$ do not intersect. Suppose that $X$ and $Y$ are isometric, and $f\:X\to Y$ is an isometry. Since each isometry is a homeomorphism, then $f$ takes the vertices of degree $1$ and $3$ to the vertices of the same degree, thus each segment in $X$ between $(m,0)$ and $\bigl(m,|\sin m|\bigr)$ is mapped by $f$ onto a segment of $Y$ between some $(n,0)$ and $\bigl(n,|\sin(n+1/2)\bigr|\bigr)$. However, since $f$ is an isometry, it has to preserve the lengths of these segments, a contradiction.
\end{rk}

As follows from Theorem~\ref{thm:GH_admiss_metrics}, the Gromov--Hausdorff distance between $X$ and $Y$ measures the least ``discrepancy'' for all possible ``alignments'' of these spaces inside metric spaces constructed on $X\sqcup Y$. A natural question arises: can all these alignments, perhaps for some special classes of metric spaces, be realized inside the same ambient metric space? The following result answers the question for the class of separable spaces.

Consider the metric space $\ell_\infty$ of all bounded sequences introduced in Example~\ref{examp:ellpANDinfty}. Recall that, by Theorem~\ref{thm:Frechet-separable}, each separable metric space is isometrically embedded into $\ell_\infty$.

\begin{prop}\label{prop:GHdistanceIn TermsOfEllInfty}
Let $X$ and $Y$ be separable metric spaces. Then
$$
d_{GH}(X,Y)=\inf d_H^{\ell_\infty}\bigl(\v(X),\psi(Y)\bigr),
$$
where the infimum is taken over all isometric embeddings $\v\:X\to\ell_\infty$ and $\psi\:Y\to\ell_\infty$.
\end{prop}

\begin{proof}
The space $X\sqcup Y$ with an admissible metric $d\in\cD(X,Y)$ is also separable, therefore, by Theorem~\ref{thm:Frechet-separable}, it can be isometrically embedded into $\ell_\infty$.
\end{proof}

For specific calculations of the Gromov--Hausdorff distance, other equivalent definitions of this distance are useful.

Recall that \emph{a relation\/} between sets $X$ and $Y$ is a subset of the Cartesian product $X\x Y$. Recall that $\cP_0(X\x Y)$ denotes the set of all nonempty subsets of $X\x Y$, i.e., the set of all nonempty relations between $X$ and $Y$. Similarly to the case of mappings, for each $\s\in\cP_0(X\x Y)$ and for every $x\in X$ and $y\in Y$, there are defined \emph{the image\/} $\s(x):=\bigl\{y\in Y:(x,y)\in\s\bigr\}$ and \emph{the preimage\/} $\s^{-1}(y)=\bigl\{x\in X:(x,y)\in\s\bigr\}$. Also, for $A\ss X$ and $B\ss Y$ their \emph{image\/} and \emph{preimage\/} are defined as the union of the images and preimages of their elements, respectively.

Let $\pi_X\:X\x Y\to X$ and $\pi_Y\:X\x Y\to Y$ be the canonical projections $\pi_X(x, y)=x$ and $\pi_Y(x,y)=y$. We denote in the same way the restrictions of these mappings to each relation $\s\ss X\x Y$. A relation $R$ between $X$ and $Y$ is called \emph{a correspondence\/} if the restrictions of the canonical projections $\pi_X$ and $\pi_Y$ to $R$ are surjective. In other words, for every $x\in X$ there exists $y\in Y$, and for every $y\in Y$ there exists $x\in X$, such that $(x,y)\in R$. Thus, the correspondence can be considered as a surjective multivalued mapping. The set of all correspondences between $X$ and $Y$ is denoted by $\cR(X,Y)$.

If $X$ and $Y$ are metric spaces, then for each relation $\s\in\cP_0(X\x Y)$ we define its \emph{distortion $\dis\s$} as follows
$$
\dis\s=\sup\Bigl\{\bigl||xx'|-|yy'|\bigr|:(x,y),\,(x',y')\in\s\Bigr\}.
$$

\begin{prb}\label{prb:monotinicityOfDistortion}
Prove that for any $\s_1,\s_2\in\cP_0(X\x Y)$ such that $\s_1\ss\s_2$, we have $\dis\s_1\le\dis\s_2$.
\end{prb}

\begin{prb}\label{prb:isometryANDcorrespondences}
Prove that for $R\in\cR(X,Y)$ it holds $\dis R=0$ if and only if $R$ is an isometry.
\end{prb}

The next two constructions establish a link between correspondences from $\cR(X,Y)$ and admissible metrics on $X\sqcup Y$.

To start with, we consider an arbitrary admissible metric and use it to construct a specific correspondence.

Let $\r\in\cD(X,Y)$ be an arbitrary admissible metric for metric spaces $X$ and $Y$, and suppose that $\r_H(X,Y)<\infty$. Choose arbitrary $r\ge\r_H(X,Y)$ for which the set $R^\r_r=\{(x,y):\r(x,y)\le r\}$ is a correspondence between $X$ and $Y$. Notice that we always can take arbitrary $r>\r_H(X,Y)$. Sometimes, for example, for compact $X$ and $Y$, we can take $r=\r_H(X,Y)$.

\begin{prop}\label{prop:DistottionForMetricCorrespondence}
Under above notations, it holds $\dis R^\r_r\le2r$.
\end{prop}

\begin{proof}
For any $(x,y),\,(x',y')\in R^\r_r$ we have
$$
\bigl||xx'|-|yy'|\bigr|=\bigl|\r(x,x')-\r(y,y')\bigr|\le\r(x,y)+\r(x',y')\le2r.
$$
It remains to pass to supremum in definition of $\dis R^\r_r$.
\end{proof}

Now we start from a correspondence and construct a specific admissible metric.

Consider arbitrary correspondence $R\in\cR(X,Y)$. Suppose that $\dis R<\infty$. Extend the metrics of $X$ and $Y$ upto a symmetric function $\r^R$ defined on pairs of points from $X\sqcup Y$: for $x\in X$ and $y\in Y$ put
$$
\r^R(x,y)=\r^R(y,x)=\inf\bigl\{|xx'|+|yy'|+\frac12\dis R:(x',y')\in R\bigr\}.
$$

\begin{prop}\label{prop:correspondence-to-pseudometric}
Under above notations, let $\dis R>0$, then $\r^R$ is an admissible metric, and $\r^R_H(X,Y)=\frac12\dis R$.
\end{prop}

\begin{proof}
To simplify notation, we put $\r:=\r^R$.

Since $\dis R>0$ then $\r$ is positively defined. Now we verify the triangle inequality. It suffices to consider the case $x_1,x_2\in X$,  $y\in Y$ and to prove the inequalities for the triangle $x_1x_2y$. Due to symmetry reasons, we prove only two inequalities $\r(x_1,y)+\r(x_2,y)\ge|x_1x_2|$ and $\r(x_2,y)+|x_1x_2|\ge\r(x_1,y)$.

Let us start from the first one. Choose arbitrary $(x_1',y_1'),\,(x_2',y_2')\in R$, then
$$
|x_1x_1'|+|y_1'y|+|x_2x_2'|+|y_2'y|+\dis R\ge|x_1x_2|-|x_1'x_2'|+|y_1'y_2'|+\dis R\ge|x_1x_2|,
$$
where the last inequality holds because $|y_1'y_2'|-|x_1'x_2'|\ge-\dis R$ by definition. Thus,
$$
\r(x_1,y)+\r(x_2,y)=\inf\bigl\{|x_1x_1'|+|y_1'y|+|x_2x_2'|+|y_2'y|+\dis R:(x_1',y_1'),\,(x_2',y_2')\in R\bigr\}\ge|x_1x_2|.
$$

Now we prove the second inequality. We have
\begin{multline*}
\r(x_2,y)+|x_1x_2|=\inf\bigl\{|x_2x_2'|+|y_2'y|+\frac12\dis R+|x_1x_2|:(x_2',y_2')\in R\bigr\}\ge\\
\ge\inf\bigl\{|x_1x_2'|+|y_2'y|+\frac12\dis R:(x_2',y_2')\in R\bigr\}=\r(x_1,y).
\end{multline*}

It remains to prove that  $\r_H(X,Y)=\frac12\dis R$. Since $R\in\cR(X,Y)$ then for each $x\in X$ there exists $y\in Y$, and for each $y\in Y$ there exists $x\in X$, such that $(x,y)\in R$ and, thus, $\r(x,y)=\frac12\dis R$, therefore $\r_H(X,Y)\le\frac12\dis R$. Besides that, for any $x\in X$ and $y\in Y$ it holds $\r(x,y)\ge\frac12\dis R$, thus $\r_H(X,Y)\ge\frac12\dis R$, and the proof is complete.
\end{proof}

\begin{thm}\label{thm:GH-metri-and-relations}
For any metric spaces $X$ and $Y$ we have
$$
d_{GH}(X,Y)=\frac12\inf\bigl\{\dis R:R\in\cR(X,Y)\bigr\}.
$$
\end{thm}

\begin{proof}
Denote by $I(X,Y)$ the right-hand side of the equality from the statement of the theorem. We prove first that $d_{GH}(X,Y)\ge I(X,Y)$. If $d_{GH}(X,Y)=\infty$, then the inequality holds.

Now, suppose that $d_{GH}(X,Y)<\infty$. Choose an arbitrary $r>d_{GH}(X,Y)$, then, by Theorem~\ref{thm:GH_admiss_metrics}, there exists $\r\in\cD(X,Y)$ for which $\r_H(X,Y)<r$. Consider the correspondence $R^\r_r$ constructed above. Then, by Proposition~\ref{prop:DistottionForMetricCorrespondence}, we have $\dis R^\r_r\le2r$, thus $I(X,Y)\le r$. Since $r$ is arbitrary, we get $I(X,Y)\le d_{GH}(X,Y)$.

We now prove that $d_{GH}(X,Y)\le I(X,Y)$. If $I(X,Y)=\infty$, then the inequality holds.

Now, suppose that $I(X,Y)<\infty$. If there exists $R\in\cR(X,Y)$ such that $\dis R=0$, then, by Problem~\ref{prb:isometryANDcorrespondences}, we have $d_{GH}(X,Y)=0$, and the equality holds. Now suppose that for all $R\in\cR(X,Y)$ we have $\dis R>0$. Choose an arbitrary $R\in\cR(X,Y)$ such that $\dis R<\infty$ and consider $\r^R\in\cD(X,Y)$ constructed above. By Proposition~\ref{prop:correspondence-to-pseudometric}, we have $\r^R_H(X,Y)=\frac12\dis R$, thus
$$
d_{GH}(X,Y)=\inf_{\r\in\cD(X,Y)}\r_H(X,Y)\le\inf\bigl\{\r^R_H(X,Y):R\in\cR(X,Y),\,\dis R<\infty\bigr\}=I(X,Y).
$$
\end{proof}

Recall that for a relation $\s$ between $X$ and $Y$, and a relation $\theta$ between $Y$ and $Z$, \emph{the composition $\theta\c\s$} is defined by the following condition: $(x,z)\in\theta\c\s$ if and only if there exists $y\in Y$ such that $(x,y)\in\s$ and $(y,z)\in\theta$.

\begin{prb}
Let $X$, $Y$, and $Z$ be metric spaces, $R_1\in\cR(X,Y)$, $R_2\in\cR(Y,Z)$. Prove that
\begin{enumerate}
\item\label{item:GH_relations:1} $R_2\c R_1\in\cR(X,Z)$;
\item\label{item:GH_relations:2} $\dis(R_2\c R_1)\le\dis R_1+\dis R_2$;
\item\label{item:GH_relations:3} derive from the previous items the triangle inequality for the Gromov--Hausdorff distance.
\end{enumerate}
\end{prb}

\begin{sol}
(\ref{item:GH_relations:1}) Let $x$ be an arbitrary point in $X$. Since $R_1$ is a correspondence, there exists $y\in Y$ such that $(x,y)\in R_1$. Since $R_2$ is a correspondence, there exists $z\in Y$ such that $(y,z)\in R_2$. Therefore, $(x,z)\in R_2\c R_1$. Similarly, for any $z\in Z$ there exists $x\in X$ such that $(x,z)\in R_2\c R_1$. Thus, $R_2\c R_1$ is a correspondence.

(\ref{item:GH_relations:2}) For any $(x,z),\,(x',z')\in R_2\c R_1$ there exist $y,\,y'\in Y$ such that $(x,y),\,(x', y')\in R_1$ and $(y,z),\,(y',z')\in R_2$, hence
$$
\bigl||zz'|-|xx'|\bigr|=\bigl||zz'|-|yy'|+|yy'|-|xx'|\bigr|\le\bigl||zz'|-|yy'|\bigr|+\bigl||yy'|-|xx'|\bigr|,
$$
and, passing to the suprema, we obtain what is required.

(\ref{item:GH_relations:3}) Denote by $R(X,Y,Z)$ the subset of $R(X,Z)$ consisting of all correspondences represented in the form $R_2\c R_1$, where $R_1\in R(X,Y)$ and $R_2\in R(Y,Z)$. Then, by Theorem~\ref{thm:GH-metri-and-relations} and the previous item, we have
\begin{multline*}
d_{GH}(X,Z)=\frac12\inf_{R\in\cR(X,Z)}\dis R\le\frac12\inf_{R\in\cR(X,Y,Z)}\dis R=
\frac12\mathop{\inf_{R_1\in\cR(X,Y)}}_{R_2\in\cR(Y,Z)}\dis(R_2\c R_1)\le\\\le
\frac12\mathop{\inf_{R_1\in\cR(X, Y)}}_{R_2\in\cR(Y, Z)}(\dis R_1+\dis R_2)\le
\frac12\inf_{R_1\in\cR(X,Y)}\dis R_1+\frac12\inf_{R_2\in\cR(Y,Z)}\dis R_2=
d_{GH}(X,Y)+d_{GH}(Y,Z),
\end{multline*}
as required.
\end{sol}

We show one more approach to the study of the Gromov--Hausdorff distance.

\begin{dfn}\label{dfn:epsilonIsometry}
A mapping $f\:X\to Y$ of metric spaces is called \emph{an $\e$-isometry\/} if $\dis f<\e$ and $f(X)$ is an $\e$-net in $Y$.
\end{dfn}

\begin{thm}\label{thm:e-isometries}
Let $X$ and $Y$ be arbitrary metric spaces and $\e>0$. Then
\begin{enumerate}
\item\label{thm:e-isometries:1} if $d_{GH}(X,Y)<\e$ then there exists a $(2\e)$-isometry $f\:X\to Y$\rom;
\item\label{thm:e-isometries:2} if there exists an $\e$-isometry $f\:X\to Y$ then it holds $d_{GH}(X,Y)<2\e$.
\end{enumerate}
\end{thm}

\begin{proof}
(\ref{thm:e-isometries:1}) By Theorem~\ref{thm:GH-metri-and-relations}, there exists a relation $R\in\cR(X,Y)$ such that $\dis R<2\e$. For each $x\in X$ we choose an arbitrary $y\in R(x)$ and put $f(x)=y$. Thus, we determined a mapping $f\:X\to Y$, and since $f\ss R$ (here we identify $f$ with its graph), by Problem~\ref{prb:monotinicityOfDistortion} it holds $\dis f\le\dis R<2\e$. Now we choose an arbitrary $y'\in Y$, an arbitrary $x\in R^{-1}(y')$, and let $y=f(x)$. Since $\diam R(x)\le\dis R<2\e$, it follows that $|yy'|<2\e$, therefore $f(X)$ is a $(2\e)$-net.

(\ref{thm:e-isometries:2}) Consider the relation $R=\bigl\{(x, y):|f(x)y|<\e\bigr\}$. Since $f(X)$ is an $\e$-net, $R$ is a correspondence. To estimate the distortion of $R$, we choose arbitrary $(x,y),(x',y')\in R$, then
$$
\bigl||xx'|-|yy'|\bigr|\le\Bigl||xx'|-\bigl|f(x)f(x')\bigr|\Bigr|+\Bigl|\bigl|f(x)f(x')\bigr|-|yy'|\Bigr|\le
\e+\bigl|f(x)y\bigr|+\bigl|y'f(x')\bigr|<3\e,
$$
therefore, $\dis R\le3\e<4\e$ and $d_{GH}(X,Y)\le\frac12\dis R<2\e$.
\end{proof}

\section{Irreducible correspondences}\label{sec:Irreducible}
\markright{\thesection.~Irreducible correspondences}
For arbitrary nonempty sets $X$ and $Y$, a correspondence $R\in\cR(X,Y)$ is called \emph{irreducible\/} if it is a minimal element of the set $\cR(X,Y)$ w.r.t\. the order given by the inclusion relation. The set of all irreducible correspondences between $X$ and $Y$ is denoted by $\cR^0(X,Y)$.

The following result is evident.

\begin{prop}\label{prop:IrredicubleAndDegrees}
A correspondence $R\in\cR(X,Y)$ is irreducible if and only if for any $(x,y)\in R$ it holds
$$
\min\bigl\{\#R(x),\#R^{-1}(y)\bigr\}=1.
$$
\end{prop}

\begin{thm}\label{thm:irreducinle-exists}
Let $X$, $Y$ be arbitrary nonempty sets. Then for every $R\in\cR(X,Y)$ there exists $R^0\in\cR^0(X,Y)$ such that $R^0\ss R$. In particular, $\cR^0(X,Y)\ne\0$.
\end{thm}

\begin{rk}
One could use the standard technique based on the Zorn lemma, but then we need to guarantee that each chain $R_1\sp R_2\sp\cdots$ has a lower bound, i.e., there exists a correspondence that belongs to all $R_i$. However, this, generally speaking, is not true. As an example, consider $X=Y=\N$ and set $R_k=\bigl\{(i,j):\max(i,j)\ge k\bigr\}$. It is clear that every $R_k$ belongs to $\cR(X,Y)$, and that these $R_k$ form a decreasing chain. However, $\cap R_k=\0$, since for any $i$ and $j$ there is $k$ for which $i<k$ and $j<k$, therefore $(i,j)\not\in R_k$.
\end{rk}

\begin{proof}[Proof of Theorem~$\ref{thm:irreducinle-exists}$.]
For each $x\in X$, choose an arbitrary $y\in R(x)$ and define a mapping $f\:X\to Y$ by setting $y=f(x)$. Let us note that $f\ss R$. Put $Y_1=f(X)$ and $Y_2=Y\sm Y_1$.

Now for each $y\in Y_2$ we choose an arbitrary $x\in R^{-1}(y)$ and define a mapping $g\:Y_1\to X$ by setting $x=g(y)$. In this case we have $g^{-1}\ss R$. Put $X_2=g(Y_2)$ and $X_1=X\sm X_2$.

Let $Y_3=f(X_2)$. It is clear that $Y_3\ss Y_1$.

Using $f$ and $g$, we define another relation: $h=f\cup g^{-1}$.

\begin{lem}
We have $h\in\cR(X,Y)$.
\end{lem}

\begin{proof}
By the definition of $f$, for each $x\in X$ it holds $\bigl(x,f(x)\bigr)\in f\ss h$.

Now consider an arbitrary $y\in Y$. If $y\in Y_1$, then, since $Y_1=\im f$, there exists $x\in X$ such that $y=f(x)$, therefore $(x,y)\in f\ss h$. If $y\in Y_2$, then, by the definition of $g$, we have $\bigl(g(y),y\bigr)\in g^{-1}\ss h$.
\end{proof}

Now we define the relation $R^0$ by removing from $h$ some $(x,y)$ for each $y\in Y_3$ according to the following rule:
\begin{enumerate}
\item if $h^{-1}(y)\cap X_1\ne\0$, then we remove $\bigl(h^{-1}(y)\cap X_2\bigr)\x\{y\}$;
\item if $h^{-1}(y)\cap X_1=\0$, i.e., $h^{-1}(y)\ss X_2$, then we remove all the elements from $h^{-1}(y)\x\{y\}$, except any one.
\end{enumerate}

\begin{lem}
We have $R^0\in\cR(X,Y)$.
\end{lem}

\begin{proof}
For every $y\in Y\sm Y_3$ we remove nothing, therefore for such $y$ there always exists $x\in X$ such that $(x,y)\in R^0$.

Now let $y\in Y_3$. If $h^{-1}(y)\cap X_1\ne\0$, then we do not remove $(x,y)$ with $x\in h^{-1}(y)\cap X_1$, thus $(x,y)\in R^0$ for such $x$. If $h^{-1}(y)\cap X_1=\0$, then we do not remove some $x\in h^{-1}(y)$, thus $(x,y)\in R^0$ for such $x$.

Now let us deal with $x\in X$. If $x\in X_1$, then we remove nothing, therefore $(x,y)\in R^0$ for some $y\in Y$. If $x\in X_2$, then, since $X_2=\im g$, there exists $y\in Y_2$ for which $x=g(y)$, but we remove nothing for such $y$, thus $(x,y)\in R^0$.
\end{proof}

\begin{lem}
We have $R^0\in\cR^0(X,Y)$.
\end{lem}

\begin{proof}
It is sufficient to show that for each pair $(x,y)\in R^0$ either $x$, or $y$ do not belong to other pairs.

If $y\in Y_2$, then $y$ is included in the only pair $\bigl(g(y),y\bigr)\in R^0$. If $y\in Y_1\sm Y_3$, then $y=f(x)$ for some $x\in X_1$, however such $x$ is included in the only pair $\bigl(x,f(x)\bigr)\in R^0$.

Finally, now let $y\in Y_3$. If $h^{-1}(y)\cap X_1\ne\0$, then we removed all the pairs of the form $(x',y)$, $x'\in X_2$, so $x\in X_1$, however, each such $x$ is included in exactly one pair, namely, in $\bigl(x,f(x)\bigr)$. If $h^{-1}(y)\cap X_1=\0$, then we removed all the pairs $(x',y)$, $x'\in h^{-1}(y)\ss X_2$, except some one, so $y$ is included in exactly one such pair.
\end{proof}

This lemma completes the proof of the theorem.
\end{proof}

Theorems~\ref{thm:irreducinle-exists} and~\ref{thm:GH-metri-and-relations}, together with Problem~\ref{prb:monotinicityOfDistortion}, implies

\begin{cor}\label{cor:GH-distance-irreducinle}
For any metric spaces $X$ and $Y$ we have
$$
d_{GH}(X,Y)=\frac12\inf\bigl\{\dis R\mid R\in\cR^0(X,Y)\bigr\}.
$$
\end{cor}

Now we give another useful description of irreducible.

\begin{prop}\label{prop:decompose_inrreducible}
For any nonempty sets $X$, $Y$, and each $R\in\cR^0(X,Y)$, there exist and unique partitions $R_X=\{X_i\}_{i\in I}$ and $R_Y=\{Y_i\}_{i\in I}$ of the sets $X$ and $Y$, respectively, such that $R=\cup_{i\in I}X_i\x Y_i$. Moreover, $R_X=\cup_{y\in Y}\bigl\{R^{-1}(y)\bigr\}$, $R_Y:=\cup_{x\in X}\bigl\{R(x)\bigr\}$,
$$
\{X_i\x Y_i\}_{i\in I}=\cup_{(x,y)\in R}\{R^{-1}(y)\x R(x)\},
$$
and for each $i$ it holds $\min\{\#X_i,\#Y_i\}=1$.

Conversely, each set $R=\cup_{i\in I}X_i\x Y_i$, where $\{X_i\}_{i\in I}$ and $\{Y_i\}_{i\in I}$ are partitions of nonempty sets $X$ and $Y$, respectively, such that for each $i$ it holds $\min\{\#X_i,\#Y_i\}=1$, is an irreducible correspondence between $X$ and $Y$.
\end{prop}

\begin{proof}
First, let $R\in\cR^0(X,Y)$. Put $R_X=\cup_{y\in Y}\bigl\{R^{-1}(y)\bigr\}$, $R_Y=\cup_{x\in X}\bigl\{R(x)\bigr\}$, and show that $R_X$ and $R_Y$ are partitions. We show it for $R_Y$ (for $R_X$ the proof is the same).

Notice first that $R_X$ and $R_Y$ are covers because $R$ is a correspondence. Now, suppose to the contrary that $R_Y$ is not a partition, i.e\., there exist two different elements of $R_Y$, say $R(x)$ and $R(x')$, such that $R(x)\cap R(x')\ne\0$. Since $R(x)\ne R(x')$, one of them, say $R(x)$, contains two different elements $y$ and $y'$, such that one of them, say $y$, belongs to $R(x)\cap R(x')$. This means that for the pair $(x,y)$ it holds $\#R(x)\ge2$ and $R^{-1}(y)\ge2$, a contradiction to Proposition~\ref{prop:IrredicubleAndDegrees}.

Now, we put $R_X=\{X_i\}_{i\in I}$. Note that for any $x,x'\in X_i$ we have $R(x)=R(x')$. Indeed, by definition, $X_i=R^{-1}(y)$ for some $y\in Y$, thus $R(x)\cap R(x')\ne\0$ and, therefore, $R(x)=R(x')$ because $R_Y$ is a partition.

Choose arbitrary $i\in I$, $x\in X_i$, and put $Y_i=R(x)$. Then, this definition is correct (does not depend on the choice of $x$). We show that the mapping $\v\:X_i\mapsto Y_i$ is a bijection between $R_X$ and $R_Y$.

If $\v$ is not injective, then there exist $x,x'\in X$ lying in different elements of the partition $R_X$ for which $R(x)=R(x')$. Thus $x,x'\in R^{-1}(y)\in R_X$ for $y\in R(x)$, a contradiction.

Finally, $\v$ is surjective, since for any $Y_i$, $y\in Y_i$, the set $R^{-1}(y)$ is an element of the partition $R_X$. Choose an arbitrary $x\in R^{-1}(y)$, then $R(x)\in R_Y$ contains $y$, therefore $\v\bigl(R^{-1}(y)\bigr)=Y_i$.

Since for any $x,x'\in X_i$ we have $R(x)=R(x')=Y_i$, then $X_i\x Y_i\ss R$. On the other hand, since $R_X$ is a partition of $X$, for any $x\in X$ there exists $X_i\in R_X$ such that $x\in X_i$, therefore each $(x,y)\in R$ is contained in some $X_i\x Y_i$.

To prove the equality $\min\{\#X_i,\#Y_i\}=1$, let us note that for any $(x,y)\in R$, $x\in X_i$, and $y\in Y_i$, we have $X_i=R^{-1}(y)$ and $Y_i=R(x)$, thus, $\{X_i\x Y_i\}_{i\in I}=\cup_{(x,y)\in R}\{R^{-1}(y)\x R(x)\}$. It suffices to apply Proposition~\ref{prop:IrredicubleAndDegrees}.

The uniqueness of the partitions is the standard fact from the set theory.

The converse is trivial.
\end{proof}

Let $X$ be an arbitrary set different from singleton, and $m$ a cardinal number, $2\le m\le\# X$. By $\cD_m(X)$ we denote the family of all possible partitions of the set $X$ into $m$ nonempty subsets.

Now let $X$ be a metric space. Then for each $D=\{X_i\}_{i\in I}\in\cD_m(X)$ we put
$$
\diam D=\sup_{i\in I}\diam X_i.
$$
Further, for any nonempty $A,B\ss X$, we have already defined $|AB|$ as $\inf\bigl\{|ab|:(a,b)\in A\x B\bigr\}$. We also need $|AB|':=\sup\bigl\{|ab|:(a,b)\in A\x B\bigr\}$. Further, for each $D=\{X_i\}_{i\in I}\in\cD_m(X)$ we put
$$
\a(D)=\inf\bigl\{|X_iX_j|:i\ne j\bigr\}\quad\text{and}\quad\b(D)=\sup\bigl\{|X_iX_j|':i\ne j\bigr\}.
$$
Also notice that $|X_iX_i|=0$, $|X_iX_i|'=\diam X_i$ and hence $\diam D=\sup_{i\in I}|X_iX_i|'$.

The next result follows easily from the definition of distortion, as well as from Proposition~\ref{prop:decompose_inrreducible}.

\begin{prop}\label{prop:disRforPartition}
Let $X$ and $Y$ be arbitrary metric spaces, $D_X=\{X_i\}_{i\in I}$, $D_Y=\{Y_i\}_{i\in I}$, $\#I\ge2$, be some partitions of the spaces $X$ and $Y$, respectively, and $R=\cup_{i\in I}X_i\x Y_i\in\cR(X,Y)$. Then
\begin{multline*}
\dis R=\sup\bigl\{|X_iX_j|'-|Y_iY_j|,\,|Y_iY_j|'-|X_iX_j|: i,j\in I\bigr\}=\\
=\sup\bigl\{\diam D_X,\,\diam D_Y,\,|X_iX_j|'-|Y_iY_j|,\,|Y_iY_j|'-|X_iX_j|: i,j\in I,\, i\ne j\bigr\}\le\\
\le\max\bigl\{\diam D_X,\diam D_Y,\b(D_X)-\a(D_Y),\b(D_Y)-\a(D_X)\bigr\}.
\end{multline*}
In particular, if $R\in\cR^0(X,Y)$, then in the previous formula we can take $R_X$ and $R_Y$ defined in Proposition~$\ref{prop:decompose_inrreducible}$ for $D_X$ and $D_Y$, respectively.
\end{prop}

It will also be convenient for us to represent a relation $\s\in\cP_0(X\x Y)$ as a bipartite graph. Then the degree $\deg$ of each vertex is defined: $\deg_\s(x)=\#\s(x)$ and $\deg_\s(y)=\#\s^{-1}(y)$ .

\begin{prb}\label{prb:nonred_property}
Let $R\in\cR^0(X,Y)$, $x\in X$, $\deg_R(x)>1$. Prove that for each $x'\in X$, $x'\ne x$, it holds $R(x)\cap R(x')=\0$.
\end{prb}

Problem~\ref{prb:nonred_property} immediately implies

\begin{cor}\label{cor:nonred_property}
Let $\#X\ge2$ and $\#Y\ge2$, then for any $R\in\cR^0(X,Y)$ there is no $x\in X$ such that $\{x\}\x Y\ss R$.
\end{cor}

\begin{examp}
Let $X=\{x_1,x_2\}$ and $Y=\{y_1,y_2\}$, then, by Corollary~\ref{cor:nonred_property}, the set $\cR^0(X,Y)$ consists only of bijections, therefore
$$
2d_{GH}(X,Y)=\bigl||x_1x_2|-|y_1y_2|\bigr|.
$$
\emph{Thus, the set of isometric classes of two-point metric spaces endowed with the Gromov--Hausdorff distance is a metric space isometric to the open ray $x>0$ of the real line $\R$ with coordinate $x$}.
\end{examp}

\begin{examp}
Consider two three-point metric spaces $X=\{x_1,x_2,x_3\}$ and $Y=\{y_1,y_2,y_3\}$. Put $\r_{ij}=|x_ix_j|=:\r_k$ and $\nu_{ij}=|y_iy_j|=:\nu_k$, where $\{i,j,k\}=\{1,2,3\}$. Without loss of generality, we assume that $\r_1\le\r_2\le\r_3$ and $\nu_1\le\nu_2\le\nu_3$. We will show that
$$
d_{GH}(X,Y)=\frac12\max\bigl\{|\r_1-\nu_1|,|\r_2-\nu_2|,|\r_3-\nu_3|\bigr\}.
$$

By Corollary~\ref{cor:GH-distance-irreducinle}, it suffices to describe all irreducible correspondences $R_m\in\cR^0(X,Y)$. We will use Proposition~\ref{prop:decompose_inrreducible}. There are three types of partitions of a three-point space: (1) into one-point subsets, (2) into one two-point subset and one-point subset, and (3) into one three-point subset. By Proposition~\ref{prop:decompose_inrreducible}, an irreducible correspondence defines a bijection between partitions, therefore both spaces have to be partitioned in the same way. Case (3) is not realized by Corollary~\ref{cor:nonred_property}, in case (1) we have a bijection, and in case (2) one-point and two-point subsets have to correspond to each other. Thus, we have two types of correspondences: if $\{i,j,k\}=\{p,q,r\}=\{1,2,3\}$, then
\begin{enumerate}
\item $R_1=\bigl\{(x_i,y_p),(x_j,y_q),(x_k,y_r)\bigr\}$,
\item $R_2=\bigl\{(x_i,y_q),(x_i,y_r),(x_j,y_p),(x_k,y_p)\bigr\}$.
\end{enumerate}
We have
\begin{align*}
&\dis R_1=\max\bigl\{|\r_{ij}-\nu_{pq}|,|\r_{ik}-\nu_{pr}|,|\r_{jk}-\nu_{qr}|\bigr\},\\
&\dis R_2=\max\bigl\{|\r_{ij}-\nu_{pq}|,|\r_{ik}-\nu_{pr}|,\r_{jk},\nu_{qr},|\r_{ij}-\nu_{pr}|,|\r_{ik}-\nu_{pq}|\bigr\}.
\end{align*}
Note that $\max\{\r_{jk},\nu_{qr}\}\ge|\r_{jk}-\nu_{qr}|$, therefore $\dis R_2\ge\dis R_1$, so it is enough to consider only the bijections $R_1$.

We also note that each bijection $R_1$ defines a bijection between the sets $\{\r_a\}$ and $\{\nu_a\}$. We show that the bijection $\cup_{i=1}^3\bigl\{(\r_i,\nu_i)\bigr\}$ is optimal, i.e., if we replace this bijection with any other bijection $\psi$, we cannot get a lower value than $M=\max_i|\r_i-\nu_i|$. Put $M'=\max_i\bigl|\r_i-\psi(\r_i)\bigr|$. We need to show that $M'\ge M$.

Let $M=|\r_1-\nu_1|$, and suppose, without loss of generality, that $\r_1\le\nu_1$, thus $M=\nu_1-\r_1$. Since $\psi(\r_1)\ge\nu_1$, then $M'\ge\psi(\r_1)-\r_1\ge M$. Similarly we deal with the case $M=|\r_3-\nu_3|$.

Now, let $M=|\r_2-\nu_2|$, and $\r_2\le\nu_2$, i.e., $M=\nu_2-\r_2$. If $\psi(\r_2)\ne\nu_1$, them $M'\ge M$. Suppose now that $\psi(\r_2)=\nu_1$, then $\psi(\r_1)$ equal either $\nu_2$, or $\nu_3$. Thus $M'\ge\psi(\r_1)-\r_1\ge\nu_2-\r_1\ge\nu_2-\r_2=M$.

Thus, the set of isometry classes of three-point metric spaces endowed with the Gromov--Hausdorff distance is a metric space isometric to the polyhedral cone
$$
\bigl\{(x,y,z):0<x\le y\le z\le x+y\bigr\}
$$
in the space $\R^3_{\infty}$, where the latter denotes the space $\R^3$ with the metric generated by the norm $\bigl\|(x,y,z)\bigr\|_{\infty}=\max\bigl\{|x|,|y|,|z|\bigr\}$. The isometry is given by the formula
$$
\{x_1,x_2,x_3\}\mapsto\frac12(\r_1,\r_2,\r_3).
$$
\end{examp}

\section{A few more examples}
\markright{\thesection.~A few more examples}

The following statement immediately follows from the definition of the Gromov--Hausdorff distance and Item~(\ref{prb:Hausdorff-evident-properties:5}) of Problem~\ref{prb:Hausdorff-evident-properties}.

\begin{examp}\label{examp:epsilon-net}
Let $Y$ be an arbitrary $\e$-net of a metric space $X$. Then $d_{GH}(X,Y)\le d_H(X,Y)\le\e$. Thus, every compact metric space is approximated (according to the Gromov-Hausdorff metric) with any accuracy by finite metric spaces.
\end{examp}

\begin{examp}\label{examp:GH_simple}
Denote by $\D_1$ a one-point metric space. Then for any metric space $X$ we have
$$
d_{GH}(\D_1,X)=\frac12\diam X.
$$
Indeed, $\cR(\D_1,X)$ consists of exactly one correspondence $R$, namely, of $R=\D_1\x X$, thus, $\dis R=\diam X$, and it remains to use Theorem~\ref{thm:GH-metri-and-relations}.
\end{examp}

\begin{examp}\label{examp:dGHbelowEstimate}
Let $X$ and $Y$ be some metric spaces, and the diameter of one of them is finite. Then
$$
d_{GH}(X,Y)\ge\frac12|\diam X-\diam Y|.
$$
Indeed, if $\diam X<\infty$ and $\diam Y=\infty$, then for any $\r\in\cD(X,Y)$ we have $\r_H(X,Y)=\infty$, otherwise $Y\ss U_r(X)$ for some finite $r$ and, therefore, $\diam Y<\infty$.

Now let the both $\diam X$ and $\diam Y$ are finite. Then it suffices to use the triangle inequality (Proposition~\ref{prop:GH-triangle-inequality}) for the triple $X$, $Y$, $\D_1$, and Example~\ref{examp:GH_simple}.
\end{examp}

\begin{examp}\label{examp:GH-ineq-max-Diam-X-Y}
Let $X$ and $Y$ be some metric spaces, then
$$
d_{GH}(X,Y)\le\frac12\max\{\diam X,\diam Y\},
$$
in particular, if $X$ and $Y$ are bounded metric spaces, then $d_{GH}(X,Y)<\infty$.

Indeed, if the diameter of one of the spaces $X$, $Y$ is infinite, then the inequality holds. If both spaces are singletons, then everything is also obvious. Now let $0<d:=\max\{\diam X,\diam Y\}<\infty$. Then for $R=X\x Y\in\cR(X,Y)$ we have $\dis R=d$, thus, by Proposition~\ref{prop:correspondence-to-pseudometric}, it holds $\r^R_H(X,Y)=\frac12\dis R=d/2$, therefore, $d_{GH}(X,Y)\le d/2$.
\end{examp}

Recall that for an arbitrary metric space $X$ and a real number $\l>0$, by $\l X$ we denote the metric space obtained from $X$ by multiplying all distances by $\l$. For $\l=0$ we set $\l X=\D_1$.

\begin{examp}\label{examp:dilatation-of-the-same-X}
For any bounded metric space $X$ and any $\l\ge0$, $\mu\ge0$, we have $d_{GH}(\l X,\mu X)=\frac12|\l-\mu|\diam X$, in particularly, for any $0\le a<b$ the curve $\g(t):=t\,X$, $t\in[a,b]$, is shortest.

Indeed, for the identity correspondence $R\in\cR(\l X,\mu X)$ we have $\dis R=|\l-\mu|\diam X$, hence $d_{GH}(\l X,\mu X)\le\frac12|\l-\mu|\diam X$ by Theorem~\ref{thm:GH-metri-and-relations}. On the other hand, by Example~\ref{examp:dGHbelowEstimate} we have
$$
d_{GH}(\l X,\mu X)\ge\frac12\bigl|\diam(\l X)-\diam(\mu X)\bigr|=\frac12|\l-\mu|\diam X.
$$
It remains to note that
$$
|\g|=\sup_{a=t_0<\cdots<t_n=b}\sum_{i=1}^nd_{GH}(t_{i-1}\,X,\,t_i\,X)=\frac12(b-a)\diam X=d_{GH}(a\,X,\,b\,X),
$$
hence $\g$ is shortest.
\end{examp}

\begin{examp}\label{examp:dilatation-and-GH}
Let $X$ and $Y$ be metric spaces, then for any $\l>0$ we have $d_{GH}(\l X,\l Y)=\l\,d_{GH}(X,Y)$. If, in addition, $d_{GH}(X,Y)<\infty$, then the equality holds for all $\l\ge0$.

Indeed, let $\l>0$. Then for each correspondence $R\in\cR (X,Y)$ and the correspondence $R_\l\in\cR(\l X,\l Y)$, which coincides with $R$ as a set, we have $\dis R_\l=\l\dis R$. It remains to use Theorem~\ref{thm:GH-metri-and-relations}. We now verify that under the condition $d_{GH}(X,Y)<\infty$ the equality holds for $\l=0$ as well. With this $\l$, we have $\l X=\l Y=\D_1$, therefore $d_{GH}(\l X,\l Y)=0$. Since $d_{GH}(X,Y)<\infty$, the value $\l\,d_{GH}(X,Y)$ vanishes too.
\end{examp}

\section{$\GH$-convergence and $\GH$-limits}
\markright{\thesection.~$\GH$-convergence and $\GH$-limits}
For a sequence $X_k$ of metric spaces that converges w.r.t\. the Gromov--Hausdorff distance to a metric space $X$, we say for short that the sequence \emph{$\GH$-converges\/} and write this as $X_k\toGH X$. We call $X$ \emph{the $\GH$-limit\/} and denote it by $\GHlim_{k\to\infty}X_k$.

Let us present a few simple observations.

\begin{enumerate}
\item Let $X=\GHlim_{k\to\infty}X_k$, $Y=\GHlim_{k\to\infty}X_k$, one of $X$ and $Y$ be compact, and the other one be complete, then $X$ and $Y$ are isometric. Indeed, due to the triangle inequality, we have $d_{GH}(X,Y)=0$. It remains to use Problem~\ref{prb:GHcompactCompleteIsometric}.
\item Hausdorff convergence implies $\GH$-convergence.
\item Each compact metric space is a $\GH$-limit of finite metric spaces (of its finite $1/k$-nets).
\end{enumerate}

\begin{thm}\label{thm:propertiesGHconverge}
Let $X_k\toGH Y$. If all $X_k$, starting from some $k$, have one of the properties listed below, then this property is inherited by the space $Y$\/\rom:
\begin{enumerate}
\item\label{thm:propertiesGHconverge:1} the diameter equals infinity\/\rom;
\item\label{thm:propertiesGHconverge:2} the diameter is bounded by a certain number $D$ \(in fact, $\diam X_k\to\diam Y$\)\rom;
\item\label{thm:propertiesGHconverge:3} separability\/\rom;
\item\label{thm:propertiesGHconverge:4} total boundedness\/\rom;
\item\label{thm:propertiesGHconverge:5} if $Y$ is complete, then bounded compactness\/\rom;
\item\label{thm:propertiesGHconverge:6} if $Y$ is complete, then the intrinsicness.
\end{enumerate}
\end{thm}

\begin{proof}
Without loss of generality, we assume that each of there properties holds for all $k$, and that $d_{GH}(X_k,Y)<1/k$. By Theorem~\ref{thm:GH-metri-and-relations}, for each $k\in\N$ there exists $R_k\in\cR(X_k,Y)$ such that $\dis R_k<2/k$. In addition, by Theorem~\ref{thm:e-isometries}, for each $k\in\N$ there exist $(2/k)$-isometries $f_k\:X_k\to Y$ and $g_k\:Y\to X_k$.

(\ref{thm:propertiesGHconverge:1}) If $\diam Y<\infty$, then, by Example~\ref{examp:dGHbelowEstimate}, we have $d_{GH}(X_k,Y)=\infty$, a contradiction.

(\ref{thm:propertiesGHconverge:2}) Since $g_k$ is $(2/k)$-isometry, then $\dis g_k<2/k$ and, hence, for any $y,y'\in Y$ we have
$\Bigl||yy'|-\bigl|g_k(y)g_k(y')\bigr|\Bigr|<2/k$, thus
$$
|yy'|<\bigl|g_k(y)g_k(y')\bigr|+2/k\le\diam X_k+2/k.
$$
Since $y,y'\in Y$ are arbitrary, we get $\diam Y\le\diam X_k+2/k$. Swapping $g_k$ and $f_k$, we get $\diam X_k\le\diam Y+2/k$, i.e., $|\diam X_k-\diam Y|<2/k$, therefore, $\diam X_k\to\diam Y$ and, thus, $\diam Y\le D$.

(\ref{thm:propertiesGHconverge:3}) In each $X_k$ we choose a countable everywhere dense subset $S_k$. Since $f_k$ is $(2/k)$-isometry, for any $y\in Y$ there exists $x_k\in X_k$ such that $\bigl|yf_k(x_k)\bigr|<2/k$. Since $S_k$ is everywhere dense in $X_k$, there exists $s_k\in S_k$ such that $|x_ks_k|<1/k$. Now the condition $\dis f_k<2/k$ implies that $|f_k(s_k)f_k(x_k)|<3/k$, therefore, $\bigl|f_k(s_k)y\bigr|<5/k$, and hence $f_k(S_k)$ is a countable $(5/k)$-net in $Y$. It remains to note that $\cup_{k=1}^\infty f_k(S_k)$ is a countable everywhere dense subset of $Y$.

(\ref{thm:propertiesGHconverge:4}) In each $X_k$ we choose a finite $(1/k)$-net $S_k$. Repeating word-by-word the reasoning from Item~(\ref{thm:propertiesGHconverge:3}), we obtain that $f_k(S_k)$ is a finite $(5/k)$-net. Since $k$ is arbitrary, we get what is required.

(\ref{thm:propertiesGHconverge:5}) Choose an arbitrary closed bounded $Z\ss Y$. We have to show that it is compact. Since $\dis f_k<2/k$, then $f_k^{-1}(Z)$ is bounded. Let $W$ be the closure of $f_k^{-1}(Z)$, then $W$ is compact, in particular, it contains a finite $(1/k)$-net $S_k$. As was discussed above, $f_k(S_k)$ is a finite $(3/k)$-net in $f(W)\sp Z$ and, thus, there exists a finite $(6/k)$-net $Z_k$ in $Z$. Since $k$ is arbitrary, we get that $Z$ is totally bounded. Since $Z$ is closed, it is compact.

(\ref{thm:propertiesGHconverge:6}) By Theorem~\ref{thm:E-MiddlePoints}, it suffices to show that for any $y,y'\in Y$ and any $\e>0$ there exists an $\e$-midpoint between $y$ and $y'$. Choose $k$ such that $4/k<\e$, and arbitrary $x_k\in R_k^{-1}(y)$, $x'_k\in R_k^{-1}(y')$. Since $X_k$ is intrinsic, there exists a $(1/k)$-midpoint $s$ between $x_k$ and $x'_k$. Choose arbitrary $w\in R_k(s)$. Since $s$ is $(1/k)$-midpoint, and $\dis R_k<2/k$, we have $\bigl||x_ks|-|x_kx'_k|/2\bigr|<1/k$, $\bigl||x'_ks|-|x_kx'_k|/2\bigr|<1/k$,  $\bigl||x_ks|-|yw|\bigr|<2/k$, $\bigl||x'_ks|-|y'w|\bigr|<2/k$, $\bigl||x_kx'_k|-|yy|\bigr|<2/k$, thus $\bigl||yw|-|yy'|/2\bigr|<4/k<\e$ and $\bigl||y'w|-|yy'|/2\bigr|<4/k<\e$, what is required.
\end{proof}

\begin{cor}
Let $X_k\toGH Y$, all $X_k$ are boundedly compact and intrinsic \(such $X_k$ are strictly intrinsic by Corollary~$\ref{cor:LocalCompactInnerCompleteIsStrictlyInner}$\), $Y$ is complete. Then $Y$ is strictly intrinsic. Moreover, if we abandon the bounded compactness property, then $Y$ may be not strictly intrinsic, namely, the GH-limit of strictly intrinsic metric spaces may be not strictly intrinsic.
\end{cor}

\begin{proof}
By Item~(\ref{thm:propertiesGHconverge:5}) and Item~(\ref{thm:propertiesGHconverge:6}) of Theorem\ref{thm:propertiesGHconverge}, the space $Y$ is boundedly compact and intrinsic. It remains to apply Corollary~\ref{cor:LocalCompactInnerCompleteIsStrictlyInner}.

To verify the second statement of the corollary, let us consider a metric graph $Y$ obtained by gluing the ends of the segments $[0,1+1/n]$, $n\in\N$ (all $0$ are glued together, and all ends $1+1/n$ are also glued together). Then $Y$ is complete, intrinsic, but not strictly intrinsic. For $X_k$ we take the space obtained from $Y$ by replacing the segment $[0,1+1/k]$ with the segment $[0,1]$. Then all $X_k$ are complete and strictly intrinsic, however, all they are not boundedly compact.
\end{proof}

\begin{prb}
Let $X=\{x^1,\ldots,x^n\}$ be a finite metric space. Prove that a sequence $X_k$ of metric spaces $\GH$-converges to $X$ if and only if for sufficiently large $k$ there exist partitions $\{X_k^i\}_{i=1}^n\in\cD_n(X_k)$ such that $|X_k^iX_k^j|\to|x^ix^j|$ for any $1\le i,\,j\le n$.
\end{prb}

Let $d_k$ be a sequence of metrics on a nonempty set $X$. We say that the metric spaces $X_k=(X,d_k)$ \emph{uniformly converges\/} to a metric space $Y=(X,d)$ if the functions $d_k$ uniformly converges to $d$, i.e., if
$$
\sup_{x,x'\in X}\bigl|d_k(x,x')-d(x,x')\bigr|\to0.
$$

\begin{prb}[Uniform convergence]
Prove that if metric spaces $X_k$ uniformly converges to a metric space $Y$ then $X_k\toGH Y$.
\end{prb}

Recall that for any Lipschitz mapping $f\:X\to Y$ of metric spaces we defined the dilatation $\dil f$ as the minimal Lipschitz constant for $f$, see Section~\ref{sec:LipschitzMappins}. Also, we defined \emph{bi-Lipschitz mapping\/} between metric spaces as a bijective mapping such that it and its inverse as Lipschitz ones. Now we define \emph{the Lipschitz distance $d_L$} between metric spaces $X$ and $Y$ as the following value:
$$
d_L(X,Y)=\inf_{f\:X\to Y}\log\bigl(\max\{\dil f,\,\dil f^{-1}\}\bigr),
$$
where infimum is taken over all bi-Lipschitz mappings (if there is no such mappings then $d_L(X,Y)=\infty$).

\begin{prb}[Lipschitz convergence]
Prove that if metric spaces $X_k$ converges to a bounded metric space $Y$ w.r.t\. the Lipschitz distance, then $X_k\toGH Y$. Does it remain true without the boundedness assumption?
\end{prb}

\begin{prb}
Let the metric space $X_k$ be obtained from the standard sphere $S^2\ss\R^3$ by removing a ball of radius $1/k$. Prove that $X_k\toGH S^2$. Prove that if we change $S^2$ to the circle $S^1$, then the similar statement does not hold. In the both cases we consider the interior metrics. 
\end{prb}


\phantomsection
\renewcommand\bibname{References to Chapter~\thechapter}
\addcontentsline{toc}{section}{\bibname}
\renewcommand{\refname}{\bibname}

\vfill\eject

\section*{\Huge Exercises to Chapter~\thechapter}
\markright{Exercises to Chapter~\thechapter.}

\begin{exe}
Prove that for any metric spaces $X$ and $Y$ there exists a realization of $(X,Y)$.
\end{exe}

\begin{exe}
Prove that for any metric spaces $X$ and $Y$ there exists at least one admissible metric, i.e., the set $\cD(X,Y)$ is not empty.
\end{exe}

\begin{exe}
Prove that for arbitrary metric spaces $X$ and $Y$, any $\s_1,\s_2\in\cP_0(X\x Y)$ such that $\s_1\ss\s_2$, we have $\dis\s_1\le\dis\s_2$.
\end{exe}

\begin{exe}
Prove that for arbitrary metric spaces $X$ and $Y$, any $R\in\cR(X,Y)$ it holds $\dis R=0$ if and only if $R$ is an isometry.
\end{exe}

\begin{exe}
Let $X$ and $Y$ be arbitrary sets, $R\in\cR^0(X,Y)$, $x\in X$, $\deg_R(x)>1$. Prove that for each $x'\in X$, $x'\ne x$, it holds $R(x)\cap R(x')=\0$.
\end{exe}

\begin{exe}
Let $X=\{x^1,\ldots,x^n\}$ be a finite metric space. Prove that a sequence $X_k$ of metric spaces $\GH$-converges to $X$ if and only if for sufficiently large $k$ there exist partitions $\{X_k^i\}_{i=1}^n\in\cD_n(X_k)$ such that $|X_k^iX_k^j|\to|x^ix^j|$ for any $1\le i,\,j\le n$.
\end{exe}

\begin{exe}[Uniform convergence]
Prove that if metric spaces $X_k$ uniformly converges to a metric space $Y$ then $X_k\toGH Y$.
\end{exe}

\begin{exe}[Lipschitz convergence]
Prove that if metric spaces $X_k$ converges to a bounded metric space $Y$ w.r.t\. the Lipschitz distance, then $X_k\toGH Y$. Does it remain true without the boundedness assumption?
\end{exe}

\begin{exe}
Let the metric space $X_k$ be obtained from the standard sphere $S^2\ss\R^3$ by removing a ball of radius $1/k$. Prove that $X_k\toGH S^2$. Prove that if we change $S^2$ to the circle $S^1$, then the similar statement does not hold. In the both cases we consider the interior metrics.
\end{exe}

\vfill\eject
 \chapter{Gromov--Hausdorff space.}
 \markboth{\chaptername~\thechapter.~Gromov--Hausdorff space.}%
          {\chaptername~\thechapter.~Gromov--Hausdorff space.}

\begin{plan}
Gromov--Hausdorff space (GH-space), distortions of a correspondence and its closure, calculating GH-distance in terms of closed correspondences, compactness of the set of all closed correspondences for compact metric spaces, continuity of distortion for compact metric spaces, optimal correspondences, existence of closed optimal correspondences for compact metric spaces, GH-space is geodesic, cover number and packing number, there relations, total boundness of families of compact metric spaces in terms of cover and packing numbers, isometric embedding of all compact metric spaces from a totally bounded family to the same compact subset of $\ell_\infty$, completeness of GH-space, separability of GH-space, $\mst$-spectrum in terms of GH-distances to simplexes, Steiner problem in GH-space.
\end{plan}

This section describes some geometrical and topological properties of the space consisting the isometry classes of compact metric spaces, endowed with the Gromov--Hausdorff metric.

We denote by $\cM$ the set of all compact metric spaces considered up to isometry (in other words, the set of isometry classes of compact metric spaces). From Propositions~\ref{prop:GH-triangle-inequality} and~\ref{prop:IsometricCompactsGHzero} it follows that the distance $d_{GH}$ is a metric on $\cM$. The metric space $(\cM,d_{GH})$ is called \emph{the Gromov--Hausdorff space}. Recall that by $\D_1$ we denoted a one-point metric space. Note that $\cM$ contains all finite metric spaces, in particular, $\D_1\in\cM$. The results presented in Examples~\ref{examp:GH_simple}--\ref{examp:dilatation-and-GH} lead to the following geometric model of $\cM$, see Figure~\ref{fig:gh-space}.

\ig{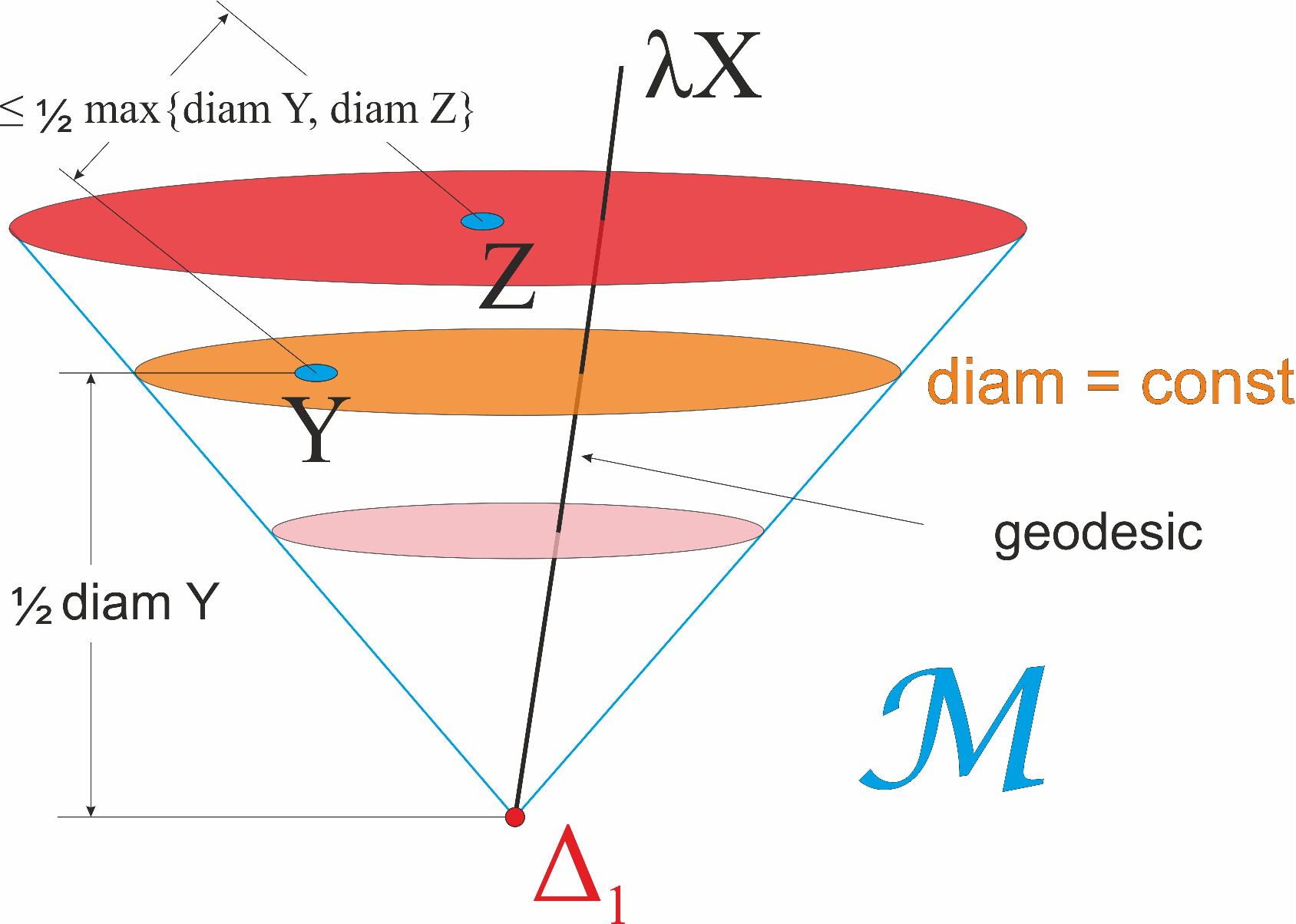}{0.3}{fig:gh-space}{Gromov--Hausdorff space: general properties.}

Indeed, according to Example~\ref{examp:dilatation-and-GH}, the operation of multiplying the metric by a number $\l>0$ is a homothety of the space $\cM$, centered in the one-point metric space $\D_1$, so $\cM$ is a cone with the vertex $\D_1$. By Example~\ref{examp:dilatation-of-the-same-X}, the curves $\g(t):=t\,X$, $X\in\cM$, $X\ne\D_1$, are shortest between any of their points, so they are the generators of the cone. Figure~\ref{fig:gh-space} also illustrates Examples~\ref{examp:GH_simple} and~\ref{examp:GH-ineq-max-Diam-X-Y}.

\section{Existence of optimal correspondences between compact metric spaces}
\markright{\thesection.~Existence of optimal correspondences between compact metric spaces}

Let $X$ and $Y$ be arbitrary metric spaces.

\begin{agree}\label{agree:MetricOnProduct}
In what follows, when we work with the metric space $X\x Y$, we always suppose that its distance function is
$$
\bigl|(x,y)(x',y')\bigr|=\max\bigl\{|xx'|,\,|yy'|\bigr\},
$$
and just this metric generates the Hausdorff distance on $\cP_0(X\x Y)$ and on all its subspaces, e.g., on $\cR(X,Y)$.
\end{agree}

\begin{prb}
Prove that the product topology on $X\x Y$ coincides with the one generated by the metric from Agreement~\ref{agree:MetricOnProduct}.
\end{prb}

\begin{prop}\label{prop:dis_closure}
For arbitrary metric spaces $X$ and $Y$, if $\bsg$ is the closure of $\s\in\cP_0(X\x Y)$, then $\dis\bsg=\dis\s$.
\end{prop}

\begin{proof}
Since $\s\ss\bsg$ then we have $\dis\s\le\dis\bsg$, hence in the case $\dis\s=\infty$ the result holds.

Now let $\dis\s<\infty$. It remains to prove that $\dis\s\ge\dis\bsg$.

By definition, for each $\e>0$ and any $(\bx,\by),\,(\bx',\by')\in\bsg$ there exist $(x,y),\,(x',y')\in\s$ such that $\max\bigl\{|\bx x|,|\by y|\bigr\}<\e/6$ and $\max\bigl\{|\bx'x'|,|\by'y'|\bigr\}<\e/6$, thus $\bigl||\bx\bx'|-|xx'|\bigr|<\e/3$ and $\bigl||\by\by'|-|yy'|\bigr|<\e/3$. Therefore,
$$
\bigl||\bx\bx'|-|\by\by'|\bigr|<\bigl||xx'|-|yy'|\bigr|+2\e/3\le\dis\s+2\e/3.
$$
Passing to the supremum in definition of $\dis\bsg$, we conclude that $\dis\bsg\le\dis\s+2\e/3$. Since $\e$ is arbitrary, we have $\dis\bsg\le\dis\s$.
\end{proof}

Denote by $\cR_c(X,Y)$ the subset of $\cR(X,Y)$ consisting of all closed correspondences $R$. Clearly that for any $R\in\cR(X,Y)$ its closure $\bR$ is also a correspondence, i.e., $\bR\in\cR_c(X,Y)$. This fact, together with Proposition~\ref{prop:dis_closure} and Theorem~\ref{thm:GH-metri-and-relations}, immediately implies

\begin{cor}
For any metric spaces $X$ and $Y$ we have
$$
d_{GH}(X,Y)=\frac12\inf\bigl\{\dis R:R\in\cR_c(X,Y)\bigr\}.
$$
\end{cor}

Now, let $X$ and $Y$ be compact metric spaces. Then $X\x Y$ is compact as well, and $\cR_c(X,Y)\ss\cH(X\x Y)$. By Corollary~\ref{cor:completenessInherit}, $\cH(X\x Y)$ is compact.

\begin{prop}\label{prop:Rc-compact}
For $X,\,Y\in\cM$ the set $\cR_c(X,Y)$ is closed in $\cH(X\x Y)$, thus, $\cR_c(X,Y)$ is compact.
\end{prop}

\begin{proof}
It suffices to show that for each $\s\in\cH(X\x Y)\sm\cR_c(X,Y)$ some its neighborhood in $\cH(X\x Y)$ does not intersect $\cR_c(X,Y)$. Notice that $\s\not\in\cR(X,Y)$ because $\cH(X\x Y)$ consists of all closed nonempty subsets of $X\x Y$, and $\cR_c(X,Y)$ equals to the set of all closed $R\in\cR(X,Y)$. Then either $\pi_X(\s)\ne X$, or $\pi_Y(\s)\ne Y$, where $\pi_X$ and $\pi_Y$ are the canonical projections. To be definite, suppose that the first condition holds, i.e., there exists $x\in X\sm\pi_X(\s)$. Since $\s$ is a closed subset of the compact $X\x Y$, then it is compact itself, and therefore $\pi_X(\s)$ is compact in $X$, thus, $\pi_X(\s)$ is closed. The latter implies that there exists an open ball $U_\e(x)$ such that $U_\e(x)\cap\pi_X(\s)=\0$. Let $U=U_{\e/2}^{\cH(X\x Y)}(\s)$, then for each $\s'\in U$ we have $d_H(\s,\s')<\e/2$, thus for any $(a',b')\in\s'$ there exists $(a,b)\in\s$ such that $\bigl|(a,b),(a',b')\bigr|<\e/2$. Since $a\in\pi_X(\s)$ then $|xa|\ge\e$. On the other hand, $|aa'|\le\bigl|(a,b),(a',b')\bigr|<\e/2$, therefore, $|xa'|>\e/2$, i.e., $a'\not\in U_{\e/2}(x)$ and hence $\pi_X(\s')\cap U_{\e/2}(x)=\0$. This implies that $\s'$ is not a correspondence, thus $\s'\not\in\cR_c(X,Y)$ and $U\cap\cR_c(X,Y)=\0$.
\end{proof}

Define a function $f\:(X\x Y)\x(X\x Y)\to\R$ as $f(x,y,x',y')=\bigl||xx'|-|yy'|\bigr|$. Clearly that $f$ is continuous. Notice that for each $\s\in\cP_0(X\x Y)$ we have
$$
\dis\s=\sup\bigl\{f(x,y,x',y'):(x,y),\,(x',y')\in\s\bigr\}=\sup f|_{\s\x\s}.
$$

\begin{prop}\label{prop:dis_continuous}
If $X,\,Y\in\cM$, then the function $\dis\:\cH(X\x Y)\to\R$ is continuous.
\end{prop}

\begin{proof}
Since $(X\x Y)\x(X\x Y)$ is compact, then the function $f$ is uniformly continuous, i.e., for any $\e>0$ there exists $\dl>0$ such that for every $(x_1,y_1,x'_1,y'_1)$ and $(x_2,y_2,x'_2,y'_2)$ with $\max\bigl\{|x_1x_2|,|y_1,y_2|,|x'_1x'_2|,|y'_1y'_2|\bigr\}<\dl$ it holds $$
\bigl|f(x_1,y_1,x'_1,y'_1)-f(x_2,y_2,x'_2,y'_2)\bigr|<\e.
$$
Thus, for any $\s\in\cH(X\x Y)$ and any $\e>0$ there exists $\dl>0$ such that for the open $\dl$-neighborhood $U=U_\dl^{X\x Y}(\s)\ss X\x Y$ of $\s$ it holds
$$
\sup f|_{U\x U}\le\sup f|_{\s\x\s}+\e.
$$
By $V$ we denote the open ball $U_\dl^{\cH(X\x Y)}(\s)\ss\cH(X\x Y)$ of radius $\dl$ centered at $\s$. Since for any $\s'\in V$ we have $\s'\ss U$, then it follows that
$$
\dis\s'=\sup f|_{\s'\x\s'}\le\sup f|_{U\x U}\le\sup f|_{\s\x\s}+\e=\dis\s+\e.
$$
Swapping $\s$ and $\s'$, we get $|\dis\s-\dis\s'|\le\e$, and hence, the function $\dis$ is continuous.
\end{proof}

\begin{dfn}
A correspondence $R\in\cR(X,Y)$ is called \emph{optimal\/} if $d_{GH}(X,Y)=\frac12\dis R$. By $\cR_\opt(X,Y)$ we denote the set of all optimal correspondences between $X$ and $Y$.
\end{dfn}

\begin{thm}\label{thm:optimal-correspondence-exists}
For any $X,\,Y\in\cM$ we have $\cR_\opt(X,Y)\cap\cR_c(X,Y)\ne\0$.
\end{thm}

\begin{proof}
By Proposition~\ref{prop:dis_continuous}, the function $\dis\:\cR_c(X,Y)\to\R$ is continuous, and by Proposition~\ref{prop:Rc-compact} the space $\cR_c(X,Y)$ is compact, thus the function $\dis$ attains its least value. The half of this value equals $d_{GH}(X,Y)$, and each correspondence $R$ which this least value is attained at is optimal.
\end{proof}

From Theorem~\ref{thm:optimal-correspondence-exists} and Propositions~\ref{prop:correspondence-to-pseudometric} and~\ref{prop:correspondence-to-pseudometric}, we immediately get

\begin{cor}\label{cor:realization}
For each nonisometric $X,\,Y\in\cM$ there exist
\begin{enumerate}
\item a correspondence $R\in\cR_c(X,Y)$ such that $d_{GH}(X,Y)=\frac12\dis R$\rom;
\item an admissible metric $\r\in\cD(X,Y)$ such that $d_{GH}(X,Y)=\r_H(X,Y)$\rom;
\item a metric space $Z$ and isometric embeddings of $X$ and $Y$ into $Z$ such that the Hausdorff distance between their images equals $d_{GH}(X,Y)$.
\end{enumerate}
\end{cor}

\section{The Gromov--Hausdorff space is geodesic}
\markright{\thesection.~The Gromov--Hausdorff space is geodesic}
Let $X$ and $Y$ be arbitrary metric spaces. We assume that the space $X\x Y$ is by default endowed with the metric
$$
\bigl|(x,y)(x',y')\bigr|=\max\bigl\{|xx'|,|yy'|\bigr\}.
$$
Also, we introduce on $X\x Y$ the following $1$-parametric family of metrics $d_t$, $t\in(0,1)$:
$$
d_t\bigl((x,y),(x',y')\bigr)=(1-t)|xx'|+t|yy'|.
$$

\begin{prb}\label{prb:GHgeodesicTopologies}
Prove that the topologies generated by all $d_t$, and by $|\cdot|$ as well, coincide with the product topology of $X\x Y$.
\end{prb}

Now we suppose that $d_{GH}(X,Y)<\infty$, and $R\in\cR(X,Y)$ is any correspondence with $\dis R<\infty$. By $R_t$, $t\in[0,1]$, we denote the metric space $(R,d_t)$ if $t\in(0,1)$, and we put $R_0=X$, $R_1=Y$.

\begin{prop}\label{prop:LipschitzGeodGH}
The mapping $t\mapsto R_t$ is Lipschitzian with the Lipschitz constant $\frac12\dis R$ \(w.r.t\. the Gromov--Hausdorff distance\/\).
\end{prop}

\begin{proof}
We have to show that $2d_{GH}(R_s,R_t)\le|s-t|\dis R$ for arbitrary $s,t\in[0,1]$. Let us start with the case $0<s,t<1$. As always, we denote by $\id$ the identical mapping. Then
\begin{multline*}
2d_{GH}(R_s,R_t)\le\dis\id=\sup\biggl\{\Bigl|d_s\bigl((x,y),(x',y')\bigr)-d_t\bigl((x,y),(x',y')\bigr)\Bigr|:(x,y),(x',y')\in R\biggr\}=\\
=|s-t|\sup\Bigl\{\bigl||xx'|-|yy'|\bigr|:(x,y),(x',y')\in R\Bigr\}=|s-t|\dis R.
\end{multline*}

Now we consider the case $s=0$, $0<t<1$. As a correspondence we take $S_t\in\cR(R_0,R_t)$ of the form
$$
S_t=\bigcup_{(x,y)\in R}\Bigl\{\bigl(x,(x,y)\bigr)\Bigr\}.
$$
Notice that such $S_t$ does not depend on $t$ as a set. We have
\begin{multline*}
2d_{GH}(R_0,R_t)\le\dis S_t=\sup\biggl\{\Bigl||xx'|-d_t\bigl((x,y),(x',y')\bigr)\Bigr|:\bigl(x,(x,y)\bigr),\bigl(x',(x',y')\bigr)\in S_t\biggr\}=\\
=t\sup\Bigl\{\bigl||xx'|-|yy'|\bigr|:\bigl(x,(x,y)\bigr),\bigl(x',(x',y')\bigr)\in S_t\Bigr\}=t\dis R=|s-t|\dis R.
\end{multline*}

All remaining cases can be proved similarly.
\end{proof}

Now, let $X$ and $Y$ be compact metric spaces, and $R\in\cR_c(X,Y)$. Then $R$ is compact as well, and by Problem~\ref{prb:GHgeodesicTopologies}, all metric spaces $R_t$ are compact too. Thus, by Proposition~\ref{prop:LipschitzGeodGH}, the mapping $t\mapsto R_t$ is a Lipschitz curve in $\cM$, and it joins $X$ and $Y$.

\begin{thm}
Given $X,Y\in\cM$ and $R\in\cR_\opt(X,Y)\cap\cR_c(X,Y)$, the curve $\g\:[0,1]\to\cM$, $\g\:t\mapsto R_t$, is a shortest geodesic with the speed $d_{GH}(X,Y)$. In particular, the length of $\g$ equals to $d_{GH}(X,Y)$, thus, the space $\cM$ is geodesic.
\end{thm}

\begin{proof}
Choose arbitrary $0\le s\le t\le 1$, then, by Proposition~\ref{prop:LipschitzGeodGH}, we have
\begin{equation}\label{eq:GHgeodesics}
d_{GH}(R_s,R_t)\le\frac{t-s}2\dis R=(t-s)\,d_{GH}(X,Y),
\end{equation}
thus, taking into account the triangle inequality, we get
$$
d_{GH}(X,Y)\le d_{GH}(X,R_s)+d_{GH}(R_s,R_t)+d_{GH}(R_t,Y)\le d_{GH}(X,Y),
$$
therefore, $d_{GH}(X,R_s)+d_{GH}(R_s,R_t)+d_{GH}(R_t,Y)=d_{GH}(X,Y)$. If for some $s$ and $t$ the inequality in Formula~(\ref{eq:GHgeodesics}) is strict, then the previous equality is not satisfied, thus $d_{GH}(R_s,R_t)=(t-s)\,d_{GH}(X,Y)$.
It remains to note that
$$
|\g|=\sup_{0=t_0<\cdots<t_n=1}\sum_{i=1}^nd_{GH}(R_{t_{i-1}},R_{t_i})=d_{GH}(X,Y),
$$
hence $\g$ is shortest.
\end{proof}

\section{Cover number and packing number}
\markright{\thesection.~Cover Number and Packing Number}

Let $X$ be an arbitrary metric space and $\e>0$. The numerical characteristics of the pair $(X,\e)$ defined below will be used by us in the study of totally bounded families of compact metric spaces, in particular, in terms of these numbers, the Gromov criterion for the precompactness of a family of compact metric spaces will be formulated.

\begin{dfn}
\emph{The cover number}
$$
\cov(X,\e)=\inf\Bigl\{n\in\N:\exists\,x_1,\ldots,x_n\in X,\ \ X=\bigcup_{i=1}^nU_\e(x_i)\Bigr\}
$$
(as usually, we put $\inf\0=\infty$). In other words, the cover number is the minimum number of open balls of radius $\e$ that cover the space $X$.

\emph{The packing number}
$$
\pack(X,\e)=\sup\bigl\{n\in\N:\exists\,x_1,\ldots,x_n\in X\,\forall\,i\ne j\ \ U_{\e/2}(x_i)\cap U_{\e/2}(x_j)=\0\bigr\}.
$$
In other words, the packing number is the maximum number of open pairwise disjoint balls of radius $\e/2$ in the space $X$.
\end{dfn}

\begin{prb}
Prove that
\begin{enumerate}
\item a metric space $X$ is bounded if only for some $\e>0$ it holds $\cov(X,\e)<\infty$ (similarly, $\pack(X,\e)<\infty$);
\item a metric space $X$ is finite if and only if there exists $n$ such that $\cov(X,\e)\le n$ for all $\e>0$ (similarly, for $\pack(X,\e)$);
\item the functions $f(\e)=\cov(X,\e)$ and $g(\e)=\pack(X,\e)$ are monotonically decreasing.
\end{enumerate}
\end{prb}

\begin{prop}\label{prop:CovParkInequality}
For any metric space $X$ and any number $\e> 0$ we have
$$
\cov(X,\e)\le\pack(X,\e)\le\cov(X,\e/4).
$$
\end{prop}

\begin{proof}
First we prove the first inequality. If $\pack(X,\e)=\infty$, then the inequality is automatically satisfied. Now let $\pack(X,\e)<\infty$ and $x_1,\ldots,x_n$, $n=\pack(X,\e)$, be a largest set of points in $X$ for which the balls $U_{\e/2}(x_i)$ are disjoint. Since this family is maximal, for any $x\in X$ there exists $x_k$ such that $U_{\e/2}(x)\cap U_{\e/2}(x_k)\ne\0$, i.e\., $|xx_k|<\e$. Then the family $\{U_\e(x_i)\}_{i=1}^n$ covers $X$, so $\cov(X,\e)\le n=\pack(X,\e)$.

Let us prove the second inequality. Again, if $\cov(X,\e/4)=\infty$, then the inequality holds. Now let $\cov(X,\e/4)<\infty$ and $x_1,\ldots, x_m$, $m=\cov (X,\e/4)$, be the smallest set of points in $X$ for which the balls $U_{\e/4}(x_i)$ cover $X$. Suppose that $\pack(X,\e)>\cov(X,\e/4)$, then there exist $x'_1,\ldots,x'_n$, $n>\cov (X,\e/4)$, such that the balls $U_{\e/2}(x'_i)$ are pairwise disjoint. On the other hand, for some $i\ne j$ there exists $k$ such that $x'_i,x'_j\in U_{\e/4}(x_k)$, therefore, $x_k\in U_{\e/2}(x'_i)\cap U_{\e/2}(x'_j)$, so this intersection is not empty, a contradiction.
\end{proof}

\begin{cor}\label{cor:CovParkFiniteness}
Let $X$ be an arbitrary metric space, then
\begin{enumerate}
\item\label{cor:CovParkFiniteness: 1} if $\pack(X,\e)<\infty$, then $\cov(X,\e)<\infty$\rom;
\item\label{cor:CovParkFiniteness: 2} if $\cov(X,\e)<\infty$, then $\pack(X,4\e)<\infty$.
\end{enumerate}
Thus, $\cov(X,\e)<\infty$ for all $\e>0$, if and only if $\pack(X,\e)<\infty$ for all $\e>0$.
\end{cor}

\begin{prop}\label{prop:CovParkTotalBound}
Let $X$ be an arbitrary metric space. Then the following statements are equivalent\/\rom:
\begin{enumerate}
\item\label{prop:CovParkTotalBound:1} $\cov(X,\e)<\infty$ for any $\e>0$\rom;
\item\label{prop:CovParkTotalBound:2} $\pack(X,\e)<\infty$ for any $\e>0$\rom;
\item\label{prop:CovParkTotalBound:3} the space $X$ is totally bounded.
\end{enumerate}
\end{prop}

\begin{proof}
$(\ref{prop:CovParkTotalBound:1})\lra(\ref{prop:CovParkTotalBound:2})$ This follows from Corollary~\ref{cor:CovParkFiniteness}.

$(\ref{prop:CovParkTotalBound:1})\lra(\ref{prop:CovParkTotalBound:3})$ The condition $\cov(X,\e)<\infty$ is equivalent to the existence of a finite cover $\{U_\e(x_i)\}_{i=1}^n$, which is equivalent to the existence of a finite $\e$-net $\{x_i\}_{i=1}^n$. Thus, the condition of Item~(\ref{prop:CovParkTotalBound:1}) is equivalent to the total boundedness of the space $X$.
\end{proof}

\begin{prop}\label{prop:CovParkEps}
Let $X$, $Y$ be metric spaces, $\dl>0$, and $d_{GH}(X,Y)<\dl$, then
\begin{enumerate}
\item\label{prop:CovParkEps:1} $\cov(X,\e)\ge\cov(Y,\e+2\dl)$,
\item\label{prop:CovParkEps:2} $\pack(X,\e)\ge\pack(Y,2\e+4\dl)$.
\end{enumerate}
\end{prop}

\begin{proof}
(\ref{prop:CovParkEps:1}) The case $\cov(X,\e)=\infty$ is obvious. Now let $m:=\cov(X,\e)<\infty$ and $\bigl\{U_\e(x_i)\bigr\}_{i=1}^m$ be a cover of $X$. By Theorem~\ref{thm:GH-metri-and-relations}, there exists $R\in\cR(X,Y)$ such that $\dis R<2\dl$. For each $i$, we choose an arbitrary $y_i\in R(x_i)$ and show that the set $\{y_i\}_{i=1}^m$ is an $(\e+2\dl)$-net, thus $\cov(Y,\e+2\dl)\le m=\cov(X,\e)$. So, we take arbitrary $y\in Y$ and choose any $x\in R^{-1}(y)$. Then for some $j$ we have $|xx_j|<\e$. Since $\dis R<2\dl$, then $|yy_j|<\e+2\dl$, as required.

(\ref{prop:CovParkEps:2}) Since the case $\pack(Y,2\e+4\dl)=\infty$ is trivial, we assume that $n:=\pack(Y,2\e+4\dl)<\infty$ and let $\bigl\{U_{\e+2\dl}(y_i)\bigr\}_{i=1}^n$ be a disjoint family of open balls in $Y$. Then for any $i\ne j$ we have $|y_iy_j|\ge\e+2\dl$. For each $i$, we choose an arbitrary $x_i\in R^{-1}(y_i)$. Since $\dis R<2\dl$, we have $|x_ix_j|>\e$, therefore the family $\bigl\{U_{\e/2}(x_i)\bigr\}_{i=1}^n$ is disjoint and, thus, $\pack(X,\e)\ge n=\pack(Y,2\e+4\dl)$.
\end{proof}

\section{Totally bounded families of compact metric spaces}
\markright{\thesection.~Totally bounded families of compact metric spaces}

We will be interested in when a particular family of compact metric spaces is totally bounded. We begin with the following auxiliary statement, which will be needed below. For $n\in\N$ we denote by $\cM_n\ss\cM$ ($\cM_{[n]}\ss\cM$) the set of all metric spaces with at most (respectively, exact) $n$ points. For $D\ge 0$, by $\cM(D)\ss\cM$ we denote the set of all compact metric spaces whose diameters do not exceed $D$. We also put $\cM_n(D)=\cM_n\cap\cM(D)$ and $\cM_{[n]}(D)=\cM_{[n]}\cap\cM(D)$. It is clear that $\cM_n=\cup_{k\le n}\cM_{[k]}$ and $\cM_n(D)=\cup_{k\le n}\cM_{[k]}(D)$.

\begin{prop}
The space $\cM_{[n]}(D)\ss\cM$ is totally bounded.
\end{prop}

\begin{proof}
Given $X\in\cM_{[n]}(D)$, we consider all possible bijections $\nu\:X\to\{1,\ldots,n\}$, and for every such $\nu$ we construct the distance matrix $f(X,\nu)=\r=(\r_{ij})$, where $\r_{ij}=\bigl|\nu^{-1}(i)\nu^{-1}(j)\bigr|$. Let $T$ be the set of all such matrices. It is clear that the mapping $g\:T\to\cM_{[n]}(D)$ such that $g\:f(X,\nu)\mapsto X$ is surjective.

We define on $T$ the distance function generated by the $\ell_\infty$-norm, so $T$ will be considered as a subset of $\R^{n^2}_\infty$. Since for every $i,j$ we have $|\r_{ij}|\le D$, the set $T$ is bounded and, therefore, totally bounded as a subset of $\R^{n^2}_\infty$.

If $X,X'\in\cM_{[n]}(D)$, $\r=f(X,\nu)$ and $\r'=f(X',\nu')$, then $R=(\nu')^{-1}\c\nu$ is a bijective correspondence between $X$ and $X'$, and $|\r\r'|_\infty=\dis R\ge2d_{GH}(X,X')$. Thus, the surjection $g$ is Lipschitzian, therefore, $\cM_{[n]}(D)=g(T)$ is also totally bounded.
\end{proof}

\begin{cor}\label{cor:cMnTotalyBound}
The space $\cM_{n}(D)\ss\cM$ is totally bounded.
\end{cor}

\begin{prb}\label{prb:cNnDIsCompact}
Prove that the set $\cM_n(D)$ is compact, while $\cM_{[n]}(D)$ for $n>1$ is not.
\end{prb}

\begin{thm}\label{thm:GromovPrecomp}
Let $\cC$ be a nonempty subset of $\cM$. Then the following statements are equivalent.
\begin{enumerate}
\item\label{thm:GromovPrecomp:1} There is a number $D\ge0$ and a function $N\:(0,\infty)\to\N$ such that for all $X\in\cC$ we have $\diam X\le D$ and $\pack(X,\e)\le N(\e)$.
\item\label{thm:GromovPrecomp:2} There is a number $D\ge0$ and a function $N\:(0,\infty)\to\N$ such that for all $X\in\cC$ we have $\diam X\le D$ and $\cov(X,\e)\le N(\e)$.
\item\label{thm:GromovPrecomp:3} The space $\cC$ with the metric $d_{GH}$ is totally bounded.
\end{enumerate}
\end{thm}

\begin{proof}
$(\ref{thm:GromovPrecomp:3})\imply(\ref {thm:GromovPrecomp:1})$. Fix an arbitrary $\e>0$. We have to find the corresponding $D$ and $N(\e)$. Since $\cC$ is totally bounded, for any $\dl>0$ there exists a finite $\dl$-net $\cC'\ss\cC$. Choose $\dl$ such that $4\dl<\e$. Since all the spaces lying in $\cC'$ are totally bounded, by Proposition~\ref{prop:CovParkTotalBound}, their packing numbers are finite. In addition, their diameters are finite. Put $D'=\max_{X'\in\cC'}\diam X'$ and $N'(\e)=\max_{X'\in\cC'}\pack(X',\e)$. For an arbitrary $X\in\cC$ there exists $X'\in\cC'$ such that $d_{GH}(X,X')<\dl$. It is easy to see that $\diam X\le\diam X'+2\dl\le D'+2\dl$, so that we can put $D=D'+2\dl$. In addition, by Proposition~\ref{prop:CovParkEps}, it holds $\pack(X,\e)\le\pack(X',\e/2-2\dl)\le N'(\e/2-2\dl)$, so we can put $N(\e)=N'(\e/2-2\dl)$.

$(\ref{thm:GromovPrecomp:1})\lra(\ref {thm:GromovPrecomp:2})$. This immediately follows from Proposition~\ref{prop:CovParkInequality}.

$(\ref {thm:GromovPrecomp:2})\imply(\ref{thm:GromovPrecomp:3})$. Fix some $\e>0$, and for each $X\in\cC$ consider a finite cover of the space $X$ by at most $n=N(\e)$ open balls of radius $\e$. By $F_X^\e$ we denote the set of centers of these balls, then $d_{GH}(X,F_X^\e)\le\e$. In addition, $F_X^\e\in\cM_n(D)$, therefore, by Corollary~\ref{cor:cMnTotalyBound}, the family $\cF^\e=\{F_X^\e\}_{X\in\cC}\ss\cM_n(D)$ is totally bounded. Since for any $X\in\cC$ and any $\e'>\e$ we have $X\in U^\cM_{\e'}(F^\e_X)$, then $\cC\ss U^\cM_{\e'}(\cF^\e)$. Since $\e$ and $\e'$ are arbitrary, we conclude that $\cC$ is also totally bounded (verify that).
\end{proof}

The following theorem allows us to realize all metric spaces from a totally bounded subset of $\cM$ as subsets of some compact subset of $\ell_\infty$.

\begin{thm}[Gromov]\label{thm:TotalBoundIntoCompactEllInfty}
For each totally bounded family $\cC\ss\cM$ there exists a compact $K\ss\ell_\infty$ such that every $X\in\cC$ is isometrically embedded into $K$.
\end{thm}

\begin{proof}
The compact $K$ is constructed as follows. By Theorem~\ref{thm:GromovPrecomp}, there exist $D\ge0$ and $N\:(0,\infty)\to\N$ such that for all $X\in\cC$ we have $\diam X\le D$ and $\cov(X,\e)\le N(\e)$. Choose an arbitrary decreasing sequence of positive numbers $E=\{\e_1,\e_2,\ldots\}$ such that $\sum_{i=1}^\infty\e_i<\infty$. This sequence and the function $N(\e)$ generate a sequence of natural numbers $N_i=N(\e_i)$. These two sequences, together with the number $D$, define the set $F_{D,E}\ss\ell_\infty$ as follows.

\begin{constr}
Put $A=\cup_{j=1}^\infty\bigl(\{1,\ldots,N_1\}\x\cdots\x\{1,\ldots,N_j\}\bigr)$. It is clear that $A$ is a countable set. Let $\ell_\infty(A)=\{f\:A\to\R:\sup|f|<\infty\}$, then $\ell_\infty(A)$ is isometric to $\ell_\infty$. For brevity, instead of $f\bigl((n_1,\ldots,n_j)\bigr)$ we will write $f(n_1,\ldots, n_j)$.

We now define the set $F_{D,E}$, composing it from all $f\:A\to\R$ that satisfy the following conditions:
\begin{enumerate}
\item\label{FDEproperties:1} $0\le f(n_1)\le D$ for all $1\le n_1\le N_1$;
\item\label{FDEproperties:2} $\bigl|f(n_1,\ldots,n_j,n_{j+1})-f (n_1,\ldots,n_j)\bigr|\le\e_j$ for all elements $(n_1,\ldots,n_j,n_{j+1})\in A$.
\end{enumerate}
\end{constr}

\begin{lem}
The set $F_{D,E}$ defined above is a compact subset of $\ell_\infty(A)$.
\end{lem}

\begin{proof}
First, note that for each function $f\in F_{D,E}$ it holds $\sup_{a\in A}\bigl|f(a)\bigr|\le D+\sum_{i=1}^\infty\e_i<\infty$, so $f\in\ell_\infty(A)$. Further, since all the inequalities defining $F_{D,E}$ are non-strict, the set $F_{D,E}$ is closed in $\ell_\infty(A)$. Since $\ell_\infty(A)$ is complete, $F_{D,E}$ is also complete. In addition, the diameter of $F_{D,E}$ is finite (it bounded by the number $D+2\sum_{i=1}^\infty\e_i$).

Put $A_{[k]}=\bigl\{(n_1,\ldots,n_k)\in A\bigr\}$ and $A_k=\cup_{j=1}^kA_j$. We denote by $\pi_k\:\ell_\infty(A)\to\ell_\infty(A_k)$ the canonical projection that maps each function $f\:A\to\R$ to its restriction on $A_k\ss A$, and let $F_k=\pi_k(F_{D,E})$. Note that $F_k$ is a closed and bounded subset of the finite-dimensional vector space $\ell_\infty(A_k)$, therefore $F_k$ is compact.

Define the mapping $\nu\:F_k\to\ell_\infty(A)$ by extending each function $f_k\in F_k$ to the entire set $A$ as follows:
$$
f_k(n_1,\ldots,n_k,n_{k+1},\ldots)=f_k(n_1,\ldots,n_k).
$$
It is clear that $\nu$ is isometric, therefore $F'_k=\nu (F_k)$ is also a compact set.

Put $e_k=\e_k+\e_{k+1}+\ldots$, then $e_k\to\0$ as $k\to\infty$. By Condition~(\ref{FDEproperties:2}), we have $F_{D,E}\ss U^{\ell_\infty(A)}_{e_k}(F'_k)$ for all $k\ge2$, which implies the total boundedness of $F_{D,E}$ (verify that).
\end{proof}

We now take the set $F_{D,2E}$ as $K$ and show that each space $X\in\cC$ can be isometrically embedded into this $K$. We consider points of the form $x_a$, $a\in A$, and again, for brevity, instead of $x_{(n_1,\ldots,n_j)}$ we write $x_{n_1\cdots n_j}$.

Take an arbitrary $X\in\cC$. Since $\cov(X,\e_1)\le N(\e_1)=N_1$, then $X$ contains an $\e_1$-net $\cup_{n_1=1}^{N_1}\bigl\{x_{n_1}\bigr\}$, i.e., the family $\bigl\{U_{\e_1}(x_{n_1})\bigr\}_{n_1=1}^{N_1}$ forms a cover of $X$. Note that some points $x_{n_1}$ may coincide.

Further, since $\cov(X,\e_2)\le N(\e_2)=N_2$, then $X$ contains an $\e_2$-net $\cup_{n_2=1}^{N_2}\bigl\{x'_{n_2}\bigr\}$, i.e., the family $\bigl\{U_{\e_2}(x'_{n_2})\bigr\}_{n_2=1}^{N_2}$ forms a cover of $X$. Fix $n_1$ and choose only those balls $U_{\e_2}(x'_{n_2})$ that satisfy $|x_{n_1}x'_{n_2}|<\e_1+\e_2$. In this way we have got at most $N_2$ balls. Enumerate them and add some copies of them to obtain exactly $N_2$ balls which we denote by $U_{\e_2}(x_{n_1n_2})$, $n_2=1,\ldots,N_2$. So, we have got a cover $\bigl\{U_{\e_2}(x_{n_1n_2})\bigr\}_{n_2=1}^{N_2}$ of the ball $U_{\e_1}(x_{n_1})$. By construction, it holds $|x_{n_1}x_{n_1n_2}|<\e_1+\e_2<2\e_1$.

Continuing this process, at the $j$-th step we get the family of balls $\bigl\{U_{\e_j}(x_a)\bigr\}_{a\in A_{[j]}}$ with $|x_{n_1\cdots n_j}x_{n_1\cdots n_jn_{j+1}}|<2\e_j$.

It is easy to see that the set $\{x_a\}_{a\in A}$ of centers of these balls is a countable everywhere dense subset of $X$ (some $x_a$ may coincide with each other). By Theorem~\ref{thm:Frechet-separable}, the space $X$ can be isometrically embedded into $\ell_\infty(A)$ by associating with each point $x$ the function $f_x\:A\to\R$ defined as follows: $f_x(a)=|xx_a|$.

\begin{lem}
For every $x\in X$ we have $f_x\in F_{D,2E}$.
\end{lem}

\begin{proof}
It is clear that $0\le f_x\le D$, so that Item~(\ref{FDEproperties:1}) from the definition of the set $F_{D,2E}$ is satisfied. Further, for each $(n_1,\ldots,n_j,n_{j+1})$, the point $x_{n_1\cdots n_jn_{j+1}}$ lies in $U_{2\e_j}(x_{n_1\cdots n_j})$, so for every $x\in X$ we have
$$
\bigl|f_x(n_1,\ldots,n_j,n_{j+1})-f_x(n_1,\ldots,n_j)\bigr|=\bigl||x_{n_1\cdots n_jn_{j+1}}x|-|x_{n_1\cdots n_j}x|\bigr|\le|x_{n_1\cdots n_jn_{j+1}}x_{n_1\cdots n_j}|<2\e_j,
$$
therefore, Item~(\ref{FDEproperties:2}) from the definition of the set $F_{D,2E}$ is also fulfilled.
\end{proof}

Thus, the mapping $x\mapsto f_x$ isometrically embeds $X$ into $K$.
\end{proof}

\section{Some other properties of Gromov--Hausdorff space}
\markright{\thesection.~Some other properties of Gromov--Hausdorff space}
In this section we apply the previous results to prove a few more properties of the Gromov--Hausdorff space $\cM$.

\subsection{Completeness of Gromov--Hausdorff space}

Theorem~\ref{thm:TotalBoundIntoCompactEllInfty} implies the following result.

\begin{thm}\label{thm:GHcomplete}
The space $\cM$ is complete.
\end{thm}

\begin{proof}
Consider an arbitrary fundamental sequence $\{X_i\}_{i=1}^\infty\ss\cM$. Then $\{X_i\}_{i=1}^\infty$ is a totally bounded subset of $\cM$. By Theorem~\ref{thm:TotalBoundIntoCompactEllInfty}, there exists a compact set $K\ss\ell_\infty$ into which all $X_i$ can be isometrically embedded. Denote by $Y_i$ the image of $X_i$. By Theorem~\ref{thm:Hausdorff-compactness}, the space $\cH(K)$ of all closed bounded subsets of $K$ is also compact, therefore the sequence $Y_i\in\cH (K)$ contains a convergent subsequence $Y_{n_i}$. Let $Y$ be the limit of this subsequence. Then $Y$ is a nonempty compact metric space and
$$
d_{GH}(X_{n_i},Y)=d_{GH}(Y_{n_i},Y)\le d_H(Y_{n_i},Y)\to0\ \ \text{as $i\to\infty$},
$$
therefore, $X_{n_i}\toGH Y$ and, since the sequence $X_i$ is fundamental, we have $X_i\toGH Y$.
\end{proof}

\subsection{Separability of Gromov--Hausdorff space}

\begin{thm}
The space $\cM$ is separable.
\end{thm}

\begin{proof}
By Corollary~\ref{cor:cMnTotalyBound}, each space $\cM_n(D)$ is totally bounded and, therefore, separable. Since $\cM_n=\cup_{k=1}^\infty\cM_n(k)$, then all $\cM_n$, as well as their union $\cup_{n=1}^\infty\cM_n$, are separable. This last union is the set of all finite metric spaces, which, as noted in Example~\ref{examp:epsilon-net}, is an everywhere dense subset of $\cM$, so that $\cM$ is separable.
\end{proof}

Recall that a complete separable metric space is called \emph{Polish}. Thus, the following result holds.

\begin{cor}
The space $\cM$ is Polish.
\end{cor}

By Problem~\ref{prb:SeparabilityAndCountableBase}, for a metric space, the separability is equivalent to having a countable base.

\begin{cor}
The space $\cM$ has a countable base.
\end{cor}

\section{Calculating $\mst$-spectrum by means of Gromov--Hausdorff distances}
\markright{\thesection.~Calculating $\mst$-spectrum by means of Gromov--Hausdorff distances}

Recall that by $\D_n$ we denoted $n$-point metric space such that all its nonzero distances equal $1$. Also, given $\l>0$ and any metric space $X$, if we multiply by $\l$ all the distances in $X$, then the resulting metric space we denote by $\l X$.

In the present section we show that the $\mst$-spectrum of an arbitrary $n$-point metric space $X$ can be represented as a linear function on the Gromov--Hausdorff distances from this space to the $\l\D_2,\ldots,\l\D_n$ for $\l\ge2\diam X$.

\begin{thm}\label{thm:spectrum-as-GH}
Let $X$ be a finite metric space, $\s(X)=(\s_1,\ldots,\s_{n-1})$, $\l\ge2\diam X$. Then
$$
\s_k=\l-2d_{GH}(\l\D_{k+1},X).
$$
\end{thm}

\begin{proof}
Choose any $1\le k\le n-1$ and arbitrary irreducible correspondence $R\in\cR^0(\l\D_{k+1},X)$. By Proposition~\ref{prop:decompose_inrreducible}, there exists partitions $R_{\l\D_{k+1}}=\{Z_i\}_{i=1}^p$ and $R_X=\{X_i\}_{i=1}^p$ of $\l\,\D_{k+1}$ and $X$, respectively, such that $R=\cup_{i=1}^pZ_i\x X_i$, and $\min\{\#Z_i,\#X_i\}=1$ for all $i$. By Proposition~\ref{prop:disRforPartition}, it holds $\dis R\ge\max\{\diam R_{\l\D_{k+1}},\diam R_X\}$. Thus, if for some $i$ we have $\#Z_i>1$, then $\dis R\ge\l\ge2\diam X$. Since $k+1\le n$, there exists $R$ such that $\#Z_i=1$ for all $i$. For such $R$, again by Proposition~\ref{prop:disRforPartition}, we have $\dis R\le\diam X$. Therefore, $\inf_{R\in\cR^0(\l\D_{k+1},X)}\dis R$ is achieved on a correspondences of the latter type. The set of these correspondences we denote by $\cR$.

Now, if $R\in\cR$, then $p=k+1$ and $R_X\in\cD_{k+1}(X)$. By Proposition~\ref{prop:disRforPartition}, we have
\begin{multline*}
\dis R=\sup\bigl\{\diam R_X,\,|X_iX_j|'-\l,\,\l-|X_iX_j|:1\le i<j\le k+1\bigr\}=\\
=\sup\bigl\{\l-|X_iX_j|:1\le i<j\le k+1\bigr\}=\l-\a(R_X),
\end{multline*}
where the second equality holds because
$$
\max\bigl\{|X_iX_j|'-\l,\diam R_X\bigr\}\le\diam X\le\l-\diam X\le\l-|X_iX_j|
$$
for any $1\le i<j\le k+1$. Corollary~\ref{cor:GH-distance-irreducinle}, together with above considerations, gives us
$$
2d_{GH}(\l\D_{k+1},X)=\min_{R\in\cR}\dis R=\min_{R\in\cR}\bigl(\l-\a(R_X)\bigr)=\l-\max_{D\in\cD_{k+1}(X)}\a(D),
$$
where the last equality holds because each $D$ generates some $R\in\cR$.

It remains to use Theorem~\ref{thm:spect-calc} which states that
$$
\s_k=\max\bigl\{\a(D):D\in\cD_{k+1}(X)\bigr\},
$$
thus, $2d_{GH}(\l\D_{k+1},X)=\l-\s_k$.
\end{proof}

\begin{cor}\label{cor:spectrum-as-GH-mod}
Let $X$ be a finite metric space and $\l\ge2\diam X$, then
$$
\mst X=\l(\#X-1)-2\sum_{k=1}^{\#X-1}d_{GH}(\l\D_{k+1},X).
$$
\end{cor}

\section{Steiner problem in Gromov--Hausdorff space}
\markright{\thesection.~Steiner problem in Gromov--Hausdorff space}
In this section we prove the following

\begin{thm}\label{thm:GH-SMT-existence}
Let $M\ss\cM$ be a finite set such that each $X\in M$ is a finite metric space. Then $\SMT(M,\cM)\ne\0$, i.e., for such $M$ the Steiner problem has a solution.
\end{thm}

\begin{rk}
For arbitrary finite $M\ss\cM$ the problem is still open.
\end{rk}

\begin{proof}
Let $n=\#M$. In Section~\ref{sec:STM} we introduced model full Steiner trees, and we have shown how to use them for calculating the length of a Steiner minimal tree. Recall the corresponding definitions in our case. A full Steiner tree has the vertices of two types only: the ones of degree $1$ which we call boundary, and the ones of degree $3$ which we call interior. In model full Steiner trees which we use to treat the problem for such $M$, the vertex set is $\{1,\ldots,2n-2\}$, where $\{1,\ldots,n\}$ are reserved for the boundary vertices. We called two model full Steiner tree equivalent if there is an isomorphism between them fixed on the boundary. By $\cB_n$ we denoted the set of all model full Steiner trees with $n$ boundary vertices considered up to this equivalence.

Enumerate the points from $M$ in an arbitrary way, i.e., we consider a bijection $\v\:\{1,\ldots,n\}\to M$. Choose an arbitrary $G\in\cB_n$, and consider a network $\G$ of the type $G$ for which $\d\G=\v$. Then all such networks for given $G$ differ from each other only in the ``positions'' of their interior vertices,  thus the set $[G,\v]$ of such networks can be identified with $\cM^{n-2}$. Then we proved (Corollary~\ref{cor:smtInTermBinTrees}) that
$$
\smt(M,\cM)=\inf\bigl\{|\G|:\G\in[G,\v],\,G\in\cB_n\bigr\}.
$$

Choose an arbitrary $G\in\cB_n$ and any $\G\in[G,\v]$. We put $X_i=\G(i)$, then $M=\{X_1,\ldots,X_n\}$. For each $ij\in E(G)$ we choose an arbitrary $R_{ij}\in\cR_\opt(X_i,X_j)$ in such a way that $R_{ji}^{-1}=R_{ij}$. Let $X_k=\{x_k^i\}_{i=1}^{n_k}$, then for any $k\in\{1,\ldots,n\}$ and any $1\le i\le n_k$ we construct a network $\G_k^i\:\{1,\ldots,2n-2\}\to\sqcup_{j=1}^{2n-2}X_j$ as follows: in each $X_j$ we choose one point $x_j^{r_j}=:\G_k^i(j)$ such that
\begin{enumerate}
\item $x_k^{r_k}=x_k^i$;
\item for any $pq\in E(G)$ we have $(x_p^{r_p},x_q^{r_q})\in R_{pq}$
\end{enumerate}
(verify that it is always possible to do). Consider the set $\cT=\{\G_k^i\}$ consisting of all $\G_k^i$ over all possible $k$ and $i$.

For any $j\in\{1,\ldots,2n-2\}$ we put $V_j=\cup_{T\in\cT}\{T(j)\}$, i.e., we gather in each $X_j$ all points that are the images of the vertices of the constructed networks. Let us note that all $V_j$ have at most $N:=\sum_{k=1}^nn_k$ points, and for each $j\in\{1,\ldots,n\}$ it holds $V_j=X_j$.

Further, for each $p$ and $q$ such that $pq\in E(G)$ we denote by $R^t_{pq}$ the set of all pairs $(x_p,x_p)$ such that for some $T\in\cT$ we have $x_p=T(p)$ and $x_q=T(q)$, i.e., we gather all pairs forming the images of the edges of the constructed networks. Thus, we obtained correspondences $R^t_{pq}\in\cR(V_p,V_q)$ such that $R^t_{pq}\ss R_{pq}$, hence
$$
d_{GH}(V_p,V_q)\le\dis R^t_{pq}\le\dis R_{pq}=d_{GH}(X_p,X_q).
$$
Denote by $\G^t$ the network $\G^t\:p\mapsto V_p$. Since $\G^t|_{\{1,\ldots,n\}}=\v$ and $\G^t$ has the type $G$, then $\G^t\in[G,\v]$. Denote by $[G,\v]^t$ the set of all such $\G^t$. Since $|\G^t|\le|\G|$, then
$$
\smt(M,\cM)=\inf\bigl\{|\G^t|:\G^t\in[G,\v]^t,\,G\in\cB_n\bigr\}.
$$
However, all $X_i$ belong to $\cM_N$, therefore, $\smt(M,\cM)=\smt(M,\cM_N)$. Moreover, if we choose $\G^t$ such that all $V_i$, $i\ge n+1$, coincide with $V_k$ for some $k\le n$, then $\smt(M,\cM_N)\le\sum_{1\le p,q\le n}d_{GH}(X_p,X_q)=:D'$. Thus, it suffices to consider only $\G^t$ with $|\G^t|\le D'$, in particular, for such $\G^t$ and any $pq\in E(G)$ we have $d_{GH}(V_p,V_q)\le D'$.

Let us put $d=\max\{\diam X_k:k=1,\ldots,n\}$. Since for any $X,Y\in\cM$ we have $d_{GH}(X,Y)\ge\frac12|\diam X-\diam Y|$, then for each $X_p$ we have $\diam X_p\le d+2(n-2)D'=:D$ (all these estimates are rather rough, however, we do not need exact ones here). In account, we proved that all $X_p$ belongs to $\cM_N(D)$. By Problem~\ref{prb:cNnDIsCompact}, the space $\cM_N(D)$ is compact, that is why the continuous function $\ell\:\cM_N(D)^{n-2}\to\R$, $\ell\:(X_{n+1},\ldots,X_{2n-2})\mapsto|\G^t|$, attains its minimum at some $\G^t_0$. It remains to notice that $\cB_n$ is finite.
\end{proof}


\phantomsection
\renewcommand\bibname{References to Chapter~\thechapter}
\addcontentsline{toc}{section}{\bibname}
\renewcommand{\refname}{\bibname}

\vfill\eject

\section*{\Huge Exercises to Chapter~\thechapter}
\markright{Exercises to Chapter~\thechapter.}

\begin{exe}
Prove that the product topology on $X\x Y$ coincides with the one generated by the metric from Agreement~\ref{agree:MetricOnProduct}.
\end{exe}

\begin{exe}
Prove that the topologies generated by all $d_t$, and by $|\cdot|$ as well, coincide with the product topology of $X\x Y$.
\end{exe}

\begin{exe}
Prove that
\begin{enumerate}
\item a metric space $X$ is bounded if only for some $\e>0$ it holds $\cov(X,\e)<\infty$ (similarly, $\pack(X,\e)<\infty$);
\item a metric space $X$ is finite if and only if there exists $n$ such that $\cov(X,\e)\le n$ for all $\e>0$ (similarly, for $\pack(X,\e)$);
\item the functions $f(\e)=\cov(X,\e)$ and $g(\e)=\pack(X,\e)$ are monotonically decreasing.
\end{enumerate}
\end{exe}

\begin{exe}
Prove that the set $\cM_n(D)$ is compact, while $\cM_{[n]}(D)$ for $n>1$ is not.
\end{exe}

\vfill\eject
 \chapter{Calculating GH-distances to simplexes and some applications.}
 \markboth{\chaptername~\thechapter.~Calculating GH-distances to simplexes and some applications.}%
          {\chaptername~\thechapter.~Calculating GH-distances to simplexes and some applications.}

\begin{plan}
GH-distance to simplexes with more points, GH-distance to simplexes with at most the same number of points, generalized Borsuk problem, solution of generalized Borsuk problem in terms of GH-distances, clique covering number and chromatic number of simple graphs, their dualities, calculating these numbers in terms of GH-distances.
\end{plan}

By \emph{simplex\/} we mean a metric space in which all non-zero distances equal to each other. If $m$ is an arbitrary cardinal number, a simplex contain $m$ points, and all its non-zero distances equal $1$, then we denote this simplex by $\D_m$. Thus, $\l\D_m$, $\l>0$, is a simplex whose non-zero distances equal $\l$. Also, for arbitrary metric space $X$ and $\l=0$, the space $\l X$ coincides with $\D_1$.

\section{Gromov--Hausdorff distance to simplexes with more points}
\markright{\thesection.~Gromov--Hausdorff distance to simplexes with more points}

The next result generalizes Theorem~4.1 from~\cite{IvaTuzSimpDist8}.

\begin{thm}\label{thm:dist-n-simplex-bigger-dim}
Let $X$ be an arbitrary metric space, $m>\#X$ a cardinal number, and $\l\ge0$, then
$$
2d_{GH}(\l\D_m,X)=\max\{\l,\diam X-\l\}.
$$
\end{thm}

\begin{proof}
If $X$ is unbounded, then $2d_{GH}(\l\D_m,X)=\infty$ by Example~\ref{examp:dGHbelowEstimate}, and we get what is required.

Now, let $\diam X<\infty$.

If $\#X=1$, then $\diam X=0$, and, by Example~\ref{examp:GH_simple}, we have
$$
2d_{GH}(\l\D,X)=\diam\l\D=\l=\max\{\l,\diam X-\l\}.
$$

If $\l=0$, then, by Example~\ref{examp:GH_simple}, we have
$$
2d_{GH}(\D_1,X)=\diam X=\max\{\l,\diam X-\l\}.
$$

Let $\#X>1$ and $\l>0$. Choose an arbitrary $R\in\cR(\l\D_m,X)$. Since $\#X<m$ and $\l>0$, then there exists $x\in X$ such that $\#R^{-1}(x)\ge2$, thus, $\dis R\ge\l$ and $2d_{GH}(\l\D_m,X)\ge\l$.

Consider an arbitrary sequence $(x_i,y_i)\in X\x X$ such that $|x_iy_i|\to\diam X$. If it contains a subsequence $(x_{i_k},y_{i_k})$ such that for each $i_k$ there exists $z_k\in\l\D$, $(z_k,x_{i_k})\in R$, $(z_k,y_{i_k})\in R$, then $\dis R\ge\diam X$ and
$$
2d_{GH}(\l\D_m,X)\ge\max\{\l,\diam X\}\ge\max\{\l,\diam X-\l\}.
$$

If such subsequence does not exist, then there exists a subsequence $(x_{i_k},y_{i_k})$ such that for any  $i_k$ there exist distinct $z_k,w_k\in\l\D_m$, $(z_k,x_{i_k})\in R$, $(w_k,y_{i_k})\in R$, and, therefore,
$$
2d_{GH}(\l\D_m,X)\ge\max\bigl\{\l,|\diam X-\l|\bigr\}\ge\max\{\l,\diam X-\l\}.
$$

Thus, in the both cases we have $2d_{GH}(\l\D,X)\ge\max\bigl\{\l,\diam X-\l\bigr\}$.

Choose an arbitrary $x_0\in X$, then, by assumption, $\#X>1$, and, thus, the set $X\sm\{x_0\}$ is not empty. Since $\#X<m$, then $\l\D_m$ contains a subset $\l\D'$ of the same cardinality with $X\sm\{x_0\}$. Let $g\:\l\D'\to X\sm\{x_0\}$ be an arbitrary bijection, and $\l\D''=\l\D_m\sm\l\D'$, then $\#\l\D''>1$. Consider the following correspondence:
$$
R_0=\Bigl\{\bigl(z',g(z')\bigr):z'\in\l\D'\Bigr\}\cup\bigl(\l\D''\x\{x_0\}\bigr).
$$
Then we can apply Proposition~\ref{prop:disRforPartition}, thus we have
$$
\dis R_0=\sup\{\l,|x_1x'_1|-\l,\l-|x_2x'_2|:x_1,x'_1,x_2,x'_2\in X,\,x_1\ne x'_1,\,x_2\ne x'_2\}=\max\{\l,\diam X-\l\},
$$
therefore,
$$
2d_{GH}(\l\D,X)=\max\{\l,\diam X-\l\},
$$
what is required.
\end{proof}

\section{Gromov--Hausdorff distance to simplexes with at most the same number of points}
\markright{\thesection.~Gromov--Hausdorff distance to simplexes with at most the same number of points}

Let $X$ be an arbitrary set different from singleton, $2\le m\le\#X$ a cardinal number, and $\l>0$. Under notations from Section~\ref{sec:Irreducible}, consider an arbitrary $D\in\cD_m(X)$, any bijection $g\:\l\D_m\to D$, and construct the correspondence $R_D\in\cR(\l\D_m,X)$ in the following way:
$$
R_D=\bigcup_{z\in\l\D_m}\{z\}\x g(z).
$$
Clearly that each correspondence $R_D$ is irreducible.

From Proposition~\ref{prop:disRforPartition} we get

\begin{prop}\label{prop:disRD}
Let $X\ne\D_1$ be an arbitrary metric space, $2\le m\le\#X$ a cardinal number, and $\l>0$. Then for any $D\in\cD_m(X)$ it holds
$$
\dis R_D=\max\{\diam D,\,\l-\a(D),\,\b(D)-\l\}.
$$
\end{prop}

\begin{proof}
If $X$ is unbounded, then $\dis R=\infty$ for any $R\in\cR(\l\D_m,X)$. Since $m\ge2$, for any $D=\{X_i\}_{i\in I}\in\cD_m(X)$ we have either $\diam D=\infty$, or $\b(D)=\infty$. Indeed,  if $\diam D<\infty$ and $\b(D)<\infty$ then for any $x,y\in X$ either $x,y\in X_i$, thus $|xy|\le\diam D$, or $x\in X_i$, $y\in X_j$, $i\ne j$, and $|xy|\le|X_iX_j|\le\b(D)$, therefore $X$ is bounded. Thus, for unbounded $X$ the right-hand side of the considered equation is infinite as well, thus we get what is required.

Now, let $\diam X<\infty$. By Proposition~\ref{prop:disRforPartition}, we have
$$
\dis R_D=\sup\bigl\{\diam D,\,\l-|X_iX_j|,\,|X_iX_j|'-\l: i,j\in I,\,i\ne j\bigr\}=\max\{\diam D,\,\l-\a(D),\,\b(D)-\l\},
$$
that completes the proof.
\end{proof}

\begin{cor}\label{cor:disRD}
Let $X\ne\D_1$ be an arbitrary metric space, $2\le m\le\#X$ a cardinal number, and $\l>0$. Then for any $D\in\cD_m(X)$ it holds
$$
\dis R_D=\max\{\diam D,\,\l-\a(D),\,\diam X-\l\}.
$$
\end{cor}

\begin{proof}
Again, for unbounded $X$ the equation evidently holds.

Consider now the case of bounded $X$. Notice that $\diam D\le\diam X$ and $\b(D)\le\diam X$. In addition, if $\diam D<\diam X$, and $(x_i,y_i)\in X\x X$ is a sequence such that $|x_iy_i|\to\diam X$, then, starting from some $i$, the points $x_i$ and $y_i$ belong to different elements of $D$, therefore, in this case we have $\b(D)=\diam X$, and the formula is proved.

Now, let $\diam D=\diam X$, then $\b(D)-\l\le\diam X$ and $\diam X-\l\le\diam X$, thus
$$
\max\{\diam D,\,\l-\a(D),\,\b(D)-\l\}=\max\{\diam X,\,\l-\a(D)\}=\max\{\diam D,\,\l-\a(D),\,\diam X-\l\},
$$
that completes the proof.
\end{proof}

\begin{prop}\label{prop:GH-dist-RD}
Let $X\ne\D_1$ be an arbitrary metric space, and $2\le m\le\#X$ a cardinal number, and $\l>0$. Then
$$
2d_{GH}(\l\D_m,X)=\inf_{D\in\cD_m(X)}\dis R_D.
$$
\end{prop}

\begin{proof}
The case of unbounded $X$ is trivial, so, let $X$ be bounded. By Corollary~\ref{cor:GH-distance-irreducinle},
$$
2d_{GH}(\l\D_m,X)=\inf_{R\in\cR^0(\l\D_m,X)}\dis R,
$$
thus it suffices to prove that for any irreducible correspondence $R\in\cR^0(\l\D_m,X)$ there exists $D\in\cD_m(X)$ such that $\dis R_D\le\dis R$.

Let us choose an arbitrary $R\in\cR^0(\l\D_m,X)$ such that it cannot be represented in the form $R_D$, then the partition $D^R_{\l\D_m}$ is not pointwise, i.e., there exists $x\in X$ such that $\#R^{-1}(x)\ge 2$, therefore, $\dis R\ge\l$.

Define a metric on the set $D^R_{\l\D_m}$ to be equal $\l$ between any its distinct elements, then this metric space is isometric to a  simplex $\l\D'_n$, $n\le m$. The correspondence $R$ generates naturally  another correspondence $R'\in\cR(\l\D'_n,X)$, namely, if $D^R_{\l\D_m}=\{\D_j\}_{j\in J}$, and $f_R\:D^R_{\l\D_m}\to D^R_X$ is the bijection generated by $R$, then
$$
R'=\bigcup_{j\in J}\{\D_j\}\x f_R(\D_j).
$$
It is easy to see that $\dis R=\max\{\l,\,\dis R'\}$. Moreover, $R'$ is generated by the partition $D'=D^R_X$, i.e., $R'=R_{D'}$, thus, by Corollary~\ref{cor:disRD}, we have
$$
\dis R'=\max\{\diam D',\,\l-\a(D'),\,\diam X-\l\},
$$
and hence,
$$
\dis R=\max\{\l, \diam D',\,\l-\a(D'),\,\diam X-\l\}=\max\{\l,\,\diam D',\,\diam X-\l\}.
$$
Since $n\le m$, the partition $D'$ has a subpartition $D\in\cD_m(X)$. Clearly, $\diam D\le\diam D'$, therefore,
$$
\dis R_D=\max\{\diam D,\,\l-\a(D),\,\diam X-\l\}\le\max\{\diam D',\,\l,\,\diam X-\l\}=\dis R,
$$
q.e.d.
\end{proof}

Considering separately the case $\l=0$, we get the following

\begin{cor}\label{cor:GH-dist-alpha-beta}
Let $X\ne\D_1$ be an arbitrary metric space, $2\le m\le\#X$ a cardinal number, and $\l\ge0$. Then
$$
2d_{GH}(\l\D_m,X)=\inf_{D\in\cD_m(X)}\max\{\diam D,\,\l-\a(D),\,\diam X-\l\}.
$$
\end{cor}

For any metric space $X$ put
$$
\e(X)=\inf\bigl\{|xy|:x,y\in X,\,x\ne y\bigr\}.
$$
Notice that $\e(X)\le\diam X$, and for a bounded $X$ the equality holds, if and only if $X$ is a simplex.

Corollary~\ref{cor:GH-dist-alpha-beta} immediately implies the following result that is proved in~\cite{IvaTuzSimpDist8}.

\begin{thm}[\cite{IvaTuzSimpDist8}]\label{thm:dist-n-simplex-same-dim}
Let $X\ne\D_1$ be a finite metric space, $m=\#X$, and $\l\ge0$, then
$$
2d_{GH}(\l\D_m,X)=\max\bigl\{\l-\e(X),\,\diam X-\l\bigr\}.
$$
\end{thm}

\section{Generalized Borsuk problem}
\markright{\thesection.~Generalized Borsuk problem}

Classical Borsuk Problem deals with partitions of subsets of Euclidean space into parts having smaller diameters. We generalize the Borsuk problem to arbitrary bounded metric spaces and partitions of arbitrary cardinality. Let $X$ be a bounded metric space, $m$ a cardinal number such that $2\le m\le\#X$, and $D=\{X_i\}_{i\in I}\in\cD_m(X)$. We say that  $D$ is a partition into subsets having \emph{strictly smaller diameters}, if there exists $\e>0$ such that $\diam X_i\le\diam X-\e$ for all $i\in I$.

By \emph{Generalized Borsuk problem\/} we call the following one: Is it possible to partition a bounded metric space $X$ into a given, probably infinite, number of subsets, each of which has a strictly smaller diameter than $X$?

We give a solution to the Generalized Borsuk problem in terms of the Gromov--Hausdorff distance.

\begin{thm}\label{thm:Borsuk}
Let $X$ be an arbitrary bounded metric space and $m$ a cardinal number such that $2\le m\le\#X$. Choose an arbitrary number $0<\l<\diam X$, then $X$ can be partitioned into $m$ subsets having strictly smaller diameters if and only if $2d_{GH}(\l\D_m,X)<\diam X$.
\end{thm}

\begin{proof}
For the $\l$ chosen, due Corollary~\ref{cor:GH-dist-alpha-beta}, we have $2d_{GH}(\l\D_m,X)\le\diam X$, and the equality holds if and only if for each $D\in\cD_m(X)$ we have $\diam D=\diam X$. The latter means that there is no partition of the space $X$ into $m$ parts having strictly smaller diameters.
\end{proof}

\begin{cor}
Let $d>0$ be a real number, and $m\le n$ cardinal numbers. By $\cM_n$ we denote the set of isometry classes of bounded metric spaces of cardinality at most $n$, endowed with the Gromov--Hausdorff distance. Choose an arbitrary $0<\l<d$. Then the intersection
$$
S_{d/2}(\D_1)\cap S_{d/2}(\l\D_m)
$$
of the spheres, considered as the spheres in $\cM_n$, does not contain spaces, whose  cardinality is less than $m$, and consists exactly of all metric spaces from $\cM_n$, whose diameters are equal to $d$ and that cannot be partitioned into $m$ subsets of strictly smaller diameters.
\end{cor}

\begin{proof}
Let $X$ belong to the intersection of the spheres, then $\diam X=d$ in accordance with Example~\ref{examp:GH_simple}. If $m>\#X$, then, due to Theorem~\ref{thm:dist-n-simplex-bigger-dim}, we have
$$
2d_{GH}(\l\D,X)=\max\{\l,\diam X-\l\}<d,
$$
therefore $X\not\in S_{d/2}(\l\D_m)$, that proves the first statement of Corollary.

Now let $m\le\#X$. Since $\diam X=d$ and $2d_{GH}(\l\D_m,X)=d$, then, due to Theorem~\ref{thm:Borsuk}, the space $X$ cannot be partitioned into $m$ subsets of strictly smaller diameters.

Conversely, each $X$ of the diameter $d$, such that $m\le\#X$ and which cannot be partitioned into $m$ subsets of strictly smaller diameter, lies in the intersection of the spheres by Theorem~\ref{thm:Borsuk}.
\end{proof}

\section{Calculating clique covering and chromatic numbers of a graph}
\markright{\thesection.~Calculating clique covering and chromatic numbers of a graph}

Recall that a subgraph of an arbitrary simple graph $G$ is called \emph{a clique}, if any its two vertices are connected by an edge, i.e., the clique is a subgraph which is a complete graph itself. Notice that each single-vertex subgraph is also a clique. For convenience, the vertex set of a clique is also referred as \emph{a clique}.

On the set of all cliques, an ordering with respect to inclusion is naturally defined, and hence, due to the above remarks, a family of maximal cliques is uniquely defined; this family forms \emph{a cover of the graph $G$} in the following sense: the union of all vertex sets of all maximal cliques coincides with the vertex set $V(G)$ of the graph $G$.

If one does not restrict himself by maximal cliques, then, generally speaking, one can find other families of cliques covering the graph $G$. One of the classical problems of the Graph Theory is to calculate the minimal possible number of cliques covering a finite simple graph $G$. This number is referred as \emph{the clique covering number\/} and is often denoted by $\theta(G)$. It is easy to see that the value $\theta(G)$ is also equal to the least number of cliques whose vertex sets form a partition of $V(G)$.

Another popular problem is to find the least possible number of colors that is necessary to color the vertices of a simple finite graph $G$ in such a way that adjacent vertices have different colors. This number is denoted by $\g(G)$ and is referred as \emph{the chromatic number of the graph $G$}.

For a simple graph $G$, by $G'$ we denote its \emph{dual graph}, i.e., the graph with the same vertex set and the complementary set of edges (two vertices of $G'$ are adjacent if and only if they are not adjacent in $G$).

\begin{prb}
For any simple finite graph $G$ it holds $\theta(G)=\g(G')$.
\end{prb}

Let $G=(V,E)$ be an arbitrary finite graph. Fix two real numbers $a<b\le2a$ and define a metric on $V$ as follows: the distance between adjacent vertices equals $a$, and nonadjacent vertices equals $b$. Then a subset $V'\ss V$ has diameter $a$ if and only if $G(V')\ss G$ is a clique. This implies that each clique covering number equals to the least cardinality of partitions of the metric space $V$ onto subsets of (strictly) smaller diameter. However, this number was calculated in Theorem~\ref{thm:Borsuk}. Thus, we get the following

\begin{cor}\label{cor:clique}
Let $G=(V,E)$ be an arbitrary finite graph. Fix two real numbers $a<b\le2a$ and define a metric on $V$ as follows: the distance between adjacent vertices equals $a$, and nonadjacent vertices equals $b$. Let $m$ be the greatest positive integer $k$ such that $2d_{GH}(a\D_k,V)=b$ \(in the case when there is no such $k$, we put $m=0$\). Then $\theta(G)=m+1$.
\end{cor}

\begin{prb}
Consider simple finite graphs $G=(V,E)$ for which the clique covering numbers $\theta(G)$ are known, and get the Gromov--Hausdorff distances between the corresponding metric spaces $V$ and simplexes $\l\D_m$ with $m\le\theta(G)$. Verify explicitly that for $k>m$ these distances are less than $\diam V$.
\end{prb}

Because of the duality between clique and chromatic numbers, we get

\begin{cor}\label{cor:chrom}
Let $G=(V,E)$ be an arbitrary finite graph. Fix two real numbers $a<b\le2a$ and define a metric on $V$ as follows: the distance between adjacent vertices equals $b$, and nonadjacent vertices equals $a$. Let $m$ be the greatest positive integer $k$ such that $2d_{GH}(a\D_k,V)=b$ \(in the case when there is no such $k$, we put $m=0$\). Then $\g(G)=m+1$.
\end{cor}

\begin{prb}
Consider simple finite graphs $G=(V,E)$ for which the chromatic numbers $\g(G)$ are known, and get the Gromov--Hausdorff distances between the corresponding metric spaces $V$ and simplexes $\l\D_m$ with $m\le\g(G)$. Verify explicitly that for $k>m$ this distances are less than $\diam V$.
\end{prb}


\phantomsection
\renewcommand\bibname{References to Chapter~\thechapter}
\addcontentsline{toc}{section}{\bibname}
\renewcommand{\refname}{\bibname}

\vfill\eject

\section*{\Huge Exercises to Chapter~\thechapter}
\markright{Exercises to Chapter~\thechapter.}

\begin{exe}
For any simple finite graph $G$ it holds $\theta(G)=\g(G')$.
\end{exe}

\begin{exe}
Consider simple finite graphs $G=(V,E)$ for which the clique covering numbers $\theta(G)$ are known, and get the Gromov--Hausdorff distances between the corresponding metric spaces $V$ and simplexes $\l\D_m$ with $m\le\theta(G)$. Verify explicitly that for $k>m$ these distances are less than $\diam V$.
\end{exe}

\begin{exe}
Consider simple finite graphs $G=(V,E)$ for which the chromatic numbers $\g(G)$ are known, and get the Gromov--Hausdorff distances between the corresponding metric spaces $V$ and simplexes $\l\D_m$ with $m\le\g(G)$. Verify explicitly that for $k>m$ this distances are less than $\diam V$.
\end{exe}


\begin{thebibliography}{9}
\bibitem{Engelking1} R.~Engelking, \emph{General Topology}, Heldermann, Berlin, 1989.
\bibitem{Nadler1}  S.~Nadler, \emph{Hyperspaces of Sets}, 1978.
\end{thebibliography}

\begin{thebibliography}{9}
\bibitem{Deza2} M.M.Deza, E.Deza, \emph{Encyclopedia of Distances}. Springer, 2009.
\end{thebibliography}

\begin{thebibliography}{9}
\bibitem{BurBurIvaEng3} D.Burago, Yu.Burago, S.Ivanov, \emph{A Course in Metric Geometry}. Graduate Studies in Mathematics, vol.33, A.M.S., Providence, RI, 2001.
\end{thebibliography}

\begin{thebibliography}{9}
\bibitem{IvaNikTuz} Ivanov A.O., Nikonov I.M., Tuzhilin A.A. \emph{Sets admitting connection by graphs of
finite length}. Sbornik: Mathematics, 2005, 196:6, 845–884.
\bibitem{IvaTuzInfTrees} Ivanov A.O., Nikonov I.M., Tuzhilin A.A. \emph{Minimal Spanning Trees on Infinite Sets}. 2014, ArXiv:1403.383.
\bibitem{ITCRC4} Ivanov A.O., Tuzhilin A.A. \emph{Minimal Networks. Steiner Problem and Its Generalizations}. CRC Press, 1994.
\bibitem{HwRiW4} Hwang F.K., Richards D.S., Winter P. \emph{The Steiner Tree Problem}. Annals of Discrete Mathematics 53, North-Holland: Elsevier, 1992. ISBN 0-444-89098-X.
\bibitem{Kruscal4} Kruskal J.B. \emph{On the shortest spanning subtree of a graph and the traveling salesman problem}, Proceedings of the American Mathematical Society, 1956, v. 7, 48--50.
\bibitem{Prim4} Prim R.C. \emph{Shortest connection networks And some generalizations}, Bell System Technical Journal, 1957, v. 36, N 6, pp. 1389--1401.
\bibitem{Fermat4} Fermat P. de, Ed. H.Tannery, ed., \emph{Oeuvres}, v. 1, Paris 1891, Supplement: Paris 1922, p. 153.
\bibitem{JK4} Jarn\'{\i}k V., K\"ossler O. \emph{O minim\'aln\'{\i}ch grafech obsahuj\'{\i}c\'{\i}ch n dan\'ych bodu}, \^Cas, P\^estov\^an\'{\i} Mat. (Essen) v. 63, pp. 223--235.
\bibitem{CourRobb4} Courant R., Robbins H. \emph{What Is Mathematics?}, Oxford University Press, 1941.
\bibitem{Melzak4} Melzak Z.A. \emph{On the problem of Steiner}, Canad. Math. Bull., 1960, N 4, pp. 143--148.
\bibitem{Hwang4} Hwang F.K. \emph{A linear time algorithm for full Steiner trees}, Oper. Res.
Letter, 1986, N 5, pp. 235--237.
\bibitem{WWZ4} Warme D.M., Winter P., Zachariasen M. \emph{Exact Algorithms for Plane Steiner Tree Problems}. A Computational Study, Technical Report DIKU-TR-98/11, Department of Computer Science, University of Copenhagen, DENMARK, 1998.
\bibitem{GGJ4} Garey M.R., Graham R.L., Johnson D.S. \emph{The complexity of computing Steiner minimal trees}, SIAM J. Appl. Math., 1977, v. 32, N 4, pp. 835--859.
\bibitem{GromovMinFil4} Gromov M. \emph{Filling Riemannian manifolds}, J. Diff. Geom., 1983, v. 18, N 1, pp. 1–147.
\bibitem{IvaTuzMinFil4} Ivanov A.O., Tuzhilin A.A. \emph{One-dimensional Gromov minimal filling}.   arXiv:1101.0106v2 [math.MG] (http://arxiv.org).
\end{thebibliography}

\begin{thebibliography}{9}
\bibitem{BurBurIvaEng5} D.Burago, Yu.Burago, S.Ivanov, \emph{A Course in Metric Geometry}. Graduate Studies in Mathematics, vol.33, A.M.S., Providence, RI, 2001.
\bibitem{Nadler5} S.Nadler, \emph{Hyperspaces of Sets}, 1978.
\bibitem{IllanesNadler5} A.Illanes, S.Nadler,  \emph{Hyperspaces}, 1999.
\bibitem{Bryant-convexity5} Bryant V.W. \emph{The convexity of the subset space of a metric space}, Compositio Mathematica, 1970, vol. 22, N 4, p. 383--385.
\end{thebibliography}

\begin{thebibliography}{9}
\bibitem{BurBurIva6} Burago D., Burago Yu., Ivanov S. \emph{A Course in Metric Geometry}. Graduate Studies in Mathematics, vol.33, A.M.S., Providence, RI, 2001.
\bibitem{Ghanaat6} Ghanaat P. ``Gromov-Hausdorff distance and applications''. In: Summer school ``Metric Geometry'', Les Diablerets, August 25--30, 2013, \url{https://math.cuso.ch/fileadmin/math/document/gromov-hausdorff.pdf}
\end{thebibliography}

\begin{thebibliography}{9}
\bibitem{BurBurIva7} Burago D., Burago Yu., Ivanov S. \emph{A Course in Metric Geometry}. Graduate Studies in Mathematics, vol.33, A.M.S., Providence, RI, 2001.
\end{thebibliography}

\begin{thebibliography}{9}
\bibitem{IvaTuzSimpDist8} A.O.Ivanov, S.Iliadis, and A.A.Tuzhilin, \emph{Geometry of Compact Metric Space in Terms of Gromov-Hausdorff Distances to Regular Simplexes}. ArXiv e-prints, {\tt arXiv:1607.06655}, 2016.
\end{thebibliography}
\end{document}